\documentclass[11pt]{amsart}
\usepackage[hmargin=3cm,vmargin=3cm]{geometry}

\usepackage[latin1]{inputenc}
\usepackage{xspace,amssymb,amsfonts,euscript}
\usepackage{amsthm,amsmath,amscd,stmaryrd,latexsym}
\usepackage{palatino, euscript}
\usepackage{enumitem}

\usepackage[small,nohug]{diagrams}
\diagramstyle[labelstyle=\scriptstyle]

\usepackage{graphicx}
\usepackage[all]{xy}
\usepackage{tikz}
\usetikzlibrary{calc, matrix, arrows, decorations.pathmorphing}
\tikzset{%
	descr/.style={fill=white},
	baseline={([yshift=-\the\dimexpr\fontdimen22\textfont2\relax]
	                    current bounding box.center)},
	->,>=angle 90
}
\usepackage[dvips]{epsfig}
\usepackage{pinlabel}

\usepackage{tensor}

\newcommand{\ig}[2]{\vcenter{\xy (0,0)*{\includegraphics[scale=#1]{#2}} \endxy}}
\newcommand{\igv}[2]{\vcenter{\xy (0,0)*{\reflectbox{\includegraphics[scale=#1, angle=180]{#2}}} \endxy}}

\newtheorem{thm}{Theorem}[section]
\newtheorem{lemma}[thm]{Lemma}

\newtheorem{prop}[thm]{Proposition}
\newtheorem{cor}[thm]{Corollary}
  
\newtheorem{conj}[thm]{Conjecture}
\newtheorem{goal}[thm]{Goal}
\newtheorem{porism}[thm]{Porism}

\newtheorem*{prop*}{Proposition}
\newtheorem*{lemma*}{Lemma}

\theoremstyle{definition}
\newtheorem{defn}[thm]{Definition}

\newtheorem{notation}[thm]{Notation}
\newtheorem{ex}[thm]{Example}
\newtheorem{example}[thm]{Example}
\newtheorem{exercise}[thm]{Exercise}

\theoremstyle{remark}
\newtheorem{remark}[thm]{Remark}

\numberwithin{equation}{section}

  \def\bg{{\mathfrak b}}  
    
    \def\DM{{\mathbb{D}}}

  \def\gg{{\mathfrak g}}  
  \def\hg{{\mathfrak h}}

    \def\KM{{\mathbb{K}}}
    
\def\MG{{\mathfrak M}}

    \def\AC{{\mathcal{A}}}
  \def\bb{{\mathbf b}}  \def\BC{{\mathcal{B}}}
    
    \def\DC{{\mathcal{D}}}
    \def\EC{{\mathcal{E}}}
    \def\FC{{\mathcal{F}}}
    \def\GC{{\mathcal{G}}}
\def\HB{{\mathbf H}}    \def\HC{{\mathcal{H}}}
    \def\IC{{\mathcal{I}}}
    \def\JC{{\mathcal{J}}}
    \def\KC{{\mathcal{K}}}
    \def\LC{{\mathcal{L}}}
    \def\MC{{\mathcal{M}}}
    \def\NC{{\mathcal{N}}}
    \def\OC{{\mathcal{O}}}
    \def\PC{{\mathcal{P}}}
    \def\QC{{\mathcal{Q}}}
    \def\RC{{\mathcal{R}}}
    \def\SC{{\mathcal{S}}}
    \def\TC{{\mathcal{T}}}
    \def\UC{{\mathcal{U}}}
    \def\VC{{\mathcal{V}}}
    \def\WC{{\mathcal{W}}}

    \def\ZC{{\mathcal{Z}}}

\def\al{\alpha}
\def\be{\beta}
\def\ga{\gamma}
\def\Ga{\Gamma}

\def\ep{\varepsilon}

\def\la{\lambda}
\def\La{\Lambda}
\def\si{\sigma}

\def\om{\omega}
\def\Om{\Omega}

\let\phi=\varphi
\let\tilde=\widetilde

\usepackage{bbm}
\def\C{{\mathbbm C}}

\def\R{{\mathbbm R}}
\def\Z{{\mathbbm Z}}

\def\one{\mathbbm{1}}

\newcommand{\un}{\underline}

\newcommand{\sqmatrix}[1]{\left[\begin{matrix} #1\end{matrix}\right]}

\newcommand{\ot}{\otimes}
\newcommand{\sot}{\dot{\otimes}}
\newcommand{\pa}{\partial}
\newcommand{\co}{\colon}

\renewcommand{\to}{\rightarrow}
\newcommand{\from}{\leftarrow}
\newcommand{\into}{\hookrightarrow}
\newcommand{\onto}{\twoheadrightarrow}

\newcommand{\sumset}{\stackrel{\scriptstyle{\oplus}}{\scriptstyle{\subset}}}

\newcommand\simto{\xrightarrow{
   \,\smash{\raisebox{-0.65ex}{\ensuremath{\scriptstyle\sim}}}\,}}

\newcommand{\define}{:=}

\renewcommand{\sl}{\mathfrak{sl}}
\newcommand{\gl}{\mathfrak{gl}}
\newcommand{\mt}{\emptyset}

\newcommand{\refequal}[1]{\xy {\ar@{=}^{#1}
(-1,0)*{};(1,0)*{}};
\endxy}

\newcommand\restr[2]{{% we make the whole thing an ordinary symbol
  \left.\kern-\nulldelimiterspace % automatically resize the bar with \right
  #1 % the function
  \vphantom{\big|} % pretend it's a little taller at normal size
  \right|_{#2} % this is the delimiter
  }}

\newcommand{\Kar}{\mathbf{Kar}}
\newcommand{\Rep}{\mathbf{Rep}}
\newcommand{\Vect}{\mathbf{Vect}}
\newcommand{\Fund}{\mathbf{Fund}}
\newcommand{\Webs}{\mathbf{Webs}}

\DeclareMathOperator{\ext}{ext}
\DeclareMathOperator{\Gr}{Gr}
\DeclareMathOperator{\Fl}{Fl}
\DeclareMathOperator{\ev}{ev}
\DeclareMathOperator{\adj}{adj}
\DeclareMathOperator{\rt}{rt}
\DeclareMathOperator{\wt}{wt}
\DeclareMathOperator{\longest}{long}
\DeclareMathOperator{\Br}{Br}
\DeclareMathOperator{\PBr}{PBr}
\DeclareMathOperator{\fin}{fin}
\DeclareMathOperator{\aff}{aff}
\DeclareMathOperator{\FT}{FT}

\DeclareMathOperator{\Hom}{Hom}
\DeclareMathOperator{\HOM}{HOM}

\DeclareMathOperator{\END}{END}

\DeclareMathOperator{\id}{id}

\newcommand{\sqot}{\boxtimes}

\newcommand{\EDiag}{\EC \Diag}

\renewcommand{\bot}{\textrm{bot}}
\newcommand{\cycle}{\copyright}
\newcommand{\Hecke}{\HB}
\newcommand{\BSBim}{\mathbb{BS}\mathrm{Bim}}
\newcommand{\SBim}{\mathbb{S}\mathrm{Bim}}
\newcommand{\Hot}{\KC}
\newcommand{\collapse}{\pi}
\DeclareMathOperator{\BS}{BS}
\newcommand{\MaxIdeal}{\MG}

\newcommand{\lineblue}{\ig{.5}{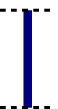}}
\newcommand{\linered}{\ig{.5}{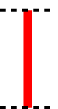}}
\newcommand{\linegreen}{\ig{.5}{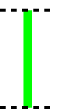}}
\newcommand{\linepurple}{\ig{.5}{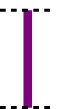}}
\newcommand{\startdotblue}{\ig{.5}{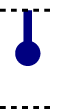}}
\newcommand{\startdotred}{\ig{.5}{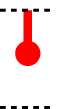}}
\newcommand{\startdotgreen}{\ig{.5}{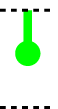}}
\newcommand{\startdotpurple}{\ig{.5}{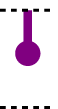}}
\newcommand{\finaldotblue}{\ig{.5}{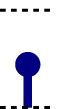}}
\newcommand{\finaldotred}{\ig{.5}{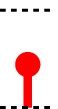}}
\newcommand{\finaldotgreen}{\ig{.5}{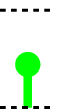}}

\newcommand{\barbblue}{\ig{.5}{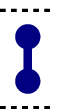}}
\newcommand{\barbred}{\ig{.5}{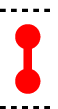}}
\newcommand{\barbgreen}{\ig{.5}{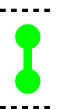}}
\newcommand{\barbpurple}{\ig{.5}{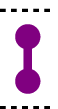}}
\newcommand{\brokenblue}{\ig{.5}{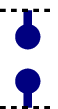}}
\newcommand{\brokenred}{\ig{.5}{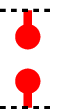}}
\newcommand{\brokengreen}{\ig{.5}{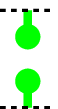}}
\newcommand{\dotfirsttoteal}{\ig{.5}{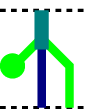}}
\newcommand{\dotfirstfromteal}{\ig{.5}{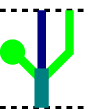}}
\newcommand{\splitred}{\ig{.5}{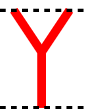}}
\newcommand{\splitblue}{\ig{.5}{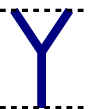}}

\newcommand{\mergered}{\ig{.5}{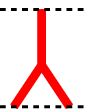}}
\newcommand{\mergeblue}{\ig{.5}{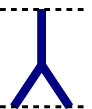}}
\newcommand{\mergegreen}{\ig{.5}{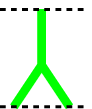}}
\newcommand{\capredcupblue}{\ig{.5}{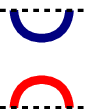}}
\newcommand{\capbluecupred}{\ig{.5}{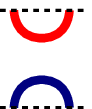}}
\newcommand{\redtoblue}{\ig{.5}{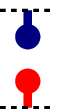}}
\newcommand{\bluetored}{\ig{.5}{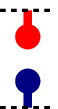}}
\newcommand{\greentored}{\ig{.5}{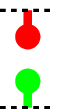}}

\newcommand{\greentoblue}{\ig{.5}{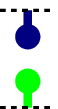}}

\newcommand{\capred}{\ig{.5}{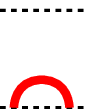}}
\newcommand{\capblue}{\ig{.5}{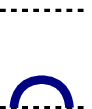}}
\newcommand{\capgreen}{\ig{.5}{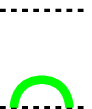}}
\newcommand{\cupred}{\ig{.5}{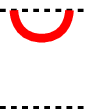}}
\newcommand{\cupblue}{\ig{.5}{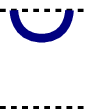}}

\renewcommand{\pitchfork}{\ig{.5}{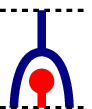}}
\newcommand{\rainbow}{\ig{.5}{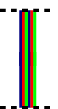}}
\newcommand{\diffstartblue}{\ig{.5}{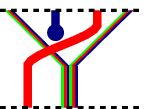}}
\newcommand{\diffstartred}{\ig{.5}{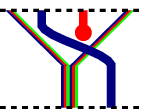}}
\newcommand{\difffinalblue}{\ig{.5}{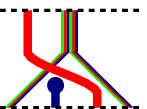}}
\newcommand{\difffinalred}{\ig{.5}{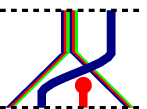}}
\newcommand{\thatterm}{\ig{.5}{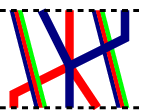}}
\newcommand{\foobara}{\ig{.5}{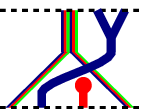}}
\newcommand{\foobarb}{\ig{.5}{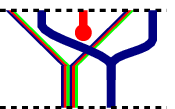}}
\newcommand{\dotdiff}{{
\labellist
\small\hair 2pt
\pinlabel {\tiny $d$} [ ] at 8 15
\endlabellist
\centering
\ig{.5}{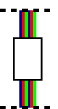}
}}
\newcommand{\thindotdiff}{{
\labellist
\small\hair 2pt
\pinlabel {\tiny $d$} [ ] at 7 15
\endlabellist
\centering
\ig{.5}{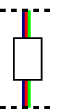}
}}
\newcommand{\cupcapbow}{\ig{.5}{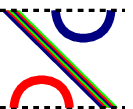}}
\newcommand{\poly}[1]{{
\labellist
\small\hair 2pt
 \pinlabel {$#1$} [ ] at 7 15
\endlabellist
\centering
\ig{.5}{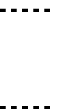}
}}

\newcommand{\Xbg}{\ig{.5}{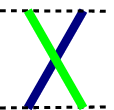}}
\newcommand{\Xgb}{\ig{.5}{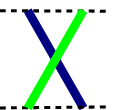}}
\newcommand{\mixedOred}{\ig{1}{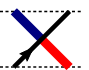}}

\newcommand{\BXto}{\vartheta}

\newcommand{\startdotsign}{\forall}
\newcommand{\finaldotsign}{A}

\newcommand{\fund}{\varpi}
\newcommand{\inv}{^{-1}}

\newcommand{\ImM}{\IC}

\newcommand{\weights}{{\rm Wt}}

\newcommand{\Diag}{\DC}

\DeclareMathOperator{\extg}{extg}
\DeclareMathOperator{\exts}{exts}
\DeclareMathOperator{\wind}{wind}
\newcommand{\unwind}{\un{\wind}}

\DeclareMathOperator{\SSat}{Satake}

\newcommand{\killit}{\nabla}
\DeclareMathOperator{\Ch}{Ch}
\DeclareMathOperator{\Wak}{Wak}
\DeclareMathOperator{\Hilb}{Hilb}
\DeclareMathOperator{\FHilb}{FHilb}
\DeclareMathOperator{\Cone}{Cone}

\newcommand{\who}{\Phi}
\newcommand{\notme}{\aleph}
\newcommand{\rotisom}{\phi}
\newcommand{\BY}{C}
\newcommand{\by}{c}
\newcommand{\idemone}{e_{\uparrow}}
\newcommand{\idemtwo}{e_{\downarrow}}
\DeclareMathOperator{\last}{last}

\newcommand{\BBB}{\boxed{B}}

\begin{document}

\begin{abstract} We initiate the study of Gaitsgory's central sheaves using complexes of Soergel bimodules, in extended affine type $A$. We conjecture that the complex associated to the standard representation can be flattened to a central complex in finite type $A$, which is the complex sought after in a conjecture of Gorsky-Negut-Rasmussen. \end{abstract}

\title{Gaitsgory's central sheaves via the diagrammatic Hecke category}

\author{Ben Elias} \address{University of Oregon, Eugene}

\maketitle

\tableofcontents

%%%%%%%%%%%%%%%%%%%%%%%%%
\section{Introduction}
\label{sec-intro}
%%%%%%%%%%%%%%%%%%%%%%%%%

%================
\subsection{The Hecke category and its center}
%================

The Hecke category $\HC_W$ is a categorification of the Iwahori-Hecke algebra $\Hecke_W$ of a Coxeter group $W$. It is a fundamental and ubiquitous object of study in categorical
representation theory, with incarnations in several different contexts. 

The Drinfeld center $\ZC(\AC)$ of a monoidal category $\AC$ is the categorical analog of the center $Z(A)$ of an algebra $A$. An object of $\ZC(\AC)$ is an object $M$ of $\AC$ equipped
with a natural isomorphism between the functors $M \ot (-)$ and $(-) \ot M$ of left and right ``multiplication'' by $M$. Note that $\ZC(\AC)$ itself is a braided monoidal category.

As the Hecke algebra often has an interesting center, one expects the Drinfeld center $\ZC(\HC_W)$ to be a category of intrinsic interest. Indeed, when $W$ is a finite Weyl group, work of
Bezrukavnikov-Finkelberg-Ostrik \cite{BFO} implies that $\ZC(\HC_W)$ is equivalent to Lusztig's category of character sheaves \cite{LuszCharactersheaves85}. There is a great deal of
recent progress towards explicitly constructing objects in $\ZC(\HC_W)$ in finite type, related to the work of \cite{EHDiag2, GNR}.

In this paper series, we focus instead on the case of affine Weyl groups, and in particular, affine type $A$. We have in mind a conjectural application to finite type $A$, see \S\ref{ssec_intro:GNR}.

%================
\subsection{Versions of the Hecke category} \label{ssec:versions}
%================

Unfortunately, it is difficult to discuss the contents of this paper without first addressing some technicalities.

When $W$ is a Weyl group, the Hecke category has classical constructions involving equivariant perverse sheaves on the flag variety, or translation functors acting on category $\OC$.
These were unified and generalized to Coxeter groups by Soergel, using the category of Soergel bimodules \cite{Soer07}. Most recently, the Hecke category was given a diagrammatic
presentation by the author and Williamson \cite{EWGr4sb}, based on the earlier works \cite{EKho}, joint with Khovanov, and \cite{ECathedral}. This diagrammatic description is the most low-tech and accessible of
the bunch, while also being more flexible and general: it agrees with the classical constructions when they are well-behaved (e.g. in characteristic zero), but it has the desired
properties even when classical constructions do not.

There are abelian, additive, and triangulated versions of the Hecke category. For example, perverse sheaves and category $\OC$ are both abelian, while Soergel bimodules are additive (they
represent projective objects in category $\OC$, or semisimple perverse sheaves). The diagrammatic category constructed by Elias-Williamson (associated with a Coxeter group $W$) will be
denoted $\Diag_W$, and it is also additive.

In this paper we are most interested in the triangulated setting. In this introduction, $\HC_W \define \Hot(\Diag_W)$ will denote the homotopy category of $\Diag_W$ (i.e. the category of
bounded complexes in $\Diag_W$, with morphisms considered modulo homotopy). This is a triangulated, graded, monoidal category, which comes equipped\footnote{This, and many other features of $\Diag_W$ we use in this paper, rely upon the Soergel conjecture, now a theorem due to \cite{EWHodge}. This requires that our category $\Diag_W$ be defined with coefficients over a field of characteristic zero.} with its \emph{perverse $t$-structure},
whose abelian core will also play a role. We will only be interested in the other versions of the Hecke category for motivational purposes, so we will not attempt to explain the precise
connections between these different categories.

\begin{remark} In the bulk of this paper we will be primarily interested in a category of pseudocomplexes over $\Diag_W$, which is another triangulated version of the Hecke category, but
we will throw this wrench at the reader at the appropriate time. \end{remark}

Each of these versions of the Hecke category has its own Drinfeld center. The Drinfeld center of the additive category $\Diag_W$ is uninteresting, because there are not enough central
objects. Both the abelian and the triangulated versions have interesting Drinfeld centers. Technically, the results of \cite{BFO} mentioned above (and \cite{Gait01} below) deal with the
Drinfeld center of the abelian version of the Hecke category. Meanwhile, this paper, and the papers \cite{EHDiag2,GNR} deal with the Drinfeld center of the triangulated version $\HC_W$.
Surprisingly and annoyingly, the literature seems to lack the foundational work needed to establish a link between the Drinfeld centers in the triangulated and abelian settings.

A second technicality is that the Hecke category $\HC_W$ depends not just on a Coxeter group $W$ but on a choice of \emph{realization} of $W$ \cite[Definition 3.1]{EWGr4sb}, which
is quite close to a choice of root datum. Thus there are several different versions of $\HC_W$ for an affine Weyl group $W$, some in which the simple roots are linearly dependent, and
some in which they are independent, etcetera. In finite and affine type $A$ there are versions associated to both $\sl_n$ and $\gl_n$. It will help to have the flexibility of discussing any of these situations. We also note that (certain) realizations of Weyl groups have associated Lie algebras, and come in Langlands dual pairs.

%================
\subsection{The center in affine type}
%================

Let us fix a finite Weyl group $W_{\fin}$, with associated affine Weyl group $W_{\aff}$. We fix a realization of $W_{\aff}$, giving by restriction a realization of $W_{\fin}$, with
associated categories $\HC_{\fin}$ and $\HC_{\aff}$. We let $\gg$ denote the complex semisimple Lie algebra associated to this realization of $W_{\fin}$. Let $\gg^\vee$ be the Langlands
dual Lie algebra, and $\Rep \gg^\vee$ denote its (semisimple monoidal) category of finite-dimensional representations. Let $\Rep_{\adj} \gg^\vee$ denote the subcategory of representations
of adjoint type, i.e. whose weights live in the root lattice.

The affine Hecke algebra $\Hecke_{\aff}$ is a deformation of the group algebra of $W_{\aff}$. In particular, it has a large commutative subalgebra of \emph{translations}, which is
naturally isomorphic to the root lattice $\La_{\rt}$ of the Langlands dual Lie algebra $\gg^\vee$. It turns out that the center $Z(\Hecke_{\aff})$ can be identified with
$W_{\fin}$-invariant sums of translations. A basis for $W_{\fin}$-invariant sums of translations can be given by the characters of irreducible representations in $\Rep_{\adj}
\gg^\vee$. Consequently, $Z(\Hecke_{\aff})$ is isomorphic to the Grothendieck group $[\Rep_{\adj} \gg^\vee]$.

A natural categorification of the isomorphism $[\Rep_{\adj} \gg^\vee] \cong Z(\Hecke_{\aff})$ would be the following.

\begin{goal} \label{goal0} There is a braided monoidal functor $\GC \co \Rep_{\adj} \gg^\vee \to \ZC(\HC_{\aff})$, which is fully faithful. \end{goal}

\begin{remark} The question of whether $\GC$ is essentially surjective or not is an important one, and to our knowledge the literature does not contain any conjectures in either direction. See Remark \ref{rmk:isthatthewholecenter} for further discussion. \end{remark}

\begin{remark} The category $\HC_{\aff}$ is both graded and triangulated, so it has a grading shift $(1)$ and a homological shift $[1]$. Between any two objects $M, N \in \HC_{\aff}$ one
could consider the usual Hom space $\Hom_{\HC}(M,N)$, or the bigraded vector space \begin{equation}\HOM_{\HC}(M,N) = \bigoplus_{k,\ell \in \Z} \Hom_{\HC}(M,N(k)[\ell]).\end{equation} This
construction is useful because the space $\HOM(M,N)$ has additional structures (e.g. it is a bimodule over a polynomial ring). Meanwhile, $\Rep_{\adj} \gg^\vee$ is not graded or
triangulated. So, when we say that the functor $\GC$ should be fully faithful, we mean that $\GC$ induces an isomorphism from $\Hom_{\Rep}(V_1,V_2)$ to $\Hom_{\HC}(\GC(V_1),\GC(V_2))$,
not to $\HOM_{\HC}(\GC(V_1),\GC(V_2))$. We can emphasize this by saying that $\GC$ is \emph{fully faithful to degree zero}. In particular, $\ZC(\HC_{\aff})$ is a much richer setting than
$\Rep_{\adj} \gg^\vee$ because of the existence of morphisms of nonzero degree. \end{remark}

Let us upgrade this goal slightly. There is also an extended affine Hecke algebra $\Hecke_{\ext}$, whose translation lattice is isomorphic to the weight lattice $\La_{\wt}$ of $\gg^\vee$.
Then $Z(\Hecke_{\ext})$ is isomorphic to the $W_{\fin}$-invariant sums of translations, so $Z(\Hecke_{\ext}) \cong [\Rep \gg^\vee]$. This extended affine Hecke algebra can be categorified
using \emph{extended Soergel bimodules} or their diagrammatic analog $\DC_{\ext}$, an enhancement of the Hecke category which adds a formal object corresponding to the standard
automorphism of the affine Dynkin diagram. We let $\HC_{\ext}$ denote the \emph{extended Hecke category}, the bounded homotopy category of $\DC_{\ext}$.

\begin{remark} This extension $\DC_{\ext}$ of Soergel bimodules is due to Mackaay-Thiel in \cite{MacThi} and its sequel, where it appears in both algebraic and diagrammatic form. However, there are aspects of the diagrammatic extension which require careful development. We give a more thorough treatment of this category in \S\ref{sec:extendeddiag}. \end{remark}

\begin{goal}\label{goal1} There is a braided monoidal functor $\GC \co \Rep \gg^\vee \to \ZC(\HC_{\ext})$, which is fully faithful. Upon restriction to $\Rep_{\adj} \gg^\vee \subset \Rep \gg^\vee$, one recovers the functor from Goal \ref{goal0}. \end{goal}

%================
\subsection{Gaitsgory's work}
%================

Goal \ref{goal1} has already been achieved by Gaitsgory \cite{Gait01}, using highly technical methods. Before discussing his work, we must shift our focus to geometry, though we will
only give a rough sketch of the goings-on.

First, the geometric Satake equivalence states that $\Rep \gg^\vee$ is equivalent to a category of equivariant perverse sheaves on the affine Grassmannian $\Gr$, with monoidal structure
given by convolution. Meanwhile, by work of Soergel and Harterich \cite{Soer90, Harterich}, one can relate (the abelian core of the perverse $t$-structure on) $\HC_{\aff}$ to a
category of equivariant perverse sheaves on an affine flag variety $\Fl$. Both $\Gr$ and $\Fl$ are infinite-dimensional ind-schemes, making their study rather technical. Using a
deformation of $\Fl$ and the nearby cycles functor, Gaitsgory constructs a monoidal functor\footnote{Gaitsgory denotes this functor $\ZC$ in \cite{Gait01}. We will use $\ZC$ for the Drinfeld center, and we name the functor $\GC$ in Gaitsgory's honor.} from perverse sheaves on $\Gr$ to perverse sheaves
on $\Fl$, whose image lies within the Drinfeld center. The perverse sheaves in the image of this functor will be referred to as \emph{Gaitsgory's central sheaves}. 

The construction used by Gaitsgory implies that these central sheaves will have a number of highly desireable properties, which we describe below in \S\ref{ssec:lovely}. The introduction
to Gaitsgory's paper is quite nice and contains additional motivation and background for this construction.

The ideas above can be encapsulated in this optimistic\footnote{This diagram mixes abelian categories and triangulated ones. In this paper we will not review how Soergel relates the
abelian category of perverse sheaves to the homotopy category of Soergel bimodules.} commutative diagram. Here, $\PC$ denotes the category of (appropriately equivariant) perverse sheaves.

\begin{equation} \label{eq:goaldiagram}
\begin{diagram}
\PC(\Gr) & \rTo^{\textrm{Satake}}_{\sim} & \Rep \gg^\vee \\
\dTo^{\textrm{Gaitsgory}} & & \dTo_{\GC} \\
\ZC(\PC(\Fl)) & \rTo^{\textrm{Soergel}}_{\sim} & \ZC(\HC_{\ext}) \end{diagram}
\end{equation}

\begin{remark} This is a diagram of braided monoidal functors, but the braiding on $\Rep \gg^\vee$ is not the standard braiding which one learns in representation theory class, where $v
\otimes w \mapsto w \otimes v$. Instead, for weight vectors $v$ and $w$, $v \otimes w \mapsto \pm w \otimes v$, for some natural choice of signs. This same braiding appears in the
geometric Satake equivalence, so we call it the \emph{Satake braiding}. We return to this in \S\ref{ssec:intro_braiding}. \end{remark}

Soergel showed that the link between geometry and Soergel bimodules (and diagrammatics) is given by the equivariant global sections functor. He gave a dictionary between geometric ideas
involving perverse sheaves, and algebraic ones involving their global sections. However, the techniques used by Gaitsgory (e.g. deformations, nearby cycles) have not (yet) been included
into Soergel's dictionary, and we do not know how to take Gaitsgory's construction and adapt it to Soergel bimodules or the diagrammatic Hecke category. In other words, the commutative
diagram \eqref{eq:goaldiagram} is still mysterious and inexplicit. While one could deduce the existence of the functor $\GC$ from Gaitsgory's work (and the horizontal equivalences), one
would not be able to directly compute with it (i.e. apply it to any representation) in the context of Soergel bimodules, using the currently known technology.

\begin{remark} We have recently learned that Achar and Rider are currently working on this very problem, with promising results. We will discuss this further below. \end{remark}

Instead, we shall attempt to reproduce Gaitsgory's functor $\GC$ from scratch in the Soergel setting, looking only at the right side of \eqref{eq:goaldiagram}, and reproducing Gaitsgory's
results by hand. Again, we do not have the technology to prove (in a satisfactory way) that any functor we construct will make \eqref{eq:goaldiagram} commute, but we hope to convince the
reader that our construction is correct nonetheless.

%================
\subsection{Our goal}
%================

This paper initiates a program to re-achieve Goals \ref{goal1} in type $A$, using purely low-tech, algebraic and diagrammatic methods. That is, we aim to make Gaitsgory's beautiful
results accessible to the common algebraist, and to give a construction of $\GC$ for the reader who\footnote{like the author} knows nothing about geometry, ind-schemes, perverse sheaves,
nearby cycles, etcetera. Of course, the reader will have to learn about Soergel bimodules or the diagrammatic Hecke category instead, but this is easier (in the author's evidently unbiased
opinion).

Because we work with complexes in a homotopy category, rather than sheaves on a geometric space, let us now excise the word ``sheaf'' from our notation. We refer to objects in the image
of $\GC$ as \emph{Gaitsgory's central complexes}.

There are advantages to the algebraic approach, making this program more significant than a mere re-proof. The diagrammatic category $\DC_{\ext}$ is a setting where one can compute
explicitly (whether by hand or with computer), allowing one direct access to the structure of these complexes, and opening up a number of combinatorial questions. We should emphasize that
the functor $\GC$ encapsulates very little of the structure of Gaitsgory's central complexes. It is fully faithful, but only to degree zero, missing out on all the morphisms of higher
degree. We believe our algebraic approach will make some of these complexes easier to study, and we include some computations in Appendix \ref{sec:computations} as a proof of concept (and
for general use). We also hope that by making Gaitsgory's functor very explicit in type $A$ and thus giving a host of examples, it will be easier to reproduce Gaitsgory's techniques directly in the algebraic context, Soergel-ifying the nearby cycles functor and related technology.

The primary accomplishment of this paper is to define the central complex $\VC = \GC(V)$ associated with the standard representation in type $A$, and to explore its properties. We also lay
the groundwork for the remainder of the program. Before discussing our approach and our results, we would like to give another motivation for why the functor $\GC$ is important, which
points to the great importance of the complex $\VC$ in particular.

%================
\subsection{Connection to knot theory and finite type $A$} \label{ssec_intro:GNR}
%================

In this section we let $W_{\fin}$ be the symmetric group $S_n$, and $\gg \cong \gg^\vee \cong \gl_n$. We adopt the standard realization of the extended affine Weyl group, see Definition
\ref{defn:standardaffine}, and use this to define $\Diag_{\ext}$ and $\HC_{\ext}$. In particular, $\HB_{\ext}$ has a central element $\det$ which is translation by the determinant
weight, and is categorified by a central object in $\DC_{\ext}$.

% \footnote{We will discuss the differences between $\sl_n$ and $\gl_n$ at length in the body of this paper, see
% \S\ref{sec:extendedbraid} and \S\ref{sec:extendeddiag}.}

Recall that the Hecke algebra $\HB_{\ext}$ has a large commutative subalgebra of translations $\La_{\wt}$, which contains the center $Z(\HB_{\ext})$. Similarly, the Hecke algebra
$\HB_{\fin}$ of the finite Weyl group has a large commutative subalgebra $J$ generated by the \emph{multiplicative Jucys-Murphy elements} or \emph{JM elements} $j_i$ for $1 \le i \le n$,
which contains the center $Z(\HB_{\fin})$. Just as $Z(\HB_{\ext})$ can be thought of as $S_n$-invariant sums of weights, $Z(\HB_{\fin})$ is spanned by symmetric polynomials in the
JM elements. We write $e_k(\overline{j})$ for the $k$-th elementary symmetric polynomial in the family $\{j_1, \ldots, j_n\}$.

\begin{remark} We refer to the commutative subalgebra $J$ as \emph{large} because it is maximal, and because it is simultaneously diagonalizable on any finite-dimensional $\HB_{\fin}$ representation, which makes it critical for representation theory, see e.g. \cite{OkoVer} or \cite{EHDiag2}. This is in analogy with the Cartan subalgebra of a semisimple Lie algebra. \end{remark}

One can think of the extended braid group $\Br_{\ext}$ as consisting of braids drawn on a cylinder. By flattening the cylinder into a planar rectangle (in the same way one flattens a
cardboard tube), one gets an ordinary braid in the finite braid group $\Br_{\fin}$. This \emph{flattening} map $\flat \co \Br_{\ext} \to \Br_{\fin}$ descends to a map $\flat \co
\HB_{\ext} \to \HB_{\fin}$, and induces the following\footnote{Technically, $\La_{\wt}$ should be replaced by the more cumbersome $\Z[v,v\inv] [\La_{\wt}]$ in this commutative diagram.} commutative diagram. \begin{equation} \label{eq:flatdiagram} \begin{diagram} Z(\HB_{\ext}) & \rInto & \La_{\wt} & \rInto & \HB_{\ext}
\\ \dTo_{\flat} & & \dTo_{\flat} & & \dTo_{\flat} \\ Z(\HB_{\fin}) & \rInto & J & \rInto & \HB_{\fin} \end{diagram} \end{equation} It is worth noting that every translation in $\La_{\wt}$ is the image of a braid, rather than a linear combination of braids, just as the JM elements are images of braids.

The translations associated with the $n$ weights inside the standard representation $V$ are flattened to the $n$ JM elements. Thus the character of $V$, living in $Z(\HB_{\ext})$, is
flattened to $e_1(\overline{j})$. The other fundamental representations $\Lambda^k V$ are sent to other elementary symmetric polynomials $e_k(\overline{j})$ in the JM elements. For
example, the determinant weight $\det$ is sent to the product of the JM elements $e_n(\overline{j}) = ft_n$, which is also known as the \emph{full twist}.

Now let us discuss the categorification of this story. Flattening does not take the Kazhdan-Lusztig basis to the Kazhdan-Lusztig basis, and one does not expect it to lift to a functor
$\DC_{\ext} \to \DC_{\fin}$. However, one might hope that it lifts to a functor on homotopy categories $\flat \co \HC_{\ext} \to \HC_{\fin}$. In the paper \cite{EFlat} we construct a functor $\DC_{\ext} \to \HC_{\fin}$ and conjecture that it extends to a functor $\HC_{\ext} \to \HC_{\fin}$. We assume this conjecture for the rest of this discussion.

Associated to any element of the braid group (finite or affine or extended) there is a \emph{Rouquier complex}, which is a (tensor-)invertible complex in $\HC$ categorifying the image of
that braid in the Hecke algebra. Inside $\HC_{\fin}$, the Jucys-Murphy elements lift to certain Rouquier complexes which we call \emph{Jucys-Murphy complexes}. Let $\JC$ denote the full
monoidal triangulated subcategory of $\HC_{\fin}$ generated by the Jucys-Murphy complexes. Meanwhile, inside $\HC_{\ext}$, there are Rouquier complexes lifting the elements of the
translation lattice, which are known as \emph{Wakimoto complexes} (they correspond to Wakimoto sheaves in geometry). Let $\Wak$ denote the full monoidal triangulated subcategory of
$\HC_{\ext}$ generated by the Wakimoto complexes. Just as JM elements $j_i$ (resp. translations) commute with each other, one has a strong categorical commutation property between
Jucys-Murphy (resp. Wakimoto) complexes, making them a ``commutative monoidal subcategory'' of $\HC$. If the conjecture in \cite{EFlat} is true, then the flattening functor is compatible with
Rouquier complexes, so one has the following commutative diagram. \begin{equation} \label{eq:flatdiagramcatfd} \begin{diagram}
\ZC(\HC_{\ext}) & \rTo^{?} & \Wak & \rInto & \HC_{\ext} \\ \dTo_{\flat} & & \dTo_{\flat} & & \dTo_{\flat} \\ \ZC(\HC_{\fin}) & \rTo^{?} & \JC & \rInto & \HC_{\fin} \end{diagram}
\end{equation}

One of the key properties of Gaitsgory's central complexes is that they have a filtration by Wakimoto complexes, corresponding to the decomposition of the corresponding representation into weight spaces. For example, let $\VC = \GC(V)$ denote the Gaitsgory central complex for the standard representation. Then $\flat(\VC)$ has a filtration by Jucys-Murphy complexes, each appearing once, which categorifies the analogous statement about the flattening of the character of $V$.

\begin{remark} \label{rmk:isthatthewholecenter} Although all Gaitsgory complexes live in $\Wak$, and all known central complexes in $\HC_{\fin}$ live in $\JC$, it is not known whether
this is true of arbitrary central objects. Thus the first, ``questionable'' horizontal arrow in each row is merely motivational. Conversely, after a conversation with Kostiantyn
Tolmachov, I now believe that every object of $\ZC(\HC_{\ext})$ which is Wakimoto filtered should be isomorphic to a Gaitsgory central complex.

Moreover, an object of $\ZC(\HC)$ is not just a complex, but comes equipped with some additional structure. Thus the questionable horizontal arrow is also forgetful. \end{remark}

Now let us explain the significance of all this to knot theory. Most of the knot polynomials in the literature, including the HOMFLYPT polynomial which specializes to many others, factor
through the quotient from the braid group $\Br_{\fin}$ to the Hecke algebra $\HB_{\fin}$. It is now understood, through work of Khovanov \cite{Khov07}, Rasmussen \cite{Rasmussen} and
others that the knot homology theories which categorify these knot polynomials can be constructed using Rouquier complexes in $\HC_{\fin}$, the categorical analog of factoring through
$\HB_{\fin}$. More recently, it was observed in the work of Gorsky-Oblomkov-Rasmussen-Shende \cite{ORS, GORS} that HOMFLYPT homology seemed to be computable using the properties of
coherent sheaves on a Hilbert scheme. This collection of ideas was placed in an overarching conjectural framework in a beautiful paper by Gorsky, Negut, and Rasmussen \cite{GNR}, which
also contains additional background on this topic.

Let us note that $\HB_{\fin}$ is equipped with a nondegenerate sesquilinear form, which restricts to a nondegenerate form on $J$. Consequently, the inclusion map $\iota \co J \into
\HB_{\fin}$ admits an adjoint map $\iota^* \co \HB_{\fin} \to J$. Now the HOMFLYPT polynomial of a braid closure can be computed by taking the braid, considering its image in
$\HB_{\fin}$, and pairing it against the identity element using the nondegenerate form. Since the identity element is in $J$, the map $(1,-)$ from $\HB_{\fin}$ actually factors through
the adjoint map $\iota^*$. Analgously, the HOMFLYPT homology of a braid closure can be computed by taking the braid, considering its Rouquier complex, and applying Hochschild
cohomology. In \cite{GNR}, they argue that HOMFLYPT homology should factor through some categorification of $\iota^*$.

We let $\Hilb_n$ and $\FHilb_n$ denote certain spaces which are closely related to the Hilbert scheme and the flag Hilbert (dg-)scheme of $n$ points in the plane, see \cite{GNR} for
details. Gorsky, Negut, and Rasmussen conjecture that $\JC$ is equivalent to $D^b(\FHilb_n)$, the monoidal derived category of coherent sheaves on the $\FHilb_n$. The subcategory
$\ZC(\HC_{\fin})$ would be equivalent to $D^b(\Hilb_n)$, the derived category of coherent sheaves on $\Hilb_n$ itself, viewed as a subcategory under pullback. The adjoint map $\iota^* \co
\HB_{\fin} \to J$ is categorified by a functor from $\HC_{\fin} \to D^b(\FHilb_n)$, obtained in some sense by taking the Hom space from all tensor powers of the full twist (this is
motivated by the theory of categorical diagonalization developed in \cite{EHDiag}, as elaborated in \cite{GNR}). To obtain the Khovanov HOMFLYPT homology of a braid, one can take the
corresponding Rouquier complex in $\HC_{\fin}$, apply $\iota^*$ to obtain a complex of coherent sheaves, and then apply some homological invariants.

Meanwhile, the map $\iota \co J \to \HB_{\fin}$ should be categorified by some functor $\iota \co D^b(\FHilb_n) \to \JC$ which sends the tautological line bundles on the flag Hilbert
scheme to the Jucys-Murphy complexes. The sheaf $\OC(1)$ is the determinant bundle of the tautological $n$-plane bundle $\TC$, which has a filtration by these line bundles. Hence $\OC(1)$
is sent to the product of the Jucys-Murphy complexes, i.e. the full twist complex $\FT_n$. Meanwhile, the image of $\TC$ under $\iota$ should be a central complex with a filtration by the
Jucys-Murphy complexes each appearing once, and a host of other properties. In \cite{GNR}, they manage to reduce all their conjectures to one conjecture, the existence of some complex
$\EC$ with the desired properties to be the image of $\TC$.

So $\EC$ should have a filtration corresponding to the first elementary symmetric polynomial in the JM elements. As observed above, so does the flattening of $\VC$, leading us to the following conjecture. It is easy to verify for $n=2$.

\begin{conj} The flattening of $\VC$, the Gaitsgory central complex for the standard representation (constructed in this paper), is a complex $\EC$ which satisfies the requirements of
Gorsky-Negut-Rasmussen's conjecture. \end{conj}

The construction of $\VC$ in this paper is completely explicit, which we hope will allow one to verify the desired properties of $\flat(\VC)$. Here is the conjectural commutative diagram, and how it should behave in the left column.

\begin{equation} \label{eq:flatdiagramcatfd2} \begin{diagram} \ZC(\HC_{\ext}) & \rTo & \Wak & \rInto & \HC_{\ext} & \qquad & \VC
\\ \dTo_{\flat} & & \dTo_{\flat} & & \dTo_{\flat} & \qquad & \dMapsto \\ \ZC(\HC_{\fin}) & \rTo & \JC & \rInto & \HC_{\fin} & \qquad & \EC \\ \dTo^{\iota^*} \uTo_{\iota} & & \dTo^{\iota^*} \uTo_{\iota} & \ldTo^{\iota^*} & & \qquad & \uMapsto \\ D^b(\Hilb_n) & \rInto & D^b(\FHilb_n) & & & \qquad & \TC \end{diagram} \end{equation} 

\begin{remark} Recent remarkable work of Oblomkov-Rozansky \cite{OblRoz} places some of these ideas in a broader algebro-geometric framework, and also explains why the flattening map from
extended affine to finite should have image related to the Hilbert scheme. The connection between their framework and the Hecke category a la Soergel is still conjectural. \end{remark}

%================
\subsection{The approach}
%================

Having motivated things, we return to the construction of the functor $\GC$. Let us outline our approach, and the reason we restrict to type $A$.

Let $\gg^\vee$ be a complex semisimple Lie algebra. The category $\Rep \gg^\vee$ is semisimple, satisfying Schur's lemma: for irreducible objects $V_\la$ and $V_{\mu}$ one has \begin{equation} \Hom(V_{\la},V_{\mu}) = \begin{cases} \C \cdot
\id & \text{if } \la = \mu,\\ 0 & \text{else}. \end{cases} \end{equation} Consequently, for any three objects $X, Y, Z$ in $\Rep \gg^\vee$, if we know their decomposition into
irreducibles, we can describe the space $\Hom(X,Y)$ explicitly using matrices, and the composition map $\Hom(Y,Z) \ot \Hom(X,Y) \to \Hom(X,Z)$ using matrix multiplication. This makes semisimple categories ``easy.''

However, $\Rep \gg^\vee$ is also a monoidal category, so we also want to understand the monoidal composition map $\Hom(X,Y) \ot \Hom(X',Y') \to \Hom(X \ot X',Y \ot Y')$. This is very
poorly understood combinatorially. For a longer exposition of this particular idea, see the introduction to \cite{ELLCC}.

Let $\Fund \gg^\vee$ denote the full monoidal subcategory whose objects are tensor products of fundamental representations. Most successful approaches to understanding $\Rep \gg^\vee$ as
a monoidal category have focused on giving a diagrammatic monoidal presentation of $\Fund \gg^\vee$ by generators and relations. This has been accomplished in rank 2 by Kuperberg
\cite{Kupe}, and in type $A$ by Cautis-Kamnitzer-Morrison \cite{CKM}. Almost nothing is known in other types.

Let $V$ denote the standard representation of $\gl_n$, and for each $1 \le k \le n$ let $V_{\fund_k} = \Lambda^k V$ denote the $k$-th exterior power of $V$. Of course, $V_{\fund_n}$ is
the invertible determinant representation. These are the \emph{fundamental representations} of $\gl_n$. Cautis-Kamnitzer-Morrison \cite{CKM} define a category $\Webs_n = \Webs(\gl_n)$
together with a fully faithful functor $\Webs_n \to \Fund \gl_n$. The objects of $\Webs_n$ are sequences of fundamental weights with orientations, which are sent to the tensor product of
the corresponding sequence of fundamental representations and their duals. Morphisms in $\Webs_n$ are encoded by planar diagrams called \emph{webs}, for which Cautis-Kamnitzer-Morrison
give a presentation. Thus to define a functor from $\Webs_n$ is straightforward: one specifies the images of the fundamental representations and their duals (the generating objects), and
the images of the generating webs (the generating morphisms), and then checks the web relations (the relations between morphisms).

Technically, \cite{CKM} contains a description of a category $\Webs(\sl_n)$, not $\Webs(\gl_n)$. It is well-known how to modify\footnote{Explcitly, allow $n$ as an object label. Omit the
``tag'' morphism and all relations which involve it.} their description to work for $\gl_n$ instead of $\sl_n$. For simplicity, we will work with $\Webs^+_n$, the positively-oriented
subcategory which corresponds to tensor products of fundamental representations without their duals.

Our approach is as follows.

\begin{goal} \label{goal1.5} Explicitly construct a braided monoidal functor $\GC$ from $\Webs^+_n$ to $\HC_{\ext}$, by \begin{enumerate} \item[(a)] constructing a complex
$\VC_{\fund_k} = \GC(V_{\fund_k})$ for each $1 \le k \le n-1$, \item[(b)] constructing a chain map for each generating web, and \item[(c)] checking the relations of the $\gl_n$-web
calculus. \end{enumerate} Then \begin{enumerate} \item[(d)] verify that the objects $\GC(V_{\fund_k})$ are central (or rather, equip them with a central structure), \item[(d')] verify that the morphisms given by generating webs are central, \item[(e)] check
that the braidings are intertwined by $\GC$, and \item[(f)] prove that $\GC$ is fully faithful to degree zero. \end{enumerate} \end{goal}

We will discuss simplifications to this process in \S\ref{sec-elaboratesimplify}.

\begin{remark} One can recover $\Webs_n$ (or an equivalent category) from $\Webs^+_n$ by inverting the determinant representation. The image under $\GC$ is already invertible in
$\HC_{\ext}$, so the functor can be extended to $\Webs_n$. In addition, various dualities on $\Rep \gl_n$ and $\Webs_n$ are intertwined with various dualities in $\HC_{\ext}$, see
\S\ref{subsec-symmetries}. This gives another way to define the functor on $\Webs_n$. Thus we continue to speak of the functor $\GC \co \Rep \gl_n \to \ZC(\HC_{\ext})$. Finally, the
extended diagrammatics for $\sl_n$ developed in \S\ref{subsec-extendedsbimquotient} should allow one to define a functor from $\Webs(\sl_n)$ as well. All these extensions are relatively
straightforward, so we do not bother to discuss them further in this paper. \end{remark}

\begin{remark} One can recover the category $\Rep \gl_n$ from $\Fund \gl_n$ or $\Webs_n$ using the Karoubi envelope, which formally adds all direct summands. To make this explicit, one
must compute the \emph{clasps}, the idempotents in $\Webs_n$ which project from a tensor product to its highest weight irreducible summand. These clasps are only known for $\gl_n$, $n \le
4$, though recent work of the author \cite{ELLCC} has provided a conjecture for all $n$. \end{remark}

\begin{remark} The flattenings of the other fundamental Gaitsgory complexes $\VC_{\fund_k}$ correspond, in the Gorsky-Negut-Rasmussen setup, to the exterior products $\Lambda^k \TC$.
These are also significant, as they (conjecturally) control higher Hochschild cohomology groups. \end{remark}

\begin{remark} There is another potential approach which only requires the construction of $\VC = \VC_{\fund_1}$, replacing the use of $\gl_n$-webs with an action of the symmetric group
on tensor powers of $\VC$ (as one expects from Schur-Weyl duality). We have not attempted this approach. \end{remark}

%================
\subsection{Properties of Gaitsgory's central complexes} \label{ssec:lovely}
%================

We now list some additional properties that one wishes Gaitsgory's central complexes to satisfy, and add them to our to-do list. Most of these were proven by Gaitsgory in the geometric setting in \cite{Gait01}. We elaborate upon these properties one by one in \S\ref{sec-elaboratesimplify}.

\begin{goal}\label{goal2} For any representation $U \in \Rep \gl_n$, the complex $\UC = \GC(U)$ in $\HC_{\ext}$ has the following properties. \begin{enumerate} \item[(g)] $\UC$ is
perverse, that is, it lies in the core of the perverse $t$-structure. \item[(h)] Stronger still, $\UC$ \emph{preserves perverses}. That is, tensoring with $\UC$ takes perverse complexes
to perverse complexes. \item[(i)] $\UC$ has a \emph{Wakimoto filtration}, a filtration whose subquotients are Wakimoto complexes, lifting its weight space filtration. \item[(j)] The maximal cell quotient of $\UC$ agrees with
the underlying graded vector space of $U$. \item[(k)] $\UC$ is equipped with a nilpotent \emph{monodromy operator} $\mu \co \UC \to \UC(-2)[2]$. \item[(l)] The functor $\GC$ is compatible
with the geometric Satake equivalence. \end{enumerate} \end{goal}

The following properties do not appear in Gaitsgory, as they deal with the triangulated structure, so we list them as conjectures instead of goals. They are supported by computation in
the cases of $n = 2, 3$.

\begin{conj}\label{conj3} If $U, U' \in \Rep \gl_n$, and $\UC = \GC(U)$, $\UC' = \GC(U')$, then: \begin{enumerate}
\item[(m)] $\Hom(\UC,\UC'(k)[c]) = 0$ for any $c < 0$ and any $k \in \Z$. This is the \emph{vanishing of negative exts}.
\item[(n)] $\Hom(\UC,\UC'(k)[c])$ is either zero for all even $c$ or for all odd $c$. (More trivially, it is also either zero for all even $k$ or for all odd $k$.) This is the \emph{(homological) parity} of Hom spaces. \end{enumerate} \end{conj}

%================
\subsection{What happens in this paper}
%================

What we accomplish is but a humble first step towards these goals, but still an illuminating one.

We explicitly provide a complex $\VC$ which should be the Gaitsgory central sheaf for the standard representation, for all $n \ge 2$. We also provide the candidate for the ``dual''
complex\footnote{These representations are only dual for $\sl_n$; over $\gl_n$ the dual of $V$ is $V_{\fund_{n-1}} \ot \det\inv$.} $\VC_{\fund_{n-1}}$. See Definitions \ref{defn:FC} and
\ref{defn:VC}.

The major technical theorem in this paper, Theorem \ref{thm:tensorBI}, is an explicit computation of (the minimal complex of) $\VC \ot B_J$ where $B_J$ is the indecomposable object of
$\Diag_{\ext}$ corresponding to the longest element of a finite parabolic subgroup $W_J$. Effectively, what we compute is the pushforward of $\VC$ from the affine flag variety to a
partial affine flag variety. In particular, $\VC \ot B_J$ is perverse. One can use this to prove that $\VC \ot B_w$ is perverse for many indecomposable objects $B_w$, but not all; the
elusive and seemingly difficult property (h) remains out of reach.

As a corollary, we can explicitly construct the isomorphism $\VC \ot B_J \to B_J \ot \VC$. This allows us, with additional computations, to construct an isomorphism $\VC \ot (-) \to (-)
\ot \VC$ as functors $\Diag_{\ext} \to \HC_{\ext}$, see Theorem \ref{thm:rotationfunctorial}. Note that this is not quite as good as proving that $\VC$ is in the Drinfeld center of
$\HC_{\ext}$, but it reduces this problem to some questions about ext vanishing, see \S\ref{ssec:intro_centrality} for more details.

When $n = 2$ or $3$, we have constructed all the fundamental representations, Goal \ref{goal1.5}(a). Thus we are able to continue the approach and construct morphisms associated with
$\gl_n$-webs. We are able to complete Goal \ref{goal1.5} for $n=2$, and most of it for $n=3$. Moreover, for $n=2$ we provide explicit constructions of $\GC(V_\la)$ for every irreducible
representation $V_\la$. In \S\ref{sec:computations} we record for posterity various structures when $n=2, 3$, such as the morphisms corresponding to the generating webs, or the braiding. The intense
computations which derived or verified these results are too long to write up, and we expect them to be superceded by more interesting general results eventually; the interested reader
is welcome to contact the author for details.

Much still remains to be done. In \S\ref{sec-elaboratesimplify} we build a framework for accomplishing Goals \ref{goal1.5} and \ref{goal2} efficiently for general $n$. This framework has
many missing steps, some large and some small. One of the purposes of this paper is to get the program out there, and make clear what still needs to be done, so that others can finish
the program.

\begin{remark} As early as 2011, Roman Bezrukavnikov suggested to me that I should look for the Soergel bimodule analog of Gaitsgory's central sheaves. In 2014 the suggestion was
repeated, with help in the interpretation from Geordie Williamson. The construction of $\VC$, the computation of $\VC \ot B_s$ for each simple reflection, and most of the proof of
centrality was completed in the summer of 2014, but in the hope that Goal \ref{goal1.5} could be completed in its entirety, I did not write it up then. Now, many years later, I have
decided to publish these results, partial though they are, to get the ball rolling, and also because of the potential applications discussed in \S\ref{ssec_intro:GNR}.  \end{remark}

\begin{remark} We have recently learned about work in progress of Achar and Rider. This work will reproduce Gaitsgory's technology (e.g. the nearby cycles functor) in the context of
parity sheaves, which is closely related to the category $\HC_{\ext}$. They too are able to explicitly produce the complex $\VC$, using a very elegant method, and their computation agrees
with ours. Moreover, they are able to work in any type, not just type $A$. It seems that they are unable to prove centrality, and our remaining results are largely complimentary.
\end{remark}

%================
\subsection{Pseudocomplexes and monodromy}
%================

Let us not conclude this introduction without correcting one major lie in the story we have told. Gaitsgory central complexes need not actually be complexes, though they are something
even better: pseudocomplexes!

As remarked in \S\ref{ssec:versions}, the Hecke category $\HC_{\ext}$ depends on a choice of realization, which is a choice of dual representations $\hg$ and $\hg^*$ of $W = W_{\ext}$, and
simple roots and coroots compatible with this action. Every morphism space in $\HC_{\ext}$ will be a bimodule for the polynomial ring $R$ whose linear terms are $\hg^*$, and any
$W$-invariant polynomial will act the same on the left and right. By convention, the degrees in $R$ are ``doubled,'' so that multiplication by a simple root has degree $2$. Let $\delta$
denote the sum of the simple roots for $W$. Depending on the choice of realization, $\delta$ could be zero, or it could be a nonzero, $W$-invariant polynomial.

When we construct $\VC$ in \S\ref{sec-construction}, we produce something which looks like a complex. However, we compute that $d^2$ is not necessarily zero; instead, $d^2 = \delta \mu$
for some chain map $\mu \co \VC \to \VC(-2)[2]$. In particular, $d^2 = 0$ if and only if $\delta=0$. This mysterious chain map $\mu$ is called the \emph{monodromy operator}, and it is a
wonderful thing.

Gaitsgory's geometric construction uses the nearby cycles functor, which automatically equips its target with a (log of) monodromy operator. In a similar context, Bezrukavnikov and Yun
\cite{BezYunMonodromy} explored a family of monodromy operators to understand the Koszul dual of the Hecke category in geometry. Makisumi \cite{MackyMonodromy} developed the Soergel-ified
analog of this idea. He considers a category of \emph{pseudocomplexes}, where $d^2$ is only equal to zero modulo positive degree polynomials, and constructs a family of monodromy
operators which measures the first order failure of $d^2$ to be zero. His category of pseudocomplexes is not monoidal in general. In our setting, we assert that $d^2 = 0$ modulo the ideal
generated by $\delta$. There is now effectively just one monodromy operator $\mu$, which measures the failure of $d^2$ to be zero in the direction of $\delta$. Because $\delta$ is
$W$-invariant, pseudocomplexes are a monoidal category.

\begin{remark} \label{rmk:monodromyconfusion} There is room for a great deal of confusion here. Following Gaitsgory's work, one might expect the monodromy to be a map $\VC \to \VC(2)$,
instead of a map $\VC \to \VC(-2)[2]$. However, the passage from geometry to Soergel bimodules often involves Kozsul duality, which swaps these grading shifts.

The conjectures of Gorsky-Negut-Rasmussen imply that their tautological complex $\EC$ has two natural endomorphisms $\mu$ and $\chi$ (called $x$ and $y$ in \cite{GNR}), where $\mu$ has
degree $(-2)[2]$ and $\chi$ has degree $[0](2)$. The map $\chi$ might also be called a monodromy operator by some. We believe that $\chi$ also lifts to an endomorphism of $\VC$. We
demonstrate $\chi$ for $n=2,3$ in \S\ref{subsec:chi}, where further discussion can be found. For the bulk of this paper, we use the term \emph{monodromy} to refer to the operator $\mu$.
\end{remark}

Thus there are two different contexts in which to study $\VC$: the context where $\delta = 0$ and $\VC$ is a genuine complex, or the context where $\delta \ne 0$ and $\VC$ is a
pseudocomplex. Note that the monodromy operator $\mu$ is a well-defined chain map in both contexts. We can perform all our homological computations within Makisumi's category of
pseudocomplexes, and will obtain equivalent results for genuine complexes in any realization where $\delta=0$. However, knowing that pseudocomplexes underlie the construction is
significant, because it gives an intrinsic construction of the monodromy map, and proves a number of facts about it for free. For example, the monodromy $\mu$ will commute with any chain map
(up to homotopy), and on tensor products, will satisfy $\mu_{F \ot G} = \mu_F \ot 1_G + 1_F \ot \mu_G$. 

\begin{remark} The essence of Makisumi's monodromy appears in a different framework in the monolithic work \cite{AMRW}, which constructs the ``free monodromic category'' which is
Koszul-dual to $\HC_W$. This more complicated homological regime subsumes the original concept of pseudocomplexes, and we could very well have worked in their setting instead. We felt
that it was best to introduce the minimal overhead, so we will stick with pseudocomplexes. \end{remark}

%================
\subsection{Structure of the paper} \label{subsec:structure}
%================

We do not recall the diagrammatic definition of $\Diag$ in this paper, and we expect the reader to be familiar with it. See \cite{EWGr4sb} for this background material.

In \cite{EFlat}, a work in preparation, we will develop the diagrammatics for the extended affine setting, by taking the category $\Diag_{\aff}$ associated to the affine Weyl group and
formally adding an invertible $1$-morphism $\Om$ which acts like the Dynkin diagram automorphism. We also give additional exposition on several topics, which will be essential for this
paper. This material has all been written already, and we temporarily include it here as \S\ref{sec:extendedbraid},
\S\ref{sec:extendeddiag}, \S\ref{sec-wakimoto}, and \S\ref{sec:extendeddiagproofs}. Eventually, \cite{EFlat} will take on an independent existence, and we will move these chapters there.

In \S\ref{sec:extendedbraid} we discuss the extended braid group and its weight lattice, and the difference between $\sl_n$ and $\gl_n$. In \S\ref{sec:extendeddiag} we discuss
realizations in the extended setting, and the extneded diagrammatics. In \S\ref{sec:extendeddiagproofs} we prove that these diagrammatics are correct. In \S\ref{sec-wakimoto} we explore
the properties of Rouquier complexes in the original and extended settings, including the perversity of Wakimoto complexes.

\begin{remark} The goal of \cite{EFlat} is to define the flattening functor, which is only important in this paper as a motivation. \end{remark}

In this paper we study the Gaitsgory central complex $\VC$ in $\HC_{\ext}$. Note that $\VC \cong \Om \ot \FC$ for some complex $\FC \in \HC_{\aff}$, which we call the \emph{twisted
(Gaitsgory complex for the) standard representation}. All the properties of $\VC$ can be explored instead in the context of $\FC$, purely within the non-extended category $\HC_{\aff}$.
Whenever possible we work with $\FC$ and $\HC_{\aff}$ in this paper, and the reader who wishes to bypass the extended Hecke category can still do a lot.

The structure of this paper is as follows.

In \S\ref{sec-pseudo} we recall the theory of pseudocomplexes. We also discuss Gaussian elimination for pseudocomplexes. We leave most of the results of this section as exercises to the
reader.

In \S\ref{sec-234} we write down the Gaitsgory central complex $\VC$ and its twisted version $\FC$ in the cases $n = 2,3,4$. These small examples allow us to highlight certain features of
the complex, which make it easier to explain the general case. They also make it easy to demonstrate the properties of these complexes concretely.

In \S\ref{sec-elaboratesimplify} we continue this introduction on a more technical level, going into detail on the various properties of Goal \ref{goal2}, which were seen in examples in \S\ref{sec-234} In \S\ref{sec-elaboratesimplify} we also provide a simplifying framework one might use to prove Goals \ref{goal1.5} and \ref{goal2}.

In Appendix \S\ref{sec:computations}, we provide the answers to a number of brute force computations, and some useful shorthand to do computations in practice. We construct
the chain maps which correspond to $\sl_n$ webs for $n = 2,3$. We also compute the braiding on $\VC$ and its complete endomorphism ring when $n=2$. We provide the Gaitsgory central
complex for any irreducible representation when $n=2$, and for $S^2 V$ when $n=3$. Finally, we discuss an endomorphism of $\VC$ which plays a role in \cite{GNR}.

We have written this paper so that the reader can choose their path. If they prefer abstractions, they can read \S\ref{sec-elaboratesimplify} first to learn about the beautiful properties satisfied by Gaitsgory central sheaves, and then can turn to \S\ref{sec-pseudo} and \S\ref{sec-234} to see them in action. The reader who, like the author, prefers to get down and dirty in examples before seeing general properties should read the paper in order, first \S\ref{sec-pseudo} and \S\ref{sec-234}, trying \S\ref{sec:computations} if they're still hungry, and then heading to \S\ref{sec-elaboratesimplify}.

%
% In \S\ref{sec-extendbraid} we lay out conventions for affine and extended Weyl groups and their braid groups. This background becomes important when we consider Wakimoto complexes. Also,
% there is both an $\sl_n$ and a $\gl_n$ version of the extended Weyl groups, so we felt it was important to make the difference very clear and precise. Ultimately this paper series will
% handle either case, though the flattening functor is only defined for the $\gl_n$ version.

% In \S\ref{sec-extenddiag} we describe extended Soergel bimodules and the corresponding diagrammatics, in both the $\sl_n$ and the $\gl_n$ versions. This could have been a stand-alone (but
% fairly unexciting and unmotivated) result, which was known to the parts of the community for some time. The $\gl_n$ version has been done previously (with different proofs) in
% \cite{MacThi}, also as a stepping stone towards a more interesting result (in their case, affine Schur-Weyl duality). Because this paper is interested in Gaitsgory's central complexes and
% not the technicalities behind the diagrammatics, we place the proofs that the diagrammatics are correct in the Appendix. These proofs are topological in nature, following closely the
% ideas in \cite{EWFenn}.
%

% In \S\ref{sec-wakimoto} we recall the theory of Rouquier complexes, certain complexes of Soergel bimodules which categorify the braid group. The purpose of this chapter is to construct
% the Wakimoto complexes, which are Rouquier complexes attached to the translation lattice inside the extended braid group.

In \S\ref{sec-construction} we write down $\VC$ and $\FC$ for arbitrary $n$. We also discuss the complex for the dual representation. We prove in this section that $\VC$ is indeed a
pseudocomplex and has the appropriate Wakimoto filtration.

In \S\ref{sec-longestfunctorial} we begin our analysis of $B_I \ot \FC$ for finitary parabolic subsets $I$. We state Theorem \ref{thm:tensorBI} describing the minimal complex of this
tensor product, and outline its proof. We use it to prove ``centrality" of $\VC$ in Theorem \ref{thm:rotationfunctorial}, modulo some further computation. In the following chapters we
prove these results, which is an immense amount of work, and can be skipped without great harm. In \S\ref{sec-descentsets} we discuss descent sets of important elements $w_I h_X$ in the
affine Weyl group. In \S\ref{sec-thicker} we provide the idempotents which project to the indecomposable object $B_{w_I h_X}$. This section uses the thick calculus of \cite{EThick}. In
\S\ref{sec-tensorBIGE} we perform a massive Gaussian elimination to describe the minimal complex of $B_I \ot \FC$. In \S\ref{sec-leftvsright} we compare $B_I \ot \FC$ with $\FC \ot B_{\tau(I)}$ and prove that they are isomorphic, finally proving Theorem \ref{thm:tensorBI}. Then in \S\ref{sec-tensorBIcomps} we provide the additional computations needed
to prove Theorem \ref{thm:rotationfunctorial}.

%================
\subsection{Acknowledgements}
%================

Many thanks to Roman Bezrukavnikov, for the original suggestion for this project, and for answering my countless questions. Thanks to Eugene Gorsky and Matt Hogancamp for explaining
their work and pushing me on with their interest on this project, and thanks to Matt for discussing the awfulness of sign conventions. Many thanks modulo 2 to Victor Ostrik for helping
me with signs. Thanks to Geordie Williamson for answering yet more stupid questions. Thanks to Pramod Achar for explaining his joint work with Laura Rider, and to Kostiantyn
Tolmachov for explaining his work. Thanks to Qiao Zhou and probably many others for their patience in the long long wait for this paper to appear. Thanks also to AIM and the organizers
for the excellent workshop ``Categorified Hecke algebras, link homology, and Hilbert schemes.''

The author was supported by NSF CAREER grant DMS-1553032, NSF FRG grant DMS-1800498, and by the Sloan Foundation.

\section{Affine and extended groups}
\label{sec:extendedbraid}
%%%%%%%%%%%%%%%%%%%%%%%%%

This chapter recalls some basic definitions, and sets a lot of notation.

%=================
\subsection{Various braid groups}
\label{subsec-affbr}
%=================

Fix $n \ge 2$. Let $(W_{\fin},S_{\fin})$ denote the Coxeter system of type $A_{n-1}$. Under the isomorphism $W_{\fin} \cong S_n$, we have $S_{\fin} = \{s_1, \ldots, s_{n-1}\}$. Here,
$s_i$ is the transposition of $i$ and $i+1$. Let $\Br_{\fin}$ denote the associated braid group, called the \emph{(finite) braid group}, with generators $f_s$ and $f_s\inv$ for $s \in
S_{\fin}$. We shorten notation and write $f_i$ for $f_{s_i}$.

Similarly, let $(W_{\aff},S_{\aff})$ denote the Coxeter system of type $\tilde{A}_{n-1}$, with $S_{\aff} = S_{\fin} \cup \{s_0\}$. Let $\Br_{\aff}$ be the associated braid group, called
the \emph{affine braid group}. We often identify $S_{\aff}$ with $\Z/n\Z$. If $k \in \Z$ and $i \in \{0, \ldots, n-1\}$ are equal modulo $n$, then $s_k$ will also denote $s_i \in
S_{\aff}$. We let $\tau$ denote the automorphism of the Coxeter system $(W_{\aff},S_{\aff})$ which sends $s_k \mapsto s_{k+1}$, and corresponds to the ``clockwise'' rotation of the affine Dynkin diagram. We let $\si$ denote the automorphism which sends $s_k \mapsto s_{-k}$, and corresponds to the flip of the affine Dynkin diagram across the affine vertex.

It is typical to visualize the finite braid group $\Br_{\fin}$ using braid diagrams in the planar strip. For example, we draw $f_2$ and $f_2^{-1}$ below, when $n=4$.
\begin{equation}
f_2 = \ig{1}{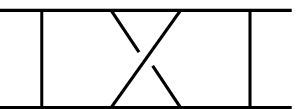}, \qquad f_2\inv = \ig{1}{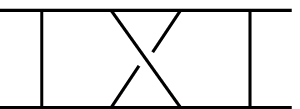}.
\end{equation}
The generator $f_i$ is called an \emph{over-crossing}, and its inverse $f_i^-$ an \emph{under-crossing}.

In similar fashion one can visualize the affine braid group $\Br_{\aff}$ using \emph{cylindrical braid diagrams}, diagrams in the cylinder which identifies the right and left sides of the planar strip. For example, $f_0$ is drawn as below.
\begin{equation}
f_0 = \ig{1}{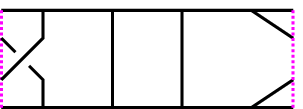}
\end{equation}
We refer to the dotted line where the gluing takes place as the \emph{seam}.

One can also visualize elements of the Coxeter groups $W_{\fin}$ and $W_{\aff}$ in the same way, with crossings $s_i$ replacing the over- and under-crossings $f_i^{\pm}$.
\begin{equation}
s_2 = s_2\inv = \ig{1}{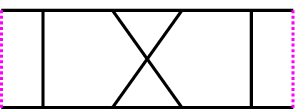}
\end{equation}
We call these \emph{cylindrical crossing diagrams}. To obtain the Coxeter group, one imposes the quadratic relation $s_i^2 = 0$ in addition to the braid relations.

Not all cylindrical braid diagrams correspond to an element of the affine braid group. To a cylindrical braid we may associate its \emph{(total) winding number} in $\Z$, which is the number of strands which go rightward across the seam, minus the number which go leftward. In particular, for any affine braid, the winding number is zero. The following diagram, which is denoted $\om$, has winding number 1.
\begin{equation}
\om = \ig{1}{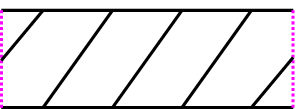}
\end{equation}
The ``inverse'' of $\om$ is pictured below.
\begin{equation}
\om\inv = \ig{1}{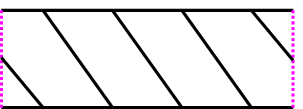}
\end{equation}
By combining affine braids with powers of $\om$, one can obtain any cylindrical braid diagram. The group obtained by taking all cylindrical braid diagrams, modulo isotopy and the braid relations, we will call the \emph{cylindrical braid group} $\Br_{\extg}$.

The subgroup $\Br_{\aff} \subset \Br_{\extg}$ is normal. Conjugation by $\om$ sends $f_k$ to $f_{k+1}$:
\begin{equation} \label{eq:omconj} \om f_k \om\inv = f_{k+1} = \tau(f_k). \end{equation}
Consequently, there is a presentation of the cylindrical braid group as
\begin{equation} \Br_{\extg} \cong \Z \ltimes \Br_{\aff} \end{equation}
where $\om$ is the generator of $\Z$. The map $\wind \co \Br_{\extg} \to \Z$ sends a cylindrical braid to its winding number, and the kernel is the subgroup $\Br_{\aff}$.

The automorphisms $\tau$ and $\si$ extend to the cylindrical braid group via $\tau(\om) = \om$ and $\si(\om) = \om\inv$. Thus $\tau$ becomes an inner automorphism, while $\si$ remains an
outer automorphism.

All the braid groups mentioned above have evaluation homomorphisms $\ev \co \Br \to S_n$, which record how a braid permutes the $n$ strands. A \emph{pure braid} (resp. a \emph{pure
affine} or \emph{cylindrical braid}) is in the kernel of this evaluation map; these form a subgroup $\PBr_{\fin}$ (resp. $\PBr_{\aff}$, $\PBr_{\extg}$). Given a pure cylindrical braid,
one can define the winding number of each strand individually, which gives a surjective homomorphism $\unwind \co \PBr_{\extg} \to \Z^n$. For example, $\unwind(\om f_{n-1} \cdots f_2 f_1) = (1,0,\ldots,0)$, as drawn below.
\begin{equation} \label{eq:y1first} \ig{1}{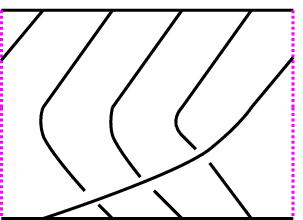}. \end{equation}
For reasons to be explained, we let $\La_{\gl_n}$ denote this lattice $\Z^n$, which is the target of the map $\unwind$.

The map $\unwind$ induces a surjective homomorphism $\unwind \co \PBr_{\aff} \to \La_{\rt}$, where $\La_{\rt} \subset \La_{\gl_n} = \Z^n$ is the sublattice $\{(x_1, \ldots, x_n) \mid
\sum x_i = 0\}$. Note that $\La_{\rt} \cong \Z^{n-1}$ is generated by the \emph{simple roots} $\al_i = (0,\ldots,0,1,-1,0,\ldots,0)$.

Let $\La_{\det} \cong \Z$ denote the sublattice of $\Z^n$ which is generated by $\det = (1,1, \ldots, 1)$. Note that $\unwind(\om^n) = \det$.

% Conjugation gives a map $\rho \co \Br_{\fin} \to \Aut(\PBr_{\extg})$, preserving the total winding number and the subgroup $\PBr_{\aff}$. We also have an action of $S_n$, and
% thereby $\Br_{\fin}$, on $\Z^n = \La_{\gl_n}$ by permuting the entries. These actions are intertwined under the map $\unwind$.

\begin{prop} The center of $\Br_{\extg}$ is generated by $\om^n$, and $\Br_{\extg}/\langle \om^n\rangle$ has trivial center. The group $\Br_{\aff}$ has trivial center. \end{prop}

\begin{remark} This proposition is often quoted without proof in the literature, and we are unaware of its attribution. \end{remark}

% \begin{proof} We begin by noting that $\om^n$ is central, and $\unwind(\om^n) = (1,1,\ldots,1)$, the generator of $\La_{\det}$. If particular, $\om^n$ can not be a commutator $[g,h]$
% of two elements $g, h \in \Br_{\extg}$, because $\La_{\gl_n}$ is abelian so any commutator would be in the kernel of $\unwind$. Consequently, if $Z(\Br_{\extg})$ is generated by
% $\om^n$, then $\Br_{\extg} / \langle \om^n\rangle$ has trivial center.
%
% Any element of the center of any of these groups must be sent to the center of $S_n$ under the map $\ev$. Of course, $Z(S_n) = 1$, so any central element is a pure cylindrical braid. Any
% central element is acted on trivially by the conjugation action of $\Br_{\fin}$, so it must be sent under $\unwind \co \Br_{\extg} \to \Z^n$ to an $S_n$-invariant vector. Thus the image of
% $\unwind$ on $Z(\Br_{\extg})$ is $\La_{\det} \subset \La_{\gl_n}$. Any element of $Z(\Br_{\extg})$ is therefore of the form $\om^{nk} z$ for $z$ in the kernel of $\unwind$. Similarly, any
% element of $Z(\Br_{\aff})$ is in the kernel of $\unwind$, since $\La_{\det} \cap \La_{\sl_n} = 0$. Note that the kernel of $\unwind$ is contained in the kernel of $\wind$, which is $\Br_{\aff}$. So we need only show that the only element of $\PBr_{\aff}$ which is in the kernel of $\unwind$ and commutes with all of $\Br_{\aff}$ is the identity element.
%
% It remains to show that any central element of $\Br_{\aff}$ in the kernel of $\unwind$ is trivial. BLAH BLAH, need a reference... \end{proof}

\begin{defn} The quotient group $\Br_{\extg} / \langle \om^n \rangle$ is the \emph{extended (affine) braid group} $\Br_{\exts}$. \end{defn}	

The map $\unwind$ descends to a map $\unwind \co \PBr_{\exts} \to \La_{\gl_n} / \La_{\det}$. Let us call this quotient lattice $\La_{\sl_n}$. We also let $\La_{\rt}$ denote the (isomorphic) image of $\La_{\rt} \subset \La_{\gl_n}$ inside this quotient.

\begin{remark} In the literature you can find both $\Br_{\extg}$ and $\Br_{\exts}$ referred to as the extended affine braid group, and rightly so. However, $\Br_{\extg}$ is the extended affine braid group with respect to $\gl_n$, while $\Br_{\exts}$ is with respect to $\sl_n$. We explain this in more detail in \S\ref{subsec-loopvscox} below. We decided to distinguish them by calling them cylindrical and extended braid groups; please complain to me if you know of better terminology! \end{remark}

\begin{remark} Since we work in type $A$, we will not bother to distinguish between the root and weight lattices, and the coroot and coweight lattices. \end{remark}

%=================
\subsection{The translation lattice}
\label{subsec-translattice}
%=================

Let us introduce a particular $\Z$-lattice inside $\Br_{\extg}$.

\begin{defn} \label{defn:ydefn} For $1 \le i \le n$, let $y_i$ denote the following pure cylindrical braid, exemplified below.
\begin{equation} y_1 = \ig{.75}{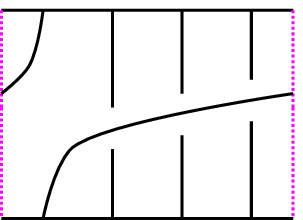} \qquad y_2 = \ig{.75}{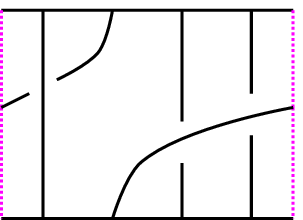} \qquad y_3 = \ig{.75}{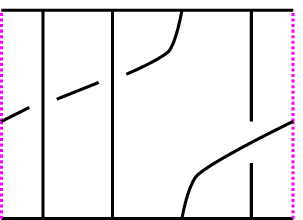} \qquad y_4 = \ig{.75}{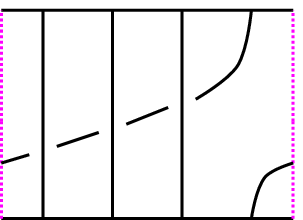} \end{equation}
One should think that the $i$-th strand is at height $i$, where $1$ is the highest and $n$ the lowest. Then $y_i$ takes the $i$-th strand around the cylinder once without changing the height, so that it goes over the strands $j$ with $i < j \le n$, across the seam, and then under the strands $k$ with $1 \le k < i$.
\end{defn}

Explicitly in terms of the Coxeter generators, we have
\begin{equation} y_i = f_{i-1}\inv \cdots f_2\inv f_1\inv \om f_{n-1} \cdots f_{i+1} f_i. \end{equation}
Equivalently, we have
\begin{equation} \label{eq:foryi} y_i = \om f_{i-2}\inv \cdots f_1\inv f_0\inv f_{n-1} \cdots f_{i+1} f_i. \end{equation}

\begin{prop} The elements $\{y_i\}$ commute with each other. The map $\unwind$ induces an isomorphism between $\La_{\gl_n}$ and the subgroup of $\Br_{\extg}$ generated by $\{y_i\}$.  One has
\begin{equation} \prod_{i=1}^n y_i = \om^n. \end{equation} \end{prop}

\begin{proof} (Sketch) Interpreting $y_i$ as a loop around the cylinder at height $i$, one can see topologically that when $i \ne j$, the loops $y_i$ and $y_j$ occur at different heights
and can be pushed past each other using isotopy. Hence $y_i y_j = y_j y_i$. The product of all the $y_i$ can be viewed as the element which takes all strands around the cylinder
simultaneously, which is $\om^n$. One can also prove these statements directly using the braid relations. The fact that $\unwind$ is injective and surjective is clear. \end{proof}

\begin{defn} For $\la \in \La_{\gl_n}$, we let $w_\la$ denote the element of $\Br_{\extg}$ living in the subgroup generated by $\{y_i\}$, for which $\unwind(w_{\la}) = \la$. Let $\ep_i = (0,\ldots,0,1,0,\ldots,0) \in \La_{\gl_n}$,
where the $1$ is in the $i$-th position. Let $\fund_i = (1,\ldots,1,0,\ldots,0)$, where the last $1$ is in the $i$-th position. The element $\fund_n = (1,1,\ldots,1)$ will also be called
$\det$. \end{defn}

For example, $y_i = w_{\ep_i}$.

\begin{defn} We say that an element $\la = (x_1, \ldots, x_n) \in \La_{\gl_n}$ is \emph{dominant} if $x_1 \ge x_2 \ge \ldots \ge x_n \ge 0$. We say that $\la$ is \emph{antidominant} if
$-\la$ is dominant. We say that an element of $\Br_{\extg}$ is \emph{positive} if it can be expressed using only $\om$ and $f_k$ for $k \in \Z$, without using $\om\inv$ or $f_k\inv$. We
say that an element is \emph{negative} if it can be expressed using only $\om\inv$ and $f_k\inv$. \end{defn}

\begin{prop} \label{prop:dominantpositive} If $\la$ is dominant, then $w_\la$ is positive. If $\la$ is antidominant then $w_\la$ is negative. \end{prop}

\begin{proof} An element of $\Br_{\extg}$ is negative if and only if its inverse is positive. Since $w_{-\la} = w_\la\inv$, it is enough to show the first statement.

Any dominant element is a sum $\sum_{k=1}^n a_k \fund_k$ for $a_k \ge 0$. Hence it is enough to show that $w_{\fund_k}$ is positive for all $k$. Note that $w_{\fund_k} = \prod_{i=1}^k
y_i$. Topologically, $w_{\fund_k}$ can be viewed as the braid which takes the first $k$ strands around the cylinder once. Since these are the highest strands, every crossing will be an
overcrossing. More explicitly, $w_{\fund_k}$ can be expressed as the positive braid in $\Br_{\fin}$ which crosses the first $k$ strands over the last $n-k$ strands, then multiplies by $\om^k$. \end{proof}

In particular since each element $\la \in \La_{\gl_n}$ can be written (not uniquely) as the sum of a dominant weight $\mu$ and an antidominant weight $\nu$, we see that each $w_{\la}$ has a \emph{positive-negative decomposition}: it can be expressed as a positive braid times a negative braid (and in this case, the positive and negative braids commute).

\begin{prop} One has $w_\la \in \Br_{\aff}$ if and only if $\la \in \La_{\rt}$. \end{prop}

\begin{proof} The winding number of $w_{\la}$ is $\sum x_i$, so $w_{\la} \in \Br_{\aff}$ if and only if $\sum x_i = 0$. \end{proof}

\begin{ex} \label{ex:longestreflection} Let $\al_{\longest} = (1, 0, \ldots, 0, -1)$ be the so-called \emph{longest root}. Then 
\begin{equation} w_{\al_{\longest}} = y_1 y_n\inv = (\om f_{n-1} \cdots f_1) (f_{n-1}\inv \cdots f_1\inv \om)\inv = f_0 f_{n-1} \cdots f_2 f_1 f_2 \cdots f_{n-1}. \end{equation}
In particular, let $f_{\longest} = f_{n-1} \cdots f_2 f_1 f_2 \cdots f_{n-1}$, which is the positive lift to $\Br_{\fin}$ of the longest reflection $s_{\longest}$ in $W_{\fin}$. Then
\begin{equation} w_{\al_{\longest}} = f_0 f_{\longest}. \end{equation}	
In particular, $w_{\al_{\longest}}$ is the standard positive lift of the usual expression $s_0 s_{\longest}$ for translation by the longest root in the affine Weyl group. \end{ex}

%=================
\subsection{Minuscule weights}
\label{subsec-translationexamples}
%=================

Because they will be important in this paper, we give examples of minuscule weights.

\begin{ex} Let $n=2$. Then
\begin{subequations}
\begin{equation} w_{(1,0)} = \om f_1, \end{equation}
\begin{equation} w_{(0,1)} = \om f_0\inv. \end{equation}
\end{subequations}
\end{ex}

\begin{ex} Let $n=3$. Then
\begin{subequations}
\begin{equation} w_{(1,0,0)} = \om f_2 f_1, \end{equation}
\begin{equation} w_{(0,1,0)} = \om f_0\inv f_2, \end{equation}
\begin{equation} w_{(0,0,1)} = \om f_1\inv f_0\inv. \end{equation}

\begin{equation} w_{(1,1,0)} = \om^2 f_1 f_2, \end{equation}
\begin{equation} w_{(1,0,1)} = \om^2 f_0\inv f_1, \end{equation}
\begin{equation} w_{(0,1,1)} = \om^2 f_2\inv f_0\inv. \end{equation}
\end{subequations}
\end{ex}

\begin{ex} \label{ex:wakimoton4} Let $n=4$. Then
\begin{subequations}
\begin{equation} w_{(1,0,0,0)} = \om f_3 f_2 f_1, \end{equation}
\begin{equation} w_{(0,1,0,0)} = \om f_0\inv f_3 f_2, \end{equation}
\begin{equation} w_{(0,0,1,0)} = \om f_1\inv f_0\inv f_3, \end{equation}
\begin{equation} w_{(0,0,0,1)} = \om f_2\inv f_1\inv f_0\inv. \end{equation}
	
\begin{equation} w_{(1,1,0,0)} = \om^2 f_2 f_1 f_3 f_2, \end{equation}
\begin{equation} w_{(1,0,1,0)} = \om^2 f_0\inv f_2 f_1 f_3, \end{equation}
\begin{equation} w_{(1,0,0,1)} = \om^2 f_1\inv f_0\inv f_2 f_1, \end{equation}
\begin{equation} w_{(0,1,1,0)} = \om^2 f_3\inv f_0\inv f_2 f_3, \end{equation}
\begin{equation} w_{(0,1,0,1)} = \om^2 f_3\inv f_1\inv f_0\inv f_2, \end{equation}
\begin{equation} w_{(0,0,1,1)} = \om^2 f_0\inv f_3\inv f_1\inv f_0\inv. \end{equation}

\begin{equation} w_{(1,1,1,0)} = \om^3 f_1 f_2 f_3, \end{equation}
\begin{equation} w_{(1,1,0,1)} = \om^3 f_0\inv f_1 f_2, \end{equation}
\begin{equation} w_{(1,0,1,1)} = \om^3 f_3\inv f_0\inv f_1, \end{equation}
\begin{equation} w_{(0,1,1,1)} = \om^3 f_2\inv f_3\inv f_0\inv. \end{equation}
\end{subequations}
\end{ex}

\begin{remark} The fastest way for me to reconstruct these reduced expressions is as follows. The formula for $w_{\fund_k}$ is straightforward. One can get any other minuscule weight
with $\sum x_i = k$ by the operation which swaps some instance of $(\ldots, 1,0, \ldots)$ with $(\ldots, 0,1, \ldots)$. If $x_i=1$ and $x_{i+1}=0$ then $f_i$ appears on the right in some
reduced expression for $w_\la$, and conjugating by $f_i$ achieves the swap to $x_i=0$ and $x_{i+1}=1$; this is because conjugation by $f_i$ sends $y_i$ to $y_{i+1}$, and commutes with
$y_j$ for $j \ne i, i+1$. In terms of the word appearing to the right of $\om^k$, conjugating by $f_i$ replaces this $f_i$ on the right with $f_{i+k}\inv$ on the left. \end{remark}

%=================
\subsection{Loop versus Coxeter}
\label{subsec-loopvscox}
%=================

We can extend the group $W_{\aff}$ by adding a formal element $\om$ which corresponds to the automorphism $\tau$, to obtain the extended affine Weyl group. This extension comes in two flavors, just as the braid group did.

\begin{defn} \label{defn:Coxpresentext} Let the \emph{cylindrical (affine) Weyl group} be $W_{\extg} = W_{\aff} \ltimes \Z$. Here, the generator of $\Z$ is denoted $\om$, and it acts on $W_{\aff}$ by
\begin{equation} \label{eq:omsi}\om s_k \om\inv = s_{k+1} = \tau(s_k). \end{equation} It is easy to see that $\om^n$ is central. Let the \emph{extended (affine) Weyl group} be $W_{\exts} = W_{\extg} / \langle \om^n \rangle$. \end{defn}

One has a presentation of $W_{\extg}$, extending the Coxeter presention of $W_{\aff}$ with a new generator $\om$ and the relation \eqref{eq:omsi}. One has a presentation of $W_{\exts}$
with the additional relation $\om^n = 1$. We call these the \emph{Coxeter presentations} of these groups.

It is well-known that $W_{\aff}$ also has a \emph{loop presentation} in terms of translations. Note that $W_{\fin} = S_n$ acts in the usual way on $\La_{\gl_n} = \Z^n$, and hence on the other lattices $\La_{\rt}$ and $\La_{\sl_n}$.

\begin{prop} There is an isomorphism $W_{\aff} \cong W_{\fin} \ltimes \La_{\rt}$. This latter group can be thought of as a subgroup of the affine linear transformations of $\La_{\gl_n}$
or $\La_{\sl_n}$, where an element $\la \in \La_{\rt}$ corresponds to translation by $\la$. \end{prop}

In Example \ref{ex:longestreflection} we have already described how translation by $\al_{\longest}$ is equal to $s_0 s_{\longest}$. This gives enough information to reconstruct the
isomorphism between the Coxeter and loop presentations.

The lattice $\La_{\sl_n}$ is a $\Z/n\Z$-extension of $\La_{\rt}$, and $\La_{\gl_n}$ is a $\Z$-extension of $\La_{\rt}$. The loop presentation of the affine Weyl group extends to
presentations of the extended and cylindrical Weyl groups.

\begin{prop} There are isomorphisms $W_{\exts} \cong W_{\fin} \ltimes \La_{\sl_n}$, and $W_{\extg} \cong W_{\fin} \ltimes \La_{\gl_n}$. \end{prop}

The isomorphisms above also deform to the Hecke algebra setting. Consider the affine Hecke algebra $\Hecke_{\aff}$ (the Hecke algebra attached to $W_{\aff}$), an extended Hecke
algebra $\Hecke_{\exts}$, and a cylindrical Hecke algebra $\Hecke_{\extg}$. Each can be defined as the quotient of the $\Z[v,v\inv]$-linear group algebra of the corresponding braid group by the quadratic relation
\begin{equation} (f_i + v)(f_i - v\inv) = 0. \end{equation}

\begin{prop} There are isomorphisms $\Hecke_{\exts} \cong \Hecke_{\fin} \ltimes \La_{\sl_n}$, and $\Hecke_{\extg} \cong \Hecke_{\fin} \ltimes \La_{\gl_n}$. \end{prop}

These isomorphisms all lift to the respective braid groups, and can be reconstructed using the correspondence $\la \mapsto w_{\la}$ for $\la \in \La_{\gl_n}$.

%=================
\subsection{The dual translation lattice}
\label{subsec-duality0}
%=================

There is an algebra homomorphism on $\Br_W$ (for any Coxeter group $W$) which sends $f_s \mapsto f_s\inv$ for every simple reflection $s$. This does not send $\be \mapsto \be\inv$ for any
braid element $\be$, since it is a homomorphism, not an antihomomorphism. Instead, if $w$ is an element of $W$ and $f_w$ its positive lift to $\Br_W$, then \begin{equation}\label{eq:bar}
f_w \mapsto f_{w\inv}\inv.\end{equation} We call this involution the \emph{bar involution}. Diagrammatically, the bar involution just turns every overcrossing into an undercrossing, and
vice versa.

It helps to visualize the bar involution as reversing the order of the height of strands. For example, in the drawing of $f_2$ it appears that the strand whose bottom label is $2$ is higher than the strand whose bottom label is $3$, and $f_2$ just permutes the strands at their respective heights. Meanwhile, $f_2\inv$ does the same when strand $2$ is lower than strand $3$.

To extend the bar involution, one is tempted to try $\om \mapsto \om\inv$, but this would not give a homomorphism as the relation \eqref{eq:omconj} would not be preserved.
The correct way to extend the bar involution to $\Br_{\extg}$ or $\Br_{\exts}$ is to send $\om \mapsto \om$. After all, $\om$ does not depend on the heights of the strands. Another reasoning is that, if we temporarily write $\om$ as $f_{\om}$ and use \eqref{eq:bar}, we obtain
\[ \om = f_{\om} \mapsto f_{\om\inv}\inv = (\om\inv)\inv = \om. \]

The bar involution preserves the winding number. Under the bar involution, a translation $w_{\la}$ is sent to the \emph{dual translation} $w_{\la}^*$. These dual translations form another
lattice inside $\Br_{\extg}$. The dual $y_i^*$ to the element $y_i$ can be drawn by simply reversing the order on the heights.

\begin{ex} When $n=4$,
\begin{equation} y_2 = \ig{1}{y2}, \qquad y_2^* = \ig{1}{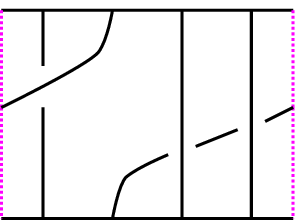}. \end{equation}
\end{ex}

\begin{ex} Let $n=2$. Then
\begin{subequations}
\begin{equation} w_{(1,0)}^* = \om f_1\inv, \end{equation}
\begin{equation} w_{(0,1)}^* = \om f_0. \end{equation}
\end{subequations}
Note that the dual translation lattice does not overlap the translation lattice except at the identity.
\end{ex}

\begin{remark} In fact, there are $n!$ different non-overlapping translation lattices, obtained by reordering the heights of the $n$ strands in the definition of $y_i$. These
lattices are in bijection with $W_{\fin}$, and the dual lattice corresponds to the longest element $w_0$. \end{remark}

One major reason to prefer the translation lattice over the dual translation lattice is Proposition \ref{prop:dominantpositive}. A dominant translation is positive, while a dominant dual
translation is neither positive nor negative.

%=================
\subsection{Flattening}
\label{subsec-flattening}
%=================

We now describe the \emph{flattening} homomorphism from $\Br_{\extg}$ to $\Br_{\fin}$. Visually, one should think that the seams of the cylinder are glued together behind the page. Then, flatten the cylinder back into a planar strip, as one flattens a cardboard tube, and consider the resulting finite braid.

\begin{defn} Let $\flat \co \Br_{\extg} \to \Br_{\fin}$ be the homomorphism defined on generators by:
\begin{subequations}
\begin{equation} \flat(f_i) = f_i \qquad \text{ for $1 \le i \le n-1$},\end{equation}
\begin{equation} \flat(\om) = f_1 f_2 \cdots f_{n-1} = \ig{1}{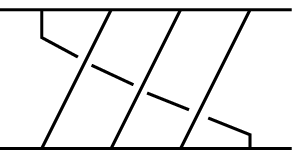},\end{equation}
\begin{equation} \flat(f_0) = \flat(\om\inv f_1 \om) = f_{n-1}\inv \cdots f_2\inv f_1 f_2 \cdots f_{n-1} = \ig{1}{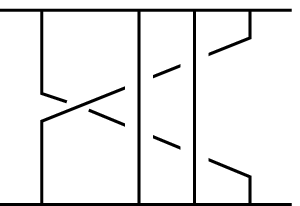}. \end{equation}
\end{subequations}
\end{defn}
We let $h$ denote $\flat(\om)$, which is the positive lift of a Coxeter element in $W_{\fin}$.

Under this flattening homomorphism, the elements $y_i$ are sent to the famous (multiplicative) Jucys-Murphy braids. More precisely, $y_n$ is sent to the identity, and for $i < n$ one has
\begin{equation} \label{eq:goestoJM} y_i \mapsto j_i = f_i f_{i+1} \cdots f_{n-1} f_{n-1} \cdots f_{i+1} f_i. \end{equation}
This is the positive braid which wraps the $i$-th strand around all the strands $k$ with $i < k \le n$. Here is the example when $i=2$ and $n=4$.
\[ j_2 = \ig{1}{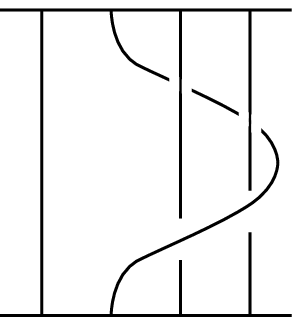} \]

\begin{remark} Often in the literature the multiplicative Jucys-Murphy braids are defined as the positive braids wrapping the $i$-th strand around the strands $k$ for $1 \le k < i$. Let
us temporarily call these the \emph{leftward} Jucys-Murphy braids, as opposed to the \emph{rightward} Jucys-Murphy braids of \eqref{eq:goestoJM}. The leftward Jucys-Murphy braids are
obtained by flattening the dual translations $y_i^*$.

There are several choices we made in our
conventions above:\begin{itemize} \item Whether $\om$ goes to the right or the left (we chose the right). \item Whether the seam of the cylinder is glued above or behind the page (we chose
behind). \item Whether the heights of the strands are increasing or decreasing, in the definition of $y_i$ (we chose decreasing). \end{itemize} We have chosen our conventions so that
$y_1$ is positive (essential for Proposition \ref{prop:dominantpositive}), and so that the flattening of $\om$ is positive. Any such choice leads to the flattening of
$y_i$ being the rightward Jucys-Murphy braid. \end{remark}

\begin{remark} If $\weights(k)$ denotes the set of weights in the $k$-th fundamental representation, then the linear combination of braids
\[ \sum_{\la \in \weights(k)} w_\la \]
flattens to the elementary symmetric polynomial
\[ e_k(j_1, j_2, \ldots, j_n) \]
of multiplicative Jucys-Murphy braids. In particular, the central element $\om^n = w_{\det}$ is flattened to the positive full twist $ft_n = \prod_{i=1}^n j_i$. \end{remark}

%%%%%%%%%%%%%%%%%%%%%%%%%
\section{Extended Soergel bimodules and diagrammatics}
\label{sec:extendeddiag}
%%%%%%%%%%%%%%%%%%%%%%%%%

The extended diagrammatic Hecke category appeared first in work of Mackaay-Thiel \cite[Chapter 2]{MacThi}. We give a different exposition of this material, in more generality and with
different proofs (focused on diagrammatic arguments, rather than using the functor to bimodules). We also treat for the first time the case when the automorphism of the Dynkin diagram has
finite order. We will assume that the reader is familiar with the usual diagrammatic Hecke category, see \cite{EWGr4sb}.

%=================
\subsection{Soergel bimodules: a reminder}
\label{subsec-sbim}
%=================

In \cite{EWGr4sb}, a \emph{realization} of a Coxeter system $(W,S)$ was defined to be a the data of:
\begin{itemize} \item a commutative base ring $\Bbbk$,
\item a free $\Bbbk$-module $\hg$, and its $\Bbbk$-dual $\hg^*$,
\item and a collection of simple roots $\{\al_s\}$ in $\hg^*$ and simple coroots $\al_s^\vee$ in $\hg$ indexed by $s \in S$, such that $\al_s(\al_s^\vee)=2$,
\end{itemize}
satisfying some conditions.
If one defines the ``reflection'' $s$ on $\hg$ to be
\begin{equation} s(x) = x - \al_s(x) \al_s^\vee,\end{equation}
then one condition states that this should produce an action of $W$ on $\hg$ (and consequently, on $\hg^*$). We will not recall the other (more technical, and usually redundant) condition, but it is necessary for the $2m$-valent vertices to be cyclic.  We also make the assumption (called Demazure surjectivity in \cite{EWGr4sb}) that $\al_s \co \hg \to \Bbbk$ and $\al_s^\vee \co \hg^* \to \Bbbk$ are both surjective, for all $s \in S$. Note that we do not assume that the roots (resp. coroots) span, or that they are linearly independent. The \emph{Cartan matrix} of the realization is given by the pairing $\al^\vee_s(\al_t)$ for all $s, t \in S$.

Here are several realizations of the affine Weyl group in type $A$, which all have the standard affine Cartan matrix.

\begin{defn} \label{defn:standardaffine} The \emph{standard affine realization} in type $\tilde{A}_{n-1}$ has the following description. The space $\hg^*$ has a basis given by $\{x_1, x_2, \ldots, x_n, \delta\}$
where $\delta$ is $W$-invariant. The finite Weyl group $S_n$ acts to permute the variables $x_i$, $1 \le i \le n$ in the usual way. The affine reflection $s_0$ acts by the formula
\begin{equation} s_0(x_1) = x_n + \delta, \quad s_0(x_n) = x_1 - \delta, \quad s_0(x_i) = x_i \textrm{ for } i \ne 1, n. \end{equation} The simple roots are $\al_i = x_i - x_{i+1}$ for
$1 \le i \le n-1$, and $\al_0 = x_n - x_1 + \delta$, so that $\delta$ is the sum of the simple roots. Meanwhile, the coroots in the dual vector space $\hg$ are given by Demazure operators. Note that the simple roots are linearly independent, but the coroots are linearly dependent. This is the realization used in \cite[Section 2.1.1]{MacThi}. \end{defn}

\begin{defn} \label{defn:rootspanaffine} The \emph{root span realization} in type $\tilde{A}_{n-1}$ is defined by setting $\hg^*$ to be the span of the roots $\al_i$ in the standard
affine realization. \end{defn}

\begin{defn} \label{defn:nondegaffine} The \emph{nondegenerate realization} in type $\tilde{A}_{n-1}$ is defined by extending the root span realization by an element $\fund_0$ for which
$\al_0^\vee(\fund_0) = 1$ and $\al_i^\vee(\fund_1) = 0$ for all $i \ne 0$. Now the simple coroots are linearly independent. Ironically, the nondegenerate realization will cause issues
for the centrality of Gaitsgory central complexes. \end{defn}

\begin{remark} For all the realizations above, one can construct another realization where the simple roots are linearly dependent by setting $\delta = 0$. All coroots already live in
the perpendicular space to $\delta$, so they still make sense in this new realization. \end{remark}

\begin{remark} The reader used to finite type may have the intuition that expanding from the span of the roots to the span of the $x_i$ is like adding ``fundamental weights'' to the
root lattice, which are dual to the coroots. However, in affine type, there is a vast difference between the standard and nondegenerate realizations; the standard realization does not
have fundamental weights, and its coroots are linearly dependent. \end{remark}

Let $R = \mathrm{Sym}(\hg^*)$ be the ring of polynomial functions on $\hg$ with doubled degree, i.e. where $\hg^*$ is given degree 2. Recall that \emph{Bott-Samelson bimodules} are
tensor products of $R$-bimodules of the form $B_s = R \ot_{R^s} R(1)$ for various $s \in S$. Here, $R^s$ is the subring of $s$-invariants, and $(1)$ is the grading shift which places
the element $1 \ot 1$ in degree $-1$. Bott-Samelson bimodules form a graded monoidal category $\BSBim$ where the space of (homogeneous) morphisms (of all degrees) forms an $R$-bimodule.

Recall that \emph{Soergel bimodules} are direct sums of grading shifts of direct summands of Bott-Samelson bimodules. The category of Soergel bimodules $\SBim$ is the smallest full
subcategory of $R$-bimodules which is closed under tensor products, grading shifts, direct sums, and direct summands, and contains the bimodules $B_s$ for each $s \in S$. Said another way, $\SBim$ is the Karoubi envelope of the additive, graded envelope of $\BSBim$.

Whenever $w = s_1 s_2 \cdots s_d$ is a reduced expression in $W$, the Bott-Samelson bimodule $B_{s_1} B_{s_2} \cdots B_{s_d}$ has a unique indecomposable direct summand in $\SBim$
which has not appeared in any shorter tensor products. This summand is called $B_w$. In \cite{Soer07} it is proven that the indecomposable objects in $\SBim$ have the form $B_w(k)$,
where $w \in W$ and $k \in \Z$, and no two such are isomorphic.

%=================
\subsection{Extended Soergel bimodules}
\label{subsec-extendedsbim}
%=================

\begin{defn} \label{defn:realizationextended} An \emph{automorphism} of a realization is the data of a bijection $\tau \co S \to S$ which gives an automorphism of $W$, and dual
automorphisms of $\hg$ and $\hg^*$ (which we also denote $\tau$) such that $\tau(\al_s) = \al_{\tau(s)}$ and $\tau(\al_s^\vee) = \al_{\tau(s)}^\vee$. We refer to a realization paired
with an automorphism $\tau$ as a \emph{realization} of the $\tau$-extended Coxeter group. \end{defn}

From this definition, one deduces that
\begin{equation} \label{eq:tauact} \tau(s)\cdot \tau(x) = \tau(s\cdot x) \end{equation}
for any $x \in \hg^*$. In particular, $x$ is $s$-invariant if and only if $\tau(x)$ is $\tau(s)$-invariant.

\begin{ex} \label{ex:standardaffinetau} On the standard affine realization of Definition \ref{defn:standardaffine}, we set
\begin{equation} \tau(\delta) = \delta, \quad \tau(x_n) = x_1 - \delta, \quad \tau(x_i) = x_{i+1} \textrm{ for } 1 \le i \le n-1. \end{equation}
Then $\tau$ defines an automorphism of the realization, corresponding to ``clockwise'' rotation of the affine Dynkin diagram. It descends to an automorphism of the root span realization. Note that $\tau^n$ acts as the identity on the simple roots, but that $\tau^n(x_i) = x_i - \delta$, so 
$\tau$ has infinite order. \end{ex}

\begin{ex} On the standard affine realization of Definition \ref{defn:standardaffine}, we set
\begin{equation} \si(\delta) = \delta, \quad \si(x_i) = -x_{n+1-i} \textrm{ for } 1 \le i \le n. \end{equation}
Then $\si$ defines an automorphism of the realization, corresponding to ``reflection through the affine root'' on the affine Dynkin diagram. It descends to an automorphism of the root span realization. Note that
\begin{equation} \si^2 = \id, \quad \si \tau \si = \tau^{-1}. \end{equation} \end{ex}

\begin{defn} Given an automorphism $\tau$ of the realization, the $R$-bimodule $R_{\tau}$ is defined as follows. It is free of rank 1 as a left or right $R$-module. The left action is
the ordinary one, while the right action of $f \in R$ is multiplication by $\tau(f)$. One defines $R_{\tau\inv}$ similarly, where the right action of $f\in R$ agrees with the left
action of $\tau\inv(f)$. \end{defn}

\begin{lemma} One has the following bimodule isomorphisms:
\begin{equation} \label{eq:tauandinv} R_{\tau} \ot R_{\tau\inv} \cong R \cong R_{\tau\inv} \ot R_{\tau},\end{equation}
\begin{equation} \label{eq:Rtauconj} R_{\tau} \ot B_s \cong B_{\tau(s)} \ot R_{\tau}. \end{equation}
\end{lemma}

\begin{proof} The isomorphism $R_{\tau} \ot R_{\tau\inv} \to R$ sends $f \ot g \mapsto f\tau(g)$, and the inverse isomorphism sends $f \mapsto f \ot 1$. It is easy to check that these are inverse bimodule isomorphisms.

The isomorphism $R_{\tau} \ot B_s \ot R_{\tau\inv} \to B_{\tau(s)}$ sends $1 \ot (f \ot_s g) \ot 1 \mapsto \tau(f) \ot_{\tau(s)} \tau(g)$. To show this map is well-defined, observe that
if $f$ is $s$-invariant then it can slide over to join $g$; this is equivalent by \eqref{eq:tauact} to $\tau(f)$ being $\tau(s)$ invariant, and sliding over to join $\tau(g)$. We leave
the remaining details to the reader. \end{proof}

We can define $R_{\tau^k}$ for $k \in \Z$ in similar fashion to $R_{\tau}$ and $R_{\tau\inv}$, and observe that
\begin{equation} R_{\tau^k} \ot R_{\tau^\ell} \cong R_{\tau^{k+\ell}}. \end{equation}

\begin{defn} The category of \emph{$\tau$-extended Soergel bimodules} is the smallest full subcategory of $R$-bimodules which is closed under tensor products, grading shifts, direct sums, and direct summands, and contains the bimodules $R_{\tau}$ and $R_{\tau\inv}$ and $B_s$ for each $s \in S$. \end{defn}

Using \eqref{eq:Rtauconj}, it is straightforward to verify that every $\tau$-extended Soergel bimodule is a direct summand of $R_{\tau^k} \ot BS$ for some Bott-Samelson bimodule $BS$,
and hence has the form $R_{\tau^k} \ot B$ for some Soergel bimodule $B$.

Note that, if $\tau$ has finite order $n$, the bimodules $R_{\tau^n}$ and $R$ are canonically isomorphic. When $S$ is a finite set, any automorphism of $S$ automatically has finite
order. However, $\tau$ need not have finite order, since $\hg^*$ (resp. $\hg$) need not be spanned by the simple roots (resp. coroots), so that the order of $\tau$ on $S$ need not agree
with its order as an automorphism on $\hg^*$.

%=================
\subsection{Extended diagrammatics}
\label{subsec-extendeddiag}
%=================

In \cite{EWGr4sb}, a diagrammatic category $\Diag$ was defined for any realization of a Coxeter group. It is a graded monoidal category where the space of (homogeneous) morphisms (of all
degrees) forms an $R$-bimodule. For some realizations, this diagrammatic category $\Diag$ is equivalent to $\BSBim$. We abusively let $B_s$ denote the generating objects of the
diagrammatic category $\Diag$; their identity maps are represented by colored strands, with one color for each $s \in S$.

In future chapters we may also use names for the particular generators of $\Diag$, such as the ``startdot'' and ``enddot,'' see \cite[page 16]{EKho}.  Let us also recall the polynomial forcing relation \cite[(5.2)]{EWGr4sb}, which we will use many times in this paper. Here $s$ is represented by the color red, and $f \in R$.

\begin{equation} \label{eq:polyforce} f \linered = \linered s(f) + \brokenred \pa_s(f). \end{equation}

Let us provide a straightforward extension of $\Diag$, which adds a new invertible object $\Om$ corresponding to the bimodule $R_{\tau}$. It is completely analogous to what is found in
\cite{MacThi}. This extension will (purposely) not include the additional isomorphism $R_{\tau^n} \cong R$ which occurs when $\tau$ has order $n$; to add this isomorphism, see
\S\ref{subsec-extendedsbimquotient}.

\begin{defn} Let $\EDiag$ be the extension of $\Diag$ defined as follows. One adds to $\Diag$ a new pair of generating objects $\Om^+$ and $\Om^-$, whose identity maps are
indicated by oriented black strands.
\begin{equation} \id_{\Om^+} = \ig{1}{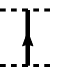}, \qquad \id_{\Om^-} = \igv{1}{Omega}. \end{equation}
There are new morphisms called \emph{oriented cups and caps}, which all have degree zero.
\begin{equation} \label{eq:blackcupscaps} \ig{1}{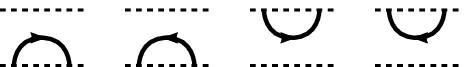} \end{equation}
They satisfy the following relations, which indicate that $\Om^+$ and $\Om^-$ are inverse and biadjoint to each other.

\begin{subequations} \label{eq:cupcaprelations}
\begin{equation} \label{eq:blackisotopy} {
\labellist
\small\hair 2pt
 \pinlabel {$=$} [ ] at 54 78
 \pinlabel {$=$} [ ] at 103 78
 \pinlabel {$=$} [ ] at 54 19
 \pinlabel {$=$} [ ] at 103 19
\endlabellist
\centering
\ig{1}{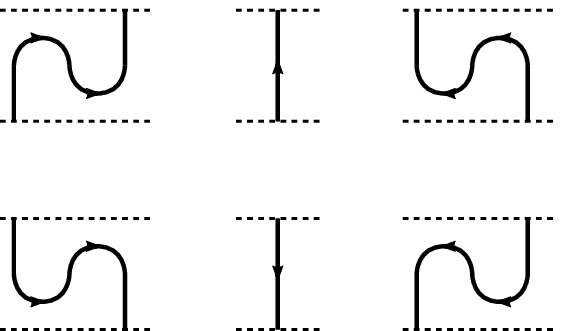}
} \end{equation}
\begin{equation} \label{eq:blackcircle} {
\labellist
\small\hair 2pt
 \pinlabel {$=$} [ ] at 54 16
 \pinlabel {$=$} [ ] at 103 16
\endlabellist
\centering
\ig{1}{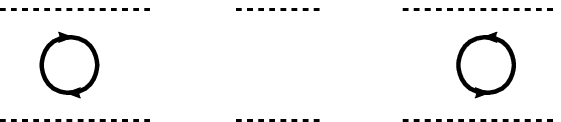}
} \end{equation}
\begin{equation} \label{eq:blackbreak} {
\labellist
\small\hair 2pt
 \pinlabel {$=$} [ ] at 38 16
\endlabellist
\centering
\ig{1}{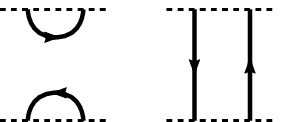}
} \end{equation}
\end{subequations}
Note: because of \eqref{eq:blackisotopy}, the rotation of any relation is another relation. Rotating \eqref{eq:blackbreak} by 90 degrees will give the orientation-reversed version of \eqref{eq:blackbreak}.

We also have new morphisms which mix black and colored strands, which we call \emph{mixed crossings}. They have degree zero.
\begin{equation} \ig{1}{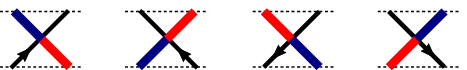} \end{equation}
Above, red represents a simple reflection $s$, and blue represents $\tau(s)$.

We impose several more relations. The first relations state that the mixed crossings are cyclic with respect to various biadjunctions.
\begin{equation} \label{eq:mixedcyclic} {
\labellist
\small\hair 2pt
 \pinlabel {$=$} [ ] at 45 14
 \pinlabel {$=$} [ ] at 88 14
\endlabellist
\centering
\ig{1}{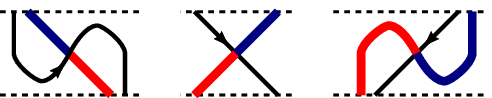}
} \end{equation}
These relations are exemplary; we also include the various horizontal and vertical flips of these.

The next set of relations allows one to ``pull'' morphisms across a black strand, at the cost of changing the colors by $\tau$.
\begin{subequations} \label{eq:pullapart}
\begin{equation} \label{eq:mixedR2} {
\labellist
\small\hair 2pt
 \pinlabel {$=$} [ ] at 40 18
\endlabellist
\centering
\ig{1}{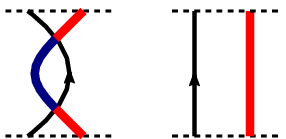}
} \end{equation}
We also include the horizontal flip of \eqref{eq:mixedR2}.
\begin{equation} \label{eq:mixeddot} {
\labellist
\small\hair 2pt
 \pinlabel {$=$} [ ] at 40 18
\endlabellist
\centering
\ig{1}{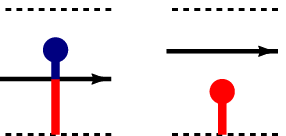}
} \end{equation}
\begin{equation} \label{eq:mixedtri} {
\labellist
\small\hair 2pt
 \pinlabel {$=$} [ ] at 40 18
\endlabellist
\centering
\ig{1}{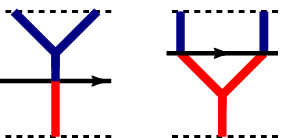}
} \end{equation}
\begin{equation} \label{eq:mixed6v} {
\labellist
\small\hair 2pt
 \pinlabel {$=$} [ ] at 40 18
\endlabellist
\centering
\ig{1}{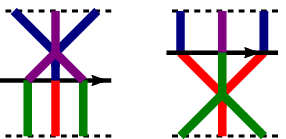}
} \end{equation}
In \eqref{eq:mixed6v}, green represents a simple reflection $t$, and purple represents $\tau(t)$. This relation is supposed to represent a host of relations, one for each pair $\{s,t\} \subset S$ with $m_{st} < \infty$, which permit one to pull the $2m_{st}$-valent vertex across the black strand.
\begin{equation} \label{eq:mixedpoly} {
\labellist
\small\hair 2pt
 \pinlabel {$=$} [ ] at 34 16
 \pinlabel {$f$} [ ] at 14 11
 \pinlabel {$\tau(f)$} [ ] at 59 18
\endlabellist
\centering
\ig{1}{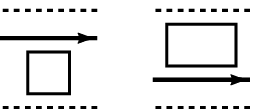}
} \end{equation}
In \eqref{eq:mixedpoly}, $f \in R$ is some homogeneous polynomial.
\end{subequations}

This concludes the definition.
\end{defn}
	
Note that $\Diag$ is a monoidal subcategory of $\EDiag$. Note also that $\EDiag$ is a monoidal category where objects have biadjoints and morphisms are cyclic, meaning that an isotopy class of diagram unambiguously represents a morphism.

\begin{defn} Fix $k \in \Z$. When $k = 0$, let $\Om^k$ denote $\one$. When $k > 0$, let $\Om^k$ denote $\Om^+ \ot \cdots \ot \Om^+$, $k$-times. When $k < 0$, let $\Om^k$ denote $\Om^- \ot \cdots \ot \Om^-$, $-k$ times. \end{defn}

\begin{lemma} One has the isomorphisms
\begin{equation} \label{eq:Omandinv} \Om^+ \ot \Om^- \cong \one \cong \Om^- \ot \Om^+.\end{equation}
Consequently, the isomorphism classes of objects constructed only from $\Om^+$ and $\Om^-$ are enumerated by $\Om^k$. \end{lemma}

\begin{proof} The isomorphism from $\Om^+ \ot \Om^- \to \one$ is given by a clockwise cap, and the inverse isomorphism by a clockwise cup. That $\Om^k$ is not isomorphic to $\Om^\ell$ for $k \ne \ell$ is clear because the generating morphisms can not change the overall exponent of $\Om^+$ in an object. \end{proof}

\begin{lemma}
\begin{equation} \label{eq:BOmconj} \Om^+ \ot B_s \cong B_{\tau(s)} \ot \Om^+. \end{equation}
\end{lemma}

\begin{proof} The isomorphism is given by a mixed crossing, with its inverse being given by another mixed crossing, thanks to \eqref{eq:mixedR2}. \end{proof}

\begin{cor} \label{cor:EDiagIndecomp} Every object of $\EDiag$ is isomorphic to an object of the form $\Om^k \ot B$, where $B$ is an object of $\Diag$, and $k \in \Z$. \qed \end{cor}
	
Henceforth we omit the tensor product notation from objects in $\EDiag$, merely using concatentation to indicate the monoidal product. We postpone the proofs of the remaining results in this section to the appendix.

\begin{lemma} \label{lem:noblackendo} Any two diagrams built entirely from black strands (and their cups and caps) and which have the same boundary are equal as morphisms in $\EDiag$. As an $R$-bimodule, $\Hom(\Om^k,\Om^k) \cong R_{\tau^k}$. \end{lemma}

In particular, given two iterated tensor products of $\Om^+$ and $\Om^-$ whose total exponent is $k$, these tensor products are canonically isomorphic to each other (and to $\Om^k$). The isomorphism is given by any diagram built entirely from black strands.

\begin{lemma} \label{lem:onlyoneiso} For any object $X$ in $\EDiag$, every degree zero diagram from $X$ to an object of the form $\Om^k B$, which is built entirely from mixed crossings and black cups and caps, is an isomorphism. Moreover, all these diagrams are equal as morphisms in $\EDiag$. \end{lemma}

Thus every object is canonically isomorphic to an object of the form $\Om^k B$ for $B \in \Diag$.

\begin{lemma} \label{lem:pullapart} For $B, C \in \Diag$ and $k, \ell \in \Z$, one has $\Hom(\Om^k B, \Om^\ell C) = 0$ unless $k = \ell$. There is an isomorphism of $R$-bimodules
\begin{equation} \label{eq:homOmB} \Hom(\Om^k B, \Om^k C) \cong \Hom(\Om^k,\Om^k) \ot_R \Hom(B,C). \end{equation}
Moreover, any morphism in $\Hom(\Om^k,\Om^k)$ can be \emph{uncolored}, that is, it is equal to a linear combination of diagrams using only black cups and caps, and polynomials. Similarly, any morphism in $\Hom(B,C)$ can be \emph{unblackened}, that is, it is equal to a linear combination of diagrams using only colored strands and polynomials.
\end{lemma}

Combining these lemmas one can prove the following. Below, $\Kar(\EDiag)$ denotes the Karoubi envelope of the additive, graded envelope of $\EDiag$.

\begin{thm} \label{thm:EDCsize} The inclusion $\Diag \to \EDiag$ is fully faithful. For $B, C$ in $\Diag$ and $k, \ell \in \Z$, one has
\begin{equation} \Hom(\Om^k B, \Om^\ell C) \cong \begin{cases} 0 & k \ne \ell \\ R_{\tau^k} \ot \Hom(B,C) & k = \ell. \end{cases} \end{equation} The Grothendieck group of $\EDiag$ is isomorphic to $\Z \ltimes \Hecke_W$, the semidirect product of $\Hecke_W$ with the group algebra of $\Z$, whose generator is the symbol $[\Om^+]$. The action of $\Z$ on $\Hecke_n$ is determined by  $[\Om^+]b_s[\Om^-] = b_{\tau(s)}$.
\end{thm}

\begin{remark} Given the existence of the functor $\Diag \to \BSBim$, it is a straightforward exercise to extend it to a functor from $\EDiag$ to extended Soergel bimodules, for which $\Om^+ \mapsto R_{\tau}$. \end{remark}

\begin{remark} In the example found in Mackaay-Thiel \cite{MacThi}, the Grothendieck group is $\Hecke_{\extg}$, as evidenced by its Coxeter presentation (see Definition
\ref{defn:Coxpresentext} and following). \end{remark}

\begin{remark} \label{rmk:stillsi} The automorphism $\si$ still acts on the category $\EDiag$, acting on the subcategory $\Diag$ by its usual automorphism, sending $\Om$ to $\Om\inv$, and reversing the
orientation on all black strands. It is straightforward to check that all the relations of $\EDiag$ are preserved by this functor. \end{remark}

%=================
\subsection{Automorphisms of finite order}
\label{subsec-extendedsbimquotient}
%=================

When the automorphism $\tau$ has finite order $n$, then $R_{\tau^n} \cong R$. One might wish to enforce an analogous property for the object $\Om$, namely that $\Om^n \cong
\one$. This is accomplished as follows.

\begin{defn} Fix $n \ge 2$. Let $\EDiag_n$ be the extension of $\EDiag$ defined as follows. One adds a pair of new morphism generators of degree zero, called \emph{$n$-valent vertices}. In the examples below, $n=4$.
\begin{equation} \label{eq:blacknv} \ig{1}{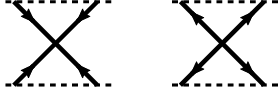} \end{equation}
Then one adds relations. The first relation states that the new morphism generators are cyclic and rotation-invariant.
\begin{subequations} \label{eq:nvalentrelns}
\begin{equation} \label{eq:blackncyclic} {
\labellist
\small\hair 2pt
 \pinlabel {$=$} [ ] at 56 21
\endlabellist
\centering
\ig{1}{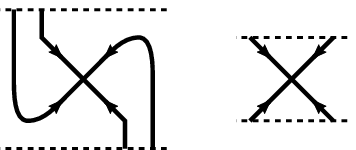}
} \end{equation}
We also impose the orientation-reversal of this relation. Now an isotopy class of diagram unambiguously represents a morphism.

The next relations effectively say that these new generators give inverse isomorphisms between $\Om^n$ and $\one$.
\begin{equation} \label{eq:blackbreakn} {
\labellist
\small\hair 2pt
 \pinlabel {$=$} [ ] at 40 18
\endlabellist
\centering
\ig{1}{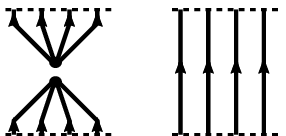}
} \end{equation}
\begin{equation} \label{eq:blackncircle} {
\labellist
\small\hair 2pt
 \pinlabel {$=$} [ ] at 40 18
\endlabellist
\centering
\ig{1}{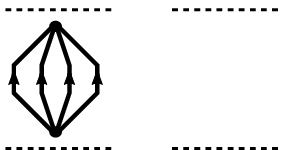}
} \end{equation}
\end{subequations}

The final relation allows the $n$-valent vertex to pull harmlessly through a colored strand.
\begin{equation} \label{eq:pullnvalent} {
\labellist
\small\hair 2pt
 \pinlabel {$=$} [ ] at 40 18
\endlabellist
\centering
\ig{1}{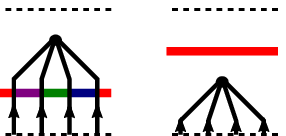}
} \end{equation}
In this picture, red is $s = \tau^4(s)$, blue is $\tau(s)$, green is $\tau^2(s)$, and purple is $\tau^3(s)$.

This ends the definition.
\end{defn}

\begin{lemma} One has $\Om^n \cong \one$ in $\EDiag_n$. \end{lemma}
\begin{proof} Left to the reader. \end{proof}
	
The following results are also proven in the appendix.
	
\begin{lemma} \label{lem:noblackendo2} Any two diagrams built entirely from black strands (and their cups and caps and $n$-valent vertices) and which have the same boundary are equal as morphisms in $\EDiag_n$. As an $R$-bimodule, $\Hom(\Om^k,\Om^k) \cong R_{\tau^k}$. \end{lemma}
	
\begin{lemma} \label{lem:pullapart2} For $B, C \in \Diag$ and $k, \ell \in \Z$, one has $\Hom(\Om^k B, \Om^\ell C) = 0$ unless $k \equiv \ell$ modulo $n$. There is an isomorphism of $R$-bimodules
\begin{equation} \label{eq:homOmB2} \Hom(\Om^k B, \Om^\ell C) \cong \Hom(\Om^k,\Om^\ell) \ot_R \Hom(B,C). \end{equation}
Moreover, any morphism in $\Hom(\Om^k,\Om^\ell)$ can be uncolored (which now permits $n$-valent vertices), and any morphism in $\Hom(B,C)$ can be unblackened.
\end{lemma}

\begin{thm} \label{thm:EDCQsize} The inclusion $\Diag \to \EDiag_n$ is fully faithful. For $B, C$ in $\Diag$ and $k, \ell \in \Z$, one has
\begin{equation} \Hom(\Om^k B, \Om^\ell C) \cong \begin{cases} 0 & k \not\equiv \ell \text{ modulo } n \\ R_{\tau^k} \ot \Hom(B,C) & k \equiv \ell \text{ modulo } n. \end{cases} \end{equation} The Grothendieck group of $\Kar(\EDiag_n)$ is isomorphic to $\Z/n\Z \ltimes \Hecke_W$, the semidirect product of $\Hecke_W$ with the group algebra of $\Z/n\Z$, whose generator is the symbol $[\Om^+]$. The action of $\Z/n\Z$ on $\Hecke_n$ is determined by $[\Om^+]b_s[\Om^-] = b_{\tau(s)}$.
\end{thm}

\begin{remark} When $\tau$ does not have order $n$, these results about morphism spaces in $\EDiag_n$ will fail. Consider multiplication by $f \in R$. Nearby, one can create two black
$n$-valent vertices using \eqref{eq:blackncircle}, pull $f$ through the middle using \eqref{eq:mixedpoly}, and then remove the black $n$-valent vertices with \eqref{eq:blackncircle}.
What remains is multiplication by $\tau^n(f)$. If multiplication by $f$ and multiplication by $\tau^n(f)$ are equal in $\EDiag_n$, but $\tau^n$ is not trivial, then this imposes
unexpected relations. \end{remark}

\begin{remark} To reiterate, when $\tau$ does have order $n$, one can still choose to work with the category $\EDiag$, which effectively ignores the fact that $\tau^n$ is trivial.
\end{remark}

\begin{example} Consider the root realization of $\tilde{A}_{n-1}$. Then $\tau$ has order $n$ on this realization. The Grothendieck group of $\EDiag_n$ will be isomorphic to
$\Hecke_{\exts}$, as evidenced by its Coxeter presentation (see Definition \ref{defn:Coxpresentext} and following). \end{example}

\begin{remark} In our constructions involving Gaitsgory's central complexes, any realization of $\tilde{A}_{n-1}$ with standard Cartan matrix will suffice. Thus it is possible (as in the
previous example) that $\tau$ may have finite order. We perform all our constructions in the category $\EDiag$, and they extend automatically to $\EDiag_n$ when $\tau$ has finite order.
So, while we prefer to work with $\gl_n$ in this paper, we include this definition of $\EDiag_n$ as a bonus for those who prefer $\sl_n$. \end{remark}

%=================
\subsection{Duality}
\label{subsec-duality}
%=================

There is a contravariant, monoidally-covariant \emph{duality} functor $\DM \co \Diag \to \Diag$, which fixes the generating objects $B_s$ for $s \in S$. It reverses grading shifts, in
that $\DM(M(1)) = \DM(M)(-1)$. It acts on diagrams by flipping them upside-down. Let us now extend this functor to $\EDiag$.

\begin{defn} Let $\DM \co \EDiag \to \EDiag$ be the contravariant, monoidally-covariant functor which extends the functor $\DM \co \Diag \to \Diag$, defined as follows. It sends $\Om^+$ to
$\Om^+$, and acts on diagrams by flipping them upside-down, and then reversing the orientation on the black strands. We define a functor $\DM \co \EDiag_n \to \EDiag_n$ in the same way.
\end{defn}

Thus, for example, the identity map of $\Om^+$ is sent to itself. It is easy to confirm that $\DM$ preserves all the relations of $\EDiag$. Note that it would be impossible to flip the
diagrams upside-down without reversing the orientation on the black strands, as the mixed crossings would be sent to mixed crossings with incorrect colorings.

Let us quickly discuss how $\DM$ descends to the Grothendieck group. The Kazhdan-Lusztig bar involution on $\Hecke_W$ is the algebra homomorphism sending $v \mapsto v\inv$, and acting on
the standard basis $\{H_w\}_{w \in W}$ by the map $H_w \mapsto H_{w\inv}\inv$. It fixes the Kazhdan-Lusztig generators $b_s = H_s + v$, so that it can also be defined as the unique
algebra homomorphism sending $v \mapsto v\inv$ and $b_s \mapsto b_s$ for each $s \in S$. Thus, it is clearly categorified by $\DM$ acting on $\Diag$. The point is that the Kazhdan-Lusztig
bar involution should extend to $\Hecke_{\extg}$ by sending $\om \mapsto \om$, as we have already noted in \S\ref{subsec-duality0}.

%=================
\subsection{Notation in affine type $A$}
\label{subsec-typeAnotation}
%=================

For the remainder of this paper, we fix a realization of the Coxeter system $(W,S)$ of type $\tilde{A}_{n-1}$ with a standard affine Cartan matrix, and an automorphism $\tau$ satisfying
$\tau(s_k) = s_{k+1}$. The main example is the standard affine realization, see Definition \ref{defn:standardaffine} and Example \ref{ex:standardaffinetau}. We always let $\delta$ denote
the sum of the simple roots, a $W$-invariant polynomial in $R$.

We let $\Diag_{\ext}$ denote the corresponding category $\EDiag$. We assume that the realization has an automorphism $\si$ with $\si(s_k) = s_{-k}$ and $\si \tau \si = \tau\inv$, so that
$\Diag_{\ext}$ is equipped with an automorphism $\si$ as in Remark \ref{rmk:stillsi}. Restricting this realization to the finite Weyl group $(W_{\fin},S_{\fin})$, we let $\Diag_{\fin}$
denote the corresponding diagrammatic category.

% \begin{notation} Consider a realization of the affine Weyl group $(W_{\aff},S_{\aff})$ in type $\tilde{A}_{n-1}$. Suppose this realization admits an automorphism $\tau$ which, on the affine Dynkin diagram $S_{\aff} = \Z/n\Z$, corresponds to rotation: $\tau(s_i) = s_{i+1}$. Whenever we are in this context, we will adopt the following notation. \begin{itemize}
% 	\item $\delta$ will denote the sum of the simple roots: $\delta = \al_0 + \al_1 + \ldots + \al_{n-1}$.
% 	\item $\Diag_{\extg}$ will denote the corresponding category $\EDiag$ for this data.
% 	\item If $\tau$ has finite order $n$ on the realization, then $\Diag_{\exts}$ will denote the corresponding category $\EDiag_n$.
% 	\item $\Diag_{\fin}$ will denote the diagrammatic category for $(W_{\fin},S_{\fin})$ associated to some (determined by context) subrealization.
% \end{itemize}
% \end{notation}
%
% In particular, the realization of Example \ref{ex:MacThiReal} is our prototypical one where $\tau$ has infinite order, and we call it the \emph{standard affine realization}. Note that $\delta = y$. Three examples where $\tau$ has finite order are: the subrealization of the standard one, spanned by the simple roots; the quotient of the standard one which sets $\delta=0$; the subquotient which both kills $\delta$ and only considers the span of the simple roots.

%%%%%%%%%%%%%%%%%%%%%%%%%
\section{Rouquier complexes and Wakimoto complexes}
\label{sec-wakimoto}
%%%%%%%%%%%%%%%%%%%%%%%%%

% Wakimoto complexes $W(\la)$
% are parametrized by arbitrary weights $\la \in \La_{\wt}$, and categorify the translation lattice in the (extended) affine Hecke algebra. As categorical lifts of standard basis vectors in
% the affine Hecke algebra, they can be realized using the technology of Rouquier complexes, see \S\ref{sec:Wakimoto}. Consequently, Wakimoto complexes are perverse. These complexes are
% still quite mysterious, and studying them (and the category $\Wak$ they generate) in detail is tantamount to categorifying the ``loop presentation'' of the affine Hecke algebra, which is
% still an open problem.

%=================
\subsection{Minimal complexes}
\label{subsec-minimal}
%=================

Let $\AC$ be an additive category. We let $\Ch(\AC)$ denote the additive category of bounded (co)chain complexes in $\AC$, and let $\Hot(\AC)$ denote the bounded homotopy category.
Namely, $\Hot(\AC)$ is the quotient of $\Ch(\AC)$ by nulhomotopic maps. We assume $\AC$ has the Krull-Schmidt property, where every object splits into indecomposable objects. Note that
splittings of objects into direct summands are not canonical in general, so that when we talk about the summands of an object we implicitly choose a splitting.

Given a complex $\MC$, the object $\MC^k$ in homological degree $k$ is a direct sum of indecomposables, as is the object $\MC^{k+1}$ in degree $k+1$, so that the differential $d \co
\MC^k \to \MC^{k+1}$ can be expressed as a matrix of morphisms between indecomposable objects. Any of these morphisms which appear as matrix coefficients we refer to as a \emph{summand
of the differential}. Similarly, for a chain map $f \co \MC \to \NC$ between two complexes, the map $f^k \co \MC^k \to \NC^k$ in homological degree $k$ can be expressed as a matrix of
morphisms, and its matrix coefficients are called \emph{summands of $f$}. To reiterate, taking a summand of a differential or a chain map makes sense only after one has chosen a
splitting of each object $\MC^k$.

We may talk about summands of a differential or a chain map for direct sum decompositions of $\MC^k$, even without the assumption that the summands are indecomposable. When $X \sumset
\MC^k$ and $Y \sumset \MC^{k+1}$, we might also refer to the matrix coefficient $X \to Y$ as \emph{the differential from $X$ to $Y$}. For example, tensor products of complexes $\MC \ot
\NC$ naturally come equipped with direct sum decompositions in each degree, $(\MC \ot \NC)^k \cong \bigoplus \MC^i \ot \NC^{k-i}$, where the summands are not necessarily indecomposable.

When the differential from $X \sumset \MC^k$ to $Y \sumset \MC^{k+1}$ is an isomorphism, there is a process called \emph{Gaussian elimination} which replaces the complex $\MC$ with a
homotopy equivalent complex $\MC'$ in which the summands $X$ and $Y$ have been removed or \emph{eliminated}. The point of Gaussian elimination is that, while it may not appear so at
first, one can alter the direct sum decompositions of $\MC^k$ and $\MC^{k+1}$ so that the two term complex $X \to Y$ is a direct summand of $\MC$, and this two term complex is
contractible.

A complex for which no summand of the differential (between indecomposable objects) is an isomorphism is called a \emph{minimal complex}. By repeated Gaussian elimination, one can prove
that any complex is homotopy equivalent to a minimal complex. In a Krull-Schmidt category this minimal complex is unique up to isomorphism, so it does not depend on the order in which
various contractible summands are eliminated. We will often refer to objects of the homotopy category $\Hot(\AC)$ abusively as complexes, but we usually have in mind the minimal complex
as our favorite representative of its homotopy equivalence class.

For more background on minimal complexes see \cite[\S 6.1]{EWHodge}. We recall the specifics of Gaussian elimination in \S\ref{sec-pseudo}, where we also discuss Gaussian elimination for
pseudocomplexes.

%=================
\subsection{Rouquier complexes and Rouquier canonicity}
\label{subsec-rouquier}
%=================

For any Coxeter system $(W,S)$, Rouquier \cite{RouqBraid-pp} defined a strict action of the braid group of $W$ on $\Hot(\SBim)$, which can be easily transported to the diagrammatic setting
$\Hot(\Diag)$.

Rouquier defined bounded complexes $F_s$ and $F_s\inv$, for each $s \in S$, as follows.
\begin{equation} F_s = \left( \begin{diagram}
\un{B}_s & \rTo^{\finaldotred} & R(1) \end{diagram} \right) \end{equation}
\begin{equation} F_s\inv = \left( \begin{diagram}
R(-1) & \rTo^{\startdotred} & \un{B}_s \end{diagram} \right) \end{equation}
The underline indicates homological degree zero.
He constructed homotopy equivalences
\begin{equation} \label{eq:rouqisom1} F_s \ot F_s\inv \cong \one \simeq F_s\inv \ot F_s \end{equation} for each $s \in S$, and \begin{equation} \label{eq:rouqisom2} F_s \ot F_t \ot \cdots \simeq F_t \ot F_s \ot \cdots,\end{equation} (with $m_{st}$ tensor factors on each side) for each $s \ne t \in S$ with $m_{st} < \infty$. In type $A$, one can see these isomorphisms written diagrammatically in \cite{EKra}.

Let $\bb$ be a braid word, a word in the generators $\{f_s^\pm\}$. Taking the corresponding tensor product of the complexes $F_s$ and $F_s\inv$, one obtains a complex $F_{\bb}$. The
homotopy equivalences above imply that when $\bb_1$ and $\bb_2$ yield the same element in the braid group, then $F_{\bb_1}$ and $F_{\bb_2}$ are isomorphic in $\Hot(\Diag)$.

However, Rouquier \cite{RouqBraid-pp} also proved a much stronger result (which is what makes this categorical action of the braid group \emph{strict}). We give our own version of the
result here.

\begin{prop} \label{prop:RouqCanon} (\emph{Rouquier Canonicity}) For each pair $\bb_1$, $\bb_2$ of braid words representing the same braid $\be$,
there exists a homotopy equivalence \[ \psi_{\bb_1, \bb_2} \co F_{\bb_1} \to F_{\bb_2}. \] These homotopy equivalences are transitive: \begin{subequations} \label{eq:rouqconsistent} \begin{equation} \label{eq:rouqtrans} \psi_{\bb_1,\bb_3} =
\psi_{\bb_2,\bb_3} \circ \psi_{\bb_1,\bb_2}. \end{equation} Moreover, if $s$ is any simple reflection, then \begin{equation} \label{eq:rouqthens} \psi_{\bb_1 f_s, \bb_2 f_s} = \psi_{\bb_1,\bb_2} \ot
\id_{F_s}, \end{equation} \begin{equation} \label{eq:rouqthensinv} \psi_{\bb_1 f_s\inv, \bb_2 f_s\inv} = \psi_{\bb_1,\bb_2} \ot
\id_{F_s\inv}. \end{equation} \end{subequations} The above are equalities in the homotopy category (i.e. the chain maps are homotopic). \end{prop}

We now discuss how one can pick out Rouquier's homotopy equivalence $\psi_{\bb_1,\bb_2}$. Let $F$ be an invertible complex in $\Hot(\Diag)$. Then $\Hom(F,F) \cong \Hom(\one,\one)$, by
tensoring with $F\inv$ and using any isomorphism $F \ot F\inv \cong \one$. The only degree zero endomorphism of $\one$ is the identity, up to scalar. Hence, Rouquier's homotopy
equivalence is already determined up to scalar (up to homotopy).

Recall that $R = \one = B_1$. It is a simple consequence of Kazhdan-Lusztig cell theory that $\one$ is never a direct summand of $B_w \ot
B_x$, unless $w = x = 1$. (Of course, $\one \ot \one \cong \one$.) Thus, $\one$ is not a direct summand of any Bott-Samelson bimodule except for $\one$ itself. Said another way, the
category $\Diag$ has a maximal two-sided monoidal ideal $\MaxIdeal$ which contains $\id_{B_x}$ for all $x \ne 1$, but not $\id_{\one}$. Taking the quotient by this maximal ideal, one
recovers the category of graded vector spaces, where $\one$ is sent to $\C$.

Consequently, the complex $F_{\bb}$ of any braid word has a unique direct summand of the form $\one$ (up to shift), as does its minimal complex. Keeping track of the shift, this summand
has the form $\one(k)[-k]$, where $[1]$ denotes the homological shift and $(1)$ the grading shift, and $k$ is the exponent of the braid. Any isomorphism $F_{\bb_1} \to F_{\bb_2}$ must
induce an isomorphism $\one(k)[-k] \to \one(k)[-k]$ between the unique copies of $\one$ in each complex. This is because the isomorphism must descend to an isomorphism modulo $\MaxIdeal$.

Summing up these last few paragraphs, we see that a homotopy equivalence $F_{\bb_1} \to F_{\bb_2}$ is already determined up to scalar, and this scalar can be interpreted as the
coefficient of the identity map $\id_{\one}$ between the unique copies of $\one$ in the chain map between these complexes. Rouquier's homotopy equivalence is uniquely determined by the
fact that it is a homotopy equivalence, and that this coefficient of the identity is $1$. This is true for the maps in \eqref{eq:rouqisom1} and \eqref{eq:rouqisom2}, and it is preserved
by the consistency requirements \eqref{eq:rouqconsistent}, so it is true in general.

\begin{remark} If desired, this argument could be modified into a reproof of Rouquier canonicity. We mention this because Rouquier canonicity is not usually stated in the literature with all the conditions mentioned above. \end{remark}

%=================
\subsection{Perversity}
\label{subsec-perverse}
%=================

Whenever $w = s_1 s_2 \cdots s_d$ is a reduced expression in $W$, $B_{s_1} B_{s_2} \cdots B_{s_d}$ has a unique indecomposable direct summand in $\Kar(\Diag)$ which has not appeared in
any shorter tensor products. This summand is called $B_w$. In \cite{EWGr4sb} it is proven that the indecomposable objects in $\Kar(\Diag)$ have the form $B_w(k)$, where $w \in W$ and $k
\in \Z$. Moreover, the objects $B_w$ are special in that they are \emph{self-dual} under a duality functor $\DM$ (in contrast, the dual of $B_w(k)$ is $B_w(-k)$).

Henceforth we assume that our realization is obtained from the geometric realization over $\R$ by base change. In \cite{EWHodge} the Soergel conjecture is proven in this case, which gives a very significant amount of control over morphism spaces between indecomposable objects. Let $\Hom^i(B_w,B_x)$ denote the space of homogeneous morphisms of degree $i$, or equivalently, the space $\Hom^0(B_w,B_x(i))$. The Soergel conjecture implies that
\begin{equation} \Hom^i(B_w, B_x) := \Hom(B_w,B_x(i)) = \begin{cases}
	0 & \text{when } i<0,\\
	0 & \text{when } i=0 \text{ and } w \ne x,\\
	\Bbbk \cdot \id_{B_w} & \text{when } i=0 \text{ and } w=x. \\ \end{cases} \end{equation}

Let $\Hot^{\ge 0}(\Diag)$ consists of all complexes homotopy equivalent to a complex where, for each $k \in \Z$, every indecomposable summand in degree $k$ has the form $B_w(\ell)$ for
$\ell \le k$.  Equivalently, one needs the minimal complex to satisfy this condition. Similarly, $\Hot^{\le 0}(\Diag)$ requires $\ell \ge k$. The intersection $\Hot^{\le
0}(\Diag) \cap \Hot^{\ge 0}(\Diag)$ consists of objects whose minimal complexes are ``supported on the diagonal:'' the summand $B_x(k)$ only appears in homological degree $k$. We call
such an object \emph{perverse}, and a complex supported on the diagonal a \emph{perverse complex}. A perverse complex is automatically minimal. For a perverse complex, all the
differentials are direct sums of maps in $\Hom^1(B_w,B_x)$ for various $w,x \in W$. There are no possible homotopies between perverse complexes, because a homotopy would be built out of
elements of $\Hom^{-1}(B_x,B_w)$, and this space is always zero.

\begin{thm} \label{thm:pervprops} The following results hold. \begin{enumerate}
	\item $F_x$ is perverse whenever $x$ is a positive lift of an element of $W$.
	\item Similarly, $F_y$ is perverse whenever $y$ is a negative lift of an element of $W$.
	\item For $s \in S$, tensoring with $F_s$ on the left or right preserves $\Hot^{\ge 0}(\Diag)$.
	\item For $s \in S$, tensoring with $F_s\inv$ on the left or right preserves $\Hot^{\le 0}(\Diag)$. \end{enumerate} \end{thm}

\begin{proof} Property (1) is \cite[Theorem 6.9]{EWHodge} and property (3) follows immediately from \cite[Lemma 6.5]{EWHodge}. Properties (2) and (4) follow from (1) and (3) by applying
the duality functor $\DM$, see \S\ref{subsec-duality2}. \end{proof}

As a consequence we have the following result, whose slick proof was explained to me by Geordie Williamson.

\begin{cor} \label{cor:posnegperv} Let $F_x F_y\inv$ be a tensor product of Rouquier complexes, where $x$ and $y$ are both positive lifts of elements of $W$. Then $F_x F_y\inv$ is perverse. \end{cor}

\begin{proof} We know that $F_x$ is perverse. Since $F_y\inv$ is a tensor product of various $F_s\inv$, we see that tensoring with $F_y\inv$ preserves $\Hot^{\le 0}(\Diag)$. Hence $F_x F_y\inv \in \Hot^{\le 0}(\Diag)$.
	
Similarly, $F_y\inv$ is perverse. Since $F_x$ is a tensor product of various $F_s$, we see that tensoring with $F_x$ preserves $\Hot^{\ge 0}(\Diag)$. Hence $F_x F_y\inv \in \Hot^{\ge 0}(\Diag)$.

Thus $F_x F_y\inv \in \Hot^{\ge 0}(\Diag) \cap \Hot^{\le 0}(\Diag)$, so $F_x F_y\inv$ is perverse. \end{proof}

%=================
\subsection{Rouquier complexes for extended Soergel bimodules}
\label{subsec-rouquierextend}
%=================

Now we fix a Coxeter system $(W,S)$, and a realization with an automorphism $\tau$. Recall the diagrammatic category $\EDiag$ from \S\ref{subsec-extendeddiag}. We claim that Rouquier complexes extend to a strict action of the extended braid group $\Br_{\ext}$ on $\Hot(\EDiag)$. First let us define the extended braid group.

\begin{defn} The \emph{extended braid group} of the data $(W,S,\tau)$ is the group $\Br_{\ext} = \Br_W \rtimes \Z$. This extends the usual braid group of $(W,S)$, adding an extra invertible element $\om$ with the relation
\begin{equation} \om f_i \om\inv = f_{\tau(i)} \end{equation}
for each $i \in S$.

When $\tau$ has finite order $n$, then $\om^n$ is central in $\Br_{\ext}$, and we can also define the quotient $\Br_{\ext,n} = \Br_W \rtimes (\Z/n\Z) \cong \Br_{\ext} / \langle \om^n
\rangle$. \end{defn}

\begin{defn} To the element $\om \in \Br_{\ext}$, we associate the complex $F_{\om} = \un{\Om}^+$, supported in degree zero. Similarly, to $\om\inv$ we associate $F_{\om\inv} =
\un{\Om}^-$. To any word in the generators of $\Br_{\ext}$ (an \emph{extended braid word}) we associate the corresponding tensor product of $F_{\om}$, $F_{\om\inv}$, $F_s$ and $F_s\inv$.
This is an object of $\Hot(\EDiag)$, which we call an \emph{extended Rouquier complex}. \end{defn}

Generalizing Rouquier's results to the extended context is rather easy, given Theorem \ref{thm:EDCsize}.

\begin{prop} If $x = \om^k \be$ and $y = \om^\ell \ga$ are elements of $\Br_{\extg}$, then there is an isomorphism of $R$-bimodules \begin{equation} \Hom_{\Hot(\EDiag)}(F_x,F_y) \cong
\begin{cases} 0 & \text{ if } k \ne \ell,\\ R_{\tau^k} \ot_R \Hom_{\Hot(\Diag)}(F_\be,F_\ga) & \text{ if } k = \ell. \end{cases} \end{equation} \end{prop}

\begin{proof} By Theorem \ref{thm:EDCsize}, there are no morphisms between any of the chain objects in $F_x$ and those in $F_y$ unless $k = \ell$. In this case, since $\Om$ is invertible,
tensoring with $\Om^{-k}$ will give the desired isomorphism. \end{proof}

\begin{cor} The extended Rouquier complexes give a strict action of $\Br_{\ext}$ on $\Hot(\EDiag)$. \end{cor}

Writing out a precise proof is somewhat annoying, but the idea is simple.

\begin{proof} (Sketch) Given two extended braid words $\bb_1$ and $\bb_2$ for the same extended braid, we define the canonical homotopy equivalence $\psi \co F_{\bb_1} \to F_{\bb_2}$ to be the unique homotopy equivalence which induces the identity map on the unique direct summand of the form $\Om^k$. Note that the category $\EDiag$ has a maximal ideal $\MaxIdeal$ which contains $\id_{B_w}$ for all $w \ne \om^k$. An isomorphism is determined by its image modulo the maximal ideal, by the same arguments as in \S\ref{subsec-rouquier}. \end{proof}

By Corollary \ref{cor:EDiagIndecomp}, the indecomposable objects in $\Kar(\EDiag)$ have the form $B_w (k)$ for $k \in \Z$ and $w \in W_{\ext}$. One can define perverse complexes just as
in the previous section.

\begin{prop} \label{thm:pervprops2} The following results hold. \begin{enumerate}
	\item $F_x$ is perverse whenever $x \in \Br_{\ext}$ is a positive lift of an element of $W_{\ext}$.
	\item Similarly, $F_y$ is perverse whenever $y \in \Br_{\ext}$ is a negative lift of an element of $W_{\ext}$.
	\item For $s \in S$, tensoring with $F_s$ on the left or right preserves $\Hot^{\ge 0}(\EDiag)$.
	\item For $s \in S$, tensoring with $F_s\inv$ on the left or right preserves $\Hot^{\le 0}(\EDiag)$. \end{enumerate}
	\item Let $x,y \in \Br_{\ext}$ be positive lifts of elements of $W_{\ext}$. Then $F_x F_y\inv$ is perverse.
 \end{prop}

\begin{proof} Since tensoring with $\Om^k$ preserves perverse complexes, the first two properties follow from the analogous properties in Theorem \ref{thm:pervprops}. Since $F_s \Om^k \cong \Om^k F_{\tau^{-k}(s)}$, the next two properties follow from the analogous properties in Theorem \ref{thm:pervprops}. The final property is proven by the same argument as Corollary \ref{cor:posnegperv}. \end{proof}

Now consider the case when $\tau$ has finite order $n$. We define Rouquier complexes for words in $\Br_{\ext,n}$ to be the images of the extended Rouquier complexes above, under the functor $\EDiag \to \EDiag_n$. It is equally straightforward to use the results of \S\ref{subsec-extendedsbimquotient} to deduce analogous properties.

\begin{prop} If $x = \om^k \be$ and $y = \om^\ell \ga$ are elements of $\Br_{\ext,n}$, then there is an isomorphism of $R$-bimodules \begin{equation} \Hom_{\Hot(\EDiag)}(F_x,F_y) \cong
\begin{cases} 0 & \text{ if } k \not\equiv \ell \text{ modulo } n,\\ R_{\tau^k} \ot_R \Hom_{\Hot(\Diag)}(F_\be,F_\ga) & \text{ if } k \equiv \ell \text{ modulo } n. \end{cases} \end{equation} \end{prop}

\begin{cor} The extended Rouquier complexes give a strict action of $\Br_{\ext,n}$ on $\Hot(\EDiag_n)$. \end{cor}
	
%=================
\subsection{Duality}
\label{subsec-duality2}
%=================

The duality functor $\DM$ from \S\ref{subsec-duality} extends to the homotopy category $\Hot(\EDiag)$. Since it is contravariant, it reverses the order of the differentials, and hence it
must also reverse the homological degree. It is easy to see that $\DM(F_s) = F_s\inv$ and $\DM(F_{\om}) = F_{\om}$. Hence $\DM(F_w) = F_{w\inv}\inv$, categorifying the action of the bar
involution $H_s \mapsto H_s\inv$ and $\om \mapsto \om$. In particular, the dual of a Rouquier complex is a Rouquier complex, and the dual of a perverse complex is a perverse complex.

%=================
\subsection{Wakimoto complexes}
\label{subsec-wakimoto}
%=================

Let us work in affine type $A$, with the conventions of \S\ref{sec:extendedbraid}. In \S\ref{subsec-translattice} we defined a collection of elements $w_\la \in \Br_{\extg}$ for each
$\la \in \La_{\gl_n}$. The (minimal complexes of the) corresponding Rouquier complexes will be denoted $W(\la)$, and called the \emph{Wakimoto complexes}. Examples will be given in the
next section.

In Proposition \ref{prop:dominantpositive} it was shown that $w_\la$ is a positive lift of an element of $W_{\extg}$ whenever $\la$ is dominant. Similarly, it is a negative lift when
$\la$ is antidominant. For any $\la$ at all, $w_\la$ is the product of a (commuting) positive lift and a negative lift. Consequently, by Corollary \ref{cor:posnegperv}, we have the
following extremely useful result.

\begin{cor} The complex $W(\la)$ is perverse, for all $\la \in \La_{\gl_n}$. \end{cor}

\begin{remark} One should think about Wakimoto complexes as the categorification of the translation lattice in the cylindrical Hecke algebra. The study of Wakimoto complexes is the first
step towards a ``loop presentation'' of the extended or affine Hecke category. However, computing the morphisms in $\Hot(\EDiag)$ between (shifts of) Wakimoto complexes is quite difficult. \end{remark}

Let us be precise about one additional feature of Wakimoto complexes: they commute. For $\la, \mu \in \La_{\gl_n}$, we know that $W(\la) \ot W(\mu) \cong W(\mu) \ot W(\la)$, because both
are isomorphic to $W(\la + \mu)$. However, Rouquier canonicity fixes a particular homotopy equivalence \begin{equation} \psi_{\la,\mu} \co W(\la) \ot W(\mu) \to W(\mu) \ot W(\la)
\end{equation} satisfying several desirable properties. It is easy to show, using Proposition \ref{prop:RouqCanon} and the discussion following it, that these maps make the full
subcategory of $\Hot(\EDiag)$ with objects $W(\la)$ into a symmetric (braided) monoidal category. In this sense, the complexes $W(\la)$ and $W(\mu)$ categorically commute.

We let $W(\la)^*$ denote $\DM(W(\la))$. These are the \emph{dual Wakimoto complexes}. They satisfy many of the same properties as the ordinary Wakimoto complexes.

Many examples of Wakimoto and dual Wakimoto complexes can be found in \S\ref{sec-234}.

%=================
\subsection{Sign conventions for Rouquier complexes}
\label{subsec-signconvention}
%=================

When tensoring two complexes $\PC$ and $\QC$, the typical convention is that the differential on $\PC^i \ot \QC^j$ is \begin{equation} \label{eq:stdsignconv} d_{\PC} \ot 1 + (-1)^i 1
\ot d_{\QC}.\end{equation} Of course, there are other choices of sign which isomorphic tensor-product complexes, such as \begin{equation}
\label{eq:altsignconv} d_{\PC} \ot 1 + (-1)^{i+1} 1 \ot d_{\QC}.\end{equation} One isomorphism between the complex $\PC \ot \QC$ defined with \eqref{eq:stdsignconv} and the one
defined with \eqref{eq:altsignconv} can be given by multiplying $\PC^i \ot \QC^j$ by $(-1)^{j+1}$.

Messing with the standard sign convention is a bad idea. Suppose one used \eqref{eq:altsignconv} consistently for all tensor products. Tensoring with the monoidal identity will introduce
a sign on all the differentials. Consequently, the unit and associator structures for the monoidal structure on complexes will need to keep track of signs, and this is a nightmare, before
one even discusses what happens to morphisms.

That said, in specific examples sometimes \eqref{eq:altsignconv} yields a complex which is easier to use. So while we will not use a non-standard sign convention, we will use non-standard
representatives of certain tensor products to make life easier.

\begin{defn} We say that $\PC \in \Ch(\Diag)$ is an \emph{even BS complex} if in homological degree $k$, $\PC^k$ is a direct sum of Bott-Samelson objects with length agreeing with
$k$ modulo $2$. We say that $\PC$ is an \emph{odd BS complex} if in homological degree $k$, $\PC^k$ is a direct sum of Bott-Samelson objects with length agreeing with $k+1$ modulo
$2$. \end{defn}

For example, the Rouquier complexes $F_s$ and $F_s\inv$ are odd BS complexes, while the monoidal identity is an even BS complex. A tensor product of two odd BS complexes is an even BS
complex, etcetera. The Gaitsgory central complex for the standard representation will be an (odd or even) BS complex, as are other complexes filtered by Rouquier complexes.

\begin{notation} Let $\QC$ be an arbitrary complex in $\Ch(\Diag)$. When $\PC$ is an even BS complex, we write $\PC \sot \QC$ for the standard tensor product, using the ordinary
sign convention \eqref{eq:stdsignconv}. When $\PC$ is an odd BS complex, we write $\PC \sot \QC$ for the complex which is isomorphic to $\PC \ot \QC$, but
whose differentials are given by \eqref{eq:altsignconv} instead. \end{notation}

\begin{lemma} If $\PC$ and $\PC'$ are (odd or even) BS complexes, then $(\PC \sot \PC) \sot \QC$ and $\PC \sot (\PC \sot \QC)$ are the same. \end{lemma}

\begin{proof} This is an easy sign chase. \end{proof}

Let us explore these sign conventions as they apply to Rouquier complexes. In both $F_s$ and $F_s\inv$, the object appearing in even homological degree is $B_s$, while the object in odd
degree is $R$. Thus each homological degree in a tensor product $F_{s_1}^{\pm 1} F_{s_2}^{\pm 1} \cdots F_{s_d}^{\pm 1}$ has a natural direct sum decomposition with Bott-Samelson objects
as summands, exemplified by $B_{s_1} \ot R \ot R \ot B_{s_4} \ot \cdots$, where each $s_k$ is replaced either by $R$ or $B_{s_k}$. Each summand of the differential is a signed dot, like
$\pm \linered \startdotblue \linepurple$ or $\pm \linered \finaldotblue \linepurple$. The question is: what is the sign? Suppose the dot appears on the $k$-th index. In the ordinary
tensor product complex, then the sign appearing is $(-1)^\ell$ where $\ell$ is the number of $R$'s which appear in indices $j$ smaller than $k$. In the isomorphic complex $F_{s_1}^{\pm 1}
\sot \cdots \sot F_{s_d}^{\pm 1}$ the sign appearing is $(-1)^m$ where $m$ is the number of $B_{s_j}$ which appear in indices $j$ smaller than $k$.

Let us illustrate with an example, where $d=4$ and the colors $s_1, s_2, s_3, s_4$ are red, green, blue, and purple respectively. Let us examine two summands of the differential:
\begin{equation} \label{eq:k3whatsign}\pm \linered \startdotblue \linepurple \co B_{s_1} \ot R \ot R \ot B_{s_4} \to B_{s_1} \ot R \ot B_{s_3} \ot B_{s_4}, \quad k=3, \end{equation}
\begin{equation} \label{eq:k2whatsign}\pm \linered \startdotgreen \linepurple \co B_{s_1} \ot R \ot R \ot B_{s_4} \to B_{s_1} \ot B_{s_2} \ot R \ot B_{s_4}, \quad k=2. \end{equation}
Using the standard convention, the sign on \eqref{eq:k3whatsign} is $-1$ because $\ell=1$, and the sign on \eqref{eq:k2whatsign} is $+1$ because $\ell=0$. What we are measuring is
invisible in the diagram for the differential; it is the number of copies of the monoidal identity (not drawn) which appear to the left of the dot. Meanwhile, using the alternate sign
convention, the sign on both \eqref{eq:k3whatsign} and \eqref{eq:k2whatsign} is $-1$ because $m=1$. This is measuring the visibly apparent statistic that one line (i.e. identity map of
some $B_s$) appears to the left of the dot. To determine the sign from the standard convention one needs to keep track of extra data: the original sequence $F_{s_1}^{\pm 1} F_{s_2}^{\pm
1} \cdots F_{s_d}^{\pm 1}$, the value of $k$, etcetera. To determine the sign from the alternate convention requires no memory: just count the number of lines to the left of the dot.

For example, here is the complex $F_s F_t F_s$, with red representing $s$ and blue representing $t$, drawn using \eqref{eq:altsignconv}.
\begin{equation} F_s F_t F_s = \left(
\begin{tikzpicture}
\node (a) at (0,0) {$B_s B_t B_s$};
\node (b1) at (4,2) {$B_s B_t(1)$};
\node (b2) at (4,0) {$B_s B_s(1)$};
\node (b3) at (4,-2) {$B_t B_s(1)$};
\node (c1) at (8,4) {$B_s(2)$};
\node (c2) at (8,0) {$B_t(2)$};
\node (c3) at (8,-4) {$B_s(2)$};
\node (d) at (12,0) {$R(3)$};
\path 
	(a) edge node[descr] {$\linered\lineblue\finaldotred$} (b1)
	(a) edge node[descr] {$-\linered\finaldotblue\linered$} (b2)
	(a) edge node[descr] {$\finaldotred\lineblue\linered$} (b3)
	(b1) edge node[descr] {$-\linered\finaldotblue$} (c1)
	(b1) edge node[descr] {$\finaldotred\lineblue$} (c2)
	(b2) edge node[descr] {$-\linered\finaldotred$} (c1)
	(b2) edge node[descr] {$\finaldotred\linered$} (c3)
	(b3) edge node[descr] {$-\lineblue\finaldotred$} (c2)
	(b3) edge node[descr] {$\finaldotblue\linered$} (c3)
	(c1) edge node[descr] {$\finaldotred$} (d)
	(c2) edge node[descr] {$\finaldotblue$} (d)
	(c3) edge node[descr] {$\finaldotred$} (d);
\end{tikzpicture} \right) \end{equation}
In all cases, the sign appearing is the number of strands to the left of the dot. This is the \emph{strand-counting sign rule for Rouquier complexes}

\begin{notation} Henceforth we change our notation for Rouquier complexes as follows. For any braid word $\bb = s_1^{\pm 1} \cdots s_d^{\pm 1}$, we let $F_{\bb}$ denote the complex
$F_{s_1}^{\pm 1} \sot \cdots \sot F_{s_d}^{\pm 1}$. When defined this way, the sign appearing on a dot in a Rouquier complex is determined by the number of strands to the left of the dot.
We identify this complex with the true tensor product $F_{s_1}^{\pm 1} \ot \cdots \ot F_{s_d}^{\pm 1}$ via the (unique) isomorphism which acts as the identity on the unique copy of
$\one$. \end{notation}

Thus we can transfer Rouquier's canonical homotopy equivalences to this setting, giving maps $\psi \co F_{\bb_1} \to F_{\bb_2}$ whenever $\bb_1$ and $\bb_2$ are equal in the braid group,
which satisfy the same compatibilities as in Proposition \ref{prop:RouqCanon}.

To reiterate, we will not be doing anything silly like using $\sot$ as a monoidal structure. We are merely choosing more convenient representatives for various tensor products. We will
always write $\sot$ when we use it (with the caveat that it is built-in to the notation for Rouquier complexes), and just $\ot$ or concatenation for the ordinary tensor product. The
advantage of the notation $\sot$ is that $\PC \sot \QC$ is automatically a complex isomorphic to the normal tensor product, while some other random sprinkling of signs on differentials
might not be.

\begin{remark} Here is one way of interpreting $\sot$ using homological algebra. Let $\gamma = \one[1]$. Then $\gamma$ lives in the Drinfeld center of the category, but its structure map
$(-) \ot \gamma \to \gamma \ot (-)$ is non-trivial, and this leads to most of the sign issues in homological algebra. Note that when people write $\QC[1]$ what they mean is the tensor
product $\gamma \QC$ (tensor product on the left, not the right). Note that $F_s \gamma$ is a complex where $R(1)$ appears in even homological degree, and $B_s$ in odd. Thus $F_s \gamma
\ot (-)$ behaves like one might desire: it puts signs on differentials when $B_s$ is to the left, and no sign when the monoidal identity is to the left. But one must also remove the shift
in the end. Ultimately, one has $F_s \sot \QC = F_s \gamma \QC \gamma\inv$, so that $F_{s_1} \sot \cdots \sot F_{s_d} \sot \QC = (F_{s_1} \gamma) \cdots (F_{s_d} \gamma) \QC \gamma^{-d}$.
Other shifts appear when dealing with $F_s\inv$. \end{remark}

\begin{remark} The reader who prefers the standard sign convention for Rouquier complexes is encouraged to use the notation from \cite{EKra}, where the authors use small circles to keep
track of the extra copies of the monoidal identity. This notation becomes cumbersome though, and we have not developed its analog for Gaitsgory central complexes. \end{remark}

%=============================
\section{Pseudocomplexes and Gaussian elimination}
\label{sec-pseudo}
%=============================

We review the homological algebra of pseudocomplexes, as found in \cite[Chapter 4]{MackyMonodromy}. These are like complexes, except that $d^2$ is nonzero in some measurable way.
As a consequence the pseudocomplex can be equipped with a homological degree 2 endomorphism which measures the failure of $d^2$ to be zero. We leave many basic facts as exercises to the reader, but outline a clear path through the material.

\subsection{Setup}

Let us describe the setting which is relevant to our paper.

Choose a realization $(\hg, \hg^*, \{\al_s\},\{\al_s^\vee\})$ of the Coxeter system $(W,S)$ of the affine Weyl group in type $A$, or its extended affine Weyl group (see Definition
\ref{defn:realizationextended}). We assume that its Cartan matrix is the standard Cartan matrix for $\tilde{A}_n$, and that the sum $\delta = \sum_{s \in S} \al_s$ of the simple roots is
nonzero. Let $R_1$ be the polynomial ring associated to this realization, and $\Diag_1$ its diagrammatic Hecke category.

Now consider another realization, where $\hg$ is replaced by the kernel of $\delta$, and $\delta$ is sent to zero in $\hg^*$. Let $R_2$ be the polynomial ring of this realization, and $\Diag_2$ its diagrammatic Hecke category. Then $R_2$ is the quotient of $R_1$ by $\delta$.

For $i = 1, 2$, the morphism spaces in $\Diag_i$ are $R_i$-bimodules which are free as left or right $R_i$-modules. As a consequence, multiplication by $\delta$ on the right is a nonzero
divisor on any morphism space in $\Diag_1$. Note that $\delta$ is $W$-invariant, so it acts the same on the right and left; this fact will not play a role until
\S\ref{subsec:pseudomonoidal}.

One has a functor $\killit \co \Diag_1 \to \Diag_2$ induced by the morphism of realizations, which kills the action of $\delta$ on the right. There is a natural isomorphism \begin{equation}
\HOM_{\Diag_1}(M,N) \ot_{R_1} R_2 \cong \HOM_{\Diag_2}(\killit M,\killit N). \end{equation} In this context, $\HOM$ denotes the graded vector space spanned by homogeneous morphisms of any
degree.

We will discuss pseudocomplexes relative to $\delta$, which is a much less general setting than found in \cite{MackyMonodromy}. However, this restrictive setting makes the theory
considerably easier to describe. Geometrically, we are considering the monodromy only from a $\C^*$ action rather than a larger torus.

\subsection{Basic definitions}

\begin{defn} A \emph{pseudocomplex} in $\Diag_1$  is the data $X = (X^i, d^i)_{i \in \Z}$ of objects $X^i$ in $\Diag_1$, and morphisms $d^i \co X^i \to X^{i+1}$ of degree zero, such that $d^{i+1} \circ d^i \in \HOM(X^i,X^{i+2}) \delta$ for all $i \in \Z$. \end{defn}

That is, a pseudo complex is almost a complex, except that instead of being zero, $d^2$ is a (right) multiple of $\delta$. The consequences of our setup above are: \begin{itemize}
\item There is a unique element $\mu^i \in \Hom(X^i,X^{i+2}(-2))$ such that $d^{i+1} \circ d^i = \mu^i \delta$ for all $i \in \Z$.
\item After applying the functor $\killit$, one obtains a genuine complex in $\Diag_2$.
\end{itemize}

\begin{remark} One can easily generalize the theory of pseudocomplexes to any operator $\delta$ in any linear category, which acts as a nonzero divisor on every $\Hom$ space. One may not necessarily have a functor which ``kills'' the action of $\delta$, but if one does, one can generalize the corresponding results about complexes in $\Diag_2$ as well. In general settings, there is no need to work in a graded category. \end{remark}

For pseudocomplexes $X$ and $Y$, we write $\Hom(X,Y[k])$ to denote the space of all linear maps $X \to Y[k]$, that is, the product \begin{equation} \Hom(X,Y[k]) = \prod_{n \in \Z}
\Hom(X^n,Y^{n+k}). \end{equation} Similarly, we write $\HOM(X,Y[k])$ to denote the direct sum $\bigoplus_{m \in \Z} \Hom(X,Y[k](m))$.

By the above discussion, any pseudocomplex $X$ is equipped with a map $\mu = \mu_X \in \Hom(X, X(-2)[2])$, such that $d^2 = \mu \delta$, which is called the \emph{monodromy} of $X$.

\begin{defn} A pseudochain map $f \co X \to Y$ of pseudocomplexes is a collection $f^i \co X^i \to Y^i$ of morphisms such that $d^i_Y f^i - f^{i+1}d^i_X \in \Hom(X^i,Y^{i+1}) \delta$ for all $i \in \Z$. \end{defn}

Again, a pseudochain map is like a chain map, but modulo $\delta$, and descends to a genuine chain map after applying the functor $\killit$. Note again that there is a unique map $\nu_f
\in \Hom(X,Y[1])$ for which $d^i_Y f^i - f^{i+1}d^i_X = \nu_f^i \delta$ for all $i \in \Z$.

\begin{exercise} Pseudocomplexes form a category $\Ch_\delta(\Diag_1)$, where the morphisms are pseudochain maps. \end{exercise}

\begin{exercise} Taking any linear map whatsoever in $\Hom(X,Y)$ and multiplying by $\delta$ will give a pseudochain map. \end{exercise}
	
\begin{exercise} Any chain map in $\Hom(\killit X, \killit Y)$ lifts to a pseudochain map in $\Hom(X,Y)$. \end{exercise}

\begin{defn} A \emph{(pseudo)homotopy} of pseudocomplexes is a morphism $h \in \Hom(X,Y[-1])$, that is, a map $h^i \in \Hom(X^i,Y^{i-1})$ for each $i \in \Z$. Given a homotopy $h$, the
corresponding \emph{nulhomotopic map} is $dh + hd \in \Hom(X,Y)$. \end{defn}

\begin{exercise} Verify that any nulhomotopic map is a pseudochain map. \end{exercise}

At first glance, homotopy looks exactly the same as for ordinary complexes. However, the set of nulhomotopic maps does not form an ideal in $\Ch_\delta(\Diag_1)$. If $f$ is a pseudochain map and $h$ is a homotopy, then the difference between $(dh + hd) f$ and $d(hf) + (hf) d$ is $h \nu_f \delta$, which can be nonzero.

\begin{defn} Let $\IC$ denote the span of the nulhomotopic maps, together with the pseudochain maps in the image of $\delta$. The category of pseudocomplexes modulo $\IC$ is denoted $K_\delta(\Diag_1)$ and is called the \emph{homotopy category of pseudocomplexes}. \end{defn}

\begin{exercise} Show that $\IC$ is an ideal in $\Ch_\delta(\Diag_1)$. \end{exercise}

\begin{prop} The category $K_\delta(\Diag_1)$ is triangulated. \end{prop}

The basic idea, just as for the ordinary homotopy category, is to define distinguished triangles as (triangles isomorphic to those) coming from cones of morphisms. We will not reproduce
the proof of this proposition here, see \cite[Lemma 4.6]{MackyMonodromy}.

\subsection{Monodromy}

We have already defined the monodromy map $\mu_X$ associated to a pseudocomplex $X$. One beautiful feature of the monodromy is that it descends to an interesting (genuine) chain map on
the (genuine) complex $\killit X$. Had one looked at $\killit X$ without the context of the pseudocomplex $X$, the existence and properties of this chain map would not be so evident. We
continue to let $\mu_X$ denote the monodromy chain map on $\killit X$.

\begin{exercise} Verify that the monodromy map $\mu_X$ is a pseudochain map $X \to X(-2)[2]$. \end{exercise}
	
\begin{lemma} (See \cite[Lemma 4.10]{MackyMonodromy}) Let $f \co X \to Y$ be any pseudochain map between pseudocomplexes. Then $f \mu_X$ and $\mu_Y f$ agree modulo nulhomotopic maps in $\Hom(X,Y(-2)[2])$. \end{lemma}

\begin{exercise} Verify that $\nu_f$ gives the homotopy between $f \mu_X$ and $\mu_Y f$. (Hint: Don't forget that $\delta$ is a nonzero divisor.) \end{exercise} 

Thus monodromy commutes with every pseudochain map (up to homotopy)! We have the following consequence in $\Diag_2$.

\begin{lemma} For any chain map $f \co \killit X \to \killit Y$, $f \mu_X \equiv \mu_Y f$ up to homotopy. \end{lemma}

\subsection{Gaussian elimination in ordinary homological algebra} \label{subsec:GE}

There is a well-known process in homological algebra, known as Gaussian elimination, which replaces complexes with homotopy equivalent complexes by removing any split maps from a
differential. We recall this first, and then in the \S\ref{subsec:GEforpseudo} state the easy generalization to pseudocomplexes. The material here is very similar to that in \cite[Appendix]{HogGor}.

Consider a complex $X$ of the following form, where $\phi$ is an isomorphism.
\begin{equation} \label{eq:GEstart} \begin{diagram} \cdots & \rTo & A & \rTo^{\left[ \begin{array}{c} a\\ b \end{array} \right]} & B \oplus F & \rTo^{\left[ \begin{array}{cc} c & e\\ d & \phi \end{array} \right]} & C \oplus F' & \rTo^{\left[ \begin{array}{cc} f & g \end{array} \right]} & D & \rTo & \cdots \end{diagram} \end{equation}
Gaussian elimination ``contracts'' away the summands $F$ and $F'$, connected by the isomorphism $\phi$, and yields the replacement complex $Y$ below.
\begin{equation} \label{eq:GEend} \begin{diagram} \cdots & \rTo & A & \rTo^{a} & B & \rTo^{c - e \phi\inv d} & C & \rTo^{f} & D & \rTo & \cdots \end{diagram} \end{equation}

\begin{prop}[Gaussian Elimination] The complexes $X$ and $Y$ are homotopy equivalent. \end{prop}
	
Note that Gaussian elimination does not change any of the differentials outside of the specific differential where a term is being contracted. It does, however, change the differential
from $B$ to $C$, altering it by the subtraction of a ``zigzag'' term $e \phi\inv d$.

\begin{exercise} Let $\alpha \co X \to Y$ be the map given below.
\begin{equation} \begin{diagram} A & \rTo & B \oplus F & \rTo & C \oplus F' & \rTo & D \\ \dTo^{1} & & \dTo^{\left[ \begin{array}{cc} 1 & 0 \end{array}\right]} & & \dTo^{\left[ \begin{array}{cc} 1 & -e\phi\inv \end{array}\right]} & & \dTo^{1} \\ A & \rTo & B & \rTo & C & \rTo & D \end{diagram} \end{equation}
Let $\beta \co Y \to X$ be the map given below.
\begin{equation} \begin{diagram} A & \rTo & B & \rTo & C & \rTo & D \\ \dTo^{1} & & \dTo^{\left[ \begin{array}{c} 1 \\ -\phi\inv d \end{array}\right]} & & \dTo^{\left[ \begin{array}{c} 1 \\ 0 \end{array}\right]} & & \dTo^{1} \\ A & \rTo & B \oplus F & \rTo & C \oplus F' & \rTo & D \end{diagram} \end{equation}
Prove that $Y$ is a complex, and that $\alpha$ and $\beta$ are chain maps. Where did you use the fact that $X$ is a complex? \end{exercise}

\begin{exercise} \label{exer:dontneedcomplex} Prove that $\alpha \circ \beta = \id_Y$. Prove that $\beta \circ \alpha + (dq + qd) = \id_X$, where $q$ is the homotopy sending $F'$ to $F$ via $\phi\inv$. Where did you use the fact that $X$ is a complex? (Hint: you didn't!) \end{exercise}
	
\subsection{Gaussian elimination for pseudocomplexes} \label{subsec:GEforpseudo}

Now we justify the statement that Gaussian elimination still holds for pseudocomplexes. Let $X$ be exactly as in \eqref{eq:GEstart}, but assume only that $X$ is a pseudocomplex. To aid in the exercises, let us give notation for the monodromy maps. Let
\begin{equation} \begin{diagram} A & \rTo^{\left[ \begin{array}{c} h \\ i \end{array} \right]} & C \oplus F' \end{diagram} \end{equation}
be the monodromy map in this particular degree, so that $h\delta = ca + eb$ and $i\delta = da + \phi b$.  Similarly, let the monodromy map in the next degree be as follows.
\begin{equation} \begin{diagram} B \oplus F & \rTo^{\left[ \begin{array}{cc} j & k \end{array} \right]} & D \end{diagram} \end{equation}
Let $Y$
be exactly as in \eqref{eq:GEend}.

\begin{exercise} Prove that $Y$ is a pseudocomplex. Compute $\mu_Y$. \end{exercise}
	
\begin{exercise} Prove that $\alpha$ and $\beta$ are pseudochain maps. Compute that $\nu_{\alpha}$ consists only of $\phi\inv k$ in one spot, and $\nu_{\beta}$ consists only of $\phi\inv i$. \end{exercise}

Because Exercise \ref{exer:dontneedcomplex} goes through verbatim for pseudocomplexes, we have the following result.

\begin{prop}[Gaussian Elimination for pseudocomplexes] The pseudocomplexes $X$ and $Y$ are homotopy equivalent. \end{prop}
	
In particular, Gaussian elimination for pseudocomplexes matches with Gaussian elimination for the ordinary complexes obtained after applying $\killit$.

\begin{exercise} The morphism $\alpha$ is supposed to intertwine the monodromy maps $\mu_X$ and $\mu_Y$ up to homotopy. Verify that this is true, and that the homotopy is given by
$\nu_{\alpha}$. Do the same for $\beta$. For this exercise, you should zoom out and look at a bigger portion of your complexes (i.e. name by $G$ the chain object in the degree before
$A$, and by $H$ the chain object in the degree after $D$. Name your chain maps, and your monodromy maps, etcetera.) \end{exercise}

\subsection{Monoidal structures} \label{subsec:pseudomonoidal}

In \cite[Chapter 4]{MackyMonodromy}, it is important to allow a version of pseudocomplexes where $d^2$ can fail to be zero not just because of multiplication by $\delta$, but because of
right multiplication by any polynomial $R$ of positive degree. Killing the right action of positive degree polynomials, one obtains the category of Soergel modules (rather than Soergel
bimodules), and this is not a monoidal category. 

Similarly, pseudocomplexes do not form a monoidal category in general. Let $X$ and $Y$ be pseudocomplexes in this more general sense, where $d^2$ is equivalent to zero modulo the action
of positive degree polynomials in the right. Define $X \ot Y$ using the usual tensor product differential (with its usual signs). Then $d^2_{X \ot Y} = d^2_X
\ot 1 + 1 \ot d^2_Y$. Now $1 \ot d^2_Y$ is in the ideal generated by positive degree polynomials acting on the right, while $d^2_X \ot 1$ has a positive degree polynomial acting in the
middle (on the right of $X$, and the left of $Y$), but not necessarily on the right of the whole tensor product. Hence $X \ot Y$ need not be a pseudocomplex.

However, in our version of pseudocomplexes, we have specified that $d^2$ is a multiple of $\delta$. Now the fact that $\delta$ is $W$-invariant implies that the left and right
actions of $\delta$ on any Hom space in $\Diag_1$ agree. Consequently, by the calculation of the previous paragraph, the tensor product of two pseudocomplexes is another pseudocomplex.

\begin{exercise} Verify (by the same calculation) that 
\begin{equation} \label{eq:monodromyontensor} \mu_{X \ot Y} = \mu_X \ot 1 + 1 \ot \mu_Y. \end{equation} \end{exercise}
	
\subsection{Simultaneous Gaussian elimination} \label{subsec:GEsimult}

We conclude with a discussion of efficiently performing simultaneous Gaussian elimination, because this is hard to find in the literature.

In a Krull-Schmidt category, one can write each homological degree as a direct sum of indecomposable objects. By iterating Gaussian elimination, one can find a homotopy-equivalent
replacement where no matrix entry of any differential is an isomorphism. This is the so-called \emph{minimal complex}, and in certain kinds of categories, it is unique up to isomorphism
of complexes.

In the process of iterated Gaussian elimination, if one only cared about the chain objects and not the differential, then one might wish to merely chop off isomorphisms
until none are left, and see which chain objects survive. The danger when doing this is that the differentials are constantly being modified by zigzag terms, and this can affect which matrix entries are isomorphisms. For example, in the
two-term complex \begin{equation} \begin{diagram} R \oplus R & \rTo^{\left[ \begin{array}{cc} 1 & 1 \\ 1 & 1 \end{array} \right]} & R \oplus R \end{diagram} \end{equation} there is an isomorphism from the first copy of $R$ in homological degree 0 to the first copy in degree 1, and from the second copy in degree 0 to the second copy in degree 1. One might hope that one can remove both isomorphisms, yielding the zero complex. However, after applying Gaussian elimination once (to either copy), what remains is
\begin{equation} \begin{diagram} R & \rTo^{0} & R \end{diagram} \end{equation}
and the differential is no longer an isomorphism. This is an example of an apparent iterated Gaussian elimination which is not actually valid.

Because of this problem, we need some additional conditions to efficiently perform iterated Gaussian elimination. For pedagogical reasons we begin by considering the case of two
isomorphisms.

\begin{defn} \label{defn:independentpair} Consider a complex of the following form, where $\phi$ and $\psi$ are isomorphisms.
\begin{equation} \label{eq:GEstart2} \begin{diagram} \cdots & \rTo & A & \rTo^{\left[ \begin{array}{c} a\\ b_1 \\ b_2 \end{array} \right]} & B \oplus F_1 \oplus F_2 & \rTo^{\left[ \begin{array}{ccc} c & e_1 & e_2\\ d_1 & \phi & \theta \\ d_2 & \tau & \psi \end{array} \right]} & C \oplus F'_1 \oplus F'_2 & \rTo^{\left[ \begin{array}{ccc} f & g_1 & g_2 \end{array} \right]} & D & \rTo & \cdots \end{diagram} \end{equation}
If $\theta \psi\inv \tau = 0$ and $\tau \phi\inv \theta = 0$ we say that $\phi$ and $\psi$ are \emph{independent}. \end{defn}

The definition of independence depends not just on the isomorphisms $\phi$ and $\psi$ but on the entire $2 \times 2$ submatrix of differentials $F_1 \oplus F_2 \to F'_1 \oplus F'_2$.

Note that after eliminating $\phi$, the new map from $F_2$ to $F'_2$ would normally be $\psi - \tau \phi\inv \theta$, and the condition of independence implies that this again equals
$\psi$. Independence implies that eliminating $\phi$ will not affect the isomophism $\psi$, and vice versa. This means that iterated Gaussian elimination is always valid, and the following exercise describes the result.

\begin{exercise} \label{exer:zigzag} Show that after eliminating both $F_1$ and $F_2$ (in either order) the remaining differential $B \to C$ has the form
\begin{equation} c - e_1 \phi\inv d_1 - e_2 \psi\inv d_2 +  e_2 \psi\inv \tau \phi\inv d_1 + e_1 \phi\inv \theta \psi\inv d_2. \end{equation}
That is, one keeps the original differential $c$, subtracts both ``length one zigzags,'' and adds both ``length two zigzags.'' \end{exercise}

\begin{remark} When $\phi$ and $\psi$ are not independent, it may still be the case that $\psi - \tau \phi\inv \theta$ is invertible, so that one can eliminate both $F_1$ and $F_2$. The
formula for the final differential $B \to C$ is not as clean. If one expands $(\psi - \tau \phi\inv \theta)\inv$ as an infinite sum using a geometric series, one can express the
differential $B \to C$ as an infinite alternating sum of ``zigzags.'' \end{remark}

Now let us state the procedure for simultaneous Gaussian elimination, which efficiently eliminates many pairwise-independent isomorphisms at once.

\begin{defn} \label{defn:independent} Consider a complex with the following form
\begin{equation} A \longrightarrow B \oplus F_1 \oplus \cdots \oplus F_m \longrightarrow C \oplus F'_1 \oplus \cdots \oplus F'_m \longrightarrow D \end{equation}
The summands of the differential will be denoted
\[ c \co B \to C, \quad d_i \co B \to F'_i, \quad e_i \co F_i \to C, \quad \theta_{ij} \co F_j \to F'_i,\]
where $\phi_i = \theta_{ii}$ is an isomorphism. Associated to a sequence $i_1, \ldots, i_r$ of indices between $1$ and $m$, and a choice of source $S$ which is either $B$ or some $F_j$, and a target $T$ which is either $C$ or some $F'_k$, we have the following \emph{length $r$ zigzag}:
\begin{equation} t \circ \phi_{i_r}\inv \circ \cdots \circ \theta_{i_3 i_2} \circ \phi_{i_2}\inv \circ \theta_{i_2 i_1} \circ \phi_{i_1}\inv \circ s. \end{equation}
Here, $s$ is $d_{i_1}$ if $S = B$ or $\theta_{i_1 j}$ if $S = F_j$, and $t$ is $e_{i_r}$ if $T = C$ or $\theta_{k i_r}$ if $T = F'_k$. We say that a zigzag is \emph{repeating} if some index appears twice in the list $i_1, \ldots, i_r$, and is \emph{non-repeating} otherwise. We say that the family $\{\phi_i\}$ of isomorphisms is \emph{independent} if for all $r \ge 1$ and all $1 \le j \le m$, all length $r$ zigzags from $F_j$ to $F'_j$ are zero. \end{defn}

\begin{exercise} Assuming that all zigzags of length $< r$ from $F_j$ to $F'_j$ are zero for each $j$, prove that any repeating zigzag of length $r$ is zero. Prove that if all 
non-repeating zigzags from $F_j$ to $F'_j$ are zero, then the isomorphisms are independent. \end{exercise}

Non-repeating zigzags have length at most $m$, so checking independence is a finite amount of work.
	
\begin{exercise} Suppose that the family of isomorphisms $\{\phi_i\}$ is independent. Prove that each pair of isomorphisms $\phi_i$ and $\phi_j$ are independent under Definition
\ref{defn:independentpair}, and prove that they remain independent after eliminating any isomorphism $\phi_k$ for $k \ne i, j$. \end{exercise}

\begin{prop}[Simultaneous Gaussian elimination] Continue the setup of Definition \ref{defn:independent}. If the family of isomorphisms $\{\phi_i\}$ is independent, then the complex is equivalent to a complex of the form
\[ A \longrightarrow B \longrightarrow C \longrightarrow D. \]
The summand of the differential from $A$ to $B$ or from $C$ to $D$ is unchanged. The summand of the differential from $B$ to $C$ is 
\begin{equation} c + \sum_{r = 1}^m (-1)^r \sum_z z \end{equation}
where the second sum is over all nonrepeating zigzags of length $r$. \end{prop}

\begin{exercise} Prove this proposition by induction. \end{exercise}

%%%%%%%%%%%%%%%%%%%%%%%%%
\section{The standard Gaitsgory complex: $n=2,3,4$}
\label{sec-234}

Let $V = V_{\fund_1}$ denote the standard representation of $\gl_n$. Our aim is to construct the Gaitsgory central complex $\VC = \GC(V)$. The complex $\VC_{\fund_{n-1}}$ for the
``dual'' representation (again, this is not quite the dual of $V$, being off by a copy of the determinant representation) can be described in essentially the same way. Before stating
this general construction, in this chapter we give examples for $n=2, 3, 4$. This will help us to fix our conventions, and will add one small wrinkle of complexity each time.

Recall our notation for the affine and finite Weyl groups from \S\ref{subsec-typeAnotation}. We set $(W,S) = (W_{\aff},S_{\aff})$. Recall that $\tau(s_k) = s_{k+1}$ and $\si(s_k) =
s_{-k}$ are two automorphisms of $(W,S)$. Note that $\si$ descends to an automorphism of $(W_{\fin}, S_{\fin})$, for which $\si(s_k) = s_{n-k}$. Both $\tau$ and $\si$ lift to the
extended Weyl group via $\tau(\om) = \om$ and $\si(\om) = \om\inv$.

We set $\Diag$ to be either $\Diag_{\aff}$ or $\Diag_{\ext}$, and $\HC$ to be its homotopy category. Whenever we discuss the (pseudo)complex $\FC$ below, one is welcome to work in either
the affine setting or the extended setting. The complex $\VC$ only exists in the extended setting, and we will explicitly write $\Diag_{\ext}$ or $\HC_{\ext}$ when the extended setting
is required. The reader willing to work with the technology of pseudocomplexes introduced in \S\ref{sec-pseudo} may assume that $\delta \ne 0$, while the reader who only likes complexes
can assume that $\delta = 0$.

%===========
\subsection{The case $n=2$}
\label{subsec-stdrep2}
%===========

Fix a realization in type $\tilde{A}_1$ with simple reflections $\{s_0,s_1\}$, and let $\delta=\al_0 + \al_1$ in $\hg^*$. We associate blue to $s_1$ and red to $s_0$.

\begin{defn} \label{def:F2} ($n=2$) Let $\FC \in \HC$ be the pseudocomplex given by
\begin{equation} \FC = \left(
\begin{tikzpicture}
\node (a) at (0,0) {$R(-1)$};
\node (a0) at (2,1) {$B_1$};
\node (a1) at (2,-1) {$B_0$};
\node (b) at (4,0) {$R(1)$};
\path
	(a) edge node[descr] {$\startdotblue$} (a0)
	(a) edge node[descr] {$\startdotred$} (a1)
	(a0) edge node[descr] {$\finaldotblue$} (b)
	(a1) edge node[descr] {$\finaldotred$} (b);
%	\draw [brown] (current bounding box.south west) rectangle (current bounding box.north east);
\end{tikzpicture} \right)
\end{equation}
\end{defn}

The columns of this diagram represent the various homological degrees of $\FC$, which will always be centered around homological degree $0$. Thus, in this example, $\FC^{-1} = R(-1)$,
$\FC^0 = B_0 \oplus B_1$, and $\FC^1 = R(1)$.

The morphisms on the individual arrows above are \emph{summands} of the differential, and we often just say ``the differential from $X$ to $Y$'' to pick out a particular summand. For
example, the differential from $R(-1)$ to $B_1$ is a blue startdot. We might also talk about summands of $d^2$; for example, ``$d^2$ from $R(-1)$ to $R(1)$" would be the morphism from
$R(-1)$ to $R(1)$ induced by $d^2$. To refer to a composition of any two arrows, we might refer to a \emph{term} of $d^2$. For example, $d^2$ from $R(-1)$ to $R(1)$ is the sum of two
terms, the term which factors through $B_1$ and the term which factors through $B_0$.

Let us compute $d^2$. There is only one nontrivial summand to check, $d^2$ from $R(-1)$ to $R(+1)$. Recall that a composition of two dots is a map called a \emph{barbell}, which is
equal to multiplication by the corresponding simple root. Thus $d^2$ from $R(-1)$ to $R(1)$ is given by $(\al_0 + \al_1) \cdot \id_R = \delta \cdot \id_R$. Thus $\FC$ is a
pseudocomplex, and $d^2=0$ if and only if $\delta = 0$.

Several informal observations of the form of $\FC$ are in order: \begin{itemize}
\item Any summand $B_x(k)$ appears in homological degree $k$, for any $x \in W$ and $k \in \Z$.
\item If $B_x(k)$ appears, then so does $B_{\tau(x)}(k)$. Moreover, the differential from $B_{\tau(x)}(k)$ to $B_{\tau(y)}(k+1)$ is just the color swap of the differential from $B_x(k)$ to $B_y(k+1)$.
\item If $B_x(k)$ appears in homological degree $k$, then $B_x(-k)$ appears in homological degree $-k$. Moreover, the differential from $B_y(-(k+1))$ to $B_x(-k)$ is just the vertical flip of the differential from $B_x(k)$ to $B_y(k+1)$.
\end{itemize}
Stated more formally, the complex $\FC$ satisfies the following properties:
\begin{itemize}
\item It is \emph{perverse}, meaning that an indecomposable bimodule in homological degree $k$ appears with shift $(k)$. Perversity implies that all summands of the differential have degree $+1$. Perversity also implies that the complex is \emph{minimal}, meaning that every summand of the differential is in the graded Jacobson radical of $\Diag$.
\item It is \emph{color rotation invariant}. More precisely, there is an obvious isomorphism from $\tau(\FC)$ to $\FC$. Consequently, there is a natural isomorphism $\Om \FC \cong \FC \Om$ in $\HC_{\ext}$.
\item It is \emph{self-dual}. More precisely, there is an obvious isomorphism from $\DM(\FC)$ to $\FC$.
\end{itemize}

\begin{defn} Define $\VC := \Om \FC$ in $\HC_{\ext}$. \end{defn}

Thus $\VC$ is also perverse, minimal, self-dual, and comes with a natural isomorphism $\VC \Om \cong \Om \VC$.

\begin{remark} The complex $\FC$ also appears to be invariant under $\si$, a special feature of $n=2$. The complex $\VC$ is not invariant under $\si$ because $\si(\Om) = \Om\inv$. In
the $\sl_2$-extended diagrammatic category, $\VC$ will in fact by $\si$-invariant, related to the fact that $V$ is self-dual. \end{remark}

Note that $\VC$ has a subcomplex
\begin{equation} \Om \left(
\begin{tikzpicture}
\node (a) at (0,0) {$ $};
\node (a0) at (2,1) {$B_1$};
\node (a1) at (2,-1) {$ $};
\node (b) at (4,0) {$R(1)$};
\path
	(a0) edge node[descr] {$\finaldotblue$} (b);
\end{tikzpicture} \right)
\end{equation}
which is the Rouquier complex for $\om f_1$, a.k.a. the Wakimoto complex $W_{(1,0)}$. The quotient by this subcomplex is
\begin{equation}\Om \left(
\begin{tikzpicture}
\node (a) at (0,0) {$R(-1)$};
\node (a0) at (2,1) {$ $};
\node (a1) at (2,-1) {$B_0$};
\node (b) at (4,0) {$ $};
\path
	(a) edge node[descr] {$\startdotred$} (a1);
\end{tikzpicture} \right)
\end{equation}
which is the Rouquier complex of $\om f_0\inv$, a.k.a. the Wakimoto complex $W_{(0,1)}$. This is the Wakimoto filtration on $\VC$.

One also has a dual Wakimoto filtration, with subcomplex
\begin{equation}W^*_{(0,1)} = \Om F_0 = \Om \left(
\begin{tikzpicture}
\node (a) at (0,0) {$ $};
\node (a0) at (2,1) {$ $};
\node (a1) at (2,-1) {$B_0$};
\node (b) at (4,0) {$R(1)$};
\path
	(a1) edge node[descr] {$\finaldotred$} (b);
\end{tikzpicture} \right)
\end{equation}
and quotient
\begin{equation}W^*_{(1,0)} = \Om F_1\inv = \Om \left(
\begin{tikzpicture}
\node (a) at (0,0) {$R(-1)$};
\node (a0) at (2,1) {$B_1$};
\node (a1) at (2,-1) {$ $};
\node (b) at (4,0) {$ $};
\path 
	(a) edge node[descr] {$\startdotblue$} (a0);
\end{tikzpicture} \right).
\end{equation}

There is a chain map $\mu \co \FC \to \FC[2](-2)$, defined as follows. (We centered this picture at homological degree $-1$.)
\begin{equation} \mu = \left(
\begin{tikzpicture}
\node (a) at (0,0) {$R(-1)$};
\node (a0) at (2,1) {$B_1$};
\node (a1) at (2,-1) {$B_0$};
\node (b) at (4,0) {$R(1)$};
\node (c) at (-4,-4) {$R(-3)$};
\node (c0) at (-2,-3) {$B_1(-2)$};
\node (c1) at (-2,-5) {$B_0(-2)$};
\node (d) at (0,-4) {$R(-1)$};

\path
	(a) edge node[descr] {$\startdotblue$} (a0)
	(a) edge node[descr] {$\startdotred$} (a1)
	(a0) edge node[descr] {$\finaldotblue$} (b)
	(a1) edge node[descr] {$\finaldotred$} (b);
	
\path
	(c) edge node[descr] {$\startdotblue$} (c0)
	(c) edge node[descr] {$\startdotred$} (c1)
	(c0) edge node[descr] {$\finaldotblue$} (d)
	(c1) edge node[descr] {$\finaldotred$} (d);
	
\path
	(a) edge node[descr] {$\id_R$} (d);

\end{tikzpicture} \right)
\end{equation}
This induces the corresponding monodromy map $\mu \co \VC \to \VC[2](-2)$. Note that $d^2 = \mu \delta$.

Let us study $\FC B_i$ for $i = 0,1$. Recall that $B_i B_i \cong B_i(-1) \oplus B_i(+1)$. Recall also that the startdot map $R(-1) \to B_i$ becomes an inclusion of a summand $B_i(-1) \sumset B_i B_i$ after tensoring with $B_i$. Similarly, the enddot map $B_i \to R(+1)$ becomes a projection to a summand $B_i B_i \to B_i(+1)$ after tensoring with $B_i$.
Using Gaussian elimination of complexes one can deduce that $\FC B_1$ is homotopy equivalent to $B_0 B_1$ in homological degree zero, because the other summands $R(-1) B_1$, $B_1 B_1$, and $R(+1) B_1$ all cancel out. Using similar arguments we see that
\begin{subequations}
\begin{equation} \FC B_1 \cong B_0 B_1 \cong B_0 \FC,\end{equation}
\begin{equation} \FC B_0 \cong B_1 B_0 \cong B_1 \FC.\end{equation}
Thus $M \ot \FC \cong \FC \ot \tau(M)$ for any object $M$ of $\Diag_{\aff}$ or $\Diag_{\ext}$. We will prove in this paper that this isomorphism 
\begin{equation} \label{eq:Ffunctorial1} (-) \ot \FC \to \FC \ot \tau(-)\end{equation} can be made functorial with respect to morphisms in $\Diag$.
	
Since $\Om B_1 \cong B_0 \Om$, we see that
\begin{equation} \VC B_1 \cong B_1 \VC, \end{equation}
\begin{equation} \VC B_0 \cong B_0 \VC. \end{equation}
\end{subequations}
We already argued that $\VC \Om \cong \Om \VC$, so it follows that $\VC M \cong M \VC$ for any object $M$ of $\Diag_{\ext}$. Using \eqref{eq:Ffunctorial1} and the functorial isomorphism \begin{equation} \label{eq:Omfunctorial1} (-) \ot \Om \to \Om \ot \tau\inv(-), \end{equation} one obtains an isomorphism
\begin{equation} \label{eq:Vfunctorial1} (-) \ot \VC \to \VC \ot (-) \end{equation} as functors $\Diag_{\ext} \to \HC_{\ext}$.

Now consider what happens when we take $\FC$ (perhaps viewed as a complex of $R$-bimodules) and restrict from $R$ to $R^{W_{\fin}}$ on the right. Here we identify $W_{\fin}$ with the
parabolic subgroup generated by $s_1$. This calculation is very similar to the calculation of $\FC B_1$. Namely, as $(R,R^{s_1})$ bimodules, $B_1$ is isomorphic to $R(+1) \oplus R(-1)$.
Gaussian elimination then cancels all terms except for $B_0$, so we conclude that \begin{equation} \FC_{R^{s_1}} \cong (B_0)_{R^{s_1}} \cong R \ot_{R^{s_0}} R_{R^{s_1}}(1).
\end{equation} In particular, after restriction to $R^{W_{\fin}}$, $\FC$ is isomorphic to a complex supported entirely in degree zero.

Moreover, applying the Soergel Satake functor from \cite{EQuantumI} to the standard representation yields the $(R^{s_0},R^{s_1})$-bimodule \begin{equation} \SSat(V) = {}_{R^{s_0}}
R_{R^{s_1}}(1). \end{equation} Thus we deduce that the induction on the left of $\SSat(V)$ agrees with the restricting on the right of $\FC$.

In $\DC_{\ext}$ there is a version of the Soergel Satake functor which yields not an $(R^{s_0},R^{s_1})$-bimodule but an $(R^{W_{\fin}},R^{W_{\fin}})$-bimodule, namely it sends $V$ to
$\Om \SSat(V)$. This extra copy of $\Om$ swaps the $R^{s_0}$ action with an $R^{s_1}$ action. Thus we conclude that
\begin{equation} R \ot_{R^{W_{\fin}}} (\Om \SSat(V)) \cong \VC_{R^{W_{\fin}}}. \end{equation}

\begin{remark} \label{rmk:restricttensorsimilar} Induction from $R^s$ to $R$ will send indecomposable singular Soergel bimodules to indecomposable Soergel bimodules. Thus, once one has
computed the restriction of $\VC$ to $R^{s_1}$ on the right, it is immediate to compute $\VC \ot_{R^{s_1}} R(1) = \VC B_1$. This explains why the computation of $\VC_{R^{s_1}}$ and $\VC
B_1$ were so similar. For the reader unfamiliar with singular Soergel bimodules, let this serve as the motto: if one is interested in computing the restriction from $R$ to the
invariants $R^J$ under a parabolic subgroup $W_J$, it is essentially equivalent to study the tensor product with $B_J$, the indecomposable associated with the longest element $w_J$ of
$w_J$. \end{remark}

\begin{remark} We use the language of bimodules above, but could also work with the diagrammatic category of singular Soergel bimodules, which is still lagging behind in the literature
but should appear relatively soon. \end{remark}

Finally, consider the quotient by the cellular ideal $\MC$, consisting of all morphisms factoring through all objects except for $R$. Let $\overline{R}$ denote the image of $R$ in the
quotient $\Diag/\MC$, and $\overline{\FC}$ denote the image of $\FC$. Then clearly $\overline{\FC}$ is just a complex with zero differential, consisting of $\overline{R}(-1)$ in degree
$-1$, and $\overline{R}(1)$ in degree $+1$. The monodromy map $\mu$ descends to $\overline{\FC}$ in the obvious way. What remains is the underlying vector space of the standard
representation $V$ of $\gl_2$, equipped with the action of its raising operator.

%===========
\subsection{The case $n=3$}
\label{subsec-stdrep3}
%===========

Fix a realization in type $\tilde{A}_2$ with simple reflections $\{s_0,s_1,s_2\}$, and let $\delta = \al_0 + \al_1 + \al_2$ in $\hg^*$. We associate blue to $s_2$, red to $s_1$, and
green to $s_0$.

\begin{defn} \label{def:F3} ($n=3$) Let $\FC \in \HC$ be the complex given by
\begin{equation}
\FC = \left(
\begin{tikzpicture}
\node (a) at (0,0) {$R(-2)$}; \node (a0) at (3,2) {$B_2(-1)$};
\node (a1) at (3,0) {$B_1(-1)$}; \node (a2) at (3,-2) {$B_0(-1)$};

\node (b) at (6,-8) {$R$};
\node (a01) at (6,6) {$B_{21}$};
\node (a12) at (6,-2) {$B_{10}$};
\node (a20) at (6,2) {$B_{02}$};

\node (b0) at (9,2) {$B_2(1)$};
\node (b1) at (9,0) {$B_1(1)$};
\node (b2) at (9,-2) {$B_0(1)$};
\node (c) at (12,0) {$R(2)$};
\path
	(a) edge node[descr] {$\startdotblue$} (a0)
	(a) edge node[descr] {$\startdotred$} (a1)
	(a) edge node[descr] {$\startdotgreen$} (a2)	
	(a0) edge node[descr] {$\finaldotblue$} (b)
	(a1) edge node[descr] {$\finaldotred$} (b)
	(a2) edge node[descr] {$\finaldotgreen$} (b)

	(b) edge node[descr] {$\startdotblue$} (b0)
	(b) edge node[descr] {$\startdotred$} (b1)
	(b) edge node[descr] {$\startdotgreen$} (b2)	
	(b0) edge node[descr] {$\finaldotblue$} (c)
	(b1) edge node[descr] {$\finaldotred$} (c)
	(b2) edge node[descr] {$\finaldotgreen$} (c)

	(a0) edge node[descr] {$- \lineblue \startdotred$} (a01)
	(a1) edge node[descr] {$- \linered \startdotgreen$} (a12)
	(a2) edge node[descr] {$- \linegreen \startdotblue$} (a20)
	(a0) edge node[descr] {$\startdotgreen \lineblue$} (a20)
	(a1) edge node[descr] {$\startdotblue \linered$} (a01)
	(a2) edge node[descr] {$\startdotred \linegreen$} (a12)

	(a01) edge node[descr] {$- \lineblue \finaldotred$} (b0)
	(a12) edge node[descr] {$- \linered \finaldotgreen$} (b1)
	(a20) edge node[descr] {$- \linegreen \finaldotblue$} (b2)
	(a20) edge node[descr] {$\finaldotgreen \lineblue$} (b0)
	(a01) edge node[descr] {$\finaldotblue \linered$} (b1)
	(a12) edge node[descr] {$\finaldotred \linegreen$} (b2);
	
\end{tikzpicture} \right)
\end{equation} \end{defn}

The differentials are all dots (possibly tensored with identity maps), with a sign which obeys the \emph{strand-counting sign rule}: it is determined by the number of identity maps appearing
to the left of the dot. Here, $B_{ij}$ is the Bott-Samelson bimodule $B_i \ot B_j$, which we identify with the indecomposable Soergel bimodule $B_{s_i s_j}$. The properties of perversity,
color rotation invariance, and self-duality are still evident.

Let us compute the summands of $d^2$. We assume $j \ne i$ below. \begin{itemize}
\item As a map $R(-2) \to B_{ij}(0)$, $d^2$ has two terms, factoring through $B_i(-1)$ and $B_j(-1)$ respectively, which cancel thanks to the sign convention.
\item The same is true for the map $B_{ij}(0) \to R(2)$.
\item As a map $R(-2) \to R(0)$, $d^2 = \al_0 + \al_1 + \al_2 = \delta \cdot \id_R$.
\item The same is true for the map $R(0) \to R(2)$.
\item As a map $B_i(-1) \to B_j(+1)$, $d^2$ has two terms, factoring through $B_{ij}(0)$ and $R(0)$ respectively, which cancel thanks to the sign convention.
\item As a map $B_i(-1) \to B_i(+1)$, $d^2$ has three terms, which we draw for $i=2$:
\begin{equation} d^2 = \lineblue \barbred + \barbgreen \lineblue + \brokenblue = \barbblue \lineblue + \barbred \lineblue + \barbgreen \lineblue. \end{equation}
Above, we used the barbell forcing relations to see that $d^2 = (\al_0 + \al_1 + \al_2) \cdot \id_{B_i} = \delta \cdot \id_{B_i}$.
\end{itemize}
Thus $\FC$ is a pseudocomplex, and is a complex if and only if $\delta=0$.

For sake of reference, the barbell forcing relations are reprinted below. For any color we have
\begin{subequations} \label{eq:barbellforce}
\begin{equation} \lineblue \barbblue + \barbblue \lineblue = 2 \brokenblue. \end{equation}
When red and blue are adjacent we have
\begin{equation} \lineblue \barbred = \barbblue \lineblue + \barbred \lineblue - \brokenblue. \end{equation}
The horizontal reflection of this relation also holds.
When blue and purple are distant (only relevant for larger $n$) we have
\begin{equation} \lineblue \barbpurple = \barbpurple \lineblue. \end{equation}
\end{subequations}
The general form of these relations states that: a polynomial $f$ to the left of an $i$-colored line is equal to $s_i(f)$ on to its right, together with the broken line (two dots) multiplied by $\pa_i(f)$. See \cite[equation (5.2)]{EWGr4sb}.

Let us now dispense with the differentials and write the complex more compactly.

\begin{equation}
\FC = \left(
\begin{tikzpicture}
\node (a) at (0,0) {$R(-2)$};
\node at (2,1) {$B_2(-1)$};
\node (b) at (2,0) {$B_1(-1)$}; 
\node at (2,-1) {$B_0(-1)$};

\node at (4,-1.5) {$R$};
\node at (4,1.5) {$B_{21}$};
\node at (4,.5) {$B_{10}$};
\node at (4,-.5) {$B_{02}$};
\node (c) at (4,0) {$\;\;\;$};

\node  at (6,1) {$B_2(1)$};
\node (d) at (6,0) {$B_1(1)$};
\node at (6,-1) {$B_0(1)$};
\node (e) at (8,0) {$R(2)$};
\path
	(a) edge (b)
	(b) edge (c)
	(c) edge (d)
	(d) edge (e);	
\end{tikzpicture} \right)
\end{equation}

Note that this complex does not have \emph{color-reversal symmetry}, i.e. $\FC \ncong \si(\FC)$. Instead, acting by $-1$ on the indices/colors in $\Z/3\Z$ will yield a new
complex $\FC^* := \si(\FC)$. In homological degree $0$, $\FC^*$ has summands $B_{12}$, $B_{20}$, $B_{01}$, and $R$.

\begin{defn} Let $\VC = \Om \FC$, and $\VC^* := \si(\VC) = \Om^{-1} \FC^*$. Let $\VC_{\fund_2} = \Om^2 \FC^*$. \end{defn}

\begin{remark} Note that $\VC^*$ and $\VC_{\fund_2}$ are isomorphic in the diagrammatic category for extended affine $\sl_2$, since $\Om^3 \cong \one$, see \S\ref{subsec-extendedsbimquotient}. \end{remark}

Here is the Wakimoto filtration: we list the subquotients in order from subcomplex to quotient.
\begin{equation}
\WC_{(1,0,0)} = \Om F_2 F_1 = \left(
\begin{tikzpicture}
\node at (0,0) {$ $};
\node at (2,1) {$ $};
\node at (2,0) {$ $}; 
\node at (2,-1) {$ $};

\node at (4,-1.5) {$ $};
\node (a) at (4,1.5) {$B_{21}$};
\node at (4,.5) {$ $};
\node at (4,-.5) {$ $};
\node at (4,0) {$\;\;\;$};

\node (b) at (6,1) {$B_2(1)$};
\node (c) at (6,0) {$B_1(1)$};
\node at (6,-1) {$ $};
\node (d) at (8,0) {$R(2)$};
\path
	(a) edge (b)
	(a) edge (c)
	(b) edge (d)
	(c) edge (d);	
\end{tikzpicture} \right)
\end{equation}
\begin{equation}
\WC_{(0,1,0)} = \Om F_0\inv F_2 = \left(
\begin{tikzpicture}
\node at (0,0) {$ $};
\node (a) at (2,1) {$B_2(-1)$};
\node at (2,0) {$ $}; 
\node at (2,-1) {$ $};

\node (b) at (4,-1.5) {$R$};
\node at (4,1.5) {$ $};
\node at (4,.5) {$ $};
\node (c) at (4,-.5) {$B_{02}$};
\node at (4,0) {$\;\;\;$};

\node at (6,1) {$ $};
\node at (6,0) {$ $};
\node (d) at (6,-1) {$B_0(1)$};
\node at (8,0) {$ $};
\path
	(a) edge (b)
	(a) edge (c)
	(b) edge (d)
	(c) edge (d);
\end{tikzpicture} \right)
\end{equation}
\begin{equation}
\WC_{(0,0,1)} = \Om F_1\inv F_0\inv = \left(
\begin{tikzpicture}
\node (a) at (0,0) {$R(-2)$};
\node at (2,1) {$ $};
\node (b) at (2,0) {$B_1(-1)$}; 
\node (c) at (2,-1) {$B_0(-1)$};

\node at (4,-1.5) {$ $};
\node at (4,1.5) {$ $};
\node (d) at (4,.5) {$B_{10}$};
\node at (4,-.5) {$ $};
\node at (4,0) {$\;\;\;$};

\node  at (6,1) {$ $};
\node at (6,0) {$ $};
\node at (6,-1) {$ $};
\node at (8,0) {$ $};
\path
	(a) edge (b)
	(a) edge (c)
	(b) edge (d)
	(c) edge (d);
\end{tikzpicture} \right)
\end{equation}
That the differentials match (between the Wakimoto complexes and the subquotients of $\VC$) follows because both complexes satisfy the strand-counting sign rule, see \S\ref{subsec-signconvention}. We leave the reader to compute the dual Wakimoto filtration of $\VC$, as well as the Wakimoto and dual Wakimoto filtrations of $\VC^*$ and $\VC_{\fund_2}$.

Let $\ImM$ denote the following pseudocomplex.
\begin{equation} \ImM = \left(
\begin{tikzpicture}
\node (a) at (0,0) {$R(-1)$};
\node (a2) at (2,1.5) {$B_2$};
\node (a1) at (2,0) {$B_1$};
\node (a0) at (2,-1.5) {$B_0$};
\node (b) at (4,0) {$R(1)$};
\path
	(a) edge node[descr] {$\startdotblue$} (a2)
	(a2) edge node[descr] {$\finaldotblue$} (b)
	(a) edge node[descr] {$\startdotred$} (a1)
	(a) edge node[descr] {$\startdotgreen$} (a0)
	(a1) edge node[descr] {$\finaldotred$} (b)
	(a0) edge node[descr] {$\finaldotgreen$} (b);
\end{tikzpicture} \right) \end{equation}
Note that $d^2 = \delta \cdot \id_R$. Then $\ImM(1)[-1]$ appears as a subcomplex of $\FC$, containing all the summands in homological degrees $2$ and $1$, together with $R$ in degree zero. Also, $\ImM(-1)[1]$ appears as a quotient complex of $\FC$, where the only summands which survive are those in homological degrees $-2$ and $-1$, together with $R$ in degree zero.

The monodromy map  $\mu \co \FC \to \FC(-2)[2]$ is the composition
\begin{equation} \FC \to \ImM(-1)[1] \to \FC(-2)[2] \end{equation}
of the quotient to $\ImM$ followed by the inclusion of $\ImM$. Explicitly, $\mu$ is given by
\begin{equation}
M = \left(
\begin{tikzpicture}
\node (a) at (4,0) {$R(-2)$};
\node at (6,1) {$B_2(-1)$};
\node (b) at (6,0) {$B_1(-1)$}; 
\node (x) at (6,-1) {$B_0(-1)$};

\node (y) at (8,-1.5) {$R$};
\node at (8,1.5) {$B_{21}$};
\node at (8,.5) {$B_{10}$};
\node at (8,-.5) {$B_{02}$};
\node (c) at (8,0) {$\;\;\;$};

\node  at (10,1) {$B_2(1)$};
\node (d) at (10,0) {$B_1(1)$};
\node at (10,-1) {$B_0(1)$};
\node (e) at (12,0) {$R(2)$};
\path
	(a) edge (b)
	(b) edge (c)
	(c) edge (d)
	(d) edge (e);

\node (a1) at (0,-4) {$R(-2)$};
\node at (2,-3) {$B_2(-1)$};
\node (b1) at (2,-4) {$B_1(-1)$}; 
\node at (2,-5) {$B_0(-1)$};

\node at (4,-5.5) {$R$};
\node (aa) at (4,-2.5) {$B_{21}$};
\node at (4,-3.5) {$B_{10}$};
\node at (4,-4.5) {$B_{02}$};
\node (c1) at (4,-4) {$\;\;\;$};

\node (xx) at (6,-3) {$B_2(1)$};
\node (d1) at (6,-4) {$B_1(1)$};
\node at (6,-5) {$B_0(1)$};
\node (e1) at (8,-4) {$R(2)$};
\path
	(a1) edge (b1)
	(b1) edge (c1)
	(c1) edge (d1)
	(d1) edge (e1);

\path
	(a) edge node[left] {$\sqmatrix{0 \\ 0 \\ 0 \\ \id_R}$} (aa)
	(x) edge node[descr] {$\id$} (xx)
	(y) edge node[right] {$\sqmatrix{0 & 0 & 0 & \id_R}$} (e1);

\end{tikzpicture} \right)
\end{equation}
We centered this picture at homological degree $-1$. This chain map in homological degree $-1$ is the identity map
\[ \id = \sqmatrix{\lineblue & 0 & 0\\ 0 & \linered & 0 \\ 0 & 0 & \linegreen}. \]
Note that $d^2_{\FC} = \mu \delta$.

\begin{remark} Of course, $\ImM$ also has a nonzero monodromy map, which the reader can readily construct. This monodromy map is compatible with the monodromy map on $\VC$, for either the subcomplex $\ImM(1)[-1]$ or the quotient complex $\ImM(-1)[1]$. This is the expected behavior whenever $\mu^2 \ne 0$. \end{remark}

Let us examine $\FC B_i$ for $i = 0,1,2$. Recall that $B_i B_{i \pm 1} B_i \cong B_w \oplus B_i$ where $w = s_i s_{i \pm 1} s_i$. So, if one only looks at the instances of (shifts of) $B_0$ inside $\FC B_0$, one will have
\begin{equation} \label{eq:FB0n3}\left(
\begin{tikzpicture}
\node (b) at (2,0) {$B_0(-2)$};
\node (c) at (4,.25) {$B_0(-2)$};
\node at (4,-.25) {$B_0(0)$};
\node (d) at (6,.25) {$B_0(0)$};
\node at (6,-.25) {$B_0(0)$};
\node (e) at (8,.25) {$B_0(0)$};
\node at (8,-.25) {$B_0(2)$};
\node (g) at (10,0) {$B_0(2)$};
\end{tikzpicture}
\right). \end{equation}
These summands cancel in pairs using Gaussian elimination. Similarly, the summands of the form $B_{10}$ will also cancel in pairs. What remains then is
\begin{equation} \FC B_0 \simeq \left(
\begin{tikzpicture}
\node (a) at (0,0) {$B_{20}(-1)$};
\node (b) at (2,2) {$B_{210}(0)$};
\node (c) at (2,-2) {$B_{w}(0)$};
\node (d) at (4,0) {$B_{20}(1)$};
\path
	(a) edge node[descr] {$- \lineblue \startdotred \linegreen$} (b)
	(a) edge node[descr] {$\dotfirsttoteal$} (c)
	(b) edge node[descr] {$- \lineblue \finaldotred \linegreen$} (d)
	(c) edge node[descr] {$\dotfirstfromteal$} (d);
\end{tikzpicture}
\right). \end{equation}
Here, $w = s_0 s_2 s_0$, and is represented by the color teal. The differentials involving $B_w$ are drawn using the thick calculus of \cite{EThick}. In particular, these differentials are dots composed with the inclusion or projection between $B_w$ and the Bott-Samelson bimodule $B_0 B_2 B_0$. The reader is encouraged to check that $d^2$ is equal to $\delta$ times the identity map of $B_{20}$. Note that we have chosen a representative of this complex which also satisfies a reasonable analog of the strand-counting sign rule.

For the reader's convenience, we include two helpful relations from the thick calculus.
\begin{equation} \ig{.5}{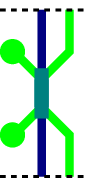} = \barbgreen \lineblue \linegreen + \brokenblue \linegreen. \end{equation}
\begin{equation} \dotfirstfromteal = \ig{.5}{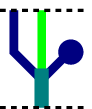}. \end{equation}
We leave it as an exercise to verify that Gaussian elimination yields the complex above.

Right multiplication of any Soergel bimodule by $B_0$ should produce only summands of the form $B_x$, where $x s_0 < x$. The key observation here is that the only summands which survive
to the minimal complex of $\FC B_0$ also satisfy $s_2 x < x$. If there is to be an isomorphism $(-) \ot \FC \to \FC \ot \tau(-)$ then $\FC B_0 \cong B_2 \FC$, which explains this phenomenon. In fact, a similar computation will show that $B_2 \FC$ is isomorphic to the same complex \eqref{eq:FB0n3}.

Thus one proves that $\FC B_k \simeq B_{k-1} \FC,$ from which it will follow that $\VC B_k \simeq B_k \VC$. Again, these can be lifted to functorial isomorphisms just as in
\eqref{eq:Ffunctorial1} and \eqref{eq:Vfunctorial1}.

We leave the expert reader as an exercise to continue this computation and prove that
\begin{equation} \label{eq:n3FB121} \FC B_{121} \simeq B_{10121}(0), \end{equation}
which is an indecomposable complex supported in homological degree zero. Similarly, $B_{101} \FC$ is isomorphic to the same complex.

% The computation of $\FC B_0$ is essentially the same as the calculation of $\FC$ after restriction on the right to $R^{s_0}$. After all, $B_0 \cong R \ot_{R^{s_0}} R(1)$, so tensoring
% with $B_0$ is the composition of (shifted) restriction and induction. One can also prove that induction on the right from $R^{s_0}$ back to $R$ is a functor which preserves
% indecomposables, and acts fully faithfully on morphisms after base change from $R^{s_0}$ to $R$, see \cite{Will-Singular}. So, just like the complex in \eqref{eq:FB0n3}, $\FC$ restricted
% to $R^{s_0}$ will have four indecomposable summands arranged in a diamond. After inducing back from $R^{s_0}$ to $R$, each indecomposable summand remains indecomposable, and one obtains
% the complex in \eqref{eq:FB0n3}.
%
% \begin{remark} This is all relatively easy to understand once one has spent the time learning singular Soergel bimodules and their diagrammatics. However, we will avoid this additional
% piece of overhead in this paper as much as possible. The discussion of the above paragraph should point to the following recourse for the mathematician unfamiliar with singular Soergel
% bimodules: if one is interested in what happens to $\FC$ after restriction to $R^{s_0}$, one can instead perform the computation $\FC B_0$ purely in the category of ordinary Soergel
% bimodules, and the result is effectively ``the same.'' \end{remark}

Let $W_{\fin}$ be generated by $\{s_1, s_2\}$, and let us abbreviate $R^{s_1,s_2}$ as $R^{12}$. In line with Remark \ref{rmk:restricttensorsimilar}, we claim that \eqref{eq:n3FB121} is equivalent to the statement that $\FC$ restricted to $R^{W_{\fin}}$ is a complex supported in homological degree zero, whose only term is the $(R,R^{W_{\fin}})$-bimodule \begin{equation} \label{eq:Vrest3} \FC_{R^{W_{\fin}}} \cong (R \ot_{R^{01}} R^{1})_{R^{12}}(5). \end{equation}
After inducing back to $R$, one obtains $B_{10121} = R \ot_{R^{01}} R^{1} \ot_{R^{12}} R(5)$. Note that the parabolic subgroup generated by $\{s_0,s_1\}$ is precisely $\tau\inv(W_{\fin})$. 

\begin{remark} For the reader interested in more details here is a step by step guide to this computation. Restricting $R$ to $R^{W_{\fin}}$ gives (unsurprisingly) one copy of $R$. Restricting $B_1$ or $B_2$ to $R^{W_{\fin}}$ gives (quantum) two copies of $R$, i.e. $R(1) \oplus R(-1)$. Restricting $B_{21}$ gives $R(-2) \oplus R(0) \oplus R(0) \oplus R(2)$. Finally, restricting $B_{10}$ to $R^{W_{\fin}}$ has a direct summand isomorphic to $R$ (this is analogous to the direct summand $B_1$ inside $B_{10} B_1$). Together, all these copies of $R$ cancel under Gaussian elimination.

Restricting $B_0$ to $R^{W_{\fin}}$ gives (unsurprisingly) one copy of $B_0$, while restricting $B_{02}$ to $R^{W_{\fin}}$ gives quantum two copies of $B_0$. These copies of $B_0$ all cancel under Gaussian elimination.

All that remains after Gaussian elimination is the other summand of $B_{10}$ restricted to $R^{W_{\fin}}$, a complex supported in homological degree $0$. We claim that this other summand is precisely the one in \eqref{eq:Vrest3}. \end{remark}

Meanwhile,
\begin{equation} \SSat(V) \cong {}_{R^{01}} R^{1}_{R^{12}}(2). \end{equation}
Thus, as for $n=2$, we have
\begin{equation} R \ot_{R^{W_{\fin}}} (\Om \SSat(V)) \cong \VC_{R^{W_{\fin}}}. \end{equation}
	
Finally, we encourage the reader to examine the quotient by $\MC$, and observe that $\overline{\VC}$ agrees with the underlying vector space of the three-dimensional standard
representation $V$, and the monodromy map descends to the principal nilpotent operator.

%===========
\subsection{The case $n=4$}
\label{subsec-stdrep4}
%===========

While we will not repeat the calculations above for the specific case $n=4$, it is important to discuss it because of an important new subtlety which arises.

Fix a realization in type $\tilde{A}_3$ with simple reflections $\{s_0,s_1,s_2,s_3\}$, and let $\al_0 + \al_1 + \al_2 + \al_3 = \delta$ in $\hg^*$. Let blue denote $s_3$, red $s_2$, green $s_1$, and purple $s_0$.

\begin{defn} \label{def:F4} ($n=4$) Let $\FC \in \KC^b(\Diag)$ be the complex given by
\begin{equation} \label{eq:F4} \FC = \left(
\scalebox{.7}{
\begin{tikzpicture}
\node (a) at (0,0) {$R(-3)$};

\node (a0) at (3,3) {$B_3(-2)$};
\node (a1) at (3,1) {$B_2(-2)$};
\node (a2) at (3,-1) {$B_1(-2)$};
\node (a3) at (3,-3) {$B_0(-2)$};

\node (b) at (6,-6) {$R(-1)$};
\node (a01) at (6,6) {$B_{32}(-1)$};
\node (a12) at (6,4) {$B_{21}(-1)$};
\node (a23) at (6,2) {$B_{10}(-1)$};
\node (a30) at (6,0) {$B_{03}(-1)$};
\node (a02) at (6,-2) {$B_{31}(-1)$};
\node (a13) at (6,-4) {$B_{20}(-1)$};

\node (a012) at (9,7) {$B_{321}$};
\node (a123) at (9,5) {$B_{210}$};
\node (a230) at (9,3) {$B_{103}$};
\node (a301) at (9,1) {$B_{032}$};
\node (b0) at (9,-1) {$B_3$};
\node (b1) at (9,-3) {$B_2$};
\node (b2) at (9,-5) {$B_1$};
\node (b3) at (9,-7) {$B_0$};

\node (c) at (12,-6) {$R(1)$};
\node (b01) at (12,6) {$B_{32}(1)$};
\node (b12) at (12,4) {$B_{21}(1)$};
\node (b23) at (12,2) {$B_{10}(1)$};
\node (b30) at (12,0) {$B_{03}(1)$};
\node (b02) at (12,-2) {$B_{31}(1)$};
\node (b13) at (12,-4) {$B_{20}(1)$};

\node (c0) at (15,3) {$B_3(2)$};
\node (c1) at (15,1) {$B_2(2)$};
\node (c2) at (15,-1) {$B_1(2)$};
\node (c3) at (15,-3) {$B_0(2)$};

\node (d) at (18,0) {$R(3)$};

\path 
	(a) edge (a0)
	(a) edge (a1)
	(a) edge (a2)
	(a) edge (a3)
	(a1) edge (b)
	(a2) edge (b)
	(a3) edge (b)
	(a0) edge (b)
	(b) edge (b0)
	(b) edge (b1)
	(b) edge (b2)
	(b) edge (b3)
	(b1) edge (c)
	(b2) edge (c)
	(b3) edge (c)
	(b0) edge (c)
	(c) edge (c0)
	(c) edge (c1)
	(c) edge (c2)
	(c) edge (c3)
	(c1) edge (d)
	(c2) edge (d)
	(c3) edge (d)
	(c0) edge (d)

	(a0) edge (a01)
	(a0) edge (a02)
	(a0) edge (a30)
	(a1) edge (a01)
	(a1) edge (a12)
	(a1) edge (a13)
	(a2) edge (a02)
	(a2) edge (a12)
	(a2) edge (a23)
	(a3) edge (a30)
	(a3) edge (a13)
	(a3) edge (a23)
	(b0) edge (b01)
	(b0) edge (b02)
	(b0) edge (b30)
	(b1) edge (b01)
	(b1) edge (b12)
	(b1) edge (b13)
	(b2) edge (b02)
	(b2) edge (b12)
	(b2) edge (b23)
	(b3) edge (b30)
	(b3) edge (b13)
	(b3) edge (b23)
	(a01) edge (b0)
	(a02) edge (b0)
	(a30) edge (b0)
	(a01) edge (b1)
	(a12) edge (b1)
	(a13) edge (b1)
	(a02) edge (b2)
	(a12) edge (b2)
	(a23) edge (b2)
	(a30) edge (b3)
	(a13) edge (b3)
	(a23) edge (b3)
	(b01) edge (c0)
	(b02) edge (c0)
	(b30) edge (c0)
	(b01) edge (c1)
	(b12) edge (c1)
	(b13) edge (c1)
	(b02) edge (c2)
	(b12) edge (c2)
	(b23) edge (c2)
	(b30) edge (c3)
	(b13) edge (c3)
	(b23) edge (c3)

	(a01) edge (a301)
	(a01) edge (a012)
	(a12) edge (a012)
	(a12) edge (a123)
	(a23) edge (a123)
	(a23) edge (a230)
	(a30) edge (a230)
	(a30) edge (a301)
	(a02) edge (a012)
	(a02) edge (a230)
	(a13) edge (a123)
	(a13) edge (a301)
	(a301) edge (b01)
	(a012) edge (b01)
	(a012) edge (b12)
	(a123) edge (b12)
	(a123) edge (b23)
	(a230) edge (b23)
	(a230) edge (b30)
	(a301) edge (b30)
	(a012) edge (b02)
	(a230) edge (b02)
	(a123) edge (b13)
	(a301) edge (b13);
	
\end{tikzpicture} } \right)
\end{equation} \end{defn}

All the summands $B_{ij}$ and $B_{ijk}$ are Bott-Samelson bimodules, which are identified with the corresponding indecomposable Soergel bimodules. Note that when the indices $2$ and $1$
both appear, $2$ is always before $1$. Similarly $1$ is always before $0$, which is always before $3$, which is always before $2$. The indices are ordered reverse-cyclically; this order makes sense for any proper subset of $S_{\aff}$. When one applies $\si$ to compute $\VC^*$ or $\VC_{\fund_{n-1}}$, the indices will be ordered cyclically instead.

For reasons of sanity, we do not draw the differentials, so we must discuss them in words. Every edge above corresponds to the addition or the removal of a single simple reflection $i$,
which is accomplished with an $i$-colored dot ($\startdotred$ or $\finaldotred$), multiplied by a sign for each identity map $\lineblue$ to the left of the dot. However, there is one
additional subtlety, which is the precise meaning of the summands $B_{31}$ and $B_{20}$ above. These are both indecomposable Soergel bimodules which can be identified with two different
(though isomorphic) Bott-Samelson bimodules.

Recall that $s_3$ and $s_1$ commute, so that $B_3 \ot B_1 \cong B_1 \ot B_3 \cong B_{s_3 s_1}$. There is supposed to be an isomorphism $\FC \cong \tau(\FC)$, and hence an isomorphism
$\phi \co \FC \to \tau^2(\FC)$. Then $\phi$ must send $B_{s_3 s_1}$ to itself by some automorphism. Said another way, suppose that we represent $B_{s_3 s_1}$ as $B_3 \ot B_1$ in $\FC$,
and as $B_1 \ot B_3$ in $\tau^2(\FC)$. This would be a natural choice, since whatever differentials appeared in $\FC$ could be color-rotated to give the differentials in $\tau^2(\FC)$.
Then $\phi$ induces some map $B_3 \ot B_1 \to B_1 \ot B_3$.

A naive choice for this isomorphism would be the colored 4-valent vertex $\Xbg$, one of the generating morphisms of $\Diag$ (not to be confused with the mixed crossing $\mixedOred$ in
$\Diag_{\ext}$). However, we shall be forced to use the 4-valent vertex multiplied by a sign! This is forced by our sign conventions on dots, which in turn are required for cancellation
in the computation of $d^2$. The differential from $B_1$ to $B_3 \ot B_1$ should be $\startdotblue\linegreen$, while the differential from $B_3$ to $B_3 \ot B_1$ is
$-\lineblue\startdotgreen$. After rotation by $\tau^2$, these two terms are swapped, so a compensating sign must be placed in the isomorphism between $B_3 \ot B_1$ and $B_1 \ot B_3$. If
one places the 4-valent vertex on top of these differentials, multiplies by $-1$, and then swaps the colors green and blue, then these differentials are exchanged with each other.
 
So, for example, let us always represent $B_{31}$ in $\FC$ as $B_3 \ot B_1$. Then the differential from $B_{31}$ to $B_{321}$ can be
efficiently encoded as the morphism $B_3 \ot B_1 \mapsto B_3 \ot B_2 \ot B_1$ which adds a signed $2$-colored dot. (The sign in this differential comes because there are an odd number of identity maps to the left of the dot.) \[ -\lineblue\startdotred\linegreen\]  However, the differential from $B_{31}$ to $B_{103}$ should be obtained as the composition $B_3 \ot B_1 \to B_1 \ot B_3 \to B_1 \ot B_0 \ot B_3$, which first applies the 4-valent vertex with a sign $(-1)$, then adds a $0$-colored dot with sign $(-1)$. \[ -\linegreen\startdotpurple\lineblue \circ - \Xbg = + \ig{.5}{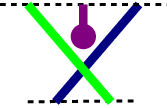}\]

With these conventions, one will obtain a pseudocomplex for which $\FC \cong \tau(\FC)$, and where $d^2=\mu \delta$, for the monodromy map $\mu$ described below.

Though the complex $\FC$ is rather large, examining each summand of $d^2$ is a reasonable exercise. As before, there are really only three observations needed:
\begin{itemize} \item Our sign conventions cause certain terms to cancel in pairs. \item For certain differentials, one must use the barbell forcing relations \eqref{eq:barbellforce} to push all polynomials to one side of the diagram. \item Once this is done, the only terms in $d^2$ which are not obviously zero are maps $B_x(k) \to B_x(k+2)$, which happen to be of the form $(\al_0 + \al_1 + \al_2 + \al_3) \cdot \id_{B_x} = \delta \cdot \id_{B_x}$.  \end{itemize}
We encourage the reader to do this exercise before reading the general proof in \S\ref{subsec-dotmaps} and \S\ref{subsec-stdrep}.

As before, we let $\VC = \Om \FC$ and $\VC^* = \si(\VC) = \Om^{-1} \si(\FC)$ and $\VC_{\fund_3} = \Om^3 \si(\FC)$. We leave the reader to compute the Wakimoto filtrations, using Example \ref{ex:wakimoton4} as a guide.

The monodromy map $\mu \co \FC \to \FC(-2)[2]$ will again be defined as a composition
\begin{equation*} \FC \to \ImM(-1)[1] \to \FC(-2)[2] \end{equation*}
where the complex $\ImM$ satisfies the strand-counting sign rule, and has objects
\begin{equation} \ImM = \left(
\scalebox{.7}{\begin{tikzpicture}
\node (a) at (0,0) {$R(-2)$};

\node at (2,1.5) {$B_3(-1)$};
\node at (2,.5) {$B_2(-1)$};
\node (b) at (2,0) {$\;\;\;\;\;$};
\node at (2,-.5) {$B_1(-1)$};
\node at (2,-1.5) {$B_0(-1)$};

\node at (4,-3) {$R$};
\node at (4,3) {$B_{32}$};
\node at (4,2) {$B_{21}$};
\node at (4,1) {$B_{10}$};
\node (c) at (4,0) {$B_{03}$};
\node at (4,-1) {$B_{31}$};
\node at (4,-2) {$B_{20}$};

\node at (6,1.5) {$B_3(1)$};
\node at (6,.5) {$B_2(1)$};
\node (d) at (6,0) {$\;\;\;\;\;$};
\node at (6,-.5) {$B_1(1)$};
\node at (6,-1.5) {$B_0(1)$};

\node (e) at (8,0) {$R(2)$};

\path
	(a) edge (b)
	(b) edge (c)
	(c) edge (d)
	(d) edge (e);
\end{tikzpicture} }
\right). \end{equation}
%%%%%%%%%%%%%%%%%%%%%%%%%
\section{Elaborations and simplifications}
\label{sec-elaboratesimplify}
%%%%%%%%%%%%%%%%%%%%%%%%%

In this chapter, an extension of the introduction, we go over in detail the desired properties of Gaitsgory's central complexes, and how one might go about proving them. In this chapter $(W,S)$ is the Coxeter system for the affine Weyl group in type $\tilde{A}_n$.

%================
\subsection{Gaussian elimination and minimal complexes} \label{subsec:someconventions}
%================

We wish to remind the reader briefly of some basic notions in the homological algebra of additive categories. We work in a (graded) Krull-Schmidt category, where every object splits into
indecomposable objects. Note that this splitting is not canonical in general, so that when we talk about the indecomposable summands of an object we implicitly choose a splitting. Given
a complex $\MC$, the object $\MC^k$ in homological degree $k$ is a direct sum of indecomposables, as is the object $\MC^{k+1}$ in degree $k+1$, so that the differential $d \co \MC^k \to
\MC^{k+1}$ can be expressed as a matrix of morphisms between indecomposable objects. Any of these morphisms which appear as matrix coefficients we refer to as a \emph{summand of the
differential}. Similarly, for a chain map $f \co \MC \to \NC$ between two complexes, the map $f^k \co \MC^k \to \NC^k$ in homological degree $k$ can be expressed as a matrix of
morphisms, and its matrix coefficients are called \emph{summands of $f$}. To reiterate, taking a summand of a differential or a chain map makes sense only after one has chosen a
splitting of each object $\MC^k$ into indecomposables.

When a summand of the differential is an isomorphism, from an indecomposable summand $X \sumset \MC^k$ to an indecomposable summand $Y \sumset \MC^{k+1}$, there is a process called
\emph{Gaussian elimination} which replaces the complex $\MC$ with a homotopy equivalent complex $\MC'$ in which the summands $X$ and $Y$ have been removed or \emph{eliminated}. The point
of Gaussian elimination is that, while it may not appear so at first, one can alter the direct sum decompositions of $\MC^k$ and $\MC^{k+1}$ so that the two term complex $X \to Y$ is a
direct summand of $\MC$, and this two term complex is contractible.

A complex for which no summand of the differential is an isomorphism is called a \emph{minimal complex}. By repeated Gaussian elimination, one can prove that any complex is homotopy
equivalent to a minimal complex, and in fact this minimal complex is unique up to isomorphism. We will often refer to objects of the homotopy category $\HC_{\ext}$ abusively as complexes,
but we usually have in mind the minimal complex as our favorite representative of its homotopy equivalence class.

For more background on minimal complexes see \cite[\S 6.1]{EWHodge}. We recall the specifics of Gaussian elimination in \S\ref{sec-pseudo}, where we also discuss Gaussian elimination for
pseudocomplexes.

%================
\subsection{Perversity}
%================

The graded monoidal category $\DC_{\ext}$ comes with indecomposable objects $\{B_w\}$ parametrized by $w \in W_{\ext}$. These are all the indecomposable objects, up to isomorphism and
grading shift. Moreover, these indecomposables are special (unlike their shifts) as they are isomorphic to their duals (under a certain duality functor). Of particular note is the fact
that \begin{equation} \label{eq:gradedschur}\Hom(B_x,B_y(k)) = \begin{cases} \C \cdot \id & \text{if } x = y \text{ and } k = 0, \\ 0 & \text{if } k < 0, \text{ or } k = 0
\textrm{ and } x \ne y \\ ?? & \text{if } k > 0. \end{cases} \end{equation} That is, the only nonzero maps between self-dual indecomposables of non-positive degree are scalar
multiples of identity maps. This is a graded analog of Schur's lemma.

A complex is called \emph{perverse} if the indecomposable objects appearing as summands in homological degree $k$ all have the form $B_w(k)$ for various $w \in W_{\ext}$. This is a very
powerful condition. For example, there are no nonzero homotopies between perverse complexes, by \eqref{eq:gradedschur}.

\begin{remark} If $\FC$ and $\GC$ are perverse then $\Hom(\FC, \GC[c])=0$ for any $c < 0$, using \eqref{eq:gradedschur}. However, it is not true that $\Hom(\FC, \GC(k)[c])=0$ for $c <
0$ and arbitrary $k \in \Z$. Thus the vanishing of negative exts, property (m) from Conjecture \ref{conj3}, is not directly related to perversity. \end{remark}

Here is another useful lemma which illustrates the utility of perversity. 

\begin{defn} A complex is \emph{multiplicity-free perverse} if it is perverse, and if for each $k \in \Z$ and $w \in W_{\ext}$ the indecomposable object $B_w(k)$ appears in homological degree $k$ with multiplicity at most $1$. \end{defn}

\begin{lemma} \label{lem:multfreejustid} Let $\MC$ be an indecomposable multiplicity-free perverse complex in $\HC_{\ext}$. Then $\Hom(\MC,\MC)$ is spanned by the identity map. \end{lemma}

\begin{proof} For each $k \in \Z$ choose a splitting of $\MC^k$ into indecomposable objects. Consider a chain map $f \co \MC \to \MC$. By \eqref{eq:gradedschur}, each summand of $f$ must
be either zero or a scalar multiple of the identity of $B_w(k)$ for some summand $B_w(k)$ in $\MC^k$. Let $\lambda_{w,k}$ be the scalar such that $f$ restricted to $B_w(k)$ is
$\lambda_{w,k} \id_{B_w(k)}$ (if such a summand does appear in $\MC$). If the summand of the differential between $B_w(k)$ and $B_{w'}(k+1)$ is nonzero, then $\lambda_{w,k} =
\lambda_{w',k+1}$ or else $f$ will not commute with the differential.

Draw a graph where the vertices are the indecomposable summands of $\MC$, and the edges are nonzero summands of the
differential. If this graph is connected, then $f$ must induce the same scalar multiple of the identity on all summands, so $f$ is a scalar multiple of the identity map of $\MC$. If this
graph is disconnected then each connected component is a direct summand of the complex $\MC$, so $\MC$ is not indecomposable. \end{proof}

The graph in the proof above does depend on the choice of splitting of each term $\MC^k$ into indecomposable objects. There are decomposable complexes with connected graphs, although any
decomposable complex has a disconnected graph for some choice of splitting. We call a complex \emph{apparently indecomposable} if, with respect to some splitting, its graph is connected.

\begin{porism} If $\MC$ is multiplicity-free perverse, then it is indecomposable if and only if it is apparently indecomposable. \end{porism}
	
\begin{proof} By the same proof as in Lemma \ref{lem:multfreejustid}, if $\MC$ is apparently indecomposable then any endomorphism of $\MC$ is a scalar multiple of the identity. But any decomposable complex has at least a two-dimensional endomorphism ring, so $\MC$ must be indecomposable. \end{proof}

An object of $\HC_{\ext}$ is called \emph{perverse} if it is isomorphic (i.e. homotopy equivalent) to a perverse complex. Equivalently, an object of $\HC_{\ext}$ is perverse if its
minimal complex is perverse. Note that any perverse complex is automatically minimal.

It is easy to decide whether a complex constructed by hand, and built out of indecomposables, is perverse: merely analyze the summands used. We will construct our complexes $\VC$ and
$\FC$ such that they are multiplicity-free perverse. However, it is quite difficult to determine whether a complex given abstractly (say, as a tensor product of other complexes) is
perverse. In our main technical result, Theorem \ref{thm:tensorBI}, we give an explicit description of the minimal complex of $B_I \ot \FC$ and $\FC \ot B_J$ thus proving that it is
multiplicity-free perverse.

%================
\subsection{The Wakimoto filtration} \label{subsec-Wakimoto}
%================

We assume that the reader has read \S\ref{subsec-wakimoto} and is familiar with Wakimoto complexes. We let $W(\la)$ denote the Wakimoto complex of $\la \in \La_{\wt}$. Let us recall two key properties of Wakimoto complexes. First, they form a categorification of the weight lattice.
\begin{equation} \label{eq:Wakwtlattice} W(\la) \ot W(\mu) \cong W(\la + \mu) \cong W(\mu) \ot W(\la). \end{equation}
Moreover, these isomorphisms are canonical, see \S\ref{subsec-wakimoto}. Consequently, Wakimoto complexes tensor-commute with each other, and are closed under taking tensor products, forming a symmetric monoidal category. We also note that the dual of a Rouquier complex is the inverse Rouquier complex, so that the dual of $W(\la)$ is $W(-\la)$.

Secondly, Wakimoto complexes are perverse.

\begin{lemma} For any $\la \in \La_{\wt}$, the complex $W(\la)$ is perverse in $\HC_{\ext}$. \end{lemma}

Here is a reminder of the proof, see \S\ref{subsec-wakimoto} for details. It is a difficult theorem, proven in \cite[Theorem 6.9]{EWHodge}, that Rouquier complexes associated with
positive lifts of reduced expressions are perverse. By duality, the same is true for negative lifts of reduced expressions, and by comparing these results, the same is true for braids
with a positive-negative decomposition. Wakimoto complexes are Rouquier complexes for such braids, so they are perverse.

We wish to say that Gaitsgory central complexes have Wakimoto filtrations, but first we should discuss what a filtration means in the homotopy category. We begin with the naive definition
of a convolution of chain complexes.

\begin{defn} \label{defn:complexfiltration} We work in the category of complexes over some additive category. Given complexes $(\NC_i,d_i)_{i \in I}$ indexed by a finite poset $I$, a \emph{convolution (thereof)} is a complex $(\MC = \oplus_{i \in I} \NC_i, d =\sum_{i,j \in I} d_{ij})$ where $d_{ij} \co \NC_j \to \NC_i$, $d_{ii} = d_i$, and $d_{ij} = 0$ for $j<i$. \end{defn}

In other words, the complex $(\MC,d)$ is like the direct sum complex $\oplus \NC_i$ but with some extra upper-triangular pieces to the differential. It has a natural $I$-indexed
filtration by subcomplexes, which in each homological degree is a filtration by direct summands, and for which the subquotient complexes are $\NC_i$. Conversely, every such filtered
complex $\MC$ can be built as a convolution of its subquotient complexes. Convolutions are just an efficient way to describe iterated mapping cones. For more details on
convolutions see \cite[\S 4]{EHDiag}.

\begin{defn} \label{defn:Wakimotofiltration} An object $\MC \in \HC_{\ext}$ is said to have a \emph{Wakimoto filtration} if its minimal complex is a convolution of the minimal complexes
of Wakimoto complexes. Objects with Wakimoto filtrations form a full subcategory $\Wak_0$ of $\HC_{\ext}$. \end{defn}

Since the minimal complexes of Wakimoto complexes are perverse, so is any object with a Wakimoto filtration. Remark \ref{rmk:whynaiveisfine} below will discuss why this naive definition
of Wakimoto filtrations is still reasonable in the homotopy category. We called the category above $\Wak_0$, so as to distinguish it from the full triangulated subcategory $\Wak$ generated by Wakimoto sheaves, which also permits homological shifts.

\begin{lemma} \label{lem:wak0closed} The subcategory $\Wak_0$ is closed under homotopy equivalence, tensor product, direct sums, duals. \end{lemma}

\begin{proof} Closure under homotopy equivalence and direct sums is obvious. A tensor product of filtered complexes is filtered by the tensor products of the subquotients; since Wakimoto
complexes are closed under tensor product, so is $\Wak_0$. Wakimoto complexes are also closed under taking duals (which reverses the order of a filtration), so $\Wak_0$ is closed under
duality. \end{proof}

\begin{remark} The subcategory $\Wak_0$ should also be closed under taking direct summands. However, this argument is not as trivial as it might seem. If $\MC$ is Wakimoto filtered and
$\NC \sumset \MC$, then the filtration on $\MC$ does descend to a filtration on $\NC$ but the subquotients in this filtration need not (a priori) be Wakimoto complexes. In other settings
there are tools for proving that direct summands of ``standardly filtered objects'' have standard filtrations. In the context of Soergel bimodules, there are functors which pick out the
various pieces in a standard filtration, so these filtrations must interact well with idempotents. In the highest weight category context, having a standard filtration is equivalent to
an Ext vanishing property, which is preserved under direct summands. Presumably some similar argument can be applied here, though we are not ready to provide it. Thanks to Geordie
Williamson and Sasha Polishchuk for discussion about this issue. \end{remark}

It was explained to me by Bezrukavnikov that Gaitsgory central complexes have Wakimoto filtrations; I do not believe this appears in Gaitsgory's paper. The filtration of $\GC(U)$ by
Wakimoto complexes is analogous to the splitting of $U$ into weight spaces. To be more precise, if the $\gl_n$-module $U$ is viewed as a module over the positive half $\bg^+$ of $\gl_n$,
then it has a filtration whose subquotients are weight spaces. We call this the \emph{weight space filtration} on $U$, and it should be lifted by the Wakimoto filtration. Thus the
multiplicity of $W(\la)$ in $\GC(U)$ agrees with the dimension of the $\la$ weight space in $U$, for $\la \in \La_{\wt}$, and the filtration respects the dominance order on weights. The following Lemma indicates an efficient way to prove these facts.

\begin{lemma} Suppose that a functor $\GC \co \Webs^+_n \to \ZC(\HC_{\ext})$ is constructed as in Goal \ref{goal1.5}. Suppose that the complexes $\VC_{\fund_i}$ for fundamental representations live in $\Wak_0$, and their Wakimoto filtration lifts their weight space filtration. Suppose one can prove that $\Wak_0$ is closed under taking direct summands.  Then for any representation $U$ of $\gl_n$, $\GC(U)$ lives in $\Wak_0$ and its Wakimoto filtration lifts its weight space filtration. \end{lemma}

\begin{proof} By Lemma \ref{lem:wak0closed}, the Gaitsgory complex for any tensor product of fundamental representations lives in $\Wak_0$. By \eqref{eq:Wakwtlattice} we can deduce that the filtration on a tensor product of fundamental representations lifts its weight space filtration. These two properties are clearly preserved under taking direct sums and duals, so it is enough to prove they hold for irreducible representations. 

Any irreducible representation $U$ appears as the top direct summand of a tensor product $T$ of fundamental representations, where all the other summands have strictly lower highest
weight. Since $\GC(U) \sumset \GC(T)$, it lives in $\Wak_0$ and its Wakimoto filtration respects the dominance order, so we must only show that the Wakimoto multiplicities agree with the
weight multiplicities. But this is true for $\GC(T)$ and true for all other direct summands by induction, so it is also true for $\GC(U)$. \end{proof}

\begin{remark} \label{rmk:whynaiveisfine} Convolutions are well-behaved under homotopy equivalence under some assumptions on the indexing poset of the filtration: replacing $\NC_i$ by
homotopy equivalent complexes $\NC_i'$, one can construct a complex $\MC'$ which is homotopy equivalent to $\MC$ and filtered by $\NC_i'$. See \cite[Proposition 4.20]{EHDiag} for a more
precise statement. When one works in the homotopy category one often discusses convolutions up to homotopy; let us say that $\MC$ is a \emph{homotopy convolution} of $\NC_i$ if $\MC$ is
homotopy equivalent to a complex $\MC'$ which is a convolution of complexes $\NC_i'$ which are homotopy equivalent to $\NC_i$. There is an important difference between convolutions and
homotopy convolutions: when $\MC$ is a true convolution of $\NC_i$ its size (in the underlying additive category) is the same as $\bigoplus N_i$; one adds the complexes together and then
adjusts the differential in an upper-triangular fashion. Meanwhile, if $\MC$ is only a homotopy convolution then the size of $\MC$ can be much smaller than $\bigoplus N_i$, as the adjusted
differential can lead to cancellation (via Gaussian elimination) between terms of $\NC_i$ and terms of $\NC_j$ in the total complex. True convolution gives control over the size of
complexes, while homotopy convolution is more natural in the context of homotopy categories.

Suppose however that $\MC$ is a homotopy convolution of perverse objects $\NC_i$. Then letting $\NC_i'$ be the (perverse) minimal complex of $\NC_i$, there is some convolution $\MC'$ of the $\NC_i'$
which is homotopy equivalent to $\MC$. However, any convolution of perverse complexes is perverse, so $\MC'$ is perverse, and perverse complexes are minimal complexes. Consequently $\MC'$ is
the minimal complex of $\MC$, and the minimal complexes of any perverse homotopy convolution are just true convolutions of the minimal complexes. Both of these notions agree with the
notion of a filtration in the abelian core of the perverse $t$-structure. \end{remark}

%================
\subsection{Preserving perversity}
%================

We say that a complex $\MC$ \emph{preserves perverses} if $\MC \ot \PC$ is perverse in $\HC_{\ext}$ whenever $\PC$ is perverse. In particular, $\MC$ must be perverse, since the monoidal
identity $\one$ is perverse.

Let us use shorthand: we write $PP(\MC,\PC)$ if $\MC \ot \PC$ is perverse (it is assumed that $\PC$ is perverse). For $U \in \Rep \gl_n$ we write $PP(\MC,U)$
for $PP(\MC,\GC(U))$. We write $PP(\MC)$ if $\MC$ preserves perverses.

\begin{lemma} The property of preserving perverses is closed under tensor products, direct sums, and direct summands. Consequently, if $PP(\VC_{\fund_i})$ holds for each fundamental
representation, then $PP(\GC(U))$ holds for any representation $U$ of $\gl_n$. \end{lemma}

\begin{proof} This is easy. Note that a direct summand of a perverse complex is perverse. \end{proof}

Preserving perverses is a very special property; for example, we do not expect any (non-identity, non-zero) complex to satisfy this property in $\HC_{\fin}$. Gaitsgory proves in the
geometric context that his central sheaves all preserve perverses. Unfortuntely, we do not know how to prove this yet. Let us indicate what we do know how to prove, for the complex $\VC$.

Any complex $\MC$ in $\Wak_0$ satisfies $PP(\MC,W(\la))$ for all $\la \in \La_{\wt}$, by Lemma \ref{lem:wak0closed}. Thus $\VC$ preserves the perversity of Wakimoto complexes. From this
one can already deduce that all tensor products $\VC^{\ot k}$ and their direct summands are perverse. Our main technical theorem, Theorem \ref{thm:tensorBI}, proves that $\VC$ preserves
the perversity of longest elements, i.e. $PP(\VC,B_J)$ for all $J \subsetneq S$. This can be bootstrapped into the following result.

\begin{prop} \label{prop:lowestcell} The property $PP(\VC,B_w)$ holds for any $w$ in the lowest two-sided Kazhdan-Lusztig cell of $W$. \end{prop} 

This result is closely related to the compatibility of $\VC$ with the geometric Satake equivalence, and with the properties of Lusztig's $a$-fuction as they apply to singular Soergel
bimodules. A proof would involve a great deal of additional exposition on these ideas, and it is not clear that this intermediate result is useful, so the reader is welcome to contact
the author for more details.

Note that when $n=2$, Proposition \ref{prop:lowestcell} implies that $\VC$ preserves perverses, because the lowest cell contains all elements except for the identity.

%================
\subsection{The underlying vector space}
%================

The diagrammatic Hecke category $\DC_{\aff}$ has a number of monoidal (two-sided cell) ideals. The maximal such ideal $\MaxIdeal$ contains the identity morphism of every indecomposable
object $B_x$ except the monoidal identity $\one = B_1$. Moreover, modulo $\MaxIdeal$, only the identity morphism of $\one$ survives (and its scalar multiples). Consequently, taking the
quotient by $\MaxIdeal$ gives a functor $\Ga \co \DC_{\aff} \to \Vect^\Z$, the category of graded vector spaces. Applying $\Ga$ to a complex in $\DC_{\aff}$ gives a complex of graded
vector spaces, and this operation preserves homotopy. Thus (taking homology) we obtain a functor from $\HC_{\aff}$ to the category of bigraded vector spaces $\Vect^{\Z\times \Z}$. When
applied to a perverse complex, the resulting bigraded vector space has grading concentrated on the (anti-)diagonal, and we identify such bigraded vector spaces with $\Vect^\Z$. We refer
to $\Ga(\PC) \in \Vect^\Z$ as the \emph{underlying graded vector space} of a perverse complex $\PC \in \HC_{\aff}$. Because $\MaxIdeal$ is a cell ideal, $\Ga$ is a monoidal functor.

We can apply this operation to $\DC_{\ext}$ as well. Now the maximal ideal $\MaxIdeal$ contains the identity morphism of every indecomposable object except for the powers of the
invertible object $\Om$. Now applying $\Ga$ to a complex yields a triply graded vector space, where the third grading records the power of $\Om$. This third grading is fairly trivial, as
$\DC_{\ext}$ already splits into blocks based on the power of $\Om$, so any indecomposable complex will be supported in a single power of $\Om$. We ignore this third grading, and continue
to identify the image of $\Ga$ with $\Vect^{\Z}$.

Any representation $U$ of $\gl_n$ has a weight space decomposition, i.e. it is graded by $\La_{\wt}$. One can collapse this $\La_{\wt}$ grading into a $\Z$-grading in such a way that
all the simple root vectors act to raise the weight by $2$. Let $\collapse(U)$ denote the $\Z$-graded vector space obtained by collapsing the grading; $\collapse$ is a faithful
(monoidal) functor from $\Rep \gl_n$ to $\Vect^\Z$. Note that, if $U$ is irreducible, then $\collapse(U)$ is supported in either even or odd degrees.

Goal \ref{goal2}(j), compatibility with the underlying vector space, is the statement that $\Ga \circ \GC \cong \collapse$ as monoidal functors $\Rep \gl_n \to \Vect^\Z$.

\begin{remark} Note that $\La_{\wt} / \La_{\rt} \cong \Z$, since we work with $\gl_n$. Thus one could keep track of an additional $\Z$-grading, measuring where each weight lives in
$\La_{\wt}/\La_{\rt}$, and this grading would be constant for all the weights in a given irreducible representation. This extra grading (which we ignore) matches the additional
$\Z$-grading coming from the power of $\Om$ (which we ignore). \end{remark}

In fact, Goal \ref{goal2}(j) is almost a consequence of the Wakimoto filtration.

\begin{lemma} Assume that Gaitsgory central sheaves have Wakimoto filtrations lifting the weight space filtration. Then $\Ga \circ \GC$ agrees with $\collapse$ on the level of objects. \end{lemma}
	
\begin{proof} For any Rouquier complex, $\Ga$ yields a one-dimensional vector space. It is trivial to verify that $\Ga(W(\la))$ agrees with the collapse of a one-dimensional space in
weight $\la$. Thus, after applying $\Ga$ to an irreducible Gaitsgory central complex, each of the Wakimotos in the filtration will contribute a one-dimensional vector space, and these
appear in either even or odd homological degrees, so the complex of vector spaces is already isomorphic to its homology. \end{proof}

For two functors which agree on objects, one can check that they are equal if they agree on generating morphisms. Thus the following lemma is obvious.

\begin{lemma} \label{lem:compatwithunderlyingvectorspace} Suppose that a functor $\GC \co \Webs^+_n \to \ZC(\HC_{\ext})$ is defined as in Goal \ref{goal1.5}. Suppose that Gaitsgory
central sheaves have Wakimoto filtrations lifting the weight space filtration. Suppose that $\Ga \circ \GC$ and $\collapse$ agree on the generating webs. Then $\Ga \circ \GC \cong
\collapse$ as functors. \end{lemma}

%================
\subsection{The great simplification} \label{subsec:greatsimple}
%================

Let us state some easy lemmas, which show the great power of Wakimoto filtrations.

\begin{lemma} Let $\MC$ be a perverse complex in $\HC_{\ext}$. Then $\Ga$ is faithful on the Hom space $\Hom(\one,\MC)$ and the Hom space $\Hom(\MC,\one)$. \end{lemma}

\begin{proof} The complex $\one$ is just the monoidal identity $\one = B_1$ of $\DC_{\ext}$ in homological degree zero. Meanwhile, by perversity, $\MC$ in homological degree zero is
some direct sum $\bigoplus B_w^{\oplus m_w}$ with no grading shifts. By \eqref{eq:gradedschur} there is only a (nonzero) map of degree zero from $B_1$ to $B_w$ if $w=1$, and this map is
just a scalar multiple of the identity map. Consequently, a chain map from $\one$ to $\MC$ is just a map $B_1 \to B_1^{\oplus m_1}$, and it is nonzero if and only if it is nonzero after
applying $\Ga$. Note that there are no potential homotopies between these complexes by perversity, so the space of chain maps agrees with $\Hom$ in $\HC$. The same argument works for
$\Hom(\MC,\one)$. \end{proof}

Note that we do not insist that $\Ga$ is faithful on the bigraded Hom space $\HOM(\one, \MC)$, or even on $\Hom(\one, \MC[k])$ for some $k \in \Z$.

\begin{lemma} Let $\MC$ be a complex in $\HC_{\ext}$ which satisfies $PP(\MC,W(\mu))$ for any Wakimoto complex $W(\mu)$, $\mu \in \La_{\wt}$. Then $\Ga$ is faithful on the Hom space $\Hom(W(\la),\MC)$ or the Hom space $\Hom(\MC,W(\la))$ for any $\la \in \La_{\wt}$. \end{lemma}

\begin{proof} The Wakimoto complex $W(\la)$ has a tensor-inverse $W(-\la)$. Thus there is an isomorphism 
\begin{equation} \Hom(W(\la),\MC) \cong \Hom(\one,\MC \ot W(-\la)) \end{equation}
obtained by tensoring with $W(-\la)$. Note that $\MC \ot W(-\la)$ is perverse by assumption. Now $\Ga$ acts faithfully on $\Hom(\one,\MC \ot W(-\la))$ by the previous lemma. Moreover, $\Ga$ commutes with the functor of tensoring with $W(-\la)$, because $\MaxIdeal$ is a tensor ideal. \end{proof}

\begin{lemma} Let $\MC$ be perverse, and let $A \to \MC' \to B \to A[1]$ be a distinguished triangle in $\HC_{\ext}$, with $A$ and $B$ perverse. If $\Ga$ is faithful on $\Hom(\MC,A)$ and $\Hom(\MC,B)$ then it is faithful on $\Hom(\MC,\MC')$. \end{lemma}
	
\begin{proof} By perversity, $\Hom(\MC,B[-1])=0$ and $\Hom(\Ga \MC, \Ga B[-1])=0$. Now comparison of the long exact sequence for $\Hom(\MC,-)$ and that for $\Hom(\Ga \MC, \Ga -)$ using the 5-lemma gives the result. \end{proof}

\begin{lemma} \label{lem:greatsimple} Let $\MC, \MC' \in \Wak_0$. Then $\Ga$ is faithful on $\Hom(\MC, \MC')$. \end{lemma}

\begin{proof} This follows by combining the two previous lemmas. \end{proof}

Lemma \ref{lem:greatsimple} can be used to simplify the construction of $\GC$.

\begin{prop} \label{prop:greatsimple} Suppose that one has achieved goals (a) and (b) and (i) and (j) of Goals \ref{goal1.5} and \ref{goal2}. Namely, suppose that one constructs
complexes $\VC_{\fund_i}$ with Wakimoto filtrations lifting their weight space filtrations, and morphisms in $\ZC(\HC_{\ext})$ associated to the generating webs which act correctly on
the underlying vector space. Then the web relations hold, so that $\GC$ is a well-defined functor. Moreover, $\GC$ is faithful. In other words, (c) and part of (f) from Goal
\ref{goal1.5} follow. \end{prop}

\begin{proof} If $U$ and $U'$ are any objects in $\Webs^+_n$ then we set $\GC(U)=\UC$ and $\GC(U')=\UC'$ to be the corresponding tensor products of the $\VC_{\fund_i}$. By Lemma
\ref{lem:wak0closed}, $\UC$ and $\UC'$ live in $\Wak_0$, so by Lemma \ref{lem:greatsimple}, $\Ga$ is faithful on $\Hom(\UC,\UC')$. We set $\GC$ of a generating web to be the
corresponding morphism in $\HC_{\ext}$. Then the web relations will hold between these morphisms, as they hold by assumption after applying $\Ga$. In particular, $\GC$ is a well-defined
functor. The composition $\Ga \circ \GC$ is isomorphic to $\collapse$ by Lemma \ref{lem:compatwithunderlyingvectorspace}, and $\collapse$ is faithful, so $\GC$ is also faithful.
\end{proof}

\begin{remark} We also checked the $\gl_n$ web relations by hand when $n=2, 3$. This takes a great deal work, but it is tractable and explicit. For reasons of sanity we did not
reproduce those computations here, since Proposition \ref{prop:greatsimple} will make them unneccesary. \end{remark}

Now let us discuss what can be done to prove that $\GC$ is full to degree zero. Suppose one has a morphism $\GC(U) \to \GC(U')$ which acts on the underlying vector spaces as a
$\gl_n$-intertwiner $U \to U'$. Then some web has the same image, and the faithfulness of $\Ga$ implies that the morphism is the image of the web. However, in theory it is possible that
there exists a morphism $\GC(U) \to \GC(U')$ which does not act on the underlying vector space as a $\gl_n$ intertwiner. We know that all morphisms must commute with the monodromy, but we
currently have no argument that all morphisms should commute with the raising and lowering operators of $\gl_n$.

Note that there are two separate questions: whether $\GC$ is full as a functor to $\ZC(\HC_{\ext})$, or whether it is full as a functor to $\HC_{\ext}$. The morphism spaces in
$\ZC(\HC_{\ext})$ are strictly smaller than those in $\HC_{\ext}$, as morphisms must commute with the central structure. We are not aware of any conjectures as to whether the forgetful
functor $\ZC(\HC_{\ext}) \to \HC_{\ext}$ is full (in degree zero) on Gaitsgory central complexes or not. Note that the forgetful functor is certainly not faithful on higher degree
endomorphisms, as multiplication by a non-invariant polynomial is not a central endomorphism.

Given that $\GC$ is a functor, the clasp idempotents in $\Webs^+_n$ are sent to idempotent endomorphisms in $\HC_{\ext}$, which should pick out the indecomposable Gaitsgory central
complexes as direct summands. Assuming that $\Wak_0$ is closed under direct summands, these summands should have Wakimoto filtrations compatible with their weight space filtrations.
Suppose that one had an independent description of the complex $\VC_{\lambda} = \GC(V_{\lambda})$ for each irreducible representation $V_{\lambda}$ of $\gl_n$. It may be possible to
prove directly that the degree zero endomorphisms of $\VC_{\lambda}$ are just scalar multiples of the identity map. If so, it would follow (by Schur's lemma) that $\GC$ is full.

The form of $\VC_{\lambda}$ is quite constrained, because it has a Wakimoto filtration and is invariant under the automorphism $\tau$; these two properties determine $\VC$ for the
standard representation, see Remark \ref{rmk:determiningVC}. For the standard representation we could prove that the endomorphism ring of $\VC$ in degree zero is spanned by the identity
map, because $\VC$ is multiplicity-free perverse. Now $\VC_{\lambda}$ is not multiplicity-free perverse in general, but a combinatorial analysis of $\VC_{\lambda}$ may still suffice to
prove that its endomorphism ring in degree zero is spanned by the identity map. This seems to me to be the most promising method to prove the fullness of $\GC$ directly. My gut feeling
is that the additional multiplicities in $\VC_{\lambda}$ will lead to a larger endomorphism ring in $\HC_{\ext}$, but that the additional constraints of $\ZC(\HC_{\ext})$ will force a
one-dimensional endomorphism ring.

%================
\subsection{Symmetries} \label{subsec-symmetries}
%================

The category $\Rep \gl_n$ has a number of symmetries, some which are not obvious until it is given its diagrammatic interpretation in terms of webs. The fact that $\Rep \gl_n$ is a symmetric monoidal category helps to obfuscate some of these symmetries, which become more distinct in the quantum deformation.

We say that a functor $F$ is \emph{monoidally covariant} if $F(M \ot N) \cong F(M) \ot F(N)$, and \emph{monoidally contravariant} if $F(M \ot N) \cong F(N) \ot F(M)$. Diagrammatic
categories (with cyclic adjunction) come naturally equipped with an action $\Z/2 \times \Z/2$, generated by the functor which flips a diagram vertically (contravariant but monoidally
covariant) and the functor which flips a diagram horizontally (covariant but monoidally contravariant). The composition of these, which rotates a diagram by 180 degrees, is the functor
sending each object to its biadjoint. However, $\Webs$ has an extra $\Z/2$ worth of symmetries, because one can also invert the orientation on a web (covariant and monoidally covariant).

There is a contravariant, monoidally covariant functor on $\Webs^+_n$ which fixes all objects and flips each diagram upside-down, and reverses the orientation so the diagram remains
positively-oriented. This functor is intertwined by the standard duality functor $\DM$ on $\Diag_{\ext}$, which also fixes all objects and flips each diagram upside-down. Note that the
extension of $\DM$ to $\HC_{\ext}$ does not fix all complexes, as it reverses differentials and homological degrees. Thus $\GC$ will not intertwine these dualities unless each
$\VC_{\fund_i}$ is self-dual under $\DM$. This is the case for the complex $\VC$ we define here.

The (contravariant, monoidally contravariant) biadjunction functor on $\Webs_n$ is the one which corresponds to the usual duality/biadjunction functor on $\Rep \gl_n$. It sends each
object to its dual, where $(U \ot U')^* \cong (U')^* \ot U^*$. It rotates all diagrams by 180 degrees. This is intertwined by the biadjunction functor on $\Diag_{\ext}$, which fixes each
$B_s$, sends $\Om$ to $\Om\inv$, and rotates all diagrams by 180 degrees.

There is another more subtle functor on $\Webs_n$, which is contravariant but monoidally covariant. It sends each fundamental representation to its dual, and reverses the orientation on
each diagram. This functor is intertwined by $\si$, the Dynkin diagram automorphism of $\Diag_{\ext}$ coming from the reflection of $\tilde{A}_{n-1}$ across the affine vertex.

Note that both $\si$ and biadjunction send $V$ to $V^*$ for the standard representation $V$. Thus it must be case that $\VC^* := \si(\VC)$ is also the biadjoint of $\VC$. By construction,
$\VC$ is built from Bott-Samelsons for collections of distinct simple reflections placed in anticlockwise order. Applying $\si$ will swap anticlockwise order with clockwise order.
Meanwhile, applying adjunction will reverse the order on a Bott-Samelson, thus also yielding clockwise order. This is why the two functors agree on $\VC$.

Assuming that all these functors intertwine appropriately on the generating objects and generating morphisms of $\Webs_n$, it will follow that they behave appropriately on the entire
category. Thus future constructions of $\VC_{\fund_i}$ are constrained to ensure that $\VC_{\fund_i}$ is self-dual under $\DM$, and that $\si(\VC_{\fund_i})$ is the adjoint of
$\VC_{\fund_i}$.

%================
\subsection{Connections to geometric Satake} \label{subsec:intro_satake}
%================

There is a pushforward map from the affine flag variety $\Fl$ to the affine Grassmanian $\Gr$, and it is clear what the analog of this functor is in the Soergel context. There is a
generalization of Soergel bimodules known as \emph{singular Soergel bimodules}, developed in \cite{WillSingular}. Instead of considering $R$-bimodules, singular Soergel bimodules work
with $(R^I,R^J)$-bimodules, where $I$ and $J$ are proper subsets of $S$, and $R^{I} = R^{W_I}$ denotes the invariants under the corresponding parabolic subgroup. The analog of the
pushforward from $\Fl$ to $\Gr$ is just restriction (on the right) from $(R,R)$-bimodules to $(R,R^{W_{\fin}})$-bimodules. We write this functor as $(-)_{R^{W_{\fin}}}$. We write the
induction functor from an $(R^I,R^J)$-bimodule to an $(R,R^J)$-bimodule as $R \ot_{R^I} (-)$. There is also a diagrammatic version of singular Soergel bimodules, and we use the same notation for the analogs of restriction and induction in this context.

The geometric Satake equivalence can also be reformulated explicitly and naively in the context of singular Soergel bimodules or their diagrammatics. This was proven\footnote{The paper
\cite{EQuantumI} only contains the proofs for $n=2, 3$. The same framework does work for general $n$, but the proofs are not yet publically available.} in type $A$ in the author's
previous work \cite{EQuantumI}, using a similar web-based approach to the one used in this paper. In particular, one constructs a 2-functor $\SC$ from a ``colored'' version of $\Webs^+_n$ to singular Soergel bimodules for the affine Weyl group. The singular Soergel bimodules in the image are $(R^I,R^J)$-bimodules, where $W_I$ and $W_J$ are isomorphic to $W_{\fin}$
(but are not necesssarily the same subset of the affine Dynkin diagram). It is straightforward to produce an analogous monoidal functor from $\Webs^+_n$ to singular Soergel bimodules
for $W_{\ext}$; because of the extra twisting object $\Om$ available in $\DC_{\ext}$ it is possible to arrange that this functor lands purely inside
$(R^{W_{\fin}},R^{W_{\fin}})$-bimodules for a fixed copy of $W_{\fin}$.

Property (l) of Goal \ref{goal2} states that, for any $U \in \Rep \sl_n$, one has \begin{equation} \label{eq:GaitsMatchesSatake} \GC(U)_{R^{W_{\fin}}} \cong R \ot_{R^{W_{\fin}}} \SC(U).\end{equation} Gaitsgory proves
this in the geometric setting. We can check this by hand for $n=2, 3$. Note that the image of $\SC$ consists of objects, not complexes of objects! Thus, this proposition implies that the
(rather complicated) complex $\GC(U)$ will, after restriction to $R^{W_{\fin}}$, contract to a complex supported entirely in homological degree zero.

The functor of restriction can be difficult to understand, but induction is relatively harmless. The induction functor $(-) \ot_{R^{W_J}} R$ from singular Soergel
$(R^I,R^J)$-bimodules to singular Soergel $(R^I,R)$-bimodules is fully faithful after base change, and preserves indecomposable objects, see \cite[Proposition 7.4.3(1)]{WillSingular}. In
particular, to compute the restriction $\GC(U)_{R^J}$, it is an equivalent problem to compute $\GC(U) \ot_{R^J} R$, which is the composition of restriction with induction. Since
$B_J \cong R \ot_{R^J} R$ up to grading shift, this latter problem is just the computation of $\GC(U) \ot B_J$ purely within the category of ordinary Soergel bimodules.

In Theorem \ref{thm:tensorBI} we compute the minimal complex of $\VC \ot B_J$ for all $J \subsetneq S$. In particular, when $J$ is a copy of the finite Weyl group, the resulting complex
is supported in a single homological degree, and matches with the expectation from geometric Satake. In fact, our proof of Theorem \ref{thm:tensorBI} can be adapted without any
significant work to the setting of singular Soergel bimodules, where it gives an explicit computation of the minimal complex of the restriction $\VC_{R^J}$. This proves
\eqref{eq:GaitsMatchesSatake} more directly for the standard representation. We do not present the singular Soergel bimodule version of Theorem \ref{thm:tensorBI} because of the
additional background required.

Of course, if one can prove \eqref{eq:GaitsMatchesSatake} for all fundamental representations and all generating webs, one proves this isomorphism of functors in general.

%================
\subsection{Centrality}\label{ssec:intro_centrality}
%================

Given a complex $\VC$ which one expects is in the Drinfeld center of a homotopy category, how does one actually prove it? How do we achieve Goal \ref{goal1.5}(d)?

Recall that $\HC_{\ext}$ is the homotopy category of an explicitly presented additive category $\DC_{\ext}$. Let us first ask: does $\VC$ commute functorially with the objects of $\DC_{\ext}$, i.e. with complexes supported in a single homological degree? Proving this is a straightforward if computationally daunting task. One must first provide a chain map $\phi_B \co B \ot \VC \to \VC \ot B$ for each generating 1-morphism $B$ of $\DC_{\ext}$. Taking tensor products of these maps will produce a map $\phi_X \co X \ot VC \to \VC \ot X$ for every 1-morphism of $\DC_{\ext}$. Then, for each of the generating 2-morphisms $f \co X \to Y$, one must show that the square
\begin{equation} \label{centralsquare}
\begin{diagram}
X \VC & \rTo^{f \ot \VC} & Y \VC \\
\dTo^{\phi_X} & & \dTo_{\phi_Y} \\
\VC X & \rTo^{\VC \ot f} & \VC Y \end{diagram}
\end{equation}
commutes up to homotopy. If we can do this, we say that \emph{$\VC$ commutes with $\DC_{\ext}$}, or that $\VC$ is in the \emph{restricted Drinfeld center} of $\HC_W$. In other words, $\phi$ gives a natural monoidal isomorphism of functors $(-) \ot \VC \to \VC \ot (-)$, when viewed as functors $\DC_{\ext} \to \HC_{\ext}$.

One of the major results of this paper is the proof that the Gaitsgory complex $\VC$ of the standard representation is in the restricted Drinfeld center, see Theorem \ref{thm:rotationfunctorial}. We reduce this to a handful of computations in Lemma \ref{lem:reducingrotation}, and then perform these computations in \S\ref{sec-tensorBIcomps}.

Once one proves that $\VC$ commutes with $\DC_{\ext}$, one still needs to prove that it commutes with any complex of objects from $\DC_{\ext}$! This is far from being an immediate
consequence. Let us give a brief summary of the problem.

Any bounded complex $\MC$ has its \emph{homological filtration}, which describes it as a convolution with subquotients $\MC^k[-k]$ supported in a single homological degree. Then $\MC \ot
\VC$ (resp. $\VC \ot \MC$) is a convolution with subquotients $\MC^k[-k] \ot \VC$ (resp. $\VC \ot \MC^k[-k]$), but these subquotients are themselves complexes, not usually supported in a
single degree. We already know how an isomorphism $\MC \ot \VC \to \VC \ot \MC$ should act on the associated graded of this filtration, but this only gives the diagonal entries of the
upper-triangular matrix which should describe a chain map on the convolution $\MC = \oplus \MC^k[-k]$. The diagonal entries will not form a chain map, since we have modified the
differential by upper-triangular terms. One must instead attempt to find a chain map $\phi_{\MC} \co \MC \ot \VC \to \VC \ot \MC$ which preserves the filtration, and acts in the
associated graded as $\oplus \phi_{\MC^k[-k]}$. The question is whether this chain map $\phi_{\MC}$ can be constructed at all, or whether it is unique.

The problem is visible already for a complex with only two nonzero degrees. Suppose one has a map $f \co X \ot Y$ (for complexes $X, Y$ both supported in degree zero) for which
\eqref{centralsquare} commutes up to homotopy. When giving an isomorphism $\Cone(f) \ot \VC \to \VC \ot \Cone(f)$ one must choose the homotopy $h$ for \eqref{centralsquare}. If two
possible homotopies $h_1$ and $h_2$ are not themselves homotopic, then their associated isomorphisms may differ. The vanishing of $\Hom(X \VC, \VC Y[-1])$ will force the homotopy $h$ to
be unique up to homotopy, solving this first potential problem. This is not the only way the problem might be resolved, but it is certainly a sufficient condition. The reader who has
never considered such issues will be well-served by carefully working through this case of a cone. When examining the situation for iterated cones, other negative exts $\Hom(FM, GN[c])$
for $c<0$ appear. This gives the gist of the problems involved in extending the natural transformation $\phi$ from objects to complexes.

In the monograph \cite{EHDrinfeld}, joint with Hogancamp, we prove a general technical lemma which allows one to upgrade objects in the restricted Drinfeld center to the true Drinfeld
center, when certain obstructions vanish. Namely, if for any objects $M, N \in \DC_{\ext}$ we have \begin{equation} \label{eq:neededforEHdrinfeld} \Hom(M \FC, \FC N[c]) = 0 \quad
\text{for all } c < 0\end{equation} then one can upgrade $\FC$.

In all the calculations we have done in small cases, \eqref{eq:neededforEHdrinfeld} holds. In particular, when $M = \one$ and $N = \one(k)$ for some $k \in \Z$, this implies property
(m). We conjecture that \eqref{eq:neededforEHdrinfeld} holds for all Gaitsgory central complexes as well. Given Theorem \ref{thm:tensorBI} it is an exercise to prove this directly for
$n=2$, and we expect that $n=3$ also has a direct proof.

% \begin{thm} When $n=2$, Conjecture \ref{conj4} holds. \end{thm}
%
% \begin{proof} (Sketch) The complex $\FC$ is supported in homological degrees between $-1$ and $+1$. It is an exercise to verify that $\Hom(\FC, \FC[c](k))=0$ for all $c < 0$ and $k \in
% \Z$. Meanwhile, for $1 \ne w \in W$, $B_w(k) \ot \FC$ is concentrated in homological degree zero. Thus \eqref{eq:conj4} holds automatically when $M = B_w$ and $N = B_{w'}(k)$ for all
% cases except when either $w$ or $w'$ is the identity element, and $c=-1$. A chain map from the complex $(0 \to C \to 0)$ supported in degree zero to $\FC[-1](k)$ would be a map $C \to R$ of some degree, such that the compositions $C \to R \to B_0$ and $C \to R \to B_1$ are both zero. But the startdot is monic, so no such chain map exists. We leave the remainder and the details to the reader. \end{proof}

Lacking a proof of \eqref{eq:neededforEHdrinfeld} in general, we are not able to prove that $\VC$ is in the true Drinfeld center. Even were we able to prove this, it would not change
the method by which one would have to compute $\phi_{\MC}$ for any given complex! One would still have to do the work of lifting the associated graded chain maps to an overall chain
map. If one can show that this lift is unique (up to homotopy) then one has effectively solved the extension problem for $\MC$. Thus in specific examples we can construct the
commutation isomorphism $\phi$, even if we currently lack the general proof of its existences. For example, we compute the braiding map $\phi_{\VC}$ when $n=2$ in
\S\ref{subsec:braidingn2}.

%================
\subsection{Braiding}\label{ssec:intro_braiding}
%================

Let us discuss what can be said about Goal \ref{goal1.5}(e). Since $\VC \ot \VC$ is Wakimoto-filtered, Lemma \ref{lem:greatsimple} will imply that morphisms can be detected on the
underlying vector space, after application of $\Ga$. Thus, to confirm that the braiding map $\phi_{\VC} \co \VC \to \VC$ agrees with the braiding on $\Rep \gl_n$, we need
only confirm that they behave the same way on the copies of $\one$ inside $\VC$.

We need to make clear some subtleties in homological algebra. In the previous section we discussed how one might try to compute how $\VC$ commutes past a complex (e.g. the morphism $\phi_{\MC}$) if one knows how $\VC$ commutes past objects (e.g. if
one knows $\phi_{\MC^k}$ for all $k$) using the homological filtration on the complex. We swept one key issue under the rug, which is that we also need to know how $\VC$ commutes past objects with a homological shift, i.e. complexes supported in a single nonzero degree (e.g. $\phi_{\MC^k[-k]}$). As always in homological algebra, shifts lead to signs.

The object $\one$ is in the Drinfeld center, and its commutation maps are always identity maps. The object $\one[1]$ is also in the Drinfeld center, but its commutation maps are not
identity maps, nor is the identity map even a chain map! That is, $\one[1] \ot \MC$ and $\MC \ot \one[1]$ are not the same complex, as the differentials differ by signs. The natural isomorphism between these complexes is the map
\begin{equation} \epsilon_{\MC} \colon \one[1] \ot \MC \to \MC \ot \one[1], \end{equation}
where $\epsilon_{\MC}$ multiples the summand $\MC^k$ in homological degree $k$ by a sign $(-1)^k$. In particular, there is a natural braiding between shifts of the monoidal identity, and it is not the identity map, instead one has the map
\begin{equation} \label{eq:braidingonone}(-1)^{k \ell} \co \one[k] \one[\ell] \to \one[\ell] \one[k]. \end{equation}
	
\begin{remark} Note that what we typically describe as $\MC[1]$ is actually $\one[1] \ot \MC$ and not $\MC \ot \one[1]$, if one uses the standard conventions for signs in tensor
products. \end{remark}

Similar constraints affect the braiding map of any object in the Drinfeld center. For example, the map $\phi_{\one}$ says how $\VC$ commutes past the monoidal identity, and it must equal
$\id_{\VC}$ or else $\phi$ will not be a monoidal natural transformation. However, $\phi_{\one[1]}$ can not equal $\id_{\VC}$, which again is not a chain map, but is equal to
$\epsilon_{\VC}$. Similarly, for any complex $\MC$, $\phi_{\MC[1]}$ and $\phi_{\MC}$ will differ by some signs.

After passage to the underlying vector space, the category is now generated by copies of the monoidal identity, which has standard braiding as given by \eqref{eq:braidingonone}. We can
ask whether $\phi_{\VC} \co \VC \to \VC$ agrees with the standard braiding on the underlying vector space. The standard braiding agrees exactly with the diagonal terms in the homological
filtration, essentially by tautology. In theory, the braiding $\VC \to \VC$ might also have upper-triangular summands \begin{equation} \one[k] \ot \one[\ell] \to \one_[k+j] \ot
\one_[\ell-j], \quad \text{for } j > 0. \end{equation}

When $n=2$ we computed directly that the braiding $\phi_{\VC}$ matches the standard braiding, and there are no upper-triangular terms. Usually the case $n=2$ is large enough to witness
any weird, unusual braidings, so this confirmation is reassuring.

\begin{remark} We believe that there should be an abstract argument why the braiding map should be diagonal after applying $\Ga$, but we have not found it. \end{remark}

Finally, we note that the standard braiding on the underlying vector space does not agree with the most naive braiding on $\Rep \gl_n$, but it does agree with what we call the
\emph{Satake braiding}. On the standard representation $V$ with basis $\{e_1, \ldots, e_n\}$, this braiding sends $e_i \ot e_j \to (-1)^{n-1} e_j \ot e_i$, because the value of $n$
determines whether these copies of the monoidal identity appear in even or odd homological degree. This is also the braiding on $\Rep \gl_n$ which arises from the geometric Satake
equivalence.

\begin{remark} To reiterate, the failure of the braiding to agree with the naive symmetric braiding $e_i \ot e_j \mapsto e_j \ot e_i$ has nothing to do with the category $\DC_{\ext}$ or
the Gaitsgory central complexes, but is entirely a feature of homological algebra itself, and its standard braiding \eqref{eq:braidingonone}. \end{remark}

%================
\subsection{Parity vanishing}\label{ssec:intro_homs}
%================

Let us briefly discuss the conjectural homological parity vanishing of Hom spaces, property (n). While the degree zero morphisms between Gaitsgory central complexes are governed by the
semisimple category $\Rep \gl_n$, the higher degree morphisms are quite interesting and mysterious. For sake of posterity, we record the higher endomorphism rings of $\VC$ when $n = 2$ in \S\ref{subsec:endosn2}, though we do not reproduce the entire computation. We have also computed the endomorphism ring for $n=3$, and is parity.

Homological parity vanishing has occurred recently in another context: the computation of $\Hom(\one, \FT_n)$ achieved in \cite{ElHog16a}. This result was extended by Mellit
\cite{Mellit}, to prove parity vanishing for $\Hom(\one, \FC)$ where $\FC$ is the Rouquier complex of a wide variety of braids. If the ideas of \S\ref{ssec_intro:GNR} are accurate,
then there should be a relationship between these morphism spaces and the ones in this paper, although we are not ready to make this connection precise. We would not have been so bold as
to have conjectured property (n) otherwise.

\begin{remark} Homological parity vanishing is a stronger and stranger condition than ordinary parity vanishing (in the sense of a \emph{parity sheaf} from \cite{JMW14}). In the category $\DC_{\ext}$, we have $\Hom(B_w,B_x(k)) = 0$ unless $k$ has the same parity as $\ell(w) - \ell(x)$. Consequently, any perverse complex splits into its
``even part,'' built from $B_w(k)$ in homological degree $k$ where $\ell(w)$ and $k$ have the same parity, and its ``odd part,'' where $\ell(w)$ and $k$ have different parity. There are
no maps from an even perverse complex to an odd perverse complex or vice versa. If $\UC$ and $\UC'$ are perverse complexes, and each is either even or odd, then $\Hom(\UC,\UC'(k)[c]) =
0$ when $k + c$ is either even or odd. However, homological parity vanishing also states that the Hom space vanishes whenever $c$ is even or odd. Instead of vanishing like a checkerboard, the Hom space vanishes like a sparse checkerboard. \end{remark}

\section{The standard Gaitsgory complex}
\label{sec-construction}
%%%%%%%%%%%%%%%%%%%%%%%%%

Having given several examples, our goal is to construct complexes $\FC$ and $\VC$ for all $n \ge 2$. This construction will take place in \S\ref{subsec-stdrep}. The construction of the
``dual" complexes $\FC^*$ and $\VC^*$ appears in \S\ref{subsec-dualrep}. First we set up notation and perform some key computations. In \S\ref{subsec-signs} we describe the objects which
appear in the complexes $\FC$, and in \S\ref{subsec-dotmaps} we describe the summands of the differential. In this chapter we also discuss the monodromy map $\mu$, the commutation
isomorphism $\Om \FC \cong \FC \Om$, and the Wakimoto filtration on $\VC$.

For this chapter fix $n \ge 2$, and let $S = S_{\aff} = \Z/n\Z$. Recall that $\tau(s_i) = s_{i+1}$ and $\si(s_i) = s_{-i}$. We work mostly in $\Diag = \Diag_{\aff}$, for some realization
of the affine Weyl group with automorphism $\tau$, though we work with $\Diag_{\ext}$ when we discuss $\VC$.

%===========
\subsection{Keeping track of signs: cyclically ordered bimodules}
\label{subsec-signs}
%===========

Our complex $\FC$ will be built from indecomposable Soergel bimodules $B_X$ attached to proper subsets $X \subset S$.

\begin{defn} \label{defn:hX} Let $X \subsetneq S$. Pick some $\notme \in S \setminus X$. Then we can place a total order on $S$, where $\notme < \notme+1 < \notme+2 < \ldots < \notme-1$. Restricting this order to the subset $X$, we write $X = \{x_1 < x_2 < \ldots < x_d\}$. Any such order on $X$ will be called a \emph{cyclic order}. Then let $h_X \in W_{\aff}$ denote the element $s_{x_d} \cdots s_{x_2} s_{x_1}$. This is a particular Coxeter element for the parabolic subgroup $W_X$, the one appearing in anticyclic order.
\end{defn}

\begin{ex} When $n=9$ and $X = \{0,1,3,4,6,8\}$, then choosing $\notme = 7$ one has $h_X = s_6 s_4 s_3 s_1 s_0 s_8$. Choosing $\notme = 2$ one has $h_X = s_1 s_0 s_8 s_6 s_4 s_3$. Of course, these two expressions for $h_X$ are equal in $W_{\aff}$, because $s_k$ commutes with $s_j$ when $j \ne k \pm 1$. \end{ex}

\begin{lemma} The element $h_X$ does not depend on the choice of cyclic order, i.e. on the choice of $\notme \in S \setminus X$.
%Similarly for $w_X^*$.
\end{lemma}

\begin{proof} One should think of $X = \coprod_a X_a$ as a disjoint union of connected components in the Dynkin diagram $\tilde{A}_{n-1}$. When $X$ is a single connected component, the
cyclic ordering on $X$ is clearly independent of the choice of $\notme$. As a consequence, in the general case, each $h_{X_a}$ is independent of the choice of $\notme$. Moreover, $h_{X_a}$
commutes with $h_{X_b}$ for all $a \ne b$, since the simple reflections in $X_a$ commute with those in $X_b$. Now $h_X$ is a product of the elements $h_{X_a}$, and while the order
in which this product appears does depend on $\notme$, the overall product $h_X$ does not. \end{proof}

\begin{lemma} Every reduced expression for $h_X$ can be described as $s_{x_d} \cdots s_{x_2} s_{x_1}$ for some order $X = \{x_1 < \ldots < x_d\}$ on $X$ which restricts to
each component $X_a$ to be the standard order. Any two reduced expressions are related by the commuting braid relations $s_k s_j = s_j s_k$ for $j \ne k \pm 1$, which reorder the
indices in $X$ without reordering the indices in any connected component $X_a$. \end{lemma}

\begin{proof} Any reduced expression for $h_X$ will have the simple reflections $s_x$, $x \in X$, appearing each once in some order. That this order restricts to each component $X_a$ to
be the standard order is straightforward. Because no simple reflection occurs twice, there can not be any braid relations of the form $s_i s_{i + 1} s_i = s_{i + 1} s_i s_{i + 1}$,
so there can only be commuting relations. \end{proof}

Recall that, for any two reduced expressions for an element $w \in W$, and for any path in the graph of reduced expressions that applies one braid relation at a time, there is a
corresponding morphism in $\Diag$. This is known as a \emph{rex move}, see \cite[\S 4.2]{EWGr4sb}. When the only braid relations involved are commutations $s_k s_j = s_j s_k$ for
$j \ne k \pm 1$, the rex move is just a product of 4-valent vertices like $\Xbg$. There are many different rex moves between two given reduced expressions, but we have the following
result (see \cite[equations (5.9) and the $m=2$ case of (5.16)]{EWGr4sb}).

\begin{prop} \label{prop:rexmovesequal} Any two rex moves built entirely from 4-valent vertices, with the same source and target, are equal. \end{prop}

Hence, if two reduced expressions are only related by commutations, then there is a canonical rex move between the corresponding Bott-Samelson objects, which is an isomorphism. One
can use this to canonically identify the Bott-Samelson objects for this reduced expression, using rex moves. Instead, as motivated by \S\ref{subsec-stdrep4}, we will canonically identify them using signed rex moves.

\begin{defn} Given a rex move $\psi$ built entirely from 4-valent vertices, the corresponding \emph{signed rex move} is just $(-1)^\ell \psi$, where $\ell$ is the number of 4-valent vertices in $\psi$. \end{defn} 

\begin{lemma} \label{lem:signedrexequal} Any two rex moves built entirely from 4-valent vertices, with the same source and target, have the same number of 4-valent vertices modulo 2. Hence, the corresponding signed rex moves are equal. \end{lemma}
	
\begin{proof} A rex move built from 4-valent vertices, between reduced expressions of length $d$, can be viewed as an expression in the Coxeter group $S_d$ which permutes the indices.
Any two expressions for the same permutation in $S_d$ have the same number of crossings modulo 2. \end{proof}

Now, for a proper subset $X \subset S$ we wish to give meaning to the symbol $B_X$, which will be isomorphic to the indecomposable Soergel bimodule $B_{h_X}$. Our goal is to choose an
isomorphism $B_X \to \BS(\un{h}_X)$ for every possible reduced expression $\un{h}_X$ of $h_X$, which will be compatible with the (signed) isomorphisms between these different
reduced expressions.

\begin{defn} \label{defn:BXto} Choose once and for all a preferred reduced expression $\un{h}^{\$}_X$ for $w_X$. We identify the symbol $B_X$ with $\BS(\un{h}^{\$}_X)$, and write $\BXto \co B_X \to \BS(\un{h}^{\$}_X)$ for the identity map. If $\un{h}_X$ is some other reduced expression for $h_X$, then we write $\BXto_X \co B_X \to \BS(\un{h}_X)$ for the composition
\begin{equation} \BXto_X = \left( B_X \to \BS(\un{h}^{\$}_X) \to \BS(\un{h}_X) \right), \end{equation}
where the first map is the identity map, and the second map is a signed rex move. We write $\BXto$ instead of $\BXto_X$ when $X$ is understood. \end{defn}

By Lemma \ref{lem:signedrexequal}, this definition of $\BXto$ does not depend on the choice of signed rex move, so $\BXto$ is well-defined.

\begin{remark}\label{rmk:changingsignglobally} If one chooses a different preferred reduced expression $\un{h}_X^{\$}$, the only effect will be that the maps $\BXto$ will be multiplied by a global sign. Ultimately this
will not affect any of the complexes we construct using $B_X$, up to isomorphism of complexes, see Remark \ref{rmk:changingsignnotchangingcomplex}. \end{remark}

\begin{remark} For a proper subset $I \subset S$ we also use the notation $B_I$ to denote the indecomposable Soergel bimodule $B_{w_I}$ associated to the longest element of the parabolic, rather than $B_{h_I}$ associated to a certain Coxeter element of thte parabolic. These notations are admittedly not consistent, but both $B_I = B_{w_I}$ and $B_X = B_{h_X}$ are used so often in this paper that it is worth having shorthand. When we use the letters $I$ and $J$ we will always be interested in $w_I$ and $w_J$, while when we use $X$, $Y$, and $Z$ we will always be interested in $h_X$, $h_Y$, and $h_Z$. If we ever need to use the element $h_I$ we will specifically write $B_{h_I}$. \end{remark}

We refer to Soergel bimodules of the form $B_X$ as \emph{cyclically-ordered Soergel bimodules}, or just \emph{cyclical bimodules} for short. We now introduce diagrammatics for dealing with cyclical bimodules. We let
\[{
\labellist
\tiny\hair 2pt
 \pinlabel {$X$} [ ] at 18 6
\endlabellist
\centering
\ig{1}{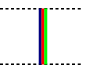}
}\] denote the identity map of $B_X$ (this picture is meant to evoke a rainbow, for which my LaTeX skills were insufficient). We let \begin{equation} \BXto = \ig{1}{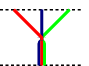} \end{equation} denote the map $\BXto$ for any given reduced expression for $\un{h}_X$ (this picture is meant to evoke a prism splitting a light beam). We let \begin{equation} \BXto\inv = \igv{1}{BXto} \end{equation} denote the inverse of $\BXto$. Then we have the following relations, where the last follows from Lemma \ref{lem:signedrexequal}.

\begin{equation} \ig{1}{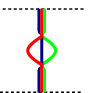} \quad = \quad \ig{1}{idBX}. \end{equation}
\begin{equation} {
\labellist
\small\hair 2pt
 \pinlabel {$=$} [ ] at 32 16
\endlabellist
\centering
\ig{1}{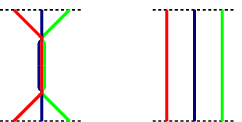}
}.\end{equation}
\begin{equation} \label{eq:BXtocross}{
\labellist
\small\hair 2pt
 \pinlabel {$=$} [ ] at 26 11
 \pinlabel {$-$} [ ] at -5 11
\endlabellist
\centering
\ig{1}{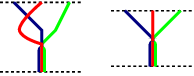}
}. \end{equation}

%===========
\subsection{Dot maps between cyclical bimodules}
\label{subsec-dotmaps}
%===========

Now we define the maps which will appear as the summands of the differential in the twisted standard complex $\FC$. The reader should remind themselves of the sign conventions for (even or odd) BS complexes from \S\ref{subsec-signconvention}.

\begin{defn} \label{defn:dotsign} Let $X$ and $Y$ be two proper subsets of $S$, with $Y = X \cup i$ for $i \notin X$. We will be defining morphisms $\startdotsign_X^Y \co B_X \to B_Y(1)$ and $\finaldotsign_Y^X \co B_Y \to B_X(1)$, which we call \emph{$XY$-dots} or simply \emph{signed dots}.
	
Choose a reduced expression $\un{h}_Y$ for $h_Y$, and let $\un{h}_X$ be the reduced expression for $h_X$ obtained by removing $s_i$. We will identify $B_X$ with $\BS(\un{h}_X)$ using
$\BXto_X$, and will identify $\BS(\un{h}_Y)$ with $B_Y$ using $\BXto_Y\inv$. To define the maps $\startdotsign_X^Y$ and $\finaldotsign_Y^X$, we need only define the corresponding maps
$\BS(\un{h}_X) \to \BS(\un{h}_Y)(1)$ and vice versa. By abuse of notation, we may also refer to these maps between Bott-Samelson objects as \emph{signed dots}.

Suppose that $s_i$ is the $k$-th simple reflection in $\un{h_Y}$. Let $\ep_{X,Y} = \ep_{Y,X} = (-1)^{k-1}$. Then the signed dot $\BS(\un{h}_X) \to \BS(\un{h}_Y)(1)$ is just $\ep_{Y,X}$ times an $i$-colored startdot in the $k$-th slot. The signed dot $\BS(\un{h}_Y) \to \BS(\un{h}_X)(1)$ is just $\ep_{X,Y}$ times an $i$-colored enddot in the $k$-th slot. \end{defn}

\begin{lemma} The maps $\startdotsign_X^Y$ and $\finaldotsign_Y^X$ are independent of the choice of $\un{h}_Y$. \end{lemma}

\begin{proof} This is a simple sign chase. \end{proof}

\begin{ex} Let $Y = \{0,2\}$ and $X = \{0\}$ (when $n \ge 3$, so that $0$ and $2$ are distant). When one chooses $\un{h}_Y = 02$, the map $\startdotsign$ is $ - \lineblue \startdotgreen$.
When one chooses $\un{h}_Y = 20$, the map is $+ \startdotgreen \lineblue$. These two maps are intertwined by $-\Xbg$ and $-\Xgb$, so they induce the same map $B_X \to B_Y(1)$. \end{ex}

\begin{ex} Let $Y = \{1,2,3,4\}$ and $X = \{1,2,4\}$. Then $w_Y$ has a unique reduced expression $4321$, so that the definition above states that one should use the reduced expression
$421$ for $w_X$. If one has decided instead that $\un{h}_X = 214$ so that $B_X \cong \BS(214)$, then to define $\startdotsign_X^Y$ from $\BS(214)$ to $\BS(4321)$, one must first apply a
signed rex move from $\BS(214)$ to $\BS(421)$, and then apply a 3-colored startdot $\BS(421) \to \BS(4321)$. In Example \ref{ex:proofofsame} we see this in practice. \end{ex}

We encourage the reader to confirm that the differentials in $\FC$ for $n=2, 3, 4$ from \S\ref{sec-234} are all of the form $\startdotsign_X^Y$ or $\finaldotsign_Y^X$. Let
us give now some lemmata which will be used to compute $d^2$. The letters $X$ and $Y$ always refer to proper subsets of $S$.

\begin{lemma} \label{lem:notsame} Let $Y$ contain distinct elements $i,j \in S$, and let $X = Y \setminus i,j$. For shorthand write $Xi$ for $X \cup \{i\}$, and $Xj$ for $X \cup \{j\}$. Then the following equations hold.
\begin{equation} \label{upup} \startdotsign_{Xi}^{Y}\startdotsign_X^{Xi} = - \startdotsign_{Xj}^Y \startdotsign_{X}^{Xj}. \end{equation}
\begin{equation} \label{downdown} \finaldotsign_{Xi}^{X} \finaldotsign_Y^{Xi} = - \finaldotsign_{Xj}^X \finaldotsign_{Y}^{Xj}. \end{equation}
\begin{equation} \label{updown} \finaldotsign_{Y}^{Xj} \startdotsign_{Xi}^Y = - \startdotsign_X^{Xj} \finaldotsign_{Xi}^{X}. \end{equation} \end{lemma}

\begin{proof} In all cases, the composition on both sides of the equality consists of two dots, one colored $i$ and one colored $j$. The only thing to check is that the signs match appropriately, which we leave to the reader. In other words, to prove \eqref{upup} one must show
\begin{equation} \label{upupep} \ep_{Y,Xi} \ep_{Xi,X} = - \ep_{Y,Xj} \ep_{Xj,X},\end{equation}
and similarly for the other equalities.  \end{proof}

\begin{lemma} \label{lem:same} Let $X \subset S$ have size at most $n-2$. Then \begin{equation} \label{selfterm} \sum_{i \in S \setminus X} \finaldotsign_{Xi}^{X} \startdotsign_{X}^{Xi} + \sum_{j \in X} \startdotsign_{X \setminus j}^{X} \finaldotsign_{X}^{X \setminus j} = \id_{B_X} \cdot (\sum_{k \in S} \al_k) = \id_{B_X} \cdot \delta. \end{equation} \end{lemma}

This lemma is best proven by example.

\begin{example} \label{ex:proofofsame} Let $n=9$ and $X = \{1,2,4,5,6\}$. We will identify $B_X$ with $\BS(21654)$ below, so that the entire computation takes place within the endomorphism ring of $\BS(21654)$.
We partition $S$ into four disjoint subsets: $\{3,2,1\}$, $\{7,6,5,4\}$, $\{8\}$, and $\{0\}$. One takes each connected component of $X$ and groups it with the missing simple reflection next in the
order (e.g. $\{4,5,6\}$ is grouped with $\{7\}$); each of the remaining simple reflections forms its own singleton set. The subsets in this partition are called the \emph{$\tau$-components} of $X$; arguments involving $\tau$-components will return in force in the technical section of this paper, starting in \S\ref{subsec-taucomponents}.

We claim that if one restricts the sums on both sides of \eqref{selfterm} to those $i, j, k$ within a single $\tau$-component, one still has an equality. If true, then summing over all
$\tau$-components proves \eqref{selfterm}.

Consider first the $\tau$-component $\{8\}$. Letting $i = 8$, we can choose the reduced expression $216548$ for $h_{Xi}$. Then $\finaldotsign_{Xi}^{X} \startdotsign_{X}^{Xi}$ will create
an 8-colored barbell on the right side of the diagram. Note that any sign appearing in $\startdotsign_{X}^{Xi}$ will cancel the sign appearing in $\finaldotsign_{Xi}^X$ (since $\ep_{X,Xi} = \ep_{Xi,X}$), so that the barbell appears
with positive sign. Since a barbell is multiplication by a simple root, we have \begin{equation} \label{eq:group8} \finaldotsign_{X8}^{X} \startdotsign_{X}^{X8} = \id_{\BS(21654)} \cdot
\al_8, \end{equation}
which matches the corresponding term $k=8$ on the RHS of \eqref{selfterm}.

An identical argument handles the group $\{0\}$. One can choose the reduced expression $216540$ for $h_{X0}$.

Consider now the $\tau$-component $\{7,6,5,4\}$. Letting $i=7$, we can choose the reduced expression $217654$ for $h_{Xi}$. Then $\finaldotsign_{Xi}^{X} \startdotsign_{X}^{Xi}$ will create a barbell $\al_7$ in the slot between $1$ and $6$ (again, with positive sign, because the sign $\ep_{X,Xi} = \ep_{Xi,X}$ appears in both $\startdotsign_{X}^{Xi}$ and $\finaldotsign_{Xi}^{X}$). Now we use the barbell forcing relations \eqref{eq:barbellforce} to force this barbell to the right. In the first step, this breaks strand $6$ with coefficient $\pa_6(\al_7) = -1$, and produces $s_6(\al_7) = \al_6 + \al_7$ in the slot between $6$ and $5$. Continuing to force this polynomial to the right, one breaks strand $5$ with coefficient $\pa_5(\al_6+\al_7) = -1$, and produces $s_5(\al_6 + \al_7) = \al_5 + \al_6 + \al_7$ in the slot between $5$ and $4$. Once more, the strand $4$ is broken with coefficient $\pa_4(\al_5 + \al_6 + \al_7) = -1$, and produces $s_4(\al_5 + \al_6 + \al_7) = \al_4 + \al_5 + \al_6 + \al_7$ on the far right of the diagram. Now, letting $j$ equal $4$, $5$, or $6$, we see that $\startdotsign_{X \setminus j}^{X} \finaldotsign_{X}^{X \setminus j}$ is exactly the broken strand with coefficient $+1$ (again, all signs cancel). Thus, adding together the whole group, all the terms with broken lines cancel, and one has
\begin{equation} \label{eq:group7654} \finaldotsign_{X7}^{X} \startdotsign_{X}^{X7} + \sum_{j = 4,5,6} \startdotsign_{X \setminus j}^{X} \finaldotsign_{X}^{X \setminus j} = \id_{\BS(21654)} \cdot (\al_4 + \al_5 + \al_6 + \al_7). \end{equation}

We claim that a similar statement holds for the $\tau$-component $\{3,2,1\}$, namely that \begin{equation} \label{eq:group321} \finaldotsign_{X3}^{X} \startdotsign_{X}^{X3} + \sum_{j =
1,2} \startdotsign_{X \setminus j}^{X} \finaldotsign_{X}^{X \setminus j} = \id_{\BS(21654)} \cdot (\al_1 + \al_2 + \al_3). \end{equation} The main wrinkle is that, letting $i = 3$, one
must choose the reduced expression $654321$ for $h_{Xi}$, since adding $3$ joined two components into one. To compute $\finaldotsign_{X3}^{X} \startdotsign_{X}^{X3}$ correctly, one must
first apply a signed rex move to go from $\BS(21654)$ to $\BS(65421)$, then multiply by $\al_3$ in the spot between $4$ and $2$, and then apply a signed rex move to go back from
$\BS(65421)$ to $\BS(21654)$. By the same polynomial forcing calculation as in the previous paragraph, we see that multiplication by $\al_3$ in the appropriate slot of $\BS(65421)$ will
break the strands $2$ and $1$ with coefficients $-1$, and will multiply on the right by $\al_1 + \al_2 + \al_3$. Pre- and post-composing with the appropriate signed rex moves (whose signs
cancel, because the rex move is paired with its inverse), we still obtain the breaking of strands $2$ and $1$, and right multiplication by $\al_1 + \al_2 + \al_3$, now as an endomorphism
of $\BS(21654)$. This proves \eqref{eq:group321}.

By analyzing \eqref{eq:group8}-\eqref{eq:group321}, we see that restricting both sides of \eqref{selfterm} to a single $\tau$-component still gives an equality. Taking the sum over all $\tau$-components, one proves \eqref{selfterm}. This example is sufficiently complicated to demonstrate the behavior which occurs for an arbitrary $X$. \end{example}

%===========
\subsection{The standard representation}
\label{subsec-stdrep}
%===========

\begin{defn} \label{def:Pk} For each $k \in \Z$ let $P_k$ denote the following set of subsets of $S = S_{\aff}$:
\begin{equation} P_k = \{ X \subsetneq S \mid |X| = m \text{ where } 0 \le m \le n-1-|k| \text{ and } n-1-k-m \text{ is even} \}. \end{equation} \end{defn}

So, for example, $P_0$ consists of all proper subsets of $S$ whose parity agrees with $n-1$, while $P_{n-1} = P_{1-n} = \{ \mt \}$, and $P_k$ is empty for $|k| \ge n$. Note that $P_k = P_{-k}$.

\begin{defn} \label{defn:FC} Let $\FC$ denote the following pseudocomplex, which we call the \emph{twisted standard complex}. The chain object in homological degree $k$ is
\begin{equation} \FC^k = \bigoplus_{X \in P_k} B_X(k). \end{equation}
The differential, restricted to a summand $B_X(k)$ in homological degree $k$, is the sum
\begin{equation} \label{def:d} d = \bigoplus_{\substack{i \notin X\\ Xi \in P_{k+1}}} \startdotsign_{X}^{Xi} \oplus \bigoplus_{\substack{j \in X\\X \setminus j \in P_{k+1}}} \finaldotsign_{X}^{X \setminus j}. \end{equation}
\end{defn}

\begin{remark} \label{rmk:onlypm1} Let us note that each summand of the differential goes from an object $B_X(k)$, with $X$ of size $m$, to an object $B_Y(k+1)$, with $Y$ of size $m \pm 1$. This differential satisfies the strand-counting sign rule, see \S\ref{subsec-signconvention}. \end{remark}
	
\begin{remark} \label{rmk:nobiggerdiff} Note that for two proper subsets $X, Y \subsetneq S$ there is a degree $+1$ map $B_X \to B_Y$ if and only if $X$ and $Y$ differ by a single index,
in which case the degree $+1$ Hom space is spanned by the $XY$-dot. In particular, every possible nonzero map between the summands of $\FC$ in adjacent homological degrees does actually
occur (up to scalar) in the differential of $\FC$. \end{remark}

\begin{remark} \label{rmk:changingsignnotchangingcomplex} In Remark \ref{rmk:changingsignglobally} we noted that changing the preferred reduced expression $\un{h}_X^{\$}$ will change the
maps $\BXto$ and $\BXto\inv$ by a global sign. This multiplies all the summands of the differential $d_{\FC}$ which go in or out of $B_X$ by $\pm 1$. The result is an isomorphic complex,
where the isomorphism is the identity everywhere except $B_X$, and on $B_X$ is $\pm \id$. Thus the choice of $\un{h}_X^{\$}$ is irrelevant, up to isomorphism. \end{remark}

\begin{lemma} \label{lem:d21} Let $X \in P_k$ and $Y \in P_{k+2}$ with $X \ne Y$.  Then the summand of $d^2$ mapping from $B_X(k)$ to $B_Y(k+2)$ is zero. \end{lemma}

\begin{proof} A term of $d^2$ which maps $B_X(k)$ to $B_Y(k+2)$ will first add or remove some index $i$ (using $\startdotsign_X^{Xi}$ or $\finaldotsign_X^{X \setminus i}$ respectively),
then add or remove some index $j$. We must have $j \ne i$ or else $Y = X$. For any such term, there is another matching term, which first adds or removes $j$, and then adds or removes
$i$. By Lemma \ref{lem:notsame}, these two terms of $d^2$ cancel. \end{proof}

\begin{lemma} \label{lem:d22} Let $X \in P_k \cap P_{k+2}$. Then the summand of $d^2$ mapping from $B_X(k)$ to $B_X(k+2)$ is $\id_{B_X} \cdot \delta$. \end{lemma}

\begin{proof} This is an immediate consequence of Lemma \ref{lem:same}. Note that if $X \in P_k \cap P_{k+2}$, then $Xi \in P_{k+1}$ for any $i \notin X$, and $X \setminus j \in P_{k+1}$ for any $j \in X$. \end{proof}

\begin{prop} $\FC$ is a well-defined pseudocomplex, which is a true complex when $\delta=0$. It is perverse, concentrated in homological degrees $-(n-1)$ through $n-1$, and is self-dual.
\end{prop}

\begin{proof} The computation of $d^2$ is accomplished by the previous lemmas. Self-duality follows because $P_k = P_{-k}$, and the dual of the morphism
$\startdotsign_{X,Y}$ is $\finaldotsign_{Y,X}$. The remaining properties are obvious from the definition. \end{proof}

\begin{defn} \label{defn:VC} Let $\VC = \VC_{\fund_1} \define \Om \FC$ inside $\HC_{\ext}$. We call it the \emph{standard (Gaitsgory) complex}. \end{defn} 

%===========
\subsection{The dual representation}
\label{subsec-dualrep}
%===========

Recall the automorphism $\si$ which flips the affine Dynkin diagram across the affine vertex, $\si(s_i) = s_{-i}$. This induces an autoequivalence $\si$ of $\Diag_{\ext}$.

\begin{lemma} \label{lem:sihX} For $X \subsetneq S$ one has $\si(h_X) = h_{\si(X)}\inv$. This is the Coxeter element $s_{x_1} s_{x_2} \cdots s_{x_d}$ where the elements $\{x_1 < x_2 < \cdots < x_d\}$ of $\si(X)$ appear in cyclic order. \end{lemma}

\begin{proof} Straightforward. \end{proof}

The entire story told above in this chapter could be repeated for the cyclic order on $X$, rather than the anticyclic order. One could define $B'_X$, identified with the Soergel bimodule
$B_{h_X\inv}$, by choosing a reduced expression of $h_X\inv$. One can define $XY$-dots, and a complex $\FC^*$ analogous to to $\FC$, whose summands in degree $k$ are $B'_X(k)$ for $X \in
P_k$. Meanwhile, applying $\si$ to $\FC$ would produce a complex whose summands in degree $k$ are $B'_{\si(X)}(k)$ for $X \in P_k$. But $\si$ preserves the set $P_k$, so these two
complexes are isomorphic. This gives two constructions of the complex $\FC^*$, from which we choose one for our definition.

\begin{defn} Let $\FC^*$ denote $\si(\FC)$. \end{defn}
	
\begin{defn} Let $\VC^*$ denote $\si(\VC) \cong \Om\inv \FC^*$. Let $\VC_{\fund_{n-1}}$ denote $\Om^{n-1} \FC^*$. \end{defn}

These will be the complexes associated with the dual representation of $V$, and the fundamental representation which is related to it by a determinant twist.  Any theorem about $\FC$ can be adapted easily to $\FC^*$ by applying the functor $\si$.

%===========
\subsection{The monodromy}
\label{subsec-monodromy}
%===========

We now construct the monodromy map, and the pseudocomplex $\ImM$ through which it factors.

\begin{defn} If $k \ge 0$, let $Q_k = P_{k+1}$. If $k \le 0$, let $Q_k = P_{k-1}$. Note that $Q_0 = P_1 = P_{-1}$. Let $\ImM$ denote the following pseudocomplex. The chain object in homological degree $k$ is
\begin{equation} \ImM^k = \bigoplus_{X \in Q_k} B_X(k). \end{equation}
The differential, restricted to a summand $B_X(k)$ in homological degree $k$, is defined analogously to \eqref{def:d} but with a sign:
\begin{equation} \label{def:dI} -d = \bigoplus_{\substack{i \notin X\\ Xi \in Q_{k+1}}} \startdotsign_{X}^{Xi} \oplus \bigoplus_{\substack{j \in X\\X \setminus j \in Q_{k+1}}} \finaldotsign_{X}^{X \setminus j}. \end{equation} \end{defn}

One can duplicate the arguments of \S\ref{subsec-stdrep} to prove the following.

\begin{lemma} $\ImM$ is a well-defined pseudocomplex, which is a true complex when $\delta=0$. It is perverse, concentrated in homological degrees $-(n-2)$ through $n-2$, and is self-dual. \end{lemma}

Now we view $\ImM$ as a sub and quotient complex of $\FC$.

\begin{lemma} There is a map $\ImM(1)[-1] \to \FC$, including into certain selected summands, which makes $\ImM(1)[-1]$ a subcomplex of $\FC$. Similarly, there is a quotient map $\FC \to \ImM(-1)[1]$. \end{lemma}

\begin{proof} For $k \ge 0$ the summands of $\ImM^k$ are indexed by $P_{k+1}$, as are the summands of $\FC^{k+1}$. For $k < 0$ the summands of $\ImM^k$ are indexed by $P_{k-1}$, which is
a subset of $P_{k+1}$, which indexes the summands of $\FC^{k+1}$. This determines the inclusion map $\ImM(1)[-1] \to \FC$. By construction, if $B_X(k)$ and $B_Y(k+1)$ are summands of
$\ImM(1)[-1]$ then the summand of the differential between them is the same in $\ImM(1)[-1]$ and in $\FC$. This is the reason for the sign in \eqref{def:dI}, because the differential is
negated by a homological shift $[-1]$. Finally, using Remark \ref{rmk:onlypm1}, one can easily argue that no summand of the differential goes from a summand $B_X(k)$ inside $\ImM(1)[-1]$
to a summand $B_Y(k+1)$ outside of it, so it is a subcomplex.

The arguments for $\ImM(-1)[1]$ are similar. \end{proof}

We can now define the \emph{monodromy map}.

\begin{defn} We let $\mu \co \FC \to \FC(2)[-2]$ denote the composition
\begin{equation} \label{eq:defnM} \FC \to \ImM(-1)[1] \to \FC(-2)[2] \end{equation}
of the quotient map and the inclusion map. \end{defn}

\begin{lemma} In the pseudocomplex $\FC$, one has $d^2 = \mu \delta$. \end{lemma}

\begin{proof} This is just a consolidation and restatement of Lemmas \ref{lem:d21} and \ref{lem:d22}. \end{proof}
	
\begin{remark} We identify $\Ga(\FC)$ or $\Ga(\VC)$ with $V$ in the obvious way. Then it is clear that $\Ga(\mu)$ agrees with the action of the regular nilpotent element in $\gl_n$. Since this is true for $\VC$, it will be true for tensor products of $\VC$ by \eqref{eq:monodromyontensor}. \end{remark}

%===========
\subsection{Rotating cyclical bimodules}
\label{subsec-rotations}
%===========

The remainder of this chapter will use the diagrammatic presentation of the extended Hecke category $\Diag_{\ext}$ from \S\ref{sec:extendeddiag}.

In \eqref{eq:Rtauconj} it is proven that $\Om B_i \cong B_{i+1} \Om$ via the mixed crossing $\mixedOred$. Composing these mixed crossings, we obtain an isomorphism $\Om B \cong
\tau(B) \Om$ for any Bott-Samelson object $B$. Let us construct the analogous isomorphism for for $B_X$.

\begin{defn} \label{defn:cycle} Let $X$ be a proper subset of $S$. Define a morphism $\cycle_X \co \Om B_X \to B_{\tau(X)} \Om$ as follows. Choose any reduced expression $\un{h}_X$ for $h_X$; we will use the reduced expression $\tau(\un{h}_X)$ for $h_{\tau(X)}$. Then let $\cycle_X$ be the composition \begin{equation} \Om B_X \simto \Om \BS(\un{w}_X) \to \BS(\tau(\un{w}_X)) \Om \simto B_{\tau(X)} \Om. \end{equation} The first map is $\id_\Om \ot \BXto_X$, the last map is $\BXto\inv_{\tau(X)} \ot \id_\Om$ ,and the middle map is a composition of mixed crossings. The overall map $\cycle_X$ is depicted below. 
\begin{equation} \label{eq:cyclemap} \cycle_X = \ig{1}{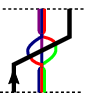} \end{equation}
\end{defn}

\begin{lemma} This map $\cycle_X$ is independent of the choice of reduced expression $\un{h}_X$. \end{lemma}

\begin{proof} It is enough to show that $\cycle_X$ is unchanged when the reduced expression $\un{h}_X$ is altered by a commuting braid relation $s_j s_k = s_k
s_j$. Using \eqref{eq:BXtocross}, one can pull a 4-valent vertex out of $\BXto$ on bottom, also creating a sign $(-1)$. Using \eqref{eq:mixed6v} (in the case of the 4-valent vertex), this
4-valent vertex can pull up through the black strand to meet $\BXto\inv$. Finally, the 4-valent vertex and the sign $(-1)$ can be absorbed into $\BXto\inv$, again using
\eqref{eq:BXtocross}. After this manipulation, one obtains a diagram which is equal to $\cycle_X$, and looks just like \eqref{eq:cyclemap} except with $\un{h}_X$ altered by this commuting
relation. \end{proof}

\begin{lemma} \label{lem:cycleit} When $X = Y \setminus j$, one has
\begin{equation} \cycle_{Y}\startdotsign_X^Y = \startdotsign_{\tau(X)}^{\tau(Y)} \cycle_X \end{equation}
and
\begin{equation} \cycle_{X}\finaldotsign_Y^X = \finaldotsign_{\tau(Y)}^{\tau(X)} \cycle_Y. \end{equation}
\end{lemma}

\begin{proof} This is an easy consequence of \eqref{eq:mixeddot}, when one uses a reduced expression $\un{h}_Y$ for $h_Y$ and the corresponding expression $\tau(\un{h}_Y)$ for $\tau(h_Y)$. \end{proof}

\begin{prop} There are isomorphisms $\cycle \co \FC \Om \to \Om \FC$ and $\cycle \co \ImM \Om \to \Om \ImM$, sending the summand $B_X(k) \Om$ to $\Om B_{\tau(X)}(k)$ by $\cycle_X$.
\end{prop}

\begin{proof} We need only check that $\cycle$ commutes with the differential, which is established in Lemma \ref{lem:cycleit}. Note that $\tau$ preserves the sets $P_k$ and $Q_k$. \end{proof}

%===========
\subsection{Wakimoto filtations}
\label{subsec-wakifilts}
%===========

\begin{thm} \label{thm:wakifilt} The complex $\VC$ has a Wakimoto filtration compatible with the weight filtration on $V$. \end{thm}

The Wakimoto complexes which should appear are the Rouquier complexes $Y_i$ of the elements $y_i \in \Br_{\ext}$, for $1 \le i \le n$, see \S\ref{subsec-translattice}. We have demonstrated the Wakimoto
filtrations in the case $n=2, 3$ in \S\ref{sec-234}, and recommend that the reader find it for $n=4$ as an exercise.

\begin{proof} The elements $y_i$ have a nice description in \eqref{eq:foryi}. In particular, $y_i = \om \tilde{h}_{X_i}$, where $X_i = S \setminus \{i-1\}$ and $\tilde{h}_{X_i}$ is some
particular negative-positive lift of the element $h_{X_i} \in W_{\ext}$. For any subset $X \subset X_i$, the subexpression of $h_{X_i}$ corresponding to this subset is just $h_X$. Thus
(ignoring the copy of $\Om$ on the left) the Rouquier complex $Y_i$ is a cube whose vertices are associated with subsets $X \subset X_i$, and where the corresponding summand is $B_{X}$.
By our sign conventions on Rouquier complexes (see \S\ref{subsec-signconvention}), the differentials in $Y_i$ are precisely the signed dot maps! Thus, to prove that $\VC$ is filtered
with subquotients $Y_i$, it is enough to prove the combinatorial statement that one can partition the summands $\coprod P_k$ of $\VC$ into cubes $Y_i$, and that there are no
differentials from the cube $Y_i$ to the cube $Y_j$ when $j > i$.

Recall that $P_k$ consists of those subsets $X \subsetneq S$ where $|k| = n - 1 - |X| - 2m$ for some positive number $m$. Let $P = \coprod_k P_k$ and $U_m = \coprod_k \{X \in P_k \mid k = n-1-|X| - 2m\}$. Then $\coprod_m U_m$ is a partition of $P$, where $m$ ranges from $0$ to $n-1$. For example, $U_0$ consists of all proper subsets $X \subsetneq S$ appearing in homological degree $k = |X|$, while $U_{n-1}$ is just the single copy of $X = \mt$ appearing in homological degree $-(n-1)$. The differentials in $\VC$ come in two kinds: maps $B_X \to B_{Xi}$ for $i \notin X$, which always send a term $X \in U_m$ to  $Xi \in U_{m-1}$, and maps $B_X \to B_{X \setminus j}$ for $j \in X$, which always send a term $X \in U_m$ to $X \setminus j \in U_m$. In particular, the terms in $U_0$ form a subcomplex of $\VC$, since $U_{-1}$ is empty so the differential preserves $U_0$. In similar fashion, the terms in $U_0 \cup U_1 \cup \cdots \cup U_d$ form a subcomplex for any $0 \le d \le n-1$.

Note that $Y_1$ is a positive Rouquier complex. The cube $Y_1$ contains all subsets $X$ of $\{1, \ldots, n-1\}$, each appearing in homological degree equal to $n-1-|X|$. It is contained in $U_0$, and contains all elements $X$ of $U_0$ for which $0 \notin X$. Any differential in $\FC$ from a term in $U_0$ must be of the form $B_X \to B_{X \setminus j}$, thus it stays within $U_0$ and preserves the condition $0 \notin X$. Hence, $Y_1$ is a subcomplex.

Note that $Y_2$ is mostly a positive Rouquier complex, except that the simple reflection $s_0$ is lifted negatively. The cube $Y_2$ contains all subsets $X \subset \{2, \ldots, n-1, 0\}$,
appearing in homological degree equal to $n-1-|X|$ if $0 \in X$, and in homological degree $n-3-|X|$ if $0 \notin X$. Thus $Y_2 \subset U_0 \cup U_1$. Note that $Y_2 \cap U_0$ consists of
those $X \in U_0$ with $0 \in X$ and $1 \notin X$, and $Y_2 \cap U_1$ consists of those $X \in U_1$ for which $0,1 \notin X$. We argue that $Y_1 \cup Y_2$ is a subcomplex. If $X \in (Y_1
\cup Y_2) \cap U_0$ then so is $X \setminus j$ for any $j \in X$, because removing an element preserves the condition that either $0 \notin X$ or $0 \in X$ and $1 \notin X$. If $X \in
(Y_1 \cup Y_2) \cap U_1$ then $0,1 \notin X$. If $j \in X$ then $X \setminus j$ still does not contain $0$ or $1$ so it remains in $(Y_1 \cup Y_2) \cap U_1$. If $i \notin X$ then $Xi \in
(Y_1 \cup Y_2) \cap U_0$: if $i \ne 0$ then $Xi \in Y_1 \cap U_0$ and if $i = 0$ then $Xi \in Y_2 \cap U_0$. Thus any possible differential preserves $Y_1 \cup Y_2$.

Similarly, $Y_3$ can be described by its intersection with each $U_m$. \begin{itemize} \item $Y_3 \cap U_0$ consists of those $X \in U_0$ with $0, 1 \in X$ and $2 \notin X$. \item $Y_3
\cap U_1$ consists of those $X \in U_1$ with either $0$ or $1$ in $X$ (but not both) and $2 \notin X$. \item $Y_3 \cap U_2$ consists of those $X \in U_2$ with $0, 1, 2 \notin X$.
\end{itemize} The reader can probably guess the rest of the proof, and is welcome to skip it.

More generally, $Y_k \cap U_m$ consists of those $X \in U_m$ for which $k-1 \notin X$ and $\{0, 1, \ldots, k-2\} \cap X$ has size $k - 1 - m$. These sets are disjoint. Otherwise, suppose
$X \in Y_k \cap Y_{\ell} \cap U_m$ with $k < \ell$. We know that $X \cap \{0, 1, \ldots, k-2\}$ has size $k-1-m$. In order for $X \cap \{0, 1, \ldots, \ell-2\}$ to also have size $\ell -
1 - m$ we must have $|X \cap \{k-1, k, \ldots, \ell-2\}| = \ell - k$, or in other words, $\{k-1, k, \ldots, \ell - 2\} \subset X$. But then $k-1 \in X$, contradicting $X \in Y_k$.

If $X \in Y_k \cap U_m$ then consider $X \setminus j$ for some $j \in X$. If $j \ge k$ then $X \setminus j \in Y_k \cap U_m$, while if $0 \le j \le k-2$ then $X \setminus
j \in Y_{k - 1} \cap U_m$. Thus the union $Y_1 \cup \cdots \cup Y_k$ is closed under differentials coming from the removal of an index.

Now consider $Xi$ for some $i \notin X$. If $0 \le i \le k-2$ then $Xi \in Y_k \cap U_{m-1}$. If $i \ge k-1$, we argue that $Xi \notin Y_{\ell} \cap U_{m-1}$ for $k < \ell$. If $X \cap \{0, 1, \ldots, k-2\}$ has size $k-1-m$ then $Xi \cap \{0, \ldots, k-2\}$ also has size $k-1-m$. If $Xi \cap \{0, 1, \ldots, \ell-2\}$ has size $\ell - 1 - (m-1)$ then $Xi \cap \{k-1, \ldots, \ell-2\}$ has size $\ell - k + 1$, which is impossible. Thus the union $Y_1 \cup \cdots \cup Y_k$ is closed under differentials coming from the addition of an index. \end{proof}

% If $i \ge k-1$ then we split into cases as follows. \begin{itemize} \item If $k-2 \notin X$
% then $Xi \in Y_{k-1} \cap U_{m-1}$. \item If $k-2 \in X$ but $k-3 \notin X$ then $Xi \in Y_{k-2} \cap U_{m-1}$. \item If $k-2, k-3 \in X$ but $k-4 \notin X$ then $Xi \in Y_{k-3} \cap
% U_{m-1}$. \item Etcetera. \item Finally, if $\{0, 1, \ldots, k-2\} \subset X$ then $k - 1 - m = k-1$ so $m=0$. In this case, we need not consider $Xi$ at all. \end{itemize}

\begin{remark} \label{rmk:determiningVC} In the complex $\VC$, some differentials live within a given subquotient complex $Y_k$, and some go between them. However, every differential is
in the $\tau$-orbit of a differential within some $Y_k$. Thus, if one knows that $\VC$ is Wakimoto-filtered as in the theorem, and that $\VC$ is $\tau$-invariant, then one can deduce the
exact form of $\VC$. Note that this process does not miss any summands of the differential, see Remark \ref{rmk:nobiggerdiff}. This is how we guessed what $\VC$ should be. \end{remark}

\section{Longest elements and functorial twisting}
\label{sec-longestfunctorial}
%%%%%%%%%%%%%%%%%%%%%%%%%

There are two main goals of the next several chapters. The first is to give an explicit description of the tensor products $B_I \FC$ and $\FC B_J$ for finitary subsets $I,J \subsetneq
S$. The second is to construct a natural isomorphism \begin{equation} \label{eq:FCtwist} \rotisom \co (-) \ot \FC \to \FC \ot \tau(-) \end{equation} of functors from $\Diag \to \Hot(\Diag)$.
That is, objects commute past $\FC$ up to a functorial twist by the Dynkin diagram automorphism. In particular, our descriptions of $B_I \FC$ and $\FC B_{\tau(I)}$ are the same. The
natural isomorphism of \eqref{eq:FCtwist} will only exist under certain restrictions on the realization.

Combined with the natural isomorphism $(-) \ot \Om \to \Om \ot \tau\inv(-)$ which is built in to the structure of $\Diag_{\ext}$, this gives a natural isomorphism \begin{equation}
\label{eq:VCcentral} (-) \ot \VC \to \VC \ot (-),\end{equation} equipping $\VC$ with the structure of an object in the restricted Drinfeld center. In this paper we focus on $\FC$ and
\eqref{eq:FCtwist} instead of $\VC$ and \eqref{eq:VCcentral}.

In this chapter we state the main results, and reduce the construction of \eqref{eq:FCtwist} to a small number of computations. We also outline the proofs. Following this, from
\S\ref{sec-descentsets} to \S\ref{sec-tensorBIcomps}, we provide the proofs. The proofs provide several insights and new technology for working with special indecomposables in the affine
Hecke category, but for most readers they are not essential; the outline gives a good idea of the contents. In \S\ref{subsec-rotisoms} the commutation maps $\rotisom_{B_s}$ are described
very explicitly.

% This is stated in Theorem \ref{thm:tensorBI}, and proven in the
% subsequent sections. In particular, this theorem will imply that $\FC \ot B_{\tau(I)} \cong B_I \ot \FC$, and will even fix this isomorphism. Moreover, the form of the answer makes it
% relatively easy to deduce how the isomorphism $\FC \ot B_{\tau(I)} \cong B_I \ot \FC$ intertwines with various morphisms in the diagrammatic Hecke category. Ultimately, this leads to the
% proof that $\FC \ot \tau(B) \cong B \ot \FC$ functorially, see Theorem \ref{thm:commutepastF}, which we also prove in this chapter.

%===========
\subsection{Statement of the theorem}
\label{subsec-statementtensorBI}
%===========

\begin{defn} Fix $I \subsetneq S = S_{\aff}$. For each $k \in \Z$ let $N^I_k$ be the subset of $P_k$ (see Definition \ref{def:Pk}) consisting of those $X \in P_k$ where $\tau(I) \subset X$. \end{defn}

So, for example, $N^I_0$ consists of all sets $\tau(I) \subset X \subsetneq S$ whose parity agrees with $n-1$, while $N^I_{n-1 - |I|} = N^I_{1-n+|I|} = \{ \tau(I) \}$, and $N^I_k$ is empty for $|k| \ge n-|I|$. Note that $N^I_k = N^I_{-k}$.

Recall the element $h_X \in W_{\aff}$ from Definition \ref{defn:hX}. Recall that $w_I$ is the longest element of the parabolic subgroup $W_I$, and $B_I$ is shorthand for $B_{w_I}$. Our
description of $B_I \FC$ will involve the indecomposable objects $B_{w_I h_X}$ for various $X \in N^I_k$. In Lemma \ref{lem:descent} we prove that
\begin{equation} \label{eq:lengthsaddIX} \ell(w_I h_X) = \ell(w_I) + \ell(h_X)\end{equation} whenever $\tau(I) \subset X$. If \eqref{eq:lengthsaddIX} holds then $B_{w_I h_X}$ is the top direct summand of $B_I B_X$.
	
% We fix\footnote{In fact, this inclusion and projection can be made canonical. To do this, one should first show that $\Hom(B_{w_I h_X},B_I B_X)$ is one-dimensional, which will follow from the results of \S\ref{sec-thicker}. Then, one can pick out the inclusion maps in this one-dimensional space by the property that it preserves the lowest degree basis element $c_{\bot}$, as discussed in \cite[\S 4]{EWHodge}. This is a fairly general phenomenon, but it has not yet been written down carefully to my knowledge, and we will not need it.} an inclusion $B_{w_I h_X} \into B_I B_X$ and a projection $B_I B_X \onto B_{w_I h_X}$. These inclusions and projections will be made more precise in \S\ref{sec-thicker}.

\begin{defn} \label{defn:N} Let $\NC_I$ denote the following pseudocomplex. The chain object in homological degree $k$ is
\begin{equation} \NC_I^k = \bigoplus_{X \in N^I_k} B_{w_I h_X}(k). \end{equation}
The differential, restricted to a summand $B_{w_I h_X}(k)$ in homological degree $k$, is the sum
\begin{equation} \label{def:dN} d = \bigoplus_{\substack{i \notin X\\ Xi \in N^I_{k+1}}} {}^I \startdotsign_{X}^{Xi} \oplus \bigoplus_{\substack{j \in X\\X \setminus j \in N^I_{k+1}}} {}^I \finaldotsign_{X}^{X \setminus j}. \end{equation}
Recall from Definition \ref{defn:dotsign} that we have signed dot maps $\startdotsign_X^Y$ or $\finaldotsign_X^Y$ from $B_X$ to $B_Y$ whenever $Y$ differs from $X$ by a single index. The maps ${}^I \startdotsign_X^Y$ and ${}^I \finaldotsign_X^Y$ in \eqref{def:dN} are obtained by taking the corresponding signed dot $B_X \to B_Y$, tensoring on the left with the identity map of $B_I$, precomposing with the inclusion $B_{w_I h_X} \into B_I B_X$, and postcomposing with the projection $B_I B_Y \onto B_{w_I h_Y}$.
\end{defn}

\begin{thm} \label{thm:tensorBI} By an explicit sequence of Gaussian eliminations, one obtains homotopy equivalences $B_I \FC \simto \NC_I$ and $\FC B_{\tau(I)} \simto \NC_I$. In
particular, $B_I \FC$ is perverse for all $I \subsetneq S$. \end{thm}

Let us put the result in some perspective. Set $J = \tau(I)$. The left tensor ideal of $B_J$ only contains indecomposable objects $B_w$ with $J$ in the right descent set of $w$, and the
right tensor ideal of $B_I$ only contains $B_w$ with $I$ in the left descent set of $I$. Thus if $\FC B_J \cong B_I \FC$ then it must be built from $B_w$ with both $J$ in its right
descent set and $I$ in its left descent set. Consequently, one should expect that all summands which do not meet this descent set criterion can be Gaussian eliminated from the complex
$B_I \FC$.

In \S\ref{sec-descentsets} we study the descent sets of elements $w_I h_X$ and $h_X w_J$. Lemma \ref{lem:descent} proves that $w_I h_X$ contains $J$ in its right descent set if and only if $J \subset X$, in which case $\ell(w_I h_X) = \ell(h_X) + \ell(w_I)$ and $w_I h_X$ also contains $I$ in its left descent set. In \S\ref{sec-thicker} we study how the tensor products $B_I B_X$ decompose. We determine that:
\begin{itemize} \item when $J \subset X$ then, aside from the top summand $B_{w_I h_X} \sumset B_I B_X$, all other direct summands $B_w$ of $B_I B_X$ will not have $J$ in the right descent set of $w$, and
\item when $J \nsubseteq X$ then all direct summands $B_w$ of $B_I B_X$ will not have $J$ in the right descent set of $w$. \end{itemize}
Thus, after removing all summands without $J$ in their right descent set, we expect that the remainder has the form $\NC_I$. Since $\NC_I$ is perverse, no further Gaussian elimination should be possible. In particular, every summand that could survive (according to the descent set criterion) will in fact survive.

In \S\ref{sec-tensorBIGE} we explicitly perform the Gaussian elimination which removes all but the desired summands from the tensor product $B_I \FC$, leaving behind the complex $\NC_I$.
This confirms the expectations above, and proves (half of) Theorem \ref{thm:tensorBI}. In order to perform the Gaussian elimination, we effectively compute the idempotents projecting to
every direct summand of $B_I B_X$ for any $X$. Thankfully, every direct summand $B_w$ of $B_I B_X$ is smooth (i.e. the Kazhdan-Lusztig polynomials of $w$ are all trivial), which
simplifies the form of these idempotents. These computations with idempotents are given in \S\ref{sec-thicker}, and use the thick calculus of \cite{EThick}.

Now consider the tensor product $\FC B_J$. By analogous arguments, the only summand $B_w$ of $B_Y B_J$ which has $I$ in its left descent set will be $B_{h_Y w_J}$ in the case when $I
\subset Y$. Thus, after Gaussian elimination, one expects $\FC B_J$ to be a complex built analogously to Definition \ref{defn:N}, except that it is built from $B_{h_Y w_J}$ when $I
\subset Y$, rather than being built from $B_{w_I h_X}$ when $J \subset X$. It is proven in Lemma \ref{lem:descent} that for any $X$ with $J \subset X$ there is a unique $Y$ of the same
size with $I \subset Y$ such that $w_I h_X = h_Y w_J$, and vice versa. By checking that the differentials match up, one proves that $\FC B_J$ and $B_I \FC$ have isomorphic descriptions; this task is accomplished in \S\ref{sec-leftvsright}.

We have already given examples of $B_I \FC$ for $n = 2, 3$ in \S\ref{sec-234}.

\begin{example} Let $n=4$. If $I = \{0\}$ and $J = \{1\}$ then
\begin{equation}
B_0 \FC \cong \FC B_1 \cong \NC_I = \left(
\begin{tikzpicture}
\node (a) at (0,0) {$B_{01}(-2)$};
\node (a0) at (2,2) {$B_{021}(-1)$};
\node (a1) at (2,0) {$B_{010}(-1)$};
\node (a2) at (2,-2) {$B_{031}(-1)$};

\node (b) at (4,-3) {$B_{01}$};
\node (a01) at (4,3) {$B_{0210}$};
\node (a12) at (4,-1) {$B_{0103}$};
\node (a20) at (4,1) {$B_{0321}$};

\node (b0) at (6,2) {$B_{021}(1)$};
\node (b1) at (6,0) {$B_{010}(1)$};
\node (b2) at (6,-2) {$B_{031}(1)$};
\node (c) at (8,0) {$B_{01}(2)$};
\path
	(a) edge (a0)
	(a) edge  (a1)
	(a) edge  (a2)	
	(a0) edge (b)
	(a1) edge  (b)
	(a2) edge  (b)

	(b) edge  (b0)
	(b) edge (b1)
	(b) edge  (b2)	
	(b0) edge (c)
	(b1) edge (c)
	(b2) edge  (c)

	(a0) edge  (a01)
	(a1) edge (a12)
	(a2) edge (a20)
	(a0) edge (a20)
	(a1) edge (a01)
	(a2) edge (a12)

	(a01) edge (b0)
	(a12) edge (b1)
	(a20) edge (b2)
	(a20) edge (b0)
	(a01) edge (b1)
	(a12) edge (b2);
	
\end{tikzpicture} \right)
\end{equation}
We encourage the reader to label each summand $B_{w_{0} h_X}$ in the complex above with the corresponding set $X \subsetneq S$. Using different reduced expressions, such as $1031$ instead of $0103$, the reader should also label each summand $B_{h_Y w_{1}}$ with the corresponding set $Y \subsetneq S$.
\end{example}

\begin{example} Let $n=4$. 
If $I = \{0,1\}$ and $J = \{1,2\}$ then
\begin{equation}
B_{010} \FC \cong \FC B_{121} \cong \NC_I = \left(
\begin{tikzpicture}
\node (a) at (0,0) {$B_{10121}(-1)$};
\node (b) at (2,2) {$B_{210121}(0)$};
\node (c) at (2,-2) {$B_{103121}(0)$};
\node (d) at (4,0) {$B_{10121}(1)$};
\path
	(a) edge (b)
	(a) edge  (c)
	(b) edge  (d)
	(c) edge (d);
\end{tikzpicture}
\right). \end{equation}
Above we used reduced expressions ending in $121$, but encourage the reader to find reduced expressions beginning with $010$.
\end{example}

Here is a crucial consequence of the description of $\FC B_J$.

\begin{prop} \label{prop:dimhomBI} For $I, I' \subsetneq S$ and $k \in \Z$, we have
\begin{equation} \label{eq:FIFJhoms} \dim \Hom(B_I \FC, B_{I'} \FC(k)) = \begin{cases} 0 & \text{if } k < 0, \\ 0 & \text{if } k=0 \text{ and } I \ne I', \\ 1 & \text{if } k=0 \text{ and } I=I', \\ ?? & \text{else}. \end{cases} \end{equation} In particular, $B_I \FC$ is indecomposable. Note that we are discussing the graded dimensions of maps in the homotopy category, with regard to the grading shift but not the homological shift. \end{prop}

\begin{proof} By Theorem \ref{thm:tensorBI}, the complex $B_I \FC$ is perverse for any $I$. There are no maps of negative graded degree between two perverse complexes, proving the first
statement.

In Lemma \ref{lem:nooverlap} we prove that $w_I h_X = w_{I'} h_{X'}$ implies that $X = X'$ and $I = I'$, under the assumption $\tau(I) \subset X$ and $\tau(I') \subset X'$.
Thus the indecomposable summands $B_{w_I h_X}$ of $B_I \FC$ and $B_{w_{I'} h_{X'}}$ of $B_{I'} \FC$ are disjoint if $I \ne I'$. There can be no nonzero chain maps of degree zero between
these perverse complexes, since there are no degree zero maps in any given homological degree.

When $I = I'$, note that $\NC_I$ is a multiplicity-free perverse complex. The remainder of the proposition follows from Lemma \ref{lem:multfreejustid} and its Porism. \end{proof} 

%===========
\subsection{Proving that twisting is functorial}
\label{subsec-provingrotation}
%===========

According to Theorem \ref{thm:tensorBI}, $B_I \FC \cong \FC B_{\tau(I)}$. By Proposition \ref{prop:dimhomBI}, the choice of this isomorphism is unique up to scalar. Thus there is a map $\rotisom_I \co B_I \FC \cong \FC B_{\tau(I)}$ which can be constructed explicitly by following the Gaussian elimination, and for which we have some freedom in a choice of scalar.

In \S\ref{subsec-rotisoms} we give an explicit description of the map $\rotisom_s$, which pins down the choice of scalar. For any sequence $\un{w}$, we define $\rotisom_{\un{w}} \co \BS(\un{w}) \FC \to \FC \BS(\tau(\un{w}))$ to be the composition of the morphisms $\rotisom_s$ over the simple reflections in the sequence. As a composition of isomorphisms, it is an isomorphism. For any morphism $f \co \BS(\un{w}) \to \BS(\un{x})$, we can consider the following square $\square_f$, where the horizontal arrows are isomorphisms.
\begin{equation} \label{fsquare} \square_f = \quad \left(
\begin{diagram}
\BS(\un{w}) \ot \FC & \rTo^{\rotisom_{\un{w}}} & \FC \ot \BS(\tau(\un{w})) \\
\dTo^{f \ot \FC} & & \dTo_{\FC \ot \tau(f)} \\
\BS(\un{x}) \ot \FC & \rTo^{\rotisom_{\un{x}}} & \FC \ot \BS(\tau(\un{x})) \end{diagram} \right)
\end{equation}
If $\square_f$ commutes for all $f$, then the morphisms $\rotisom_{\un{w}}$ define a natural isomorphism of functors $(-) \ot \FC \to \FC \ot \tau(-)$.

\begin{remark} In this section, the functors we discuss are functors $\Diag \to \Hot(\Diag)$. We have already discussed in \S\ref{ssec:intro_centrality} what it would take to extend the
isomorphism \eqref{eq:FCtwist} to an isomorphism of endofunctors of $\Hot(\Diag)$. \end{remark}

\begin{thm} \label{thm:rotationfunctorial} Assume that $\Diag$ is constructed using the root realization or the standard affine realization. The morphisms $\rotisom_{\un{w}}$ define a natural isomorphism $(-) \ot \FC \to \FC \ot \tau(-)$ as functors $\Diag \to \Hot(\Diag)$. That is, $\square_f$ commutes for all $f$. \end{thm}

We prove this theorem, modulo some crucial computations, in Lemma \ref{lem:reducingrotation} below.

It is easy to observe that if $\square_f$ and $\square_g$ commute, then $\square_{f \circ g}$ commutes. Clearly $\square_{\id}$ commutes for any identity map. Thus one need only check
$\square_f$ for the generating morphisms of $\Diag$.

It will be temporarily useful to state a slightly weaker condition. For a morphism $f \co \BS(\un{w}) \to \BS(\un{x})$ we say that $\square_f$ \emph{commutes up to scalar} if there is some scalar $\lambda$ such that the following diagram commutes.
\begin{equation} \label{fsquarelambda} \left(
\begin{diagram}
\BS(\un{w}) \ot \FC & \rTo^{\rotisom_{\un{w}}} & \FC \ot \BS(\tau(\un{w})) \\
\dTo^{f \ot \FC} & & \dTo_{\lambda \FC \ot \tau(f)} \\
\BS(\un{x}) \ot \FC & \rTo^{\rotisom_{\un{x}}} & \FC \ot \BS(\tau(\un{x})) \end{diagram} \right)
\end{equation}
Of course, if $\lambda=1$ then $\square_f$ commutes honestly. Note that if $\square_f$ and $\square_g$ each commute up to a scalar, then $\square_{f \circ g}$ also commutes up to a scalar, and the scalar for $f \circ g$ is the product of the scalars
for $f$ and for $g$.  

It is a relatively straightforward exercise when $n=2$ to compute $\square_f$ when $f$ is multiplication by any linear polynomial. In particular, $\square_{\al_0}$ and $\square_{\al_1}$
commute. However, consider a realization in type $\tilde{A}_1$ which possesses a fundamental weight $\om_0$ with $\pa_0(\om_0) = 1$ and $\pa_1(\om_0) = 0$. One can compute that left
multiplication by $\om_0$ on $\FC$ is not homotopic to right multiplication by any polynomial, even accounting for possible monodromy. Consequently, $\square_{\om_0}$ does not commute,
not even up to scalar, and there is no isomorphism of functors $(-) \ot \FC \to \FC \ot \tau(-)$ for this realization! For this reason we restrict our attention momentarily to the root
realization (or more generally, to the (non-full) subcategory of $\Diag$ where the polynomial ring is generated by the simple roots) where there is still hope of such an isomorphism.
Afterwards we will state the generalization to other realizations.

We now analyze the generators of $\Diag_S$. Consider the square \eqref{fsquare} associated with the merging trivalent vertex.
\begin{equation} \label{mergesquare}
\begin{diagram}
B_s B_s \FC & \rTo^{\rotisom_{ss}} & \FC B_{\tau(s)} B_{\tau(s)} \\
\dTo^{\mergered \FC} & & \dTo_{\FC \mergeblue} \\
B_s \FC & \rTo^{\rotisom_s} & \FC B_{\tau(s)} \end{diagram}
\end{equation}

\begin{lemma} \label{lem:bsbsbs-1} For $s \in S$, we have $\dim \Hom^{-1}(B_s B_s \FC, \FC B_{\tau(s)}) = 1$. \end{lemma}

\begin{proof} We use Proposition \ref{prop:dimhomBI}, and the fact that $\FC B_{\tau(s)} \cong B_s \FC$. We have \[ \dim \Hom^{-1}(B_s B_s \FC ,  B_s \FC) = \dim \Hom^{-2}(B_s \FC, B_s \FC) + \dim \Hom^{0}(B_s \FC, B_s \FC) = 1,\] since $B_s B_s \cong B_s(1) \oplus B_s(-1)$. 
\end{proof}

Thus the two paths around \eqref{mergesquare} live in a one-dimensional hom space, so \eqref{mergesquare} commutes up to a scalar $\lambda$. In fact, $\lambda$ is invertible, since both
paths are nonzero: following the proof of Lemma \ref{lem:bsbsbs-1}, the space $\Hom^{-1}(B_s B_s \FC, B_s \FC)$ is isomorphic to $\Hom^0(B_s \FC, B_s \FC)$ via the trivalent vertex. By
similar arguments, the square associated to the splitting trivalent vertex also commutes up to an invertible scalar.

Let us argue that the square associated to the $2m_{st}$-valent vertex commutes up to an scalar. Suppose for sake of example that $m_{st}=3$. The $2m_{st}$-valent vertex factors as a composition $B_s B_t B_s \to B_{sts} \to B_t B_s B_t$ of a projection and an inclusion to the longest element of the dihedral parabolic subgroup. We need only show that the square associated to the projection (resp. the inclusion) commutes up to scalar, which follows from the following lemma.

\begin{lemma} \label{lem:bsbtbs0} If $m_{st} = 3$ then $\dim \Hom^0(\FC B_s B_t B_s,\FC B_{sts}) = 1$. \end{lemma}

\begin{proof} We know that $B_s B_t B_s \cong B_{sts} \oplus B_s$. By Proposition \ref{prop:dimhomBI}, $\Hom^0(\FC B_s, \FC B_{sts}) = 0$, and $\dim \Hom^0(\FC B_{sts}, \FC B_{sts}) = 1$.
\end{proof}

To deal with the case $m_{st} = 2$, we have the analogous lemma. \begin{lemma} If $m_{su} = 2$ then $\dim \Hom^0(\FC B_s B_u, \FC B_{su}) = 1$. \end{lemma}

The remaining generators of $\Diag_S$ are the dots, which have the following squares.
\begin{equation} \label{startdotsquare}
\begin{diagram}
R \FC & \rTo^{\rotisom_{\mt}} & \FC R \\
\dTo^{\startdotred \FC} & & \dTo_{\FC \startdotblue} \\
B_s \FC & \rTo^{\rotisom_s} & \FC B_{\tau(s)} \end{diagram}
\end{equation}
\begin{equation} \label{enddotsquare}
\begin{diagram}
B_s \FC & \rTo^{\rotisom_s} & \FC B_{\tau(s)} \\
\dTo^{\finaldotred \FC} & & \dTo_{\FC \finaldotblue} \\
R \FC & \rTo^{\rotisom_{\mt}} & \FC R \end{diagram}
\end{equation}
We have not yet computed the dimension of $\Hom^1(\FC, \FC B_s)$, so we can not immediately deduce that their squares commute up to scalar. Checking directly that these squares commute is our crucial computation.

\begin{lemma} \label{lem:reducingrotation} If \eqref{startdotsquare} and \eqref{enddotsquare} commute then $\square_f$ commutes for all $f \in \Diag$ for the root realization. If $\square_{x_1}$ also commutes, then $\square_f$ will commute for all $f \in \Diag$ for the standard affine realization.  \end{lemma}

\begin{proof} Note that $\Diag$ for the root realization is generated by diagrams (without polynomials), since the roots can be expressed as barbells. Since all the generators of $\Diag$
produce squares which commute up to scalar, then every square for any morphism in $\Diag$ will commute up to scalar. Let $\Theta$ be the autoequivalence of $\Diag$ defined as follows.
It fixes all objects. It rescales all morphisms, sending $f$ to $\lambda f$ where $\lambda$ is the scalar making \eqref{fsquarelambda} commute. This is a functor, because the scalars
multiply over composition. Then the morphisms $\rotisom_{\un{w}}$ define a natural isomorphism $(-) \ot \FC \to \FC \ot \tau(\Theta(-))$.

In \cite{EHBraidToolkit}, we classified\footnote{The paper \cite{EHBraidToolkit} is still work in preparation, though the relevant portion has been written up, and is available upon
request.} the automorphisms of the category $\Diag$ which fix objects. We do this under the assumption that the automorphism rescales polynomials according to their degree, but we prove
this assumption always holds for the root realization. Any such automorphism will act by an invertible scalar on the one-dimensional degree 1 hom space $\Hom^1(R,B_s)$, and similarly on
$\Hom^1(B_s,R)$. It was proven that these two invertible scalars (for each simple reflection $s$) determine the automorphism uniquely. In particular, if an automorphism of $\Diag$ fixes
the startdots and enddots of each color, then it is the identity functor! This is true for $\Theta$ by assumption, so every $\square_f$ commutes honestly.

This proof will extend from the root realization to a larger realization if and only if $\square_f$ commutes for any (linear) polynomial $f$. Since the polynomial ring of the standard
affine realization is generated over the simple roots by $x_1$, it is enough to check $\square_{x_1}$. Note that this proof will also work for some other realizations, such as any
extension of the root realization or the standard realization by $W$-invariant polynomials. \end{proof}

In \S\ref{sec-tensorBIcomps} we will confirm that \eqref{startdotsquare} and \eqref{enddotsquare} and $\square_{x_1}$ commute, thus proving Theorem \ref{thm:rotationfunctorial}. The reader can read \S\ref{sec-tensorBIcomps} now if desired.

\begin{remark} If we knew a priori that \eqref{startdotsquare} and \eqref{enddotsquare} commuted up to scalar, there would be a further simplification. Namely, we have the freedom of
choosing the scalar on the maps $\rotisom_s$. We could choose each $\rotisom_s$ such that \eqref{startdotsquare} commuted with scalar $1$. We know that $\square_{\al_s}$ commutes with
scalar $1$ (this requires no choices, since $\rotisom_{\mt}$ is the identity map). Thus \eqref{enddotsquare} also must commute with scalar $1$. However, the further content of Theorem
\ref{thm:rotationfunctorial} is that the map $\rotisom_s$ given in \S\ref{subsec-rotisoms} is the correct rescaling which makes \eqref{startdotsquare} commute. \end{remark}

\section{Preliminaries on descent sets}
\label{sec-descentsets}
%%%%%%%%%%%%%%%%%%%%%%%%%

\begin{notation} We write $\LC(w)$ for the left descent set of an element $w$ in a Coxeter group, and $\RC(w)$ for the right descent set. \end{notation}

In this section we prove some results about the elements of the form $w_I h_X$. In particular, Lemmas \ref{lem:rewriteaffine} and \ref{lem:rewriteaffine2} give a fairly useful way to
rewrite the expression $w_I h_X$. In the course of our analysis we will set up notation and introduce some ideas which will play a role in the more difficult computations to follow.

%===========
\subsection{Three settings}
\label{subsec-settings}
%===========

We will very frequently be using the (finite) symmetric group as a testing ground for ideas involving the affine Weyl group.

\begin{notation} This is the \emph{(special) finite setting}. Let $W_T = S_{N+1}$ be the symmetric group with simple reflections indexed by the ordered set $T = \{1 < 2 < \ldots < N\}$.
For a subset $X \subset T$, let $h_X = s_{x_1} s_{x_2} \cdots s_{x_d}$ where $X = \{x_1 > x_2 > \ldots > x_d\}$. When $X$ is fixed, then for $i \in T$ we set $\epsilon_i=1$ if $i \in X$
and $\epsilon_i = 0$ if $i \notin X$. We set $I = \{1, 2, \ldots, N-1\} \subset T$, and $J = \{2, \ldots, N-1, N\}$. If $i \in I$ is a simple reflection, then $\tau(i)$ denotes the simple reflection $i+1 \in J$. \end{notation}

\begin{notation} This is the \emph{(general) finite setting}. The notation is the same as above, except that $I$ is now some subset of $\{1, \ldots, N-1\}$, and $J = \tau(I)$.
\end{notation}

\begin{notation} This is the \emph{affine setting}. Let $(W,S)$ denote the affine Weyl group in type $\tilde{A}_n$, with rotation automorphism $\tau$. Recall from Definition \ref{defn:hX}
that for $X \subsetneq S$ the element $h_X$ is the product of the simple reflections in $X$ in decreasing order, with respect to the cyclic order on $X$ induced by some element $\notme
\notin X$. When $X$ is fixed, then for $i \in S$ we set $\epsilon_i=1$ if $i \in X$ and $\epsilon_i = 0$ if $i \notin X$. We fix some $I \subsetneq S$, and let $J = \tau(I)$.
\end{notation}

In the discussions below we will prove a fact first for the special finite setting, and then use this as a lemma to prove the result for the affine setting. Usually the process of going
from the special finite setting to the affine setting introduces two additional wrinkles: the first wrinkle arises when (the Dynkin diagram of) $I$ has multiple connected components, and
the second wrinkle arises when $I \cup J \cup X = S$. This first wrinkle can be isolated if one considers the general finite setting as an intermediate case; we chose not to do this to
eliminate some repetition. The reader is well-served by attempting the general finite setting for themselves.

%===========
\subsection{Descent sets: special finite setting}
\label{subsec-descentsfinite}
%===========

\begin{lemma} \label{lem:rightdescentfinite} We work in the special finite setting. For any subset $X \subset T$, the following are equivalent: \begin{itemize}
\item $J \subset X$.
\item $J \subset \RC(w_I h_X)$. \end{itemize} Moreover, in this case, one has: \begin{enumerate} \item $\ell(w_I h_X) = \ell(w_I) + \ell(h_X)$. \item $w_I h_X = h_Y w_J$ for some uniquely determined $Y \subset T$. \item $Y$ is the conjugation of $X$ by the longest element of $W_T$. \item $\RC(w_I h_X) = X$. \item $\LC(h_Y w_J) = Y$. \end{enumerate} \end{lemma}

\begin{proof} Let us write $h_X$ as \begin{equation} h_X = s_N^{\epsilon_N} \cdots s_2^{\epsilon_2} s_1^{\epsilon_1}. \end{equation} Consider the string diagram for the permutation $w_I
h_X$. Our favorite strand is the strand which is rightmost on top. This strand is fixed by $w_I$, so it is also rightmost below $w_I$ and above $h_X$. If some $\epsilon_i = 0$ for $i \in J$ then this strand crosses fewer than $N-1$ other strands in $h_X$, and will not be the first strand on bottom. But if $J \subset \RC(w_I h_X)$ then the last $N$ strands on the bottom must all cross each other, and our favorite strand is one of these. This is impossible, since our favorite strand does not cross $N-1$ other strands.

Conversely, if $J \subset X$ there are two cases: either $X = J$ or $X = T$. We observe that
\begin{equation} w_I h_J = h_I w_J \end{equation} and
\begin{equation} w_I h_T = w_T = h_T w_J, \end{equation}
and the lengths add in these expressions. In both cases everything is as desired. \end{proof}

%===========
\subsection{Statement in the affine setting}
\label{subsec-descentsstatement}
%===========

Our goal is to prove the affine analog of Lemma \ref{lem:rightdescentfinite}. 

\begin{lemma} \label{lem:descent} We work in the affine setting. For any proper subset $X \subsetneq S$, the following are equivalent: \begin{itemize}
\item $J \subset X$.
\item $J \subset \RC(w_I h_X)$. \end{itemize}  Moreover, in this case, one has: \begin{enumerate} \item $\ell(w_I h_X) = \ell(w_I) + \ell(h_X)$. \item $w_I h_X = h_Y w_J$ for some uniquely determined $Y \subsetneq S$. \item $Y$ is obtained from $X$ by conjugation by the longest element within each $\tau$-component, see below. \item $\RC(w_I h_X) = J \cup \{x \in X \mid x-1 \notin X\}$. \item $\LC(h_Y w_J) = I \cup \{y \in Y \mid y+1 \notin Y\}$. \end{enumerate} \end{lemma}

We will prove this lemma in \S\ref{subsec-descentaffine}, after introducing notation and several other ideas.

%===========
\subsection{Uniqueness}
\label{subsec-uniqueness}
%===========

Here is an easier fact about the descent set of $w_I h_X$.

\begin{lemma} Both $\LC(w_I h_X)$ and $\RC(w_I h_X)$ are subsets of $X \cup I$. \end{lemma}

\begin{proof} The only simple reflections appearing in $h_X$ or $w_I$ are within $X \cup I$. \end{proof}

Assuming Lemma \ref{lem:descent} we have a uniqueness result about elements of the form $w_I h_X$.

\begin{lemma} \label{lem:nooverlap} Let $I, I' \subsetneq S$. Let $X, X' \subsetneq S$ with $\tau(I) \subset X$ and $\tau(I') \subset X'$. If $w_I h_X = w_{I'} h_{X'}$ then $I = I'$ and $X = X'$. \end{lemma}

\begin{proof} 
Clearly if $I = I'$ then dividing by $w_I$ we get that $X = X'$. By symmetry, it is enough to prove that $I \subset I'$.

Suppose to the contrary that $i \in I, i \notin I'$. Then $i+1 \in \tau(I)$ so $i+1 \in \RC(w_I h_X) = \RC(w_{I'} h_{X'}) = \tau(I') \cup \{x \in X' \mid x-1 \notin X'\}$. But $i+1 \notin
\tau(I')$ so $i+1 \in X'$ and $i \notin X'$. Meanwhile, $i \in \LC(w_I h_X) = \LC(w_{I'} h_{X'}) \subset I' \cup X'$. But $i \notin I'$ and $i \notin X'$, a contradiction. \end{proof}

%===========
\subsection{$\tau$-components}
\label{subsec-taucomponents}
%===========

\begin{defn} \label{defn:taucomponent} We work in the affine setting. A \emph{$\tau$-component of $S$} is a subset $A$ such that whenever $j \in A$ then
\begin{itemize} \item $j+1 \in A$ if and only if $j \in I$, and \item $j-1 \in A$ if and only if $j \in J$. \end{itemize} Note that $j \in I$ if and only if $j+1 \in J = \tau(I)$. \end{defn}

Another way to state this definition is as follows. For each connected component $B \subset I$, then $\tau(B)$ is a connected component of $J$, and $A = B \cup \tau(B)$ is a
$\tau$-component. We have $B = A \cap I$ and $\tau(B) = A \cap J$. Meanwhile, each element of $S \setminus (I \cup J)$ forms its own $\tau$-component, which is a singleton, and any singleton $\tau$-component lives in $S \setminus (I \cup J)$. The following lemma is easy.

\begin{lemma} The set $S$ is a disjoint union of its $\tau$-components. \end{lemma}

\begin{example} \label{ex:main} Suppose that $I = \{1, 2, 3, 5, 9, 10\}$ and $J = \{2,3,4, 6, 10, 11\}$ and $n = 12$. Then the $\tau$-components are \[ A_0 = \{1,2,3,4\}, \quad A_1 =
\{5,6\}, \quad A_2 = \{7\}, \quad A_3 = \{8\}, \quad A_4 = \{9,10,11\}, \quad A_5 = \{0\}.\] \end{example}

\begin{defn} \label{defn:taucomporder} We will place an order on each $\tau$-component $A$. If $A \subsetneq S$ then we use the natural cyclic order on $A$. In the special case where $|I| = n-1$ so that there is a unique $\tau$-component $A$ and $A=S$, let $i$ be the unique element of $S \setminus I$, and use the order $\{i+1 < i+2 < \ldots < i-1 < i\}$ on $A$. \end{defn}

\begin{lemma} For this order on a $\tau$-component $A$, one has $I \cap A = A \setminus \max(A)$ and $J \cap A = A \setminus \min(A)$. \end{lemma}

\begin{proof} This is obvious by construction. \end{proof}

Recall that to define the element $h_X$ we need to choose an element $\notme \notin X$, and use this to define a cyclic order on $X$. When $\notme$ is fixed, we will label our
$\tau$-components in a particular way.

\begin{notation} We work in the affine setting. Fix some $\notme \in S$. Then let $A_0$ be the $\tau$-component containing $\notme$, and label the remaining $\tau$-components $A_1, A_2,
\ldots, A_d$ in increasing order. Moreover, split $A_0 \setminus \notme$ into two pieces, $A_{> \notme}$ and $A_{< \notme}$, using the order on $A_0$ from Definition \ref{defn:taucomporder}. \end{notation}

\begin{example} Continue Example \ref{ex:main}. Our labels are consistent so long as $\notme \in \{1, 2, 3, 4\}$. If $\notme = 2$ then $A_{< \notme} = \{1\}$ and $A_{> \notme} =
\{3,4\}$. If $\notme = 1$ then $A_{< \notme} = \mt$ and $A_{> \notme} = \{2,3,4\}$. \end{example}

\begin{remark} Suppose that $X \cup I \cup J \ne S$. We may as well choose $\notme \in S \setminus (X \cup I \cup J)$, meaning that $A_0 = \{\notme\}$. Now $A_{>\notme} = \mt$, $A_{<
\notme} = \mt$, $I \cap A_0 = \mt$. This situation reduces to the general finite setting, for $T = S \setminus \{\notme\}$, and the reader who wants practice should always attempt this
case first. \end{remark}

%===========
\subsection{Rewriting $w_I h_X$}
\label{subsec-rewritingwIhX}
%===========

Our next step is to rewrite $w_I h_X$ in a more useful form.

If $\notme \notin X \subsetneq S$, then we can use the cyclic order on $X$ induced from $\notme$. Decomposing $X$ based on the $\tau$-components of $S$, we have
\begin{equation} h_X = h_{X \cap A_{< \notme}} h_{X \cap A_d} \cdots h_{X \cap A_1} h_{X \cap A_{> \notme}}. \end{equation}
Meanwhile, $I$ is union of its connected components $I \cap A_j$ (some of which are empty). Thus
\begin{equation} \label{eq:decomposewI} w_I = w_{I \cap A_0} w_{I \cap A_d} \cdots w_{I \cap A_1} \end{equation}
and all of these terms commute. Moreover, any element of $I \cap A_j = A_j \setminus \max(A_j)$ commutes with any element of $A_{j'}$ unless $j' = j$ or $j' = j - 1$. Finally, any element of $I \cap A_1$ commutes with $A_{< \notme}$. Thus in the product $w_I h_X$,  the factor $w_{I \cap A_1}$ commutes past everything in $h_X$ up until $h_{X \cap A_1}$, the factor $w_{I \cap A_2}$ commutes past everything up until $h_{X \cap A_2}$, and so forth. From this we have proven the following lemma.

\begin{lemma} \label{lem:rewriteaffine} We work in the affine setting, and fix $\notme \in S$. For any $X \subsetneq S$ with $\notme
\notin X$, we can apply only commutation relations to obtain the equality
\begin{equation} \label{eq:rewriteaffine} w_I h_X = (w_{I \cap A_0} h_{X \cap A_{< \notme}}) (w_{I \cap A_d} h_{X \cap A_d}) \cdots (w_{I \cap A_1} h_{X \cap A_1}) h_{X \cap A_{> \notme}}. \end{equation} \end{lemma}

\begin{example} \label{ex:exemplarybad} Continue Example \ref{ex:main}. Suppose that $X = \{0, 1, 3, 4, 6, 7, 8, 10, 11\}$ and $\notme = 2$. Then
\begin{equation} w_I h_X = w_{123} w_{5} w_{9\; 10} s_1 s_0 s_{11} s_{10} s_8 s_7 s_6 s_4 s_3 = (w_{123} s_1) (s_0) (w_{9\; 10} s_{11} s_{10}) (s_8) (s_7) (w_5 s_6) (s_4 s_3). \end{equation}
\end{example}

\begin{remark} When we say that the equality \eqref{eq:rewriteaffine} only uses commutation relations, we are also tacitly permitting ourselves to use any desired reduced expression for $w_I$, in this case a reduced expression compatible with \eqref{eq:decomposewI}. That we only use commutation relations will become relevant in the categorification, where using arbitrary reduced expressions for $w_I$ is also permissible, see the next chapter. \end{remark}
	
The expression \eqref{eq:rewriteaffine} simplifies when $X \cap A_{< \notme}$ is empty. If $X \cap A_{< \notme}$ is empty then $\min(A_0) \notin X$, so we may as well assume instead that $\notme = \min(A_0)$.

\begin{defn} For a $\tau$-component $A$, let the \emph{interior} of $A$ be $\mathring{A} = A \cap I \cap J$. Said another way, $\mathring{A} = A \setminus \{\max(A), \min(A)\}$. \end{defn}
	
Crucially, any simple reflection in $\mathring{A}$ will commute with any simple reflection in any other $\tau$-component.

\begin{lemma} \label{lem:rewriteaffine2} We work in the affine setting, and set $\notme = \min(A_0)$. For any $X \subsetneq S$ with $\notme \notin X$, we can apply only commutation relations to obtain the equality
\begin{equation} \label{eq:rewriteaffine2} w_I h_X = h_{I \cap A_0} (w_{I \cap A_d} h_{X \cap A_d}) \cdots (w_{I \cap A_1} h_{X \cap A_1}) (w_{\mathring{A_0}} h_{X \cap A_0}). \end{equation}
\end{lemma}

\begin{proof} We rewrite $w_{I \cap A_0} = h_{I \cap A_0} w_{\mathring{A_0}}$, and then use \eqref{eq:rewriteaffine}, observing that $w_{\mathring{A_0}}$ commutes past all the terms until $h_{X \cap A_{> \notme}} = h_{X \cap A_0}$. \end{proof}

%===========
\subsection{Reduction to components}
\label{subsec-reductiontocomponents}
%===========

The rest of this chapter focuses on proving Lemma \ref{lem:descent}, and the reader can skip it if desired.

\begin{lemma} \label{lem:reducetocomponent} We work in the affine setting. For any $X \subsetneq S$ and any $\tau$-component $A$, we have $J \cap A \subset X$ if and only if $J \cap A \subset \RC(w_{I \cap A} h_{X \cap A})$ if and only if $J \cap A \subset \RC(w_I h_X)$. \end{lemma}
	
\begin{proof} The fact that $J \cap A \subset X$ if and only if $J \cap A \subset \RC(w_{I \cap A} h_{X \cap A})$ follows immediately from Lemma \ref{lem:rightdescentfinite} applied to $T = A$, except in the special case where $A = S$. We will treat this special case below.

Suppose that we can choose some $\notme \notin X$ with $\notme \notin A$. Then using $\notme$ to number the $\tau$-components, $A = A_i$ for some $i \ne 0$. Using \eqref{eq:rewriteaffine} we can write
\begin{equation} w_I h_X = g f g' \end{equation}
where $f = w_{I \cap A} h_{X \cap A}$ is in the parabolic subgroup generated by $A$, $g$ is in the parabolic subgroup generated by $A_0$ and $A_j$ for $j > i$, and $g'$ is in the parabolic subgroup generated by $A_{> \notme}$ and $A_j$ for $j < i$. In particular, the simple reflections in $J \cap A$ commutes with the simple reflections in $g'$. A general result in Coxeter theory (see the next lemma) will imply that, for any simple reflection $s \in J \cap A$, we have $s \in \RC(f)$ if and only if $s \in \RC(g f g')$, which gives us the desired result.

Suppose that we are forced to label $A = A_0$, that is, $S \setminus X \subset A$. Now we prove that $J \cap A \subset X$ if and only if $J \cap A \subset \RC(w_I h_X)$ by a different
and subtler argument, which does not look at $J$ one simple reflection at a time, but all at once. If we can do this then we also handle the exceptional case $A=S$ from above.

Suppose first that $J \cap A_0 \subset X$. Then we must have $\notme = \min(A_0)$ and $X \cap A_0 = J \cap A_0$. Let $f = h_{I \cap A_0}$ and $f' = w_{\mathring{A}} h_{J cap A_0} = w_{J \cap A_0}$. Using \eqref{eq:rewriteaffine2} we see that 
\begin{equation} w_I h_X = f g f' \end{equation} where $f \in W_{I \cap A_0}$, $f' \in W_{J \cap A_0}$, and $g$ is in the parabolic subgroup complementary to $A_0$. In Lemma \ref{lem:rightdescentfinite} we have already proven that $J \cap A_0 \subset \RC(f')$. That is containment is unchanged by left multipication by $g$ and then $f'$ is not surprising. After all, $g$ is in a disjoint parabolic subgroup. Meanwhile, $f$ satisfies the property that $\ell(f w) = \ell(f) + \ell(w)$ for all $w \in W_{J \cap A_0}$ which one can see by counting strand crossings (this argument recurs in the next section, after further related exposition, so don't worry if that didn't make sense for you now). The fact that $J \cap A_0 \subset \RC(f')$ implies $J \cap A_0 \subset \RC(fgf')$ is another consequence of general results in Coxeter theory (see the next lemma).

Now suppose that $J \cap A_0 \nsubseteq X$. We wish to deduce that $J \cap A_0 \nsubseteq \RC(w_I h_X)$. A typical example is Example \ref{ex:exemplarybad}. When we wanted to prove the analogous statement in the finite setting in Lemma \ref{lem:rightdescentfinite}, we made an argument counting strand crossings. In that proof, we argued that
if $J \nsubseteq X$ then a particular strand did not cross enough other strands for $J$ to be in the descent set of $w_I h_X$. The Coxeter-theoretic version of this argument examines the
positive roots sent to negative roots, and we now sketch a similar argument in the affine setting.

Consider the cylindrical crossing diagram of $w_I h_X$, and our favorite strand labeled $\notme$, where $\notme$ is assumed to live in $J \cap A_0$. (The crossing labeled $i$ crosses the
strands labeled $i$ and $i+1$.) Using \eqref{eq:rewriteaffine} we write $w_I h_X = fg$ where $f = (w_{I \cap A_0} h_{X \cap A_{< \notme}})$ and $g$ is the rest. Note that $g$ fixes the
strand $\notme$. Since $f$ lives in the parabolic subgroup $I \cap A_0$, isomorphic to a symmetric group on $|I \cap A_0|+1$ letters, our favorite strand has the potential to cross at
most $|I \cap A_0|$ other strands. In order for $J \cap A_0$ to be in the right descent set, it must cross $|J \cap A_0|$ particular strands, including the strand $\notme+1$. Since $|I
\cap A_0| = |J \cap A_0|$, we know exactly which strands it must cross. However, in the affine setting, the positive root sent to a negative root by an element of $W$ is not just
determined by which strands cross but by how many times they wind around the cylinder before they cross. In order for our favorite strand $\notme$ to cross $\notme+1$ within the final
term $f$, $\notme+1$ must have been pulled by $g$ all the way around the cylinder to become part of the block of strands permutated by $I \cap A_0$. Thus even if our favorite strand
crosses $\notme+1$, it does not correspond to the correct positive root, and $J \cap A_0$ will not be in the right descent set. For the reader unfamiliar with the affine root system, it
is a good exercise to make this argument more rigorous.
\end{proof}

In the above proof we made use of two general facts about disjoint parabolic subgroups of Coxeter groups.

\begin{lemma} \label{lem:descentsetlemma} For any Coxeter system, partition the simple reflections into two disjoint subsets $B$ and $C$. Let $f, f' \in W_B$ and $g, g' \in W_C$. Then
$\ell(fg) = \ell(gf) = \ell(f) + \ell(g)$. Let $s \in B$ be a simple reflection commuting with $C$, and $t \in B$ be an arbitrary simple reflection. \begin{enumerate} \item $s \in \RC(f)$ if and only if $s \in \RC(gfg')$. \item If $t \notin \RC(f)$ then $t \notin \RC(fg)$. \item $t \in \RC(f)$ if and only if
$t \in \RC(gf)$. \item If $t \notin \RC(f)$ then $t \notin \RC(gfg')$. \item Suppose that $\ell(f f') = \ell(f) + \ell(f')$. Then $t \in \RC(f')$ implies that $t \in \RC(fgf')$. \end{enumerate} \end{lemma}

\begin{proof} These results are easiest to prove using basic facts about positive and negative roots in the geometric representation, which has a basis given by simple roots, see
\cite[Chapter 5]{HumpCox}. For a simple reflection $t$ we know that $wt < w$ if and only if $w(\alpha_t)$ is a negative root. Let the \emph{$B$-span} be the span of $\{\alpha_u \mid
u \in B\}$, and a \emph{$B$-root} be a root in the $B$-span. The action of $W_B$ on any root can only change the coefficients in the $B$-span. Thus if $W_B$ sends a positive root to a
negative root, then that positive root is a $B$-root.
	
Let $\un{w}$ be an expression for an element $w$. The mechanics of the deletion condition states that, if $wt < w$, then there is a rightmost simple reflection in $\un{w}$
whose removal gives an expression for $wt$. Moreover, this simple reflection can be determined by looking at where the expression sends the positive root $\alpha_t$: acting from right to
left one reflection at a time, some simple reflection first produces a negative root, and this is the one to be removed.

Let us prove (1). The element $g'$ (or any subexpression thereof) will fix $\alpha_s$, since $s$ commutes with $W_C$. Since $g'(\alpha_s) = \alpha_s$, the process of determining which simple reflection first sends $\alpha_s$ to a negative root (if any) is the same for $gf$ and for $gfg'$. Thus $s \in \RC(gfg')$ if and only if $s \in \RC(gf)$, which reduces the problem to part (3).

Let us prove (3). Clearly $f$ sends $\alpha_t$ to a $B$-root. Whether this is positive or negative is unchanged by $g \in W_C$. Thus whether $\alpha_t$ is sent to a positive or a negative root is determined by $f(\alpha_t)$, whence the result.

Let us prove (5). We know that $f'$ sends $\alpha_t$ to a negative $B$-root. Then $g$ sends this root to another negative root, so the coefficients of this root in the $C$-span are negative, or zero if $gf'(\alpha_t) = f'(\alpha_t)$. In the former case, $f$ will not change the coefficients in the $C$-span, so $fgf'(\alpha_t)$ is still negative. In the latter case, the result follows from the assumption that $\ell(f f') = \ell(f) + \ell(f')$.

Let us prove (2). If $g(\alpha_t)$ has any nonzero coefficients in the $C$-span then they must be positive, and $gf(\alpha_t)$ must be positive too. If $g(\alpha_t) = \alpha_t$ then $fg(\alpha_t)$ is positive because $t \notin \RC(f)$. Thus $t \notin \RC(fg)$.

Applying (2) and (3) immediately yields (4). \end{proof}

%===========
\subsection{Proof of the descent set lemma}
\label{subsec-descentaffine}
%===========

Let us recall the statement of the lemma.

\begin{lemma*} We work in the affine setting. For any proper subset $X \subsetneq S$, the following are equivalent: \begin{itemize}
\item $J \subset X$.
\item $J \subset \RC(w_I h_X)$. \end{itemize}  Moreover, in this case, one has: \begin{enumerate} \item $\ell(w_I h_X) = \ell(w_I) + \ell(h_X)$. \item $w_I h_X = h_Y w_J$ for some uniquely determined $Y \subsetneq S$. \item $Y$ is obtained from $X$ by conjugation by the longest element within each $\tau$-component. \item $\RC(w_I h_X) = J \cup \{x \in X \mid x-1 \notin X\}$. \item $\LC(h_Y w_J) = I \cup \{y \in Y \mid y+1 \notin Y\}$. \end{enumerate} \end{lemma*}

Let us restate condition (3). Consider a $\tau$-component $A$. If $J \subset X$ then there are two cases: either $X \cap A = J \cap A$ and $Y \cap A = I \cap A$, or $X \cap A = A$ and $Y
\cap A = A$. This is true as well when $A = \{i\}$ is a singleton: then $i \in X$ if and only if $i \in Y$.

\begin{example} Continue Example \ref{ex:main}. If $J \subset X$ then $X$ must contain $\{2, 3, 4\}$ and $Y$ must contain $\{1, 2, 3\}$. Also, $X$ may or may not contain $1$, and $1 \in
X$ if and only if $4 \in Y$. Similarly, $X$ may or may not contain $7$, and $7 \in X$ if and only if $7 \in Y$. \end{example}

\begin{proof} The equivalence between $J \subset X$ and $J \subset \RC(w_I h_X)$ follows from Lemma \ref{lem:reducetocomponent}, by intersecting each set with each $\tau$-component $A$. Let us now assume that $J \subset X$, and prove the remaining statements. Note that $\notme = \min(A_0)$ by necessity, so that we may use Lemma \ref{lem:rewriteaffine2}. 

To show that the lengths add in the expression $w_I h_X$, rewrite this expression using \eqref{eq:rewriteaffine2}. With the exception of the initial factor $h_{I \cap A_0}$, the remainder
is a product of terms $w_{I \cap A} h_{X \cap A}$ for various disjoint parabolic subgroups $A$. The lengths add in each of these individual expressions by Lemma
\ref{lem:rightdescentfinite}, and the lengths add between them since they are in disjoint parabolic subgroups. Finally, the initial factor $h_{I \cap A_0}$ only adds crossings involving
the strand $\notme$, which has not let crossed anything, so it can only increase the length.

Now let us prove the statement about right descent sets. We know that $J \subset \RC(w_I h_X)$. If $x \in X$ and $x-1 \notin X$ then $x \in
\RC(h_X)$ so $x \in \RC(w_I h_X)$ since the lengths add in this expression. Thus $J \cup \{x \in X \mid x-1 \notin X\} \subset \RC(w_I h_X)$.

To show the reverse containment we must rule out several cases. Let $x \in S$ with $x \notin J$. \begin{itemize} \item Suppose that $x \notin X$ and $x \notin I$. Then the simple
reflection $x$ does not appear in the product $w_I h_X$, so it can not be in the right descent set. \item Suppose that $x \notin X$ and $x \in I$. In particular, $x = \min(A)$ for some
$\tau$-component $A$, and the only place where $x$ appears in the product $w_I h_X$ is inside $w_{I \cap A}$. From Lemma \ref{lem:rightdescentfinite} we know that $x \notin \RC(w_{I \cap
A} h_{X \cap A})$. Now we apply Lemma \ref{lem:descentsetlemma}(4) to deduce that $x \notin \RC(w_I h_X)$. \item Suppose that $x \in X$ and $x-1 \in X$. Several cases are possible: $x \in I$
or $x \in S \setminus (I \cup J)$, $x-1 \in J$ or $x-1 \in S \setminus(I \cup J)$. Regardless, $x$ is not in the right descent set of some appropriately small product (like $s_x s_{x-1}$
when both are in $S \setminus (I \cup J)$, or $w_{I \cap A} h_{X \cap A} s_{x-1}$ when $x \in I \cap A$, etcetera). Applying Lemma \ref{lem:descentsetlemma}(4) again, we see that $x \notin \RC(w_I
h_X)$. \end{itemize}

Assuming that $J \subset X$ we can use \eqref{eq:rewriteaffine2} to rewrite $w_I h_X$. Assuming that $I \subset Y$, we can use an entirely analogous argument to rewrite $h_Y w_J$. If $Y$
and $X$ are related as in the Lemma, then Lemma \ref{lem:rightdescentfinite} states that $w_{I \cap A} h_{X \cap A} = h_{Y \cap A} w_{J \cap A}$ for any $\tau$-component $A$. (One must
also handle the edge case where $A = S$, but this is easily dealt with.) Applying this to each term gives the desired equality $w_I h_X = h_Y w_J$. We elaborate upon this in \S\ref{sec-leftvsright}.

All the arguments used above for the right descent set of $w_I h_X$ have analogues for the left descent set of $h_Y w_J$. \end{proof}
%%%%%%%%%%%%%%%%%%%%%%%%%
\section{Tensoring with the longest element: a thicker calculus}
\label{sec-thicker}
%%%%%%%%%%%%%%%%%%%%%%%%%

En route to proving Theorem \ref{thm:tensorBI}, we will need to study the indecomposable objects $B_{w_I h_X}$ in great detail. In this chapter, we compute the idempotent which picks out
the top summand $B_{w_I h_X}$ inside the tensor product $B_I B_X$ when $\ell(w_I) + \ell(h_X) = \ell(w_I h_X)$. We provide some diagrammatic tools for working with these idempotents.

The motivation for this chapter is really the Gaussian elimination of the next chapter. If motivation is required, one can read \S\ref{subsec-GEobjects} now.

%===========
\subsection{Recollection of the thick calculus}
\label{subsec-thickrecall}
%===========

Let $I$ be a parabolic subgroup inside a symmetric group. The \emph{thicker calculus} developed in \cite{EThick} gives a diagrammatic presentation of the category which is monoidal
generated not just by $B_s$, but also by $B_I$ for all $I$.  Let us recall this thicker calculus.
 
We will use a \emph{thick strand} to denote the object $B_I$. For any reduced expression $\un{w}$ for $w_I$, $B_I$ is a direct summand of the Bott-Samelson object $\BS(\un{w})$.
We draw the (degree zero) inclusion and projection maps as follows, and call them \emph{thick splitters and mergers}.
\begin{equation} \ig{1}{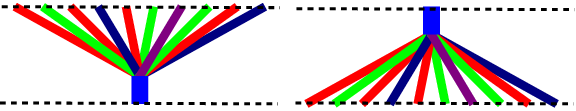} \end{equation}
(The inclusion and projection maps are only defined up to scalar, but there is a preferred normalization which preserves the $1$-tensor, see \cite{EThick}.) They compose in one direction to the identity of $B_I$, and in the other direction to an idempotent $e_I$ in $\BS(\un{w})$ which is described explicitly in \cite{EThick}; the details will not be relevant.

When $I' \subset I$ one can choose a reduced expression for $w_I$ that ends with a reduced expression for $w_{I'}$. By applying a thick splitter from $I$ and then a thick merger to $I'$ on the subexpression, we can effectively split off a copy of $I'$ and leave behind a reduced expression for $w_I w_{I'}\inv$.
\begin{equation} \label{eq:splitoffI} \ig{1}{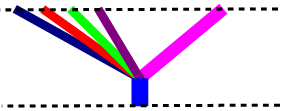} \end{equation}

For any $s \in I$ there is a \emph{thick trivalent vertex} of degree $-1$, defined by choosing any reduced expression $\un{w}$ for $w_I$ ending in $s$ (the choice of reduced expression will end up being irrelevant).
\begin{equation} \label{thicktridefn}
{
\labellist
\small\hair 2pt
 \pinlabel {$=$} [ ] at 95 30
\endlabellist
\centering
\ig{1}{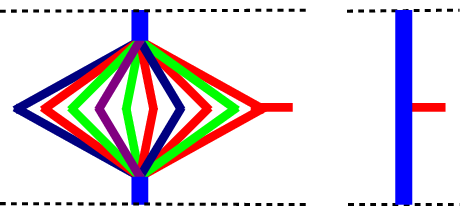}
}
\end{equation}
Pulling the $s$-colored strand up to down yields degree $-1$ maps $B_I B_s \to B_I$ and $B_I \to B_I B_s$, analogous to the trivalent vertices $\splitred$ and $\mergered$ in the case $I = \{s\}$. These morphisms are related by the self-adjunction of $B_s$, so that we may consider diagrams up to isotopy.

\begin{remark} This is not exactly the definition of the thick trivalent vertex given in \cite[(4.10)]{EThick}, which uses a fixed reduced expression not necessarily ending in $s$, and a more complicated morphism $a_s$ in the middle. However, it is equivalent to this one using \cite[(4.2) and (4.3)]{EThick}. \end{remark}

These new trivalent vertex satisfies an analog of the unit and associativity axiom, just as ordinary trivalent vertices do in the case $I = \{s\}$. The unit axiom is
\begin{equation} \label{eq:thickunit} {
\labellist
\small\hair 2pt
 \pinlabel {$=$} [ ] at 38 18
\endlabellist
\centering
\ig{1}{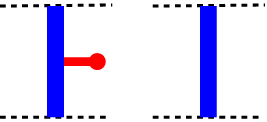}
} \end{equation}
and the associativity axiom is
\begin{equation} \label{eq:thickassoc} {
\labellist
\small\hair 2pt
 \pinlabel {$=$} [ ] at 36 18
\endlabellist
\centering
\ig{1}{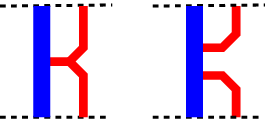}
}. \end{equation}
A consequence is that ``an empty eyehole is zero,''
\begin{equation} \label{eq:thickneedle} \ig{1}{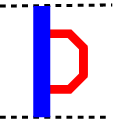} = 0. \end{equation}
One also can absorb $2m_{st}$-valent vertices into the thick strand, for $s,t \in I$.
\begin{equation} \label{eq:thick2m} {
\labellist
\small\hair 2pt
 \pinlabel {$=$} [ ] at 36 18
\endlabellist
\centering
\ig{1}{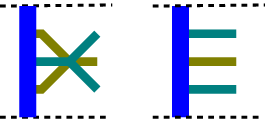}
} \end{equation}

The idempotent decomposition $B_I B_s \cong B_I(1) \oplus B_I(-1)$ is given by the following factored idempotent decomposition.
\begin{equation} \label{eq:thicksdecomp} {
\labellist
\small\hair 2pt
 \pinlabel {$=$} [ ] at 35 18
 \pinlabel {$+$} [ ] at 82 18
 \pinlabel {$\frac{\al_s}{2}$} [ ] at 62 6
 \pinlabel {$\frac{\al_s}{2}$} [ ] at 106 28
\endlabellist
\centering
\ig{2}{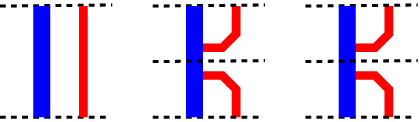}
} \end{equation}
The first idempotent $e_+$ projects to $B_I(1)$, and the second $e_-$ projects to $B_I(-1)$.
	
If $I$ and $I'$ are two parabolic subgroups which commute with each other, then $B_I B_{I'} \cong B_{I'} B_I$, and the isomorphism is drawn as a new kind of 4-valent vertex, whose definition is given explicitly in \cite{EThick}.
\begin{equation} \label{eq:thicksupercrossing} \ig{1}{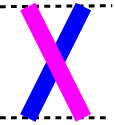} \end{equation}
It satisfies properties analogous to the ordinary 4-valent vertex, which we will not bother to explicitly recall.

%===========
\subsection{Smoothness and plethysm: finite case}
\label{subsec-prelimsmoothBI}
%===========

Recall that smoothness of an element $w$ in a Coxeter group $W$ is the statement that all Kazhdan-Lusztig polynomials are trivial. If $\{H_y\}_{y \in W}$ is the standard basis of $\Hecke(W)$, let $\Sigma_w = \sum_{y \le w} v^{\ell(w) - \ell(y)} H_y$. Then $w$ is \emph{smooth} if $\Sigma_w$ is equal to the Kazhdan-Lusztig basis element $b_w$. An equivalent condition is that $\Sigma_w$ is self-dual under the bar involution. Smoothness is the exception, but there are many nice families of elements which are smooth.

\begin{lemma} \label{lem:smoothIJ} Let $(W,S)$ be a Coxeter system, and $I, J$ two disjoint subsets of $S$. If $w \in W_I$ and $x \in W_J$ are both smooth, then $wx$ is smooth, and $b_{wx} = b_w b_x$. \end{lemma}

\begin{proof} For any $y \in W_I$ and $z \in W_J$, multiplying the standard basis elements $H_y$ and $H_z$ in the Hecke algebra will yield $H_{yz}$, and $\ell(yz) = \ell(y) + \ell(z)$.
Thus $\Sigma_w \Sigma_x = \Sigma_{wx}$. Moreover, $\Sigma_w$ and $\Sigma_x$ are self-dual, and thus so is their product. \end{proof}

\begin{notation} \label{not:suffixfinite} In the special finite setting, we have $I = \{1, \ldots, N-1\} \subset T = \{1, \ldots, N\}$. If $X$ has the form $\{m, \ldots,N-1, N\}$ for some $1 \le m \le N+1$, we call $X$ a \emph{suffix} of $T$. Note that $X$ has size $k = N+1-m$, and by convention $X = \mt$ when $m = N+1$. We denote $w_I h_X$ by $\who_k$, thus
\begin{equation} \who_k = w_I h_X = w_I s_N s_{N-1} \cdots s_m. \end{equation} \end{notation}
	
We sometimes index suffixes by $m$, which helps keep track of simple reflections, and sometimes by $k$, because it helps keep track of size.

We claim that $w_I h_X$ is smooth when $X$ is a suffix of $T$. This is a well-known fact, with many proofs. For example, it is proven in \cite[\S 6.7]{EHDiag2}, using the
technology of pattern avoidance. One can also prove it geometrically, observing that these elements correspond to Schubert varieties inside projective space. We sketch yet another proof
here, because are interested in the plethysm which is involved in the proof.

If $X$ is a suffix of $T$ (e.g. $X = \{m+1, \ldots, N\}$) then there is exactly one simple reflection (e.g. $s_m$) that one can add to $X$ to get another suffix of $T$. Meanwhile, every
other element of $T$ is in the right descent set $\RC(w_I h_X)$. Let us distinguish between two cases. If $X$ is empty so that $s_N$ is the unique element of $T$ not in $\RC(w_I h_X)$,
we say that $N$ \emph{initiates} the suffix. Meanwhile, if $X$ is a nonempty suffix of $T$, and $s_m$ is the unique element of $T$ not in $\RC(w_I h_X)$, we say that $m$ \emph{grows}
the suffix $X$. Whenever $m$ grows the suffix $X$, it is clear that $m+1 \in X$. Either way, we write $Xm$ for $X \cup \{m\}$ when $m \notin X$.

\begin{lemma} \label{lem:finitepleth} Let $X$ be a suffix of $T$. When $m$ grows the suffix $X$, we have \begin{subequations} \label{eq:wIhXsmfinite}
\begin{equation}\label{eq:wIhXsmfinitecase1} b_{w_I h_X} b_{s_m} = b_{w_I h_{Xm}} + b_{w_I h_{X \setminus m+1}} \end{equation}
or in other words
\begin{equation} b_{\who_k} b_{s_{N-k}} = b_{\who_{k+1}} + b_{\who_{k-1}}. \end{equation}
When $m = N$ initiates the suffix $X = \mt$, we have
\begin{equation} \label{eq:wIhXsmfinitecase2} b_{w_I} b_{s_N} = b_{w_I s_N} \end{equation}
or in other words
\begin{equation} b_{\who_0} b_{s_N} = b_{\who_1}. \end{equation}
In all other cases we have
\begin{equation} \label{eq:wIhXsmfinitecase3} b_{w_I h_X} b_{s_j} = (v+v\inv) b_{w_I h_X}. \end{equation} \end{subequations}
In particular, $w_I h_X$ is smooth whenever $X$ is a suffix. \end{lemma}

\begin{proof} (Sketch) One has the general fact that $b_w b_s = (v + v\inv) b_w$ whenever $ws < w$. The equation \eqref{eq:wIhXsmfinitecase3} holds because $s_j \in \RC(w_I h_X)$ in all
other cases.

The equation \eqref{eq:wIhXsmfinitecase2} holds by Lemma \ref{lem:smoothIJ}, and serves as the base case for our proof of \eqref{eq:wIhXsmfinitecase1}. We induct on the size of $X$. For
$w \in W_T$ set $W_{\le w} = \{y \in W_T \mid y \le w\}$, an ideal in the Bruhat order.

The action of the parabolic subgroup $\{1,s_m\}$ on the right will split $W_T$ into cosets of size $2$. Consider the set $W_{\le w_I h_{X \setminus m+1}}$. We claim that $y \in W_{\le
w_I h_{X \setminus m+1}}$ implies $y s_m \in W_{\le w_I h_{X \setminus m+1}}$, and thus $W_{\le w_I h_{X \setminus m+1}}$ is a union of its right cosets for $s_m$. On the other hand, we
claim that all the elements of $W_{\le w_I h_X} \setminus W_{\le w_I h_{X \setminus m+1}}$ are minimal in their right $s_m$ cosets. We leave these claims as exercises to the reader.

From this observation, one computes that \begin{equation} \Sigma_{w_I h_X} b_{s_m} = \Sigma_{w_I h_{Xm}} + \Sigma_{w_I h_{X \setminus m+1}}. \end{equation} If we assume inductively that $w_I h_X$ and $w_I h_{X \setminus m+1}$ are smooth, then corresponding $\Sigma$ elements are self-dual, and thus so is $\Sigma_{w_I h_{Xm}}$. Thus $w_I h_{Xm}$ is smooth. Replacing the $\Sigma$s above with $b$s, we get \eqref{eq:wIhXsmfinitecase1}. \end{proof}

%===========
\subsection{Idempotents: finite case}
\label{subsec-idempsfinite}
%===========

We continue to work in the special finite setting. Our goal is to explicitly find the idempotents in the direct sum decompositions which categorify \eqref{eq:wIhXsmfinite}.

For $1 \le m \le N+1$, with $k = N+1-m$, let $X$ be the suffix $\{m, \ldots, N\} \subset T$. We want to compute the idempotent $e_k$ projecting to the top summand $B_{\who_k}$ inside the tensor product $B_I B_X$. We will denote $e_k$ with an oval labeled by $k$.
\begin{equation} {
\labellist
\small\hair 2pt
 \pinlabel {$k$} [ ] at 26 17
\endlabellist
\centering
\ig{1}{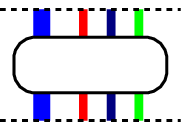}
} \end{equation}
Of course $e_0$ is just the identity map of $B_I$. By \eqref{eq:wIhXsmfinitecase2}, $B_I B_{s_N}$ is indecomposable, so $e_1$ is just the identity map of $B_I B_{s_N}$.

By \eqref{eq:wIhXsmfinitecase1} we know that
\begin{equation} B_{\who_k} B_{s_{N-k}} = B_{\who_{k+1}} + B_{\who_{k-1}}. \end{equation}
The next proposition computes $e_m$ recursively, based on this known decomposition.

\begin{remark} These plethysm rules may remind some readers of the tensor product decomposition rules for representations of $\sl_2$. From this one might
mistakenly guess that the Jones-Wenzl recursion should be used to construct these idempotents, as happens for dihedral groups, see \cite{ECathedral}. In fact, the recursion for idempotents is much simpler, and unrelated to Jones-Wenzl projectors. \end{remark}

\begin{prop} \label{prop:recursionetc} The equations \eqref{eq:emrecursion} through \eqref{eq:commutespastem} hold. \end{prop} 

The recursive formula for $e_{k+1}$ in terms of $e_k$ and $e_{k-1}$ is
\begin{equation} \label{eq:emrecursion} {
\labellist
\small\hair 2pt
 \pinlabel {$=$} [ ] at 65 49
 \pinlabel {$+$} [ ] at 147 46
 \pinlabel {$(-1)$} [ ] at 158 20
 \pinlabel {$k$} [ ] at 25 49
 \pinlabel {$k+1$} [ ] at 101 49
 \pinlabel {$k$} [ ] at 197 15
 \pinlabel {$k-1$} [ ] at 194 51
 \pinlabel {$k$} [ ] at 197 87
\endlabellist
\centering
\ig{1}{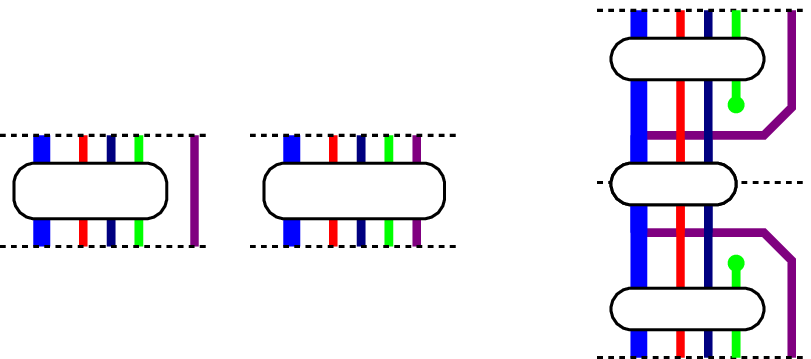}
} \end{equation}
when $k \ge 1$. The last (purple) strand is colored $s_{N-k}$, and the penultimate (green) strand is colored $s_{N-k+1}$. The formula \eqref{eq:emrecursion} expresses the identity of $B_{\who_{k}} B_{s_{N-k}}$ as a sum of two orthogonal idempotents, $\idemone$ which projects to $B_{\who_{k+1}}$, and $\idemtwo$ which projects to $B_{\who_{k-1}}$. There is a sign on $\idemtwo$, which we have placed lower down because we will always include it in the projection map. Note that the central idempotent (oval) in $\idemtwo$ is redundant and can be ignored, though we draw it to emphasize how the morphism factors.

This is supposed to define an idempotent.
\begin{equation} \label{eq:ekidemp} e_k^2 = e_k. \end{equation}

The second idempotent $\idemtwo$ factors as a projection and an inclusion. Composing them the other direction should give the identity of $B_{\who_{k-1}}$.
\begin{equation} \label{eq:LIFdown} {
\labellist
\small\hair 2pt
 \pinlabel {$= (-1)$} [ ] at 83 53
 \pinlabel {$k$} [ ] at 25 51
 \pinlabel {$k-1$} [ ] at 21 87
 \pinlabel {$k-1$} [ ] at 21 15
 \pinlabel {$k-1$} [ ] at 133 53
\endlabellist
\centering
\ig{1}{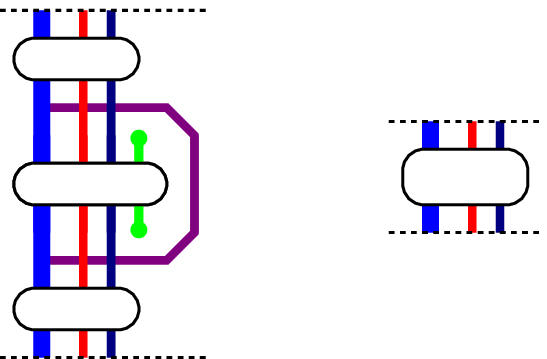}
} \end{equation}
Hence $\idemtwo$ is actually an idempotent.

The idempotents $\idemone$ and $\idemtwo$ are orthogonal thanks to
\begin{equation} \label{eq:emortho} {
\labellist
\small\hair 2pt
 \pinlabel {$k+1$} [ ] at 32 15
 \pinlabel {$k-1$} [ ] at 21 51
\endlabellist
\centering
\ig{1}{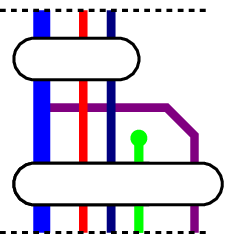}
} \quad = 0, \end{equation}
which also expresses the fact that $\Hom^0(B_{\who_{k+1}}, B_{\who_{k-1}}) = 0$.

Most of the occurrences of $e_{k-1}$ in the formulas above are actually redundant, thanks to an iteration of
\begin{equation} \label{eq:emabsorb} {
\labellist
\small\hair 2pt
 \pinlabel {$=$} [ ] at 102 30
 \pinlabel {$k+1$} [ ] at 61 17
 \pinlabel {$k+1$} [ ] at 157 29
 \pinlabel {$k$} [ ] at 54 42
\endlabellist
\centering
\ig{1}{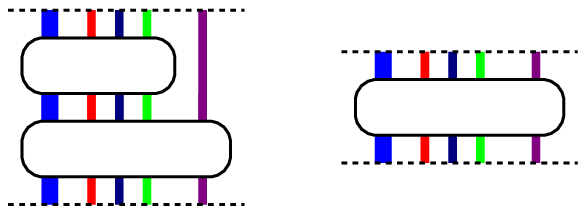}
}. \end{equation}
This formula expresses the general fact that top idempotents factor through top idempotents.

Finally, let $j \in \{1, \ldots, N-k-1\}$ so that $j$ commutes with $X = \{N-k+1,\ldots,N\}$. Then a $j$-colored thick trivalent vertex also commutes with the idempotent $e_k$.
\begin{equation} \label{eq:commutespastem} {
\labellist
\small\hair 2pt
 \pinlabel {$=$} [ ] at 111 32
 \pinlabel {$k$} [ ] at 54 42
 \pinlabel {$k$} [ ] at 150 17
 \pinlabel {$j$} [ ] at 84 14
 \pinlabel {$j$} [ ] at 184 30
\endlabellist
\centering
\ig{1}{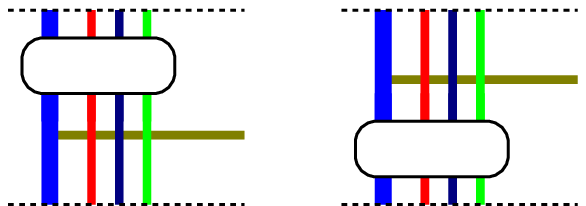}
} \end{equation}

\begin{proof} We will use \eqref{eq:emrecursion} to define the endomorphisms $e_k$ for each $k \ge 2$, given that $e_0$ and $e_1$ are identity maps. Then we prove the remaining formulae, together with the fact that $e_k$ is an idempotent and projects to $B_{\who_k}$, by a simultaneous induction.
	
The base case of our induction are all the pictures that only use $e_0$ and $e_1$, namely $k=1$ for \eqref{eq:LIFdown} and \eqref{eq:ekidemp} and \eqref{eq:commutespastem}, and $k=0$ for
\eqref{eq:emabsorb}. Of these only \eqref{eq:LIFdown} at $k=1$ is nontrivial. The LHS is just like \eqref{eq:thickneedle} but with a green barbell inside. Forcing the barbell out using
the polynomial forcing relation, we break the purple strand with coefficient $-1$, and this gives the RHS by \eqref{eq:thickunit}. The diagram with a polynomial outside is zero by
\eqref{eq:thickneedle}.

We now assume all the equations hold when the only involve idempotents $e_\ell$ for $\ell \le k$, and check the equations which involve the idempotent $e_{k+1}$. We use \eqref{eq:emrecursion} to define $e_{k+1}$, so obviously this equation holds.

The formula \eqref{eq:emrecursion} writes $e_{k+1}$ in terms of two diagrams which each have $e_k$ on top. By \eqref{eq:ekidemp}, $e_k$ is idempotent, so we deduce \eqref{eq:emabsorb}. Similarly, choose some $j \in \{1, \ldots, N-k-2\}$. Then \eqref{eq:emrecursion} writes $e_{k+1}$ in terms of two diagrams which all commute with the $j$-colored thick trivalent vertex (by induction) so we deduce \eqref{eq:commutespastem} for $k+1$.

Also, note that $\idemtwo$, the rightmost diagram in \eqref{eq:emrecursion}, is an idempotent by \eqref{eq:LIFdown}. Hence $e_{k+1}$ is the difference of two idempotents $e_k \ot \id$ and $\idemtwo$, and the latter absorbs the former, so $e_{k+1}$ is also an idempotent. This proves \eqref{eq:ekidemp} for $k+1$.

Now resolve the copy of $e_{k+1}$ inside \eqref{eq:emortho}. The result is
\begin{equation} {
\labellist
\small\hair 2pt
 \pinlabel {$=$} [ ] at 67 64
 \pinlabel {$+ (-1)$} [ ] at 155 68
 \pinlabel {$k+1$} [ ] at 107 47
 \pinlabel {$k$} [ ] at 25 47
 \pinlabel {$k$} [ ] at 201 87
 \pinlabel {$k$} [ ] at 201 14
 \pinlabel {$k-1$} [ ] at 22 83
 \pinlabel {$k-1$} [ ] at 99 83
 \pinlabel {$k-1$} [ ] at 198 123
 \pinlabel {$k-1$} [ ] at 198 51
\endlabellist
\centering
\ig{1}{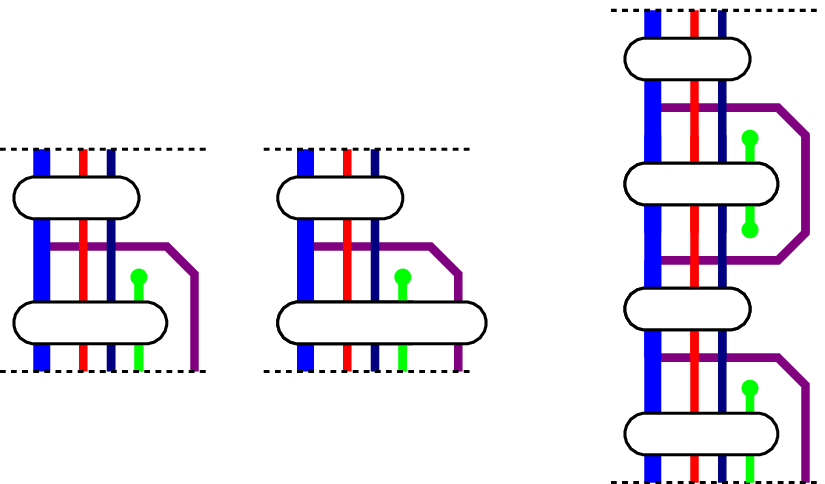}
} \end{equation}
Resolving the rightmost diagram with \eqref{eq:LIFdown}, this diagram becomes equal to the LHS. Thus the first diagram on the RHS is zero, proving \eqref{eq:emortho}.

Finally we must show the most interesting computation, \eqref{eq:LIFdown} for $k+1$. Now we let olive denote the next color in line, $s_{N-k-1}$. Applying the recursive definition we have
\begin{equation} \label{eq:LIFdowncomp}{
\labellist
\small\hair 2pt
 \pinlabel {$=$} [ ] at 74 88
 \pinlabel {$+$} [ ] at 164 90
 \pinlabel {$k+1$} [ ] at 28 86
 \pinlabel {$k$} [ ] at 27 123
 \pinlabel {$k$} [ ] at 26 50
 \pinlabel {$k$} [ ] at 114 50
 \pinlabel {$k$} [ ] at 114 86
 \pinlabel {$k$} [ ] at 114 123
 \pinlabel {$k$} [ ] at 201 15
 \pinlabel {$k$} [ ] at 201 50
 \pinlabel {$k$} [ ] at 201 123
 \pinlabel {$k$} [ ] at 201 160
 \pinlabel {$k-1$} [ ] at 199 86
\endlabellist
\centering
\ig{1}{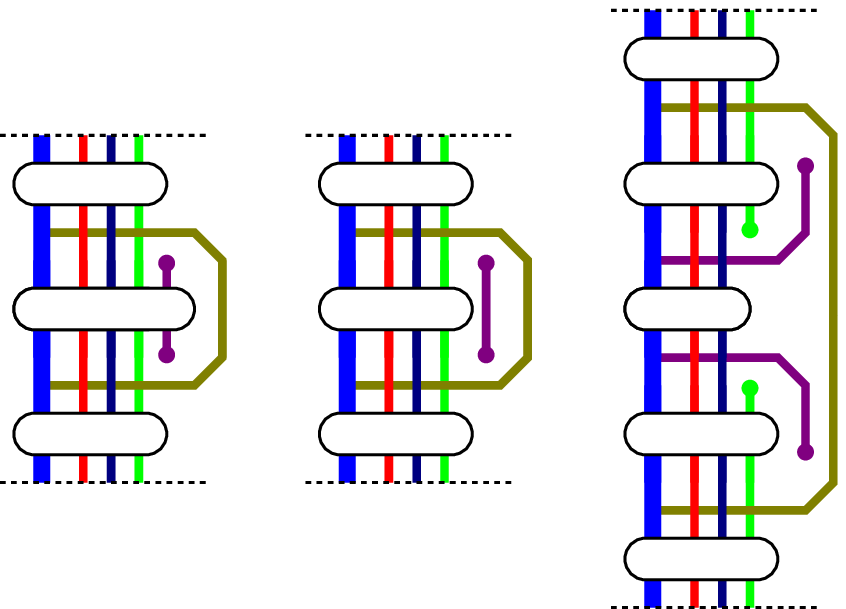}
} \end{equation}
In the second diagram on the RHS of \eqref{eq:LIFdowncomp}, the purple strands will disappear using \eqref{eq:thickunit}. Then the olive strands can slide through the rest by \eqref{eq:commutespastem} until they meet in the middle, where the result is zero by \eqref{eq:thickneedle}. Hence only the first diagram on the RHS of \eqref{eq:LIFdowncomp} contributes. Let us force the purple barbell to the right using \eqref{eq:barbellforce}. In one term the olive strand is broken with coefficient $\partial_{N-k-1}(\alpha_{N-k}) = -1$, and this term gives the identity map (with coefficient $-1$) thanks to \eqref{eq:thickunit} and \eqref{eq:ekidemp}.
\begin{equation} {
\labellist
\small\hair 2pt
 \pinlabel {$=$} [ ] at 76 52
 \pinlabel {$k$} [ ] at 25 52
 \pinlabel {$k$} [ ] at 25 87
 \pinlabel {$k$} [ ] at 25 15
 \pinlabel {$k$} [ ] at 115 52
\endlabellist
\centering
\ig{1}{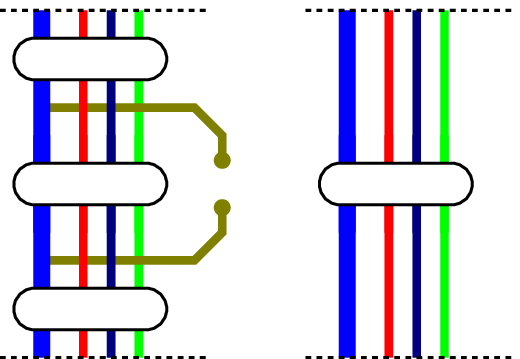}
} \end{equation} 
In the other term, the olive strand is unbroken and some polynomial appears outside it. Again, the olive strands slide through the rest by \eqref{eq:commutespastem} and the result is zero by \eqref{eq:thickneedle}. Thus the final result is minus the identity of $e_k$, as desired. \end{proof}

\begin{cor} \label{cor:ekiswhok} The image of $e_k$ is isomorphic to the indecomposable object $B_{\who_k}$, whose symbol in the Grothendieck group is equal to $b_{\who_k}$. \end{cor}
	
\begin{proof} This follows since the images of these idempotents also categorify \eqref{eq:wIhXsmfinite}. \end{proof}

\begin{remark} This proof does not rely on the Soergel conjecture, only on the fact that $B_I$ is smooth, so it generalizes to arbitrary characteristic. \end{remark}

%===========
\subsection{Funky trivalent vertices}
\label{subsec-funkytri}
%===========

Now we plan to categorify \eqref{eq:wIhXsmfinitecase3}. We want idempotents $e_+$ and $e_-$ which give the direct sum decomposition
\begin{equation} \label{whoksjpm} B_{\who_k} B_{s_j} \cong B_{\who_k}(+1) \oplus B_{\who_k}(-1) \end{equation}
whenever $j \ne N-k$. There are two separate cases, when $j < N-k$ and when $j > N-k$.

Suppose $j < N-k$, and draw it using the color olive as before. Henceforth, we will encode the diagram on either side of \eqref{eq:commutespastem} as a new \emph{funky trivalent vertex}.
\begin{equation}
{
\labellist
\small\hair 2pt
 \pinlabel {$=$} [ ] at 108 31
 \pinlabel {$k$} [ ] at 53 41
 \pinlabel {$k$} [ ] at 148 30
\endlabellist
\centering
\ig{1}{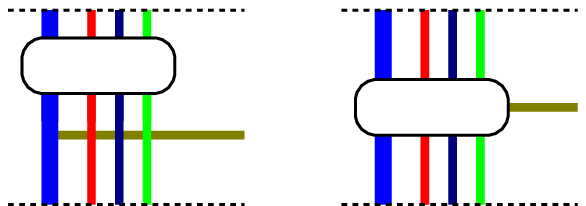}
} \end{equation}
This is the ``easy'' case of the funky trivalent vertex.

Now suppose that $j > N-k$, so it is one of the colors already present in $X = \{N-k+1, \ldots, N\}$.  We seek to define a \emph{funky trivalent vertex} here as well. Let us illustrate them by example when the size of $X$ is three.
\begin{equation} \label{eq:funkytrihard} {
\labellist
\small\hair 2pt
 \pinlabel {$=$} [ ] at 85 242
 \pinlabel {$=$} [ ] at 85 143
 \pinlabel {$=$} [ ] at 85 42
\endlabellist
\centering
\ig{1}{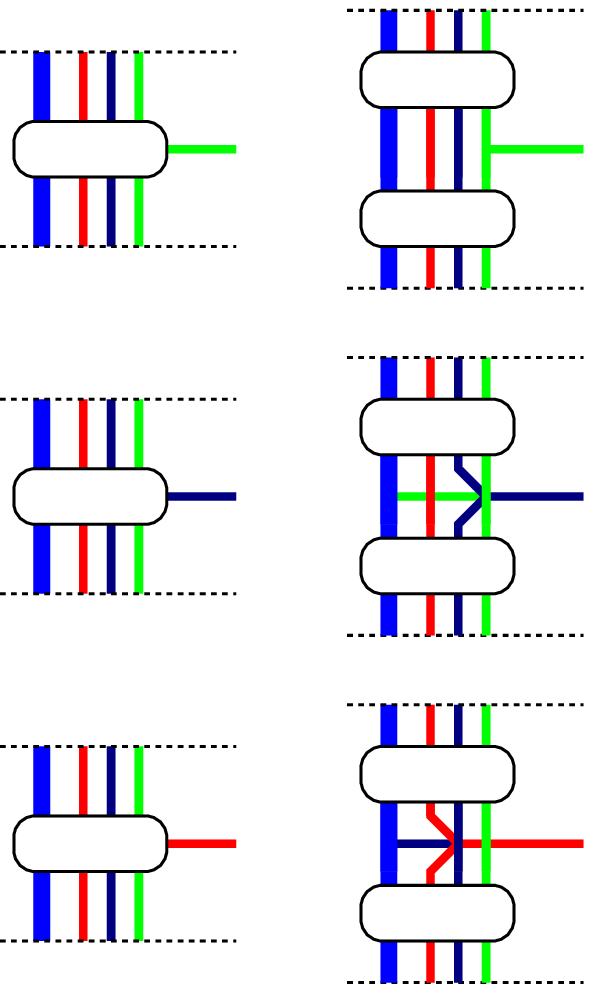}
} \end{equation}
When $j = N-k+1$ is the last strand in $X$, we merely use a trivalent vertex. Otherwise, we commute $j$ past the strands in $X$ until it meets the color $j-1$, then we apply a $6$-valent vertex to replace $j$ with $j-1$, and then we pull $j-1$ past the remaining colors until it can do a thick trivalent vertex with $I$.

\begin{remark} In this paper, we will only ever need to use the decomposition \eqref{whoksjpm} when $j < N-k$, so we will only ever need to use easy funky trivalent vertices. We include the more complicated cases $j > N-k$ because we feel they belong here, but the reader is welcome to only read the easy parts of this section. \end{remark}

\begin{remark} The details of this construction are highly reminiscent to the morphisms $a_j$ which appear in the definition of the thick trivalent vertices, see \cite[Definition
3.34]{EThick}. Suppose that $X = T$ so that $w_I h_X = \who_N = w_T$. Then $e_N$ does actually project to $B_{w_T}$, the longest element for $W_T$, and the funky trivalent vertices
defined in \eqref{eq:funkytrihard} agree with the thick trivalent vertices defined in \cite{EThick}. In fact, the thick calculus of \cite{EThick} could be reproven using the formulas just
computed here, although this would not reprove the alternate description of $e_N$ in \cite{EThick} which uses Manin-Schechtmann theory. \end{remark}

Let us prove that the funky trivalent vertices also satisfy their own analogs of the unit and associativity axioms.
\begin{equation} \label{eq:funkyunit} {
\labellist
\small\hair 2pt
 \pinlabel {$=$} [ ] at 82 30
 \pinlabel {$j$} [ ] at 57 38
\endlabellist
\centering
\ig{1}{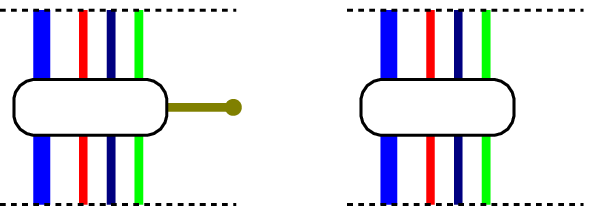}
} \end{equation}
\begin{equation} \label{eq:funkyassoc} {
\labellist
\small\hair 2pt
 \pinlabel {$=$} [ ] at 75 22
 \pinlabel {$j$} [ ] at 54 30
 \pinlabel {$j$} [ ] at 142 33
\endlabellist
\centering
\ig{1}{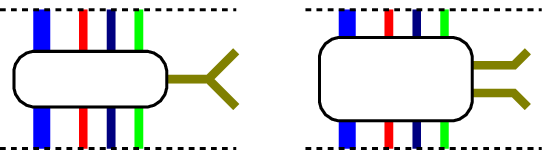}
} \end{equation}

The RHS of \eqref{eq:funkyassoc} uses an abuse of notation, which we now explain. For each $j \ne N-k$, the funky trivalent vertex can be described as ``apply $e_k$, then do some interesting diagram where the action is happening, and then apply $e_k$ again.'' Thus if you square this operation, and use the fact that $e_k$ is an idempotent, the diagram reads, ``$e_k$ then action then $e_k$ then action then $e_k$.'' Meanwhile, there is another thing one might do, which is ``$e_k$ then action then action then $e_k$,'' which removes the middle copy of $e_k$. In fact, these operations are equal. Here is what we mean explicitly when $X$ has size $4$ and $j = N-1$ is red; this is the schematic picture for any $j > N-k+1$.
\begin{equation} \label{eq:funkyaux} {
\labellist
\small\hair 2pt
 \pinlabel {$=$} [ ] at 82 61
 \pinlabel {$k$} [ ] at 28 62
\endlabellist
\centering
\ig{1}{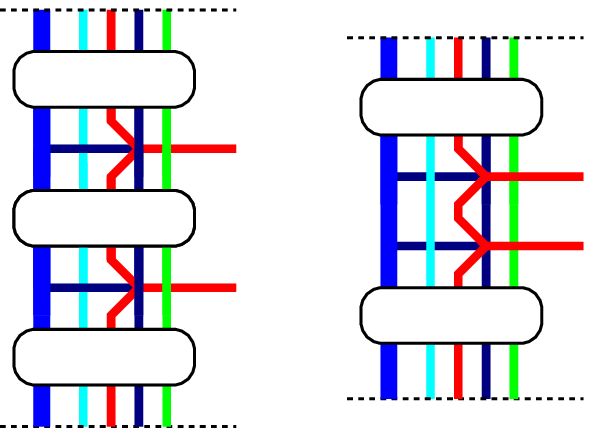}
} \end{equation}
So the square of the funky trivalent vertex has two descriptions, and the second description looks like the one used abusively in \eqref{eq:funkyassoc}.  When $j = N-k+1$ or $j < N-k$, it is easy to prove this fact about the square of the funky trivalent vertex, but the proof is nontrivial for $j > N-k+1$, and appears in the next proposition.

If \eqref{eq:funkyunit} and \eqref{eq:funkyassoc} hold, then an empty eyehole is zero as in \eqref{eq:thicksdecomp}. This enables a factored idempotent decomposition just as in \eqref{eq:thicksdecomp}, which is really our goal.
\begin{equation} \label{eq:whokjdecomp} {
\labellist
\small\hair 2pt
 \pinlabel {$=$} [ ] at 69 22
 \pinlabel {$+$} [ ] at 154 22
 \pinlabel {$k$} [ ] at 27 22
 \pinlabel {$k$} [ ] at 107 22
 \pinlabel {$k$} [ ] at 189 22
 \pinlabel {$\frac{\al_j}{2}$} [ ] at 133 8
 \pinlabel {$\frac{\al_j}{2}$} [ ] at 217 33
\endlabellist
\centering
\ig{1}{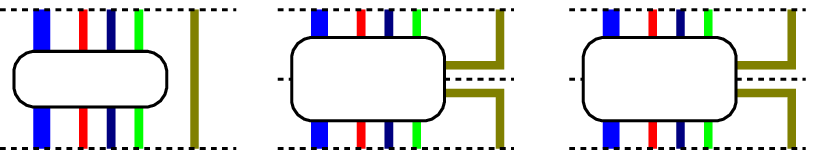}
} \end{equation}

\begin{prop} Equations \eqref{eq:funkyunit} and \eqref{eq:funkyassoc} hold for all $j \ne N-k$. Equation \eqref{eq:funkyaux} holds for $j > N-k+1$. \end{prop}
	
\begin{proof} We check the cases $j < N-k$, $j = N-k+1$, and $j > N-k+1$ separately, since the funky trivalent vertex is defined in different ways each time. For $j < N-k$, the funky
trivalent vertex is essentially just a thick trivalent vertex with some 4-valent crossings, so \eqref{eq:funkyunit} follows from \eqref{eq:thickunit} and \eqref{eq:funkyassoc} follows
from \eqref{eq:thickassoc}. When $j = N-k+1$, the funky trivalent vertex is just an ordinary trivalent vertex, and so \eqref{eq:funkyunit} follows from the ordinary unit relation and
\eqref{eq:funkyassoc} from the ordinary associativity relation.

Hence we assume that $j > N-k+1$. We illustrate the proof in the example where $X$ has size $4$ and $j = N-1$, as the calculation is effectively the same in all other cases.
\begin{equation} {
\labellist
\small\hair 2pt
 \pinlabel {$=$} [ ] at 83 43
 \pinlabel {$+$} [ ] at 173 43
\endlabellist
\centering
\ig{1}{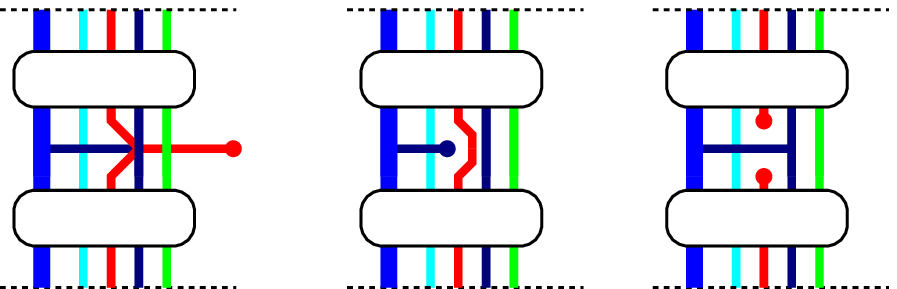}
} \end{equation}
The first diagram on the RHS simplifies by \eqref{eq:thickunit} and becomes the idempotent $e_k$. The second diagram on the RHS is zero by \eqref{eq:emortho} (and if necessary \eqref{eq:emabsorb}). Thus we have proven \eqref{eq:funkyunit}.

Let us prove \eqref{eq:funkyassoc} assuming \eqref{eq:funkyaux}. We have
\begin{equation} {
\labellist
\small\hair 2pt
 \pinlabel {$=$} [ ] at 79 55
 \pinlabel {$=$} [ ] at 167 55
\endlabellist
\centering
\ig{1}{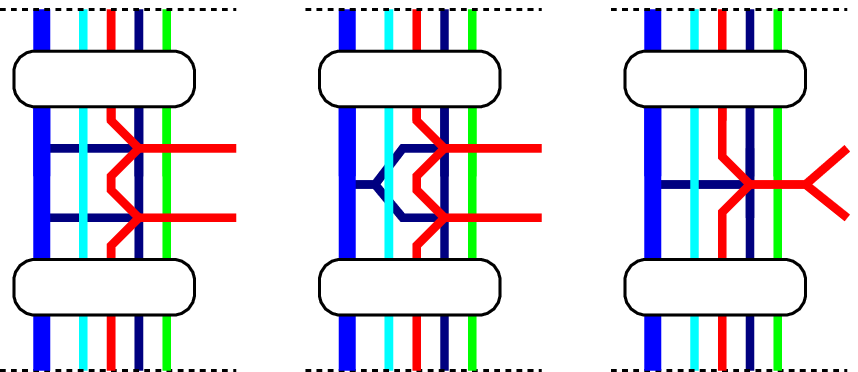}
} \end{equation}
where the first equality was \eqref{eq:thickassoc} and the second came from two-color associativity.

Now let us prove \eqref{eq:funkyaux}. We use induction as follows: instead of placing the idempotent $e_k$ in the middle (the labeled oval in \eqref{eq:funkyaux}), we place $e_{\ell}$
for some $\ell \le k$. We then prove by induction that this idempotent $e_{\ell}$ can be removed. The base cases are $e_0$ and $e_1$, which are both equal to their respective identity
maps. For the induction step, \eqref{eq:emrecursion} writes $e_{\ell}$ in terms of $e_{\ell-1} \ot \id$ and another idempotent $\idemtwo$, so if we can show that putting $\idemtwo$ in
place of $e_{\ell}$ yields zero, then we finish the proofs.

The most interesting cases are $j = N-\ell, N-\ell+1, N-\ell+2$ which look like 
\begin{equation} \label{eq:funkyauxcomp1} \ig{1}{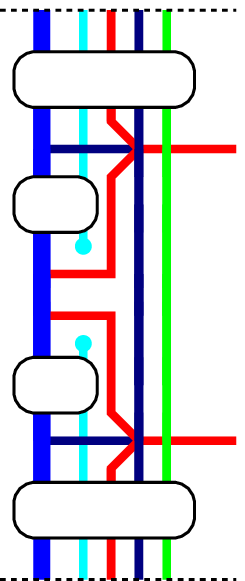}, \quad \ig{1}{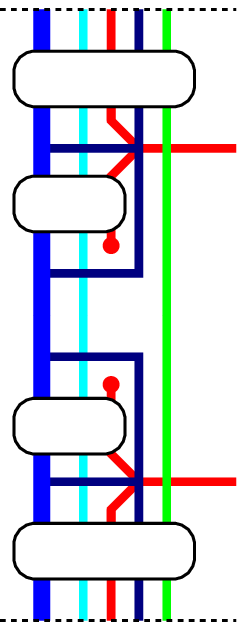}, \quad \ig{1}{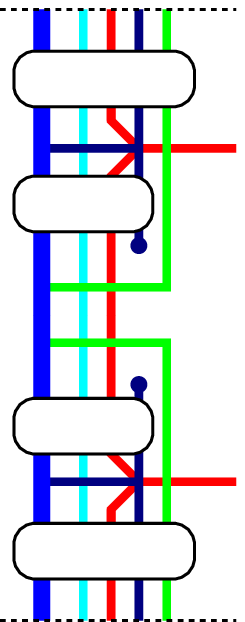}. \end{equation}
In each diagram there are two ``small idempotents'' (the two smaller ovals in the middle) which get somewhat in the way. In the first diagram of \eqref{eq:funkyauxcomp1}, the small idempotents can commute past the blue strand by \eqref{eq:commutespastem} and then get absorbed by $e_k$ on the outside by \eqref{eq:emabsorb}. Thus we can pretend the small idempotents were not there in this first diagram. It is not so easy to prove that the small idempotents are redundant in the other two diagrams (or for $\ell > N-\ell+2$ either) but in fact they are! Let us assume this for now, and discuss how to remedy the situation after seeing the arguments.

Consider the first diagram of \eqref{eq:funkyauxcomp1}, ignoring the small idempotents. A combination of \eqref{eq:thickassoc} and \eqref{eq:thick2m} yield
\begin{equation} \ig{1}{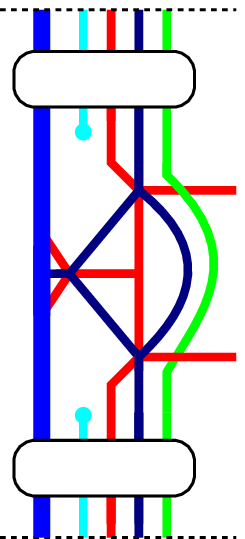}. \end{equation}
In each half of this diagram there are two red strands close to the aqua dot, which could be fused together using the first relation in \eqref{eq:barbellforce}. When fusing two strands, a polynomial appears on one side or the other. If the polynomial is not next to the aqua dot, the result is zero by \eqref{eq:emortho}. So for both the top and bottom half the polynomial is next to the aqua dot. Now we have a subdiagram
\begin{equation} \ig{1}{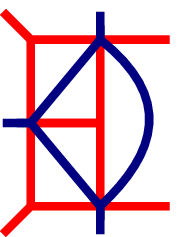} \end{equation}
which is a degree $-3$ diagram in a Hom space for which the lowest degree is $-1$, so it is zero. It is important in the above calculation that the entire diagram was zero, while either half need not be; otherwise a simpler argument using only the top or bottom half of the diagram would suffice.

Consider the second diagram of \eqref{eq:funkyauxcomp1}, ignoring the small idempotents. We can resolve the red dot into the 6-valent and the two resolutions vanish, either by \eqref{eq:emortho} or by \eqref{eq:thickneedle}.

Consider the third diagram of \eqref{eq:funkyauxcomp1}, ignoring the small idempotents. One can resolve the blue dot into the 6-valent vertex. The two resolutions vanish by \eqref{eq:emortho}, either for a red dot under a blue cup, or a blue dot under a green cup.

For $\ell$ larger or smaller, assuming one can deal with the small idempotents appropriately, the result will be zero by \eqref{eq:emortho}.

To deal with the small idempotents, we must add the statement that they are redundant (in pictures like \eqref{eq:funkyauxcomp1}) to our inductive hypothesis. Again, this is proven by
applying the recursive definition to the small idempotents - if the reader attempts this, they will see that one recovers exactly the same intermediate diagrams as in the computations
above, which vanish for the same reasons. \end{proof}

\begin{remark} We will only need the easy funky trivalent vertices (for $j < N-k$) in the rest of this paper, but felt it was appropriate to include them anyway. \end{remark}

%===========
\subsection{Checkpoint}
\label{subsec-minisummary}
%===========

Let us summarize what was just accomplished. In \cite{EThick} a thick calculus is developed which computes the idempotents projecting to the indecomposable objects $B_I$ corresponding
to longest elements $w_I$ of parabolic subgroups, and adds them as new objects. Complicated morphisms are encoded as thick trivalent vertices, and it is shown that these satisfy analogs
of the unit and associativity axiom, making them easy to work with.

In the last few sections, we extended this thick calculus to the elements $\who_k$, or equivalently $w_I h_X$ when $X$ is a suffix of $T$. We computed the funky trivalent vertices, and
showed that they also satisfy the unit and associativity axioms. (They also satisfy an analog of \eqref{eq:thick2m}, and the proof is quite similar and slightly easier than the
corresponding proof in \cite{EThick}.) All the plethysm relations of \eqref{eq:wIhXsmfinite} were lifted to explicit idempotent decompositions.

Now we have the ability to fully decompose the identity of the tensor product $B_I B_X$ into orthogonal idempotents for any $X \subset T$ (not necessarily a suffix). Because we will be
doing this in a more complicated setup soon, let us elaborate upon the process.

Suppose that $X=T$ so that $B_I B_X = B_I B_N B_{N-1} \cdots B_2 B_1$, and begin with the identity map, which is $e_0 \ot \id \ot \cdots \ot \id$. We ``absorb'' one index from $X$ into
the idempotent at a time. In the first stage, we absorb $B_N$ and replace $e_0 \ot \id$ with $e_1$ (which is also the identity map). In the next stage, we replace $e_1 \ot \id$ with
$e_2$ and with $\idemtwo$, the two idempotents from \eqref{eq:emrecursion}. Now we have two idempotents which project to $B_{\who_2} B_{N-2} \cdots B_2 B_1$ and $B_{\who_0} B_{N-2}
\cdots B_2 B_1$ respectively. For the second summand we can repeatedly use \eqref{eq:thicksdecomp} (or equivalently, \eqref{eq:whokjdecomp}) to absorb the remaining indices of $X$,
producing $2^{N-2}$ orthogonal idempotents, since $s_j \in \RC(\who_0)$ for $1 \le j \le N-2$. The other idempotent is dealt with in the next stage. At the $k$-th stage we replace $e_k
\to \id$ with $e_{k+1}$ and $\idemtwo$, yielding two orthogonal idempotents projecting to $B_{\who_{k+1}} B_{N-k-1} \cdots B_2 B_1$ and $B_{\who_{k-1}} B_{N-k-1} \cdots B_2 B_1$
respectively. The latter can be dealt with by repeated use of \eqref{eq:whokjdecomp} (which only uses the ``easy'' funky trivalent vertices), and the former is the input to the next
stage.

Suppose that $X = T \setminus i$ for some $i$ in the middle of $T$. Then $X$ has two connected components, which commute with each other; one is a suffix and one is contained entirely
in $I$. The suffix should be dealt with using the process above. The subset of $I$ should be dealt with purely using thick trivalent vertices and the decomposition
\eqref{eq:thicksdecomp} (or equivalently, \eqref{eq:whokjdecomp}). Note that these two subsets can be dealt with in either order, and one recovers the same idempotents! The pieces of
the idempotents belonging to these two components commute with each other, using \eqref{eq:commutespastem} and other relations involving 4-valent vertices. Similarly, an arbitrary $X
\subset T$ is a disjoint union of its maximal suffix and the maximal union of connected components contained in $I$. These two parabolic subgroups commute with each other, and can be dealt with independently as above.

For the rest of this chapter, we will be repeating this process to decompose the tensor product $B_I B_X$ in the affine setting. Within each $\tau$-component, we want to use the process above to decompose the identity map of some tensor factor in $B_I B_X$ into orthogonal idempotents. The main complication is
that the $\tau$-component in question need not be at the beginning or the end of the expression $w_I h_X$, but could be somewhere in the middle, or could even be split between the
beginning and the end! In our Gaussian elimination procedure to come, we will attack a single $\tau$-component for all $X$ simultaneously, so we do not have the freedom of attacking the
$\tau$-components in cyclic order (this depends on $X$). The process is no longer entirely local.

The rest of this chapter is devoted to handling the complications that can occur, and introducing notation to keep track of idempotents. The main notational confusion is that one can
not work with elements of the affine Weyl group $W_{\aff}$, but must work in the affine Hecke algebra $\Hecke_{\aff}$. Above we could use $\who_k \in W_{\aff}$ to index our idempotents,
but this only labels the local idempotents, not the global summands in the decomposition process. For example, we are really interested in the product $b_{\who_k} b_{N-k} \cdots b_2
b_1$, as opposed to the element $\who_k s_{N-k} \cdots s_2 s_1$ which does not depend on $k$.

%===========
\subsection{Smoothness and plethysm: affine case}
\label{subsec-prelimsmoothBIaffine}
%===========

Let $A$ be a $\tau$-component. We seek a procedure which takes the identity map of $B_I B_X$ (or some other idempotent, wait a moment) and start absorbing the indices of $X \cap A$
one by one, decomposing the identity into a sum of orthogonal idempotents. We will absorb the indices of $X$ from $\max(A)$ down to $\min(A)$. When we
are done, we will have a number of idempotents which involve all the strands in $X \cap A$ and the thick strand for $I \cap A$, but should not involve the strands in the rest of $X$,
beyond the 4-valent vertices needed to commute strands past each other. Having completed this, we move on to the next $\tau$-component.

We should not assume that we begin with the identity map of $B_I B_X$, but that we begin having already absorbed some of the other $\tau$-components, and proceed to absorb one
more. We should hope that the idempotent decomposition one obtains at the end does not depend on the order in which the $\tau$-components were absorbed.

Recall \eqref{eq:rewriteaffine} which rewrites $w_I h_X$ as
\[ w_I h_X = (w_{I \cap A_0} h_{X \cap A_{< \notme}}) (w_{I \cap A_d} h_{X \cap A_d}) \cdots (w_{I \cap A_1} h_{X \cap A_1}) h_{X \cap A_{> \notme}}. \]
The final complication will be whether our component in question $A$ is the special component $A_0$, or one of the middle components $A_i$, $i \ne 0$.

\begin{notation} We work in the affine setting, and fix $I \subsetneq S$. A subset $X \subsetneq S$ is a \emph{$\tau$-suffix} if, for each $\tau$-component $A$, $X \cap A$ is a suffix of
$A$. If $m \notin X$ we speak of $m$ as \emph{initiating} or \emph{growing} a $\tau$-suffix based on what it does to the $\tau$-component containin $A$, see the paragraph before Lemma \ref{lem:finitepleth}. \end{notation}

\begin{example} If $I = \{1,2,3,5,9,10\}$ and $X$ is a $\tau$-suffix, then $2 \in X$ implies $3, 4 \in X$. \end{example}

In the definition below, analogous to Lemma \ref{lem:finitepleth}, one should think of $Y$ as the strands in a $\tau$-suffix which have been already absorbed at some intermediate stage
of the process.

\begin{defn} \label{defn:byXY} Fix $I \subsetneq S$ as usual. Let $X \subsetneq S$ be an arbitrary subset, and $Y \subset X$ be a $\tau$-suffix with respect to $I$. We will recursively define elements
$\by_{I,X,Y} \in \Hecke_{\aff}$.

We begin with \begin{equation} \label{Yempty} \by_{I,X,\mt} = b_{w_I} b_{h_X}. \end{equation}
Suppose that $m \in X$, $m \notin Y$, and $Ym = Y \cup \{m\}$ is also a $\tau$-suffix. If $m$ initiates a $\tau$-suffix then we define
\begin{equation}\label{initiate} \by_{I,X,Y} = \by_{I,X,Ym}. \end{equation}
If $m$ grows a $\tau$-suffix, so that $m+1 \in Y$ and $m, m+1$ are in the same $\tau$-component, then
\begin{equation}\label{grow} \by_{I,X, Y} = \by_{I,X,Ym} + \by_{I,X,Y \setminus m+1}. \end{equation}
\end{defn}

This recursive definition is concise, but it obfuscates the fact that these elements $\by$ can really be defined $\tau$-component by $\tau$-component, as we will see in the proof below.

\begin{prop} \label{prop:byXY} This definition is consistent, giving rise to a well-defined element $\by_{I,X,Y}$ for all $Y \subset X \subsetneq S$ with $Y$ a $\tau$-suffix. Moreover, when $Y = X$ is a $\tau$-suffix we have
\begin{equation} \by_{I,X,X} = b_{w_I h_X}, \end{equation}
and this Kazhdan-Lusztig basis element is smooth.

If $j \in X$, $j \notin Y$ and $Yj$ is not a $\tau$-suffix, then necessarily $j \in I$ and
\begin{equation} \label{godown} \by_{I,X,Y} = (v+v\inv) \by_{I,(X\setminus j),Y}. \end{equation}
\end{prop}

\begin{proof} The formula \eqref{eq:rewriteaffine} only used commutation braid relations, so it continues to hold in the Hecke algebra with the Kazhdan-Lusztig generators.
\begin{equation} \label{eq:rewriteaffine3} b_{w_I} b_{h_X} = (b_{w_{I \cap A_0}} b_{h_{X \cap A_{< \notme}}}) (b_{w_{I \cap A_d}} b_{h_{X \cap A_d}}) \cdots (b_{w_{I \cap A_1}} b_{h_{X \cap A_1}}) b_{h_{X \cap A_{> \notme}}}. \end{equation}
Consider what happens within a single $\tau$-component $A$ for $A \ne A_0$. The rules \eqref{initiate} and \eqref{grow} agree with \eqref{eq:wIhXsmfinitecase2} and \eqref{eq:wIhXsmfinitecase1} respectively, when adding reflections $m \in A$ to $Y$. Hence, the process of going from $b_{w_I} b_{h_X}$ to $\by_{I,X,Y \cap A}$ is the same as replacing
\begin{equation} b_{w_{I \cap A}} b_{h_{X \cap A}} \rightsquigarrow b_{w_I h_{Y \cap A}} b_{h_{(X \setminus Y) \cap A}}. \end{equation}
That is, it absorbs the suffix of $X \cap A$ into a Kazhdan-Lusztig basis element, and leaves behind the rest of $b_{h_X}$. Since this operation is local to the $\tau$-component, it commutes with what happens in other $\tau$-components.

Let us examine what happens for the special $\tau$-component $A_0$. In this case, any suffix in $X \cap A_0$ is contained within $X \cap A_{> \notme}$. Note also that any simple reflection in $I \cap A_{> \notme}$ commutes with any simple reflection in $A_{< \notme}$ and any simple reflection in $A_i$ for $i \ne 0$. The Kazhdan-Lusztig basis element $b_{w_{I \cap A_0}}$ is in the left ideal of $b_{w_{I \cap A_{> \notme}}}$, so we can write it as
\begin{equation} b_{w_{I \cap A_0}} = p b_{w_{I \cap A_{> \notme}}} \end{equation}
for some $p$. Then we can rewrite \eqref{eq:rewriteaffine3} as
\begin{equation} \label{eq:rewriteaffine4} b_{w_I} b_{h_X} = p b_{h_{X \cap A_{< \notme}}} (b_{w_{I \cap A_d}} b_{h_{X \cap A_d}}) \cdots (b_{w_{I \cap A_1}} b_{h_{X \cap A_1}}) (b_{w_{I \cap A_{> \notme}}} b_{h_{X \cap A_{> \notme}}}). \end{equation}
Again, what happens with suffixes is local in the last factor involving the parabolic subgroup $A_{> \notme}$.

Putting this together, we see that $\by_{I,X,Y}$ does not depend on what order the indices were added to $Y$, because it can be defined without induction by the wordy formula
\begin{equation} \label{eq:rewriteaffine5} \by_{I,X,Y} = p b_{h_{X \cap A_{< \notme}}} (b_{w_{I \cap A_d} h_{Y \cap A_d}} b_{h_{(X \setminus Y) \cap A_d}}) \cdots (b_{w_{I \cap A_{> \notme}} h_{Y \cap A_{> \notme}}} b_{h_{(X \setminus Y) \cap A_{> \notme}}}). \end{equation}

Working componentwise suffices to prove \eqref{godown}, using \eqref{eq:wIhXsmfinitecase3}. The only thing worth mentioning is that when $j \in A_0$ there are two cases. If $j \in A_{>
\notme}$ then the situation follows from \eqref{eq:wIhXsmfinitecase3} applied to the parabolic subgroup $A_{> \notme}$. If $j \in A_{< \notme}$ then $j \in I$ and it is clear that $b_j$
can be absorbed into $b_{w_{I \cap A_0}}$ at the cost of a factor of $(v+v\inv)$.

Now consider the case when $X$ is a $\tau$-suffix. Note that $X \cap A_{< \notme} = \mt$. Probably there is a faster proof, but for safety's sake we will treat the $\tau$-component $A_0$ separately. Let $Y = X \setminus A_0$, which is also a $\tau$-suffix. Then we have
By \eqref{eq:rewriteaffine3} we have
\begin{equation} \label{eq:rewriteaffine6} \by_{I,X,Y} =  b_{w_{I \cap A_0}} \left( \prod_{i \ne 0} (b_{w_{I \cap A_i} h_{X \cap A_i}}) \right) b_{h_{X \cap A_0}} = q b_{h_{X \cap A_0}}. \end{equation}
Now each of $w_{I \cap A_i} h_{X \cap A_i}$ is smooth, as is $w_{I \cap A_0}$, and so is their product by Lemma \ref{lem:smoothIJ}. Thus $q$ is a smooth Kazhdan-Lusztig basis element, $q = b_a$ for some $a \in W_{\aff}$.  We wish to argue that initiating and then growing the suffix in the component $A_0$ will preserve smoothness and eventually yield the Kazhdan-Lusztig basis element. Thankfully this is accomplished entirely by right multiplication of simple reflections on a Kazhdan-Lusztig basis element, which brings the problem back down to earth. The argument is almost identical to the sketch in Lemma \ref{lem:finitepleth}. We give this in example below, and leave the reader to formulate the general proof. \end{proof}

\begin{example} Suppose that $I = \{1,2,3, 5, 9,10\}$ and $n = 12$. Suppose that $X = S \setminus \{1\}$, with $\notme = 1$. We are considering the product
\begin{equation} w_{123} (s_0) (w_{9,10} s_{11} s_{10} s_9) (s_8) (s_7) (w_5 s_6 s_5) (s_4 s_3 s_2) = a (s_4 s_3 s_2). \end{equation}
The element $a$ is smooth, and $s_4$ does not appear in it, so $a s_4$ is smooth by Lemma \ref{lem:smoothIJ}; this handles the initiation of the suffix. For $w \in W$ write $K_w = \{ x \in W \mid x \le w\}$. Now the set $K_a$ is a union of $s_3$ cosets, since $s_3$ commutes with everything in $a$ except $w_{123}$, which is the longest element of a larger parabolic subgroup. Meanwhile, the set $K_{a s_4} \setminus K_a$ consists only of elements minimal in their $s_3$ coset. By the same argument as before, we deduce that
\begin{equation} \Sigma_{a s_4} b_{s_3} = \Sigma_{a s_4 s_3} + \Sigma_{a} \end{equation}
from which one deduces that $a s_4 s_3$ is smooth. Now we repeat the argument: $K_{a s_4}$ is a union of $s_2$ cosets, $K_{a s_4 s_3} \setminus K_{a s_4}$ contains only elements minimal in their $s_2$ coset, so
\begin{equation} \Sigma_{a s_4 s_3} b_{s_2} = \Sigma_{a s_4 s_3 s_2} + \Sigma_{a s_4} \end{equation}
from which one deduces that $a s_4 s_3 s_2$ is smooth. \end{example}

%===========
\subsection{Idempotents: affine case}
\label{subsec-idempsaffine}
%===========

Now we categorify the elements $\by_{I,X,Y}$ with objects $\BY_{I,X,Y}$, using explicit idempotent decompositions. We can do this recursively by categorifying the formulas
\eqref{initiate}, \eqref{grow}, and \eqref{godown}. We can also do this directly, categorifying \eqref{eq:rewriteaffine5}. These will agree, proving that the recursive approach is
well-defined.

Let us draw the idempotent $e_{I,X,Y}$ in $B_I B_X$ which projects to the direct summand $\BY_{I,X,Y}$ in three different ways. Here is the first style, which requires much explanation.
\begin{equation} \label{style1} e_{I,X,Y} = \ig{1}{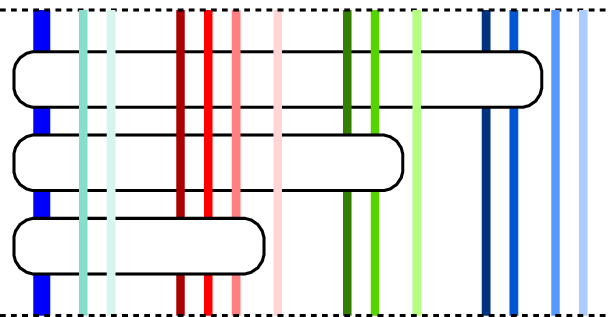} \end{equation}
	
In this picture there are three $\tau$-components $\{A_0, A_1, A_2\}$, and the intersection of ecah $\tau$-component with $X$ is depicted with a different color palette. That is, $X \cap
A_2$ is colored with shades of red, $X \cap A_1$ is colored with shades of green, and $X \cap A_0$ is colored with shades of blue. The shades get lighter as the index decreases in cyclic
order. Note that $X \cap A_0$ is split into two pieces, $A_{> \notme}$ at the end and $A_{< \notme}$ at the start. We continue to draw $B_I$ as a thick blue strand. The $\tau$-suffix $Y$
contains the first 3 red strands, the first 2 green strands, and the first 2 blue strands (which are both in $A_{> \notme}$). The remaining strands in each color are in $X \setminus Y$.

Now let us explain the ovals in \eqref{style1}, starting with the one on the bottom. Let $T = A_2$. In \S\ref{subsec-idempsfinite} we defined an idempotent $e_3$ within the parabolic
subgroup $W_T$, that projects from $B_{I \cap T} B_{Y \cap T}$ to its top summand, when $Y \cap T$ is a suffix of size $3$. This oval roughly represents the idempotent $e_3$, except that
the parabolic on the left is $B_I$ rather than $B_{I \cap T}$, and there are some strands in $A_{< \notme}$ in the way (the blue strands going over the oval). Instead of being defined by the recursion in \eqref{eq:emrecursion} we can use the similar recursion
\begin{equation} \label{recursionstyle1}
{
\labellist
\small\hair 2pt
 \pinlabel {$=$} [ ] at 90 54
 \pinlabel {$+$} [ ] at 195 54
\endlabellist
\centering
\ig{1}{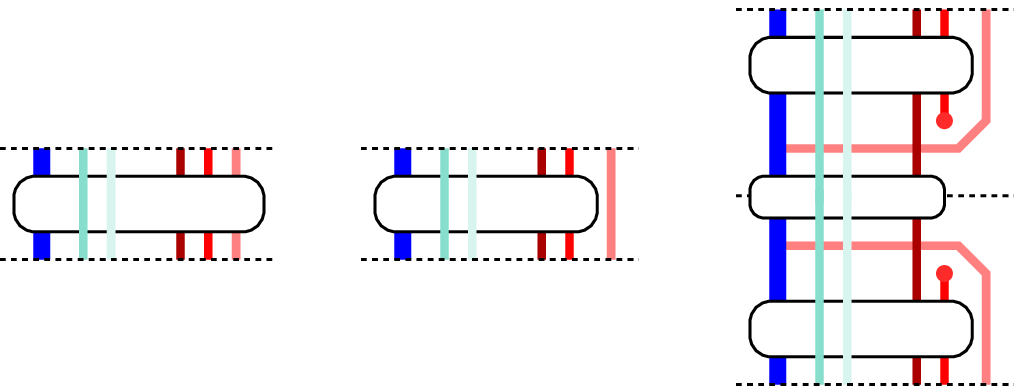}
} \end{equation}
with base case
\begin{equation} \label{recursionstyle1basecase} {
\labellist
\small\hair 2pt
 \pinlabel {$=$} [ ] at 73 17
\endlabellist
\centering
\ig{1}{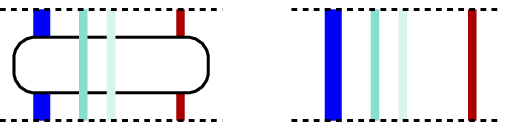}
}. \end{equation}
In other words, we do the same exact recursion but with some extra blue strands that interact with nothing, passing over everything with 4-valent vertices. The rightmost diagram of \eqref{recursionstyle1} makes sense: every strand in $Y \cap A_2$ but the first strand will be contained in $I$ and will commute with $A_{< \notme}$, so the 4-valent vertices and the thick trivalent vertices used are all valid. Proposition \ref{prop:recursionetc} applies to these idempotents as well with essentially the same proof: neither the extra blue strands nor the fact that $I$ is larger than just $I \cap A_2$ affect any part of the calculations.

The other two ovals are defined similarly, only with more extraneous strands to the left. They are also idempotents. Moreover, all these ovals commute with each other, using
\eqref{eq:commutespastem} and other distant color sliding relations.

Here is a second way to draw the idempotent. Note that $I$ is a disjoint union of its connected components, and $B_I \cong B_{I \cap A_0} B_{I \cap A_d} \cdots B_{I \cap A_1}$. This
isomorphism can be realized using thick splitters, see \eqref{eq:splitoffI}. Recall also the thick 4-valent vertices which can commpute $B_{I \cap A_i}$ past $B_{I \cap A_{i'}}$ or past
$B_j$ when $j$ commutes with $B_{I \cap A_i}$, see \eqref{eq:thicksupercrossing}. In this diagram we use the isomorphic object $B_{I \cap A_0} \cdots B_{I \cap A_1}$ instead of $B_I$.
\begin{equation} \label{style2} e_{I,X,Y} = \ig{1}{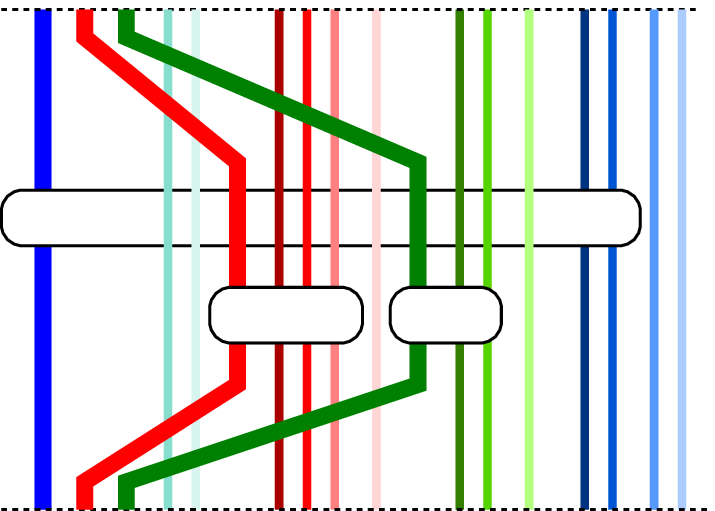} \end{equation}

In this picture, we bring each longest element $B_{I \cap A_i}$ close to the strands of each color, entirely analogously to the expression \eqref{eq:rewriteaffine3}. Now the two small
ovals really are exactly the idempotents $e_k$ defined in \S\ref{subsec-idempsfinite}. However, to avoid splitting up $B_{I \cap A_0}$, we can not bring this thick strand all the way
across to be next to $A_{> \notme}$; it may not commute with the red strands in $X$. For the blue idempotent we continue to use the alternate version \eqref{recursionstyle1} defined
above.

\begin{remark} To account for the signed isomorphisms between different reduced expressions of $X$, we should also used signed crossings when we pull thick strands across thin strands.
In the diagram \eqref{style2}, each thick crossing appears twice so the signs cancel. \end{remark}

Our third way of drawing the idempotent will be most analogous to \eqref{eq:rewriteaffine4}. We take a reduced expression for $w_{I \cap A_0}$ and write it as $r w_{I \cap A_{> \notme}}$
for some $r \in W_{\aff}$. We can use a thick splitter, see \eqref{eq:splitoffI}, to give a morphism from $B_I$ to $\BS(\un{r}) B_{I \cap A_{>\notme}}$ for any reduced expression
$\un{r}$ for $r$ (which we draw in purple). Then $B_{I \cap A_{>\notme}}$ can be pulled across the rest of the diagram.
\begin{equation} \label{style3} e_{I,X,Y} = \ig{1}{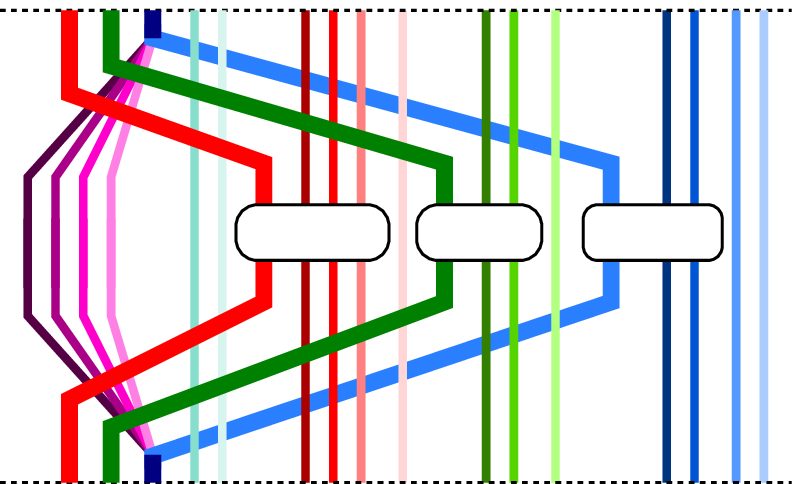} \end{equation}

At this point, the following things are actually quite straightforward exercises, mostly following from distant sliding relations and the work done in \S\ref{subsec-idempsfinite}. 

\begin{prop} \label{prop:easystuffnow} These definitions of $e_{I,X,Y}$ all agree, and $e_{I,X,Y}$ is an idempotent. These idempotents satisfy recursive formulas which categorify \eqref{initiate} and \eqref{grow}. Consequently, the image of $e_{I,X,Y}$, which we denote $\BY_{I,X,Y}$, has symbol $\by_{I,X,Y}$ in the Grothendieck group. Finally, funky trivalent vertices and an analog of \eqref{eq:whokjdecomp} can be used to categorify \eqref{godown}. \end{prop}

%%%%%%%%%%%%%%%%%%%%%%%%%
\section{Tensoring with the longest element: Gaussian elimination}
\label{sec-tensorBIGE}
%%%%%%%%%%%%%%%%%%%%%%%%%

We work in the affine setting, with $I \subsetneq S$ and $J = \tau(I)$. We aim to apply Gaussian elimination to the complex $B_I \FC$ until it has the form $\NC_I$. That is, we must
eliminate all summands of $B_I B_X$ when $J \nsubseteq X$, and all summands except the top summand when $J \subset X$. There are many procedures by which one could accomplish this: we
present here the one which seemed easiest to explain.

In \S\ref{subsec-GEobjects} we outline which summands we will eliminate and in which order. In \S\ref{subsec-GEmorphisms} we will work carefully with morphisms and examine the
differentials, ensuring that our iterated Gaussian eliminations are valid, and keeping track of what remains.

This entire chapter should be considered as the proof of the part of Theorem \ref{thm:tensorBI} which says that $B_I \FC \cong \NC_I$.

%===========
\subsection{Working with objects: a walkthrough}
\label{subsec-GEobjects}
%===========

We will formally state our algorithm for Gaussian elimination later, but let us walk carefully through the process first.

Let $A$ be a $\tau$-component of $I$. Let us assume for sake of an example that $A = \{1, 2, 3, 4\}$ so that $I \cap A = \{1, 2, 3\}$ and $J \cap A = \{2, 3, 4\}$. Recall from \S\ref{subsec-idempsfinite} the direct sum decompositions
\begin{subequations} \label{eq:BIB3andor4}
\begin{equation} \label{eq:BIB3} B_I B_3  \cong B_I(-1) \oplus B_I(+1), \end{equation}
\begin{equation} \label{eq:BIB4} B_I B_4 \cong B_{w_I h_4}, \end{equation}
and
\begin{equation} \label{eq:BIB43} B_I B_4 B_3 \cong B_{w_I h_{43}} \oplus B_I. \end{equation}
\end{subequations}

Our goal is to examine the product $B_I B_Z$ for various $Z \subsetneq S$. It turns out that the behavior of $B_I B_Z$ depends on $Z \cap \{3,4\}$, while grouping the summands $B_I B_Z$
into blocks based on the set $X = Z \cap (S \setminus \{3,4\})$ is a useful organizational tool. So for that purpose, let us introduce some notation.

When $X \subset S$ and $3, 4 \notin X$, we let $X4$ denote $X \cup \{4\}$, and $3X$ denote $X \cup \{3\}$. Earlier in this paper we might have denoted $X \cup \{3\}$ by $X3$, but we
will now pay close attention to the anticyclic order. Since $4 \notin 3X$, there is an anticyclic order where $3$ is first and the rest of $X$ follows. Similarly, since $3 \notin X4$
there is an anticyclic order where $4$ is last.

In an anticyclic order on the set $X \cup \{3,4\}$ (when this is a proper subset of $S$), it must be the case that $4$ precedes $3$, but it may not be possible to place these indices at the start or end of the order. However, for any given $\notme \notin X \cup \{3,4\}$, we can write $X = X_1 \sqcup X_2$ where $x_1 > 4 > 3 > x_2$ for any $x_1 \in X_1$ and $x_2 \in X_2$. We will denote $X \cup \{3,4\}$ as $X_1 43 X_2$.

Together, the four subsets $X$, $3X$, $X4$, and $X_1 43 X_2$ will be called the \emph{block} associated to $X$.

Now we examine the analogs of \eqref{eq:BIB3andor4} for products $B_I B_Z$, with $Z$ in the block associated to $X$. We claim that
\begin{subequations} \label{eq:BIBX3andor4}
\begin{equation} \label{eq:BIB3X} B_I B_{3X} \cong B_I B_X(-1) \oplus B_I B_X(+1). \end{equation}
This is now evident since $B_{3X} \cong B_3 B_X$, so we can apply \eqref{eq:BIB3}. Meanwhile, we have the exciting decomposition
\begin{equation} \label{eq:BIBX4} B_I B_{X4} \cong B_I B_{X4}. \end{equation}
We do not claim, as in \eqref{eq:BIB4}, that $B_I B_{X4}$ is indecomposable, only that we do not choose to decompose it further at this time. For that matter, we also do not claim that either summand $B_I B_X(\pm 1)$ of the RHS of \eqref{eq:BIB3X} is indecomposable.

It takes slightly more work to explain the next decomposition. We aim to define an object $\BY_{I,X_1 43 X_2,43}$ in the Karoubi envelope of $\Diag$ such that
\begin{equation} \label{eq:BIBX43} B_I B_{X_1 43 X_2} \cong \BY_{I,X_1 43 X_2,43} \oplus B_I B_X. \end{equation}
\end{subequations}
Again, neither summand need be indecomposable.

For reasons of sanity, let us shorten $\BY_{I,X_1 43 X_2,43}$ to $\BY_{X,43}$.

First suppose for simplicity that $X_1$ is empty, which is possible if and only if $5 \notin X$. Then $B_{X_1 43 X_2} \cong B_4 B_3 B_{X_2}$, and \eqref{eq:BIB43} implies the decomposition
\begin{equation} B_I B_{X_1 43 X_2} \cong B_{w_I h_{43}} B_{X_2} \oplus B_I B_{X_2}. \end{equation}
Setting $\BY_{X,43} = B_{w_I h_{43}} B_{X_2}$ gives the desired decomposition.

When $X_1$ is nonempty it is more annoying to write down $\BY_{X,43}$ formulaically, though we will describe it as the image of a specific idempotent. Note that $w_I$ and $h_{X_1}$ need not commute! For example, one might have $X_1 = \{1 > 0 > n-1 > \ldots > 5\}$. However, $s_3$ is guaranteed to commute with any simple reflection in $X_1$, and the decomposition \eqref{eq:BIB43} really only relies upon the fact that $3 \in I$.

This allows us to construct the following idempotent decomposition of $B_I B_{X_1 43 X_2}$.
\begin{equation} \label{eq:twoidempotents} {
\labellist
\small\hair 2pt
 \pinlabel {$=$} [ ] at 117 227
 \pinlabel {$=$} [ ] at 293 229
 \pinlabel {$+ (-1)$} [ ] at 154 72
\endlabellist
\centering
\ig{1}{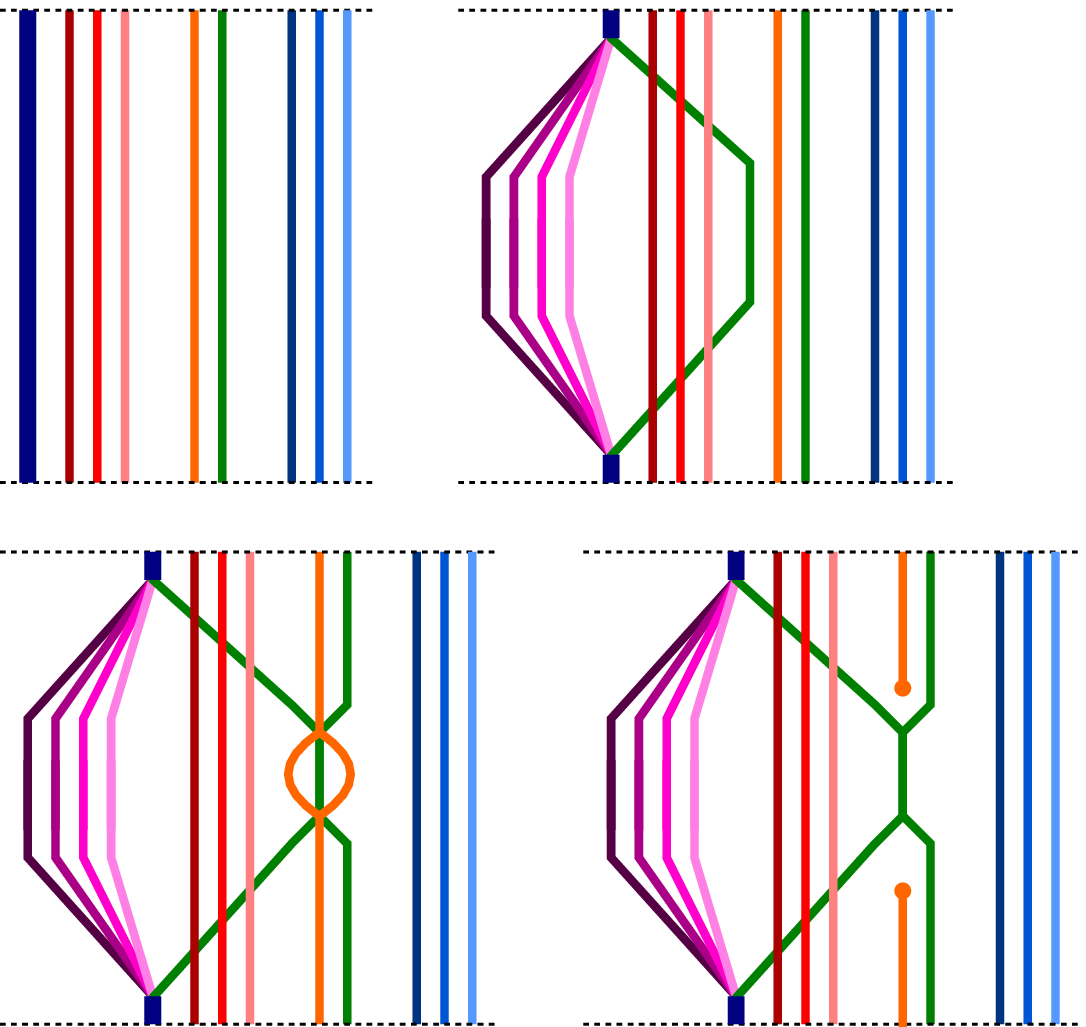}
} \end{equation}
In this picture, $B_I$ is split into a reduced expression ending in $3$ (most of which is purple, while $3$ is green), the $3$-colored strand is pulled through $X_1$ (which is different shades of red), used to decompose $B_3 B_4 B_3 \cong B_{343} \oplus B_3$ in the usual way ($4$ is orange), and the $3$-colored strand is returned to its source in $B_I$. Meanwhile $X_2$ (different shades of blue) hangs out on the side doing nothing. We refer to the idempotents on the bottom of \eqref{eq:twoidempotents} as $\idemone$ and $\idemtwo$ respectively. We define $\BY_{X,43}$ as the image of $\idemone$.

Using associativity, thick associativity \eqref{eq:thickassoc}, and the definition of the thick trivalent vertex, we can rewrite $\idemtwo$ as
\begin{equation} \label{eq:idemtwofactor} (-1) \ig{1}{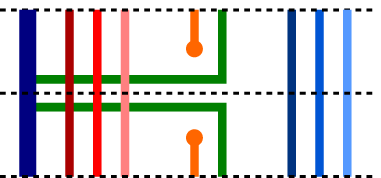}, \end{equation}
so that $\idemtwo$ factors through the summand $B_I B_{X_1 X_2} = B_I B_X$, as desired.

What do the direct sum decompositions \eqref{eq:BIBX3andor4} allow us to do? We have seen that, for various $Z$ in the block associated to $X$, the direct summand $B_I B_X$ will appear
(possibly with shift) in $B_I B_Z$. Our first simultaneous Gaussian elimination that we apply to $B_I \FC$ will remove all summands of the form $B_I B_X$ from this block. In fact, we
can do this for all blocks at once, simultaneously! What remains are the summands $B_I B_{X4}$ from \eqref{eq:BIBX4} and $\BY_{X,43}$ from \eqref{eq:BIBX43}, both of which have a $4$
in them somewhere. Thus we call this process \emph{eliminating the absence of $4$}.

Let us enumerate the summands we will eliminate, within the block of $X$. Suppose that $3X \in P_k$, so that $B_{3X}(k)$ appears in homological degree $k$ of $\FC$. Then $B_I B_{3X}(k) \cong B_I B_X(k-1) \oplus B_I B_X(k+1)$ appears in homological degree $k$ in $B_I \FC$. We refer to these two summands as the images of two idempotents, $e_-$ and $e_+$ respectively. Whenever $3X \in P_k$ it is always the case that $X \in P_{k-1}$, so that $B_I B_X(k-1)$ appears in homological degree $k-1$ and can cancel against the image of $e_-$. Conversely, almost every copy of $B_I B_X$ is cancelled by the image of $e_-$ inside $B_I B_{3X}$, with one exception! If $\ell$ is the maximal homological degree in which $X$ appears in $P = \sqcup_k P_k$, then $B_I B_X(\ell)$ in degree $\ell$ has not yet been cancelled, because $3X \notin P_{\ell + 1}$. Note that $\ell = n - 1 - |X| \ge 0$.

Whenever $3X \in P_k$ it is almost always the case that $X_1 43 X_2 \in P_{k+1}$, so that $B_I B_{X_1 43 X_2}(k+1)$ appears in homological degree $k+1$. The direct summand $B_I B_X(k+1) \sumset B_I B_{X_1 43 X_2}(k+1)$, the image of $\idemtwo$, can cancel against the image of $e_+$. Again, there is a unique homological degree $\ell - 1$  where $3X \in P_{\ell-1}$ but $X_1 43 X_2 \notin P_{\ell}$, which is again\footnote{This is even true when $|X| = n-2$ so that $X \cup \{3,4\} = S$. Then $\ell = 1$. While $B_{X_1 43 X_2}$ does not make sense, it is still true that $3X \in P_{\ell-1}$ but $X \cup \{3,4\} \notin P_{\ell}$,
and everything we say applies equally to this case.} the maximal degree in which $3X$ appears in $P$. The number $\ell$ matches the $\ell$ from the previous paragraph.

Now there are two copies of $B_I B_X(\ell)$ which have not been cancelled yet, the image of $e_+$ in homological degree $\ell-1$, and the vanilla summand $B_I B_X(\ell)$ in homological degree $\ell$, so we will cancel them against each other!

We note quickly that the summands of the form $B_I B_{X4}$ have been entirely ignored in this process.

Thus we have our family of Gaussian eliminations. For each $X$ with $3,4 \notin X$: \begin{itemize}
\item For each $k$ with $3X \in P_k$ we cancel the image of $e_-$ in $B_I B_{3X}(k)$ against $B_I B_X(k-1)$ in degree $k-1$.
\item For each $k$ with $3X \in P_k$ we cancel the image of $e_+$ in $B_I B_{3X}(k)$ against a summand in degree $k+1$, which is either \begin{itemize}
	\item $B_I B_X(k+1)$ when $k+1 = \ell$ is the maximal degree in which $X$ appears, or
	\item the image of $\idemtwo$ in $B_I B_{X_1 43 X_2}(k+1)$ when $k+1 < \ell$.
	\end{itemize}
\end{itemize}
In the following section we prove that this family of isomorphisms, for all such $X$ at once, can be simultaneously Gaussian eliminated, because the family is independent (see Definition \ref{defn:independent}). Moreover, this family of isomorphisms produces no nonzero zigzag terms between the survivors of the elimination.

Having eliminated the absence of $4$, we can now eliminate the absence of $3$. Suppose that $2, 3, 4 \notin Y$. Both $Y$ and $Y \cup \{2\}$ are valid choices for $X$ with $3, 4 \notin X$. Thus the survivors of the previous elimination will have the form
\begin{equation} \label{eq:Yblock} B_I B_{Y4}, \quad B_I B_{2Y4}, \quad \BY_{Y,43}, \quad \BY_{Y2,43}. \end{equation}
We think of these four terms as being the \emph{block} associated to $Y$, completely analogous to the four terms $X$, $3X$, $X4$, and $X_1 43 X_2$ in the previous argument. Now the relevant feature is the intersection with $\{2,3\}$ rather than with $\{3,4\}$, and we think that $4$ is always present.

Our process will be entirely analogous. We start by observing that the second term $B_I B_{2Y4}$ (by which we mean, the second kind of survivor in \eqref{eq:Yblock}) will split into two copies of $B_I B_{Y4}$ with grading shifts $-1$ and $+1$, corresponding to idempotents $e_-$ and $e_+$. This is because $B_{2Y4} \cong B_2 B_{Y4}$, and $B_I B_2 \cong B_I(-1) \oplus B_I(+1)$ since $2 \in I$. The image of $e_-$ will cancel against some copy of the first term $B_I B_{Y4}$ in the previous homological degree. The image of $e_+$ will cancel against  either against a summand $B_I B_{Y4}$ inside the fourth term $\BY_{Y2,43}$, or in the final homological degree, it will
cancel against the only copy of the first term $B_I B_{Y4}$ which hasn't already been cancelled. What survives will be the terms of the form $\BY_{Y,43}$ (which are ignored entirely in
this process) and the other summand of $\BY_{Y2,43}$, which would be denoted $\BY_{I,Y432,432}$ in \S\ref{subsec-idempsaffine}, and which we shorten to $\BY_{Y,432}$ here.

In the next step we can eliminate the absence of $2$, treating analogously the four survivors $\BY_{Z,43}$, $\BY_{1Z,43}$, $\BY_{Z,432}$, and $\BY_{Z1,432}$ for $1, 2, 3, 4 \notin Z$.
Again crucially, $\BY_{1Z,43} \cong \BY_{Z,43}(-1) \oplus B_{Z,43}(+1)$, coming from $B_I B_1 \cong B_I(-1) \oplus B_I(+1)$, which holds since $1 \in I$.

The pattern stops here though. We can not eliminate the absence of $1$ because $0 \notin I$, so $B_I B_0$ is not just two shifted copies of $B_I$. Thus this process will only eliminate
the absence of $J \cap A = \{2,3,4\}$, as desired.

Next we should eliminate the absence of $j$, for some other $j \in J \cap A$ maximal in its $\tau$-component (i.e. initiating the $\tau$-suffix). The arguments are entirely similar, only
we must work from the start with funky summands (e.g. $\BY_{Xj,432}$ and $\BY_{X(j-1),4321}$ and others, for $1,2, 3, 4, j-1, j \notin X$). This makes the notation much more painful, but
the arguments are almost identical. The precise definition of all the idempotents involved was the focus of \S\ref{subsec-idempsaffine}, which we assume the reader has read henceforth. We
continue in this fashion, eliminating the absence of a $\tau$-suffix at each step, until one has eliminated the absence of every $j \in J$. The result will be a complex with the same objects as $\NC_I$.

%===========
\subsection{Working with morphisms}
\label{subsec-GEmorphisms}
%===========

The reader should recall the definition of independent isomorphisms for Gaussian elimination, and of zigzag morphisms, from Definition \ref{defn:independent}. Our aim in this section is
to prove that the Gaussian elimination procedure from the previous section can be achieved with simultaneous elimination of independent isomorphisms, and that no zigzag morphisms
interfere with the surviving differentials.

Let us continue the example where $A = \{1,2,3,4\}$ that was worked in previously. Our first step was to eliminate the absence of $4$. Until otherwise stated we assume that $3,4 \notin X$. As a reminder, we wish to eliminate the following.
\begin{itemize}
\item For each $k$ with $3X \in P_k$ we cancel the image of $e_-$ in $B_I B_{3X}(k)$ against $B_I B_X(k-1)$ in degree $k-1$.
\item For each $k$ with $3X \in P_k$ we cancel the image of $e_+$ in $B_I B_{3X}(k)$ against a summand in degree $k+1$, which is either \begin{itemize}
	\item $B_I B_X(k+1)$ when $k = \ell$ is the maximal degree in which $3X$ appears, or
	\item the image of $\idemtwo$ in $B_I B_{X_1 43 X_2}(k+1)$ when $k \ne \ell$.
	\end{itemize}
\end{itemize}
To quickly state that a term $B_I B_X(\ell+1)$ appears in its final homological degree, we might refer to this copy of $X$ as $X_{\last}$. Given a direct summand $F$ of some $B_I B_Z$, we say that $Z$ is the \emph{indexset} of $F$.

The idempotents $e_-$ and $e_+$ for the direct sum decomposition $B_I B_{3X} \cong B_I B_X(-1) \oplus B_I B_X(+1)$ are obtained by taking the idempotents from \eqref{eq:thicksdecomp} and tensoring with $\id_{B_X}$. These idempotents are factored, with the projection and inclusion maps being denoted $p_+, p_-, i_+,i_-$. 
It is worth mentioning that these idempotents are not canonical, but that $p_-$ and $i_+$ are canonical (up to scalar). The fact that $p_- i_+ = 0$ will play a significant role below. Similarly, we can factor the idempotent $\idemtwo$ as in \eqref{eq:idemtwofactor}, and refer to the projection and inclusion as $p_{\downarrow}$ and $i_{\downarrow}$ (both are canonical up to scalar), and we include the sign in $p_{\downarrow}$.

First, let us confirm that the differentials we plan to contract are actually isomorphisms.

The differential $B_I B_X(k-1) \to B_I B_{3X}(k)$ is a (signed) $3$-colored startdot. By postcomposing with $p_-$ we get the induced map between summands,
\begin{equation} \label{itseasy} {
\labellist
\small\hair 2pt
 \pinlabel {$=$} [ ] at 68 26
\endlabellist
\centering
\ig{1}{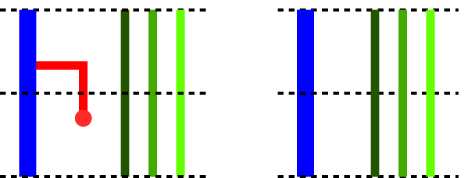}
} \end{equation}
which is equal to the identity map by \eqref{eq:thickunit}. There might be a sign involved in the differential if one chooses a different reduced expression for $B_{3X}$, although one still gets an isomorphism.

Similarly, the differential $B_I B_{3 X_{\last}}(\ell) \to B_I B_{X_{\last}}(\ell+1)$ is a (signed) $3$-colored enddot. By precomposing with $i_+$ we get the induced map between summands, and this is exactly the vertical flip of \eqref{itseasy}.

The differential $B_I B_{3X}(k) \to B_I B_{X_1 43 X_2}(k+1)$ is a (signed) $4$-colored startdot, with a change in reduced expression. By precomposing with $i_+$ and postcomposing with $p_{\downarrow}$  we get the induced map between summands
\begin{equation} \label{itseasy2}{
\labellist
\small\hair 2pt
 \pinlabel {$\pm$} [ ] at -15 39
 \pinlabel {$= \mp$} [ ] at 103 38
\endlabellist
\centering
\ig{1}{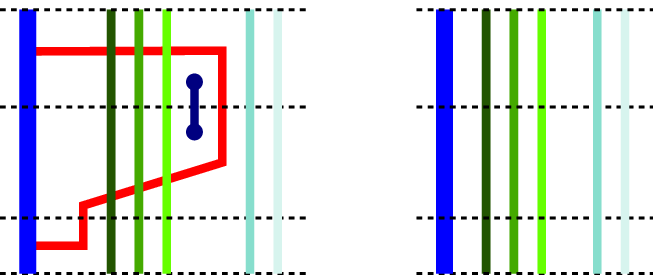}
} \end{equation}
which is again (up to sign) just the identity map.

Next, we begin to examine the possible zigzags. Recall that a non-repeating length $d$ zigzag has the form $S \to F'_1 \from F_1 \to F'_2 \from F_2 \to \ldots \to F'_d \from F_d \to T$,
where $F_i \to F_i'$ is one of the isomorphisms being contracted, and $F_i' \from F_i$ is its inverse. The other maps $F_i \to F'_j$ are differentials in the complex, and the $\{F_i\}$
are assumed to be distinct. In order for there to be a differential from $F_i$ to $F'_j$, the indexsets of $F_i$ and $F'_j$ must differ by a single index.

Let us just shorten $X_1 43 X_2$ to $X43$ in the discussion below.

Let us consider all possible differentials $F_1 \to F'_2$. Let $Z$ be the indexset of $F_1$ and $Z'$ the indexset of $F'_2$. One constraint is that $Z$ and $Z'$ differ by a single index. Another is that $Z$ is either $3X$ or $X$, and $Z'$ is either $3Y$, $Y43$, or $Y_{\last}$, for some $X, Y$ with $3, 4 \notin X, Y$, possibly with $X = Y$. The possibilities which remain are: \begin{itemize}
\item $Z = 3X$, $F_1$ is the image of $e_+$, $Z' = X43$, $F'_2$ is the image of $\idemtwo$. But then $F'_1$ must also have indexset $X43$ (it can't be have indexset $X_{\last}$ because $X_{\last}$ and $X43$ can not occur in the same homological degree). This is forbidden since $F'_1 = F'_2$ is not allowed.
\item $Z = 3X$, $F_1$ is the image of $e_+$, $Z' = X_{\last}$, $F'_2 = B_I B_{X_{\last}}$. By the same argument, $F'_1 = F'_2$ which is forbidden.
\item $Z = 3X$, $F_1$ is the image of $e_+$, $Z' = 3Y$ for some $Y$ differing from $X$ by a single index, and $F'_2$ is the image of $e_-$ in $B_I B_{3Y}$. This is permitted. But then the differential from $F_1$ to $F'_2$ is a signed dot, precomposed with $i_+$ and postcomposed with $p_-$.
\begin{equation} \label{itseasy3} {
\labellist
\small\hair 2pt
 \pinlabel {or} [ ] at 104 26
\endlabellist
\centering
\ig{1}{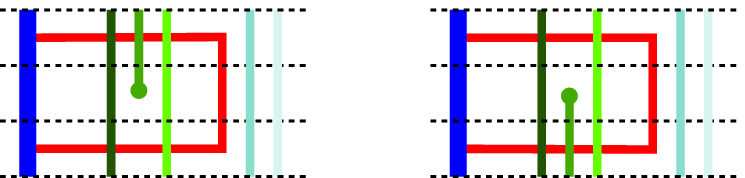}
} \end{equation}
This is zero by \eqref{eq:thickneedle}. Even though $i_+$ and $p_-$ are inclusions/projections for different idempotent decompositions (one for $B_I B_{3X}$, the other for $B_I B_{3Y}$), they are still orthogonal (the extra signed dot commutes past and plays no role).
\item $Z = X$, $F_1$ is $B_I B_X$, $Z' = 3X$, and $F'_2$ is the image of $e_-$. But this is forbidden, because $F'_1 = F'_2$.
\item $Z = X$, $F_1$ is $B_I B_X$, $Z' = Y_{\last}$ for some $Y$ differing from $X$ by a single index, and $F'_2 = B_I B_{Y_{\last}}$. This is permitted, and the differential is just a signed dot. Finally, we have found a nonzero differential.
\end{itemize}
Thus, any nonzero map $F_1 \to F'_2$ must have the form $B_I B_X \to B_I B_{Y_{\last}}$, with $F'_1$ the image of $e_-$ in $B_I B_{X3}$, and $F_2$ the image of $e_+$ in $B_I B_{Y3}$.

As a consequence, two consecutive forward maps $F_i \to F'_{i+1}$ and $F_{i+1} \to F'_{i+2}$ between contracting objects in a zigzag can not both be nonzero. It can not be the case that
both $F_{i+1}$ has indexset $Y3$ for some $Y$, and $X$ for some $X$. In particular, no zigzag can have length $3$ or more, lest it have two consecutive forward maps. Zigzags of length $2$
are very constrained.

The isomorphisms are independent if every non-repeating zigzag of length $d \ge 1$, whose source and target are contracting terms (i.e. $S = F_0$ and $T = F_{d+1}$, still distinct), is
zero. But any such zigzag has two consecutive forward maps $S \to F'_1$ and $F_1 \to T$, so it is zero. Thus the isomorphisms are independent.

Now we ask whether there are any nonzero zigzags whose source and target survive the elimination. In this case $S = B_I B_{U4}$ or $S$ is the image of $\idemone$ in $B_I B_{U43}$,
and similarly $T = B_I B_{V4}$ or $T$ is the image of $\idemone$ in $B_I B_{V43}$, for some $U, V$ with $3, 4 \notin U, V$, possibly with $U = V$. We know that all zigzags of
length $3$ or more are zero. If there were a zigzag of length $2$ then $F_2$ has indexset $3Y$ for some $Y$, and $F'_2$ has indexset $Y_{\last}$. For there to be a nonzero differential
from $F_2$ to $T$ we need either $V4$ or $V43$ to differ by one index from $3Y$. The only possibility is that $V = Y$ and $T$ has indexset $Y43$. But this is impossible since $Y_{\last}$
and $Y43$ can not appear in the same homological degree.

Consider a zigzag of length $1$. Let $Z$ be the indexset of $F_1$; we know that either $Z = X$ or $Z = 3X$ for some $X$ with $3, 4 \notin X$. Both $S$ and $T$ have a $4$ somewhere, so the zigzag must remove the $4$ (either in the differential $S \to F'_1$ or in the contraction $F'_1 \from F_1$) and then add it back (in the differential $F_1 \to T$). In particular, $X = V$. The possibilities are:
\begin{itemize} \item $F_1$ has indexset $3V$, and $T$ has indexset $V43$. The differential from $F_1$ to $T$ is a $4$-colored startdot, precomposed with $i_+$, and postcomposed with $\idemone$.
\begin{equation} \label{itseasy4} \ig{1}{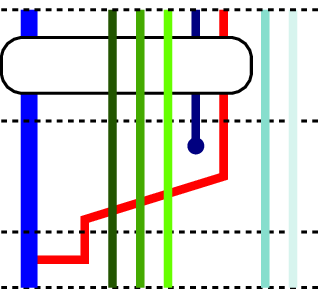} \end{equation}
But this is zero by \eqref{eq:emortho}. Said another way, the composition of $i_+$ and the $4$-colored startdot is exactly $i_{\downarrow}$, which is orthogonal to $\idemone$.
\item $F_1$ has indexset $V$, and $T$ has indexset $V4$. Then $F'_1$ has indexset $3V$, so $S$ has indexset $V43$. Now the differential from $S$ to $F'_1$ is the $4$-colored enddot, precomposed with $\idemone$ and postcomposed with $p_-$. This is just \eqref{itseasy4} flipped upside-down, and vanishes for the same reason.
\end{itemize}
Thus all zigzags between surviving terms are zero.

To summarize, we have analyzed the possible zigzags. This analysis involved a combinatorial aspect, which limited the indexsets based on the choice of our contracting isomorphisms, and a
diagrammatic part, where we determined that certain idempotents are orthogonal to certain maps. The result was that we can simultaneously Gaussian eliminate the absence of $4$ from $B_I
\FC$. What remains is a complex built from terms $B_I B_{X4}$ and $\BY_{X,43}$, where the differentials between surviving terms remain ``unchanged'': they are signed dots, pre- and post-composed with the appropriate idempotents.

Now we wish to simultaneously Gaussian eliminate the absence of $3$. As in \S\ref{subsec-GEobjects}, the objects $B_I B_{Y4}$, $B_I B_{2Y4}$, $\BY_{Y,43}$, and $\BY_{Y2,43}$, now play the
role that $B_I B_X$, $B_I B_{3X}$, $B_I B_{X4}$, and $B_I B_{X43}$ used to play. Said another way, all our indexsets now contain $4$, and our combinatorial analysis now checks the intersection of the indexset with $\{2,3\}$ rather than $\{3,4\}$. Thus the combinatorial analysis is completely analogous to the previous argument.

The diagrammatic analysis is also similar. For example, we need to show that the composition $B_I B_{X4}(+1) \to B_I B_{2X4} \to B_I B_{2Y4} \to B_I B_{Y4}(-1)$ is zero. This composition is (up to sign)
\begin{equation} \label{itseasy5} {
\labellist
\small\hair 2pt
 \pinlabel {or} [ ] at 104 26
\endlabellist
\centering
\ig{1}{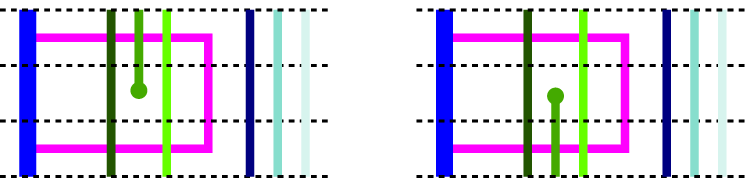}
}.\end{equation}
This is zero by \eqref{eq:thickneedle}. Effectively, this is the same calculation as \eqref{itseasy3}, except that the empty eyehole is colored $2$ insetad of $3$, and there is an extra $4$-colored strand hanging out to the right.

Similarly, we need to show that the composition $B_I B_{X4}(+1) \to B_I B_{2X4} \to \BY_{X2,43} \to \BY_{X,432}$ is zero.
\begin{equation} \label{itseasy6} \ig{1}{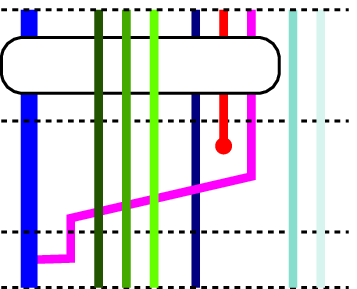} \end{equation}
We have to use a more complicated idempotent (absorbing $432$ rather than just $43$), but the result is still zero by \eqref{eq:emortho}. Again, the composition of $i_+$ with the signed dot is the inclusion $i_{\downarrow}$, so it is orthogonal to $\idemone$. 

Having done these examples, we claim that our work is effectively done. Using the factored idempotent decompositions developed in \S\ref{subsec-idempsfinite} and
\S\ref{subsec-idempsaffine}, one can repeat the calculation for any index one hopes to eliminate the absence of, whether it initiates a new suffix or grows a suffix, and regardless of
which other indices have already been eliminated. One must merely use the idempotents $e_{I,X,Y}$ and the corresponding funky trivalent vertices, rather than the simpler idempotents
(special cases) above. The proof is almost identical, only the notation is worse and the diagrammatics are more annoying to draw. The reader will hopefully forgive us for not formally
writing up the tedious general case.

In conclusion, we have proven that $B_I \FC$ is homotopy equivalent to $\NC_I$, by confirming that the differentials after Gaussian elimination are as simple as possible: they agree with the differentials before elimination, pre- and post-composed with the appropriate inclusions and projections from the surviving direct summands.

%%%%%%%%%%%%%%%%%%%%%%%%%
\section{Left versus right}
\label{sec-leftvsright}
%%%%%%%%%%%%%%%%%%%%%%%%%

Theorem \ref{thm:tensorBI} states that both $B_I \FC$ and $\FC B_{\tau(I)}$ are homotopy equivalent to $\NC_I$. In the previous chapter we performed the explicit Gaussian elimination
which gave the contraction $B_I \FC \simto \NC_I$. In this short chapter we ask about $\FC B_{\tau(I)}$. We continue to write $J = \tau(I)$.

Consider the following symmetry $\Theta$ of the category $\Diag_{\aff}$, a composition of two easier symmetries. The first is the color reflection $\si$, and the second (monoidally
contravariant) flips all diagrams left-to-right. Composing these two symmetries, one obtains a symmetry $\Theta$ which preserves cyclically ordered sets, with $\Theta(B_X) \cong
B_{\si(X)}$. One also has $\Theta(B_I) \cong B_{\si(I)}$. Applying $\Theta$ to $\FC$ yields a pseudocomplex isomorphic to $\FC$, but not quite equal on the nose. The difference is that
the signs on the differentials in $\FC$ are determined by counting strands to the left of the dot, while the signs in $\Theta(\FC)$ are determined by counting strands to the right of the
dot. This changes the sign on certain differentials: the maps $B_X \to B_Y$ when the pair $\{|X|,|Y|\}$ is equal to either $\{1,2\}$ or $\{3,4\}$ or $\{5,6\}$, etcetera. One choice of
isomorphism $\FC \to \Theta(\FC)$ will send $B_X \to B_X$ with a sign which is $-1$ when $|X| \{2,3,6,7, 10, 11, \ldots\}$ and is $+1$ otherwise.

The point then is that, applying $\Theta$ to all the calculations of the preceding chapters, as they apply to $B_{\si(J)} \FC$, we obtain an analgous Gaussian elimination for the complex
$\FC B_J$. The independence of the contracting isomorphisms and the lack of nonzero zigzags between surviving terms implies that the differential on the minimal complex is induced from
the differential on $\FC B_J$ by pre- and post-composing with the appropriate inclusions and projections. Thus, if we don't want to keep track of the sign flip in the isomorphism between $\FC$ and $\Theta(\FC)$, we don't have to; choosing the normal signs on differentials in $\FC B_J$, the same result gives the normal sign rule for the minimal complex.

Thus, for free by symmetry, we can deduce that $\FC B_J$ is homotopy equivalent to the complex $\MC_J$ defined below (c.f. Definition \ref{defn:N}).

\begin{defn} Fix $I \subsetneq S = S_{\aff}$. For each $k \in \Z$ let $M^J_k$ be the subset of $P_k$ (see Definition \ref{def:Pk}) consisting of those $Y \in P_k$ where $I \subset Y$. \end{defn}

\begin{defn} \label{defn:M}  Let $\MC_J$ denote the following pseudocomplex. The chain object in homological degree $k$ is
\begin{equation} \MC_J^k = \bigoplus_{Y \in M^J_k} B_{h_Y w_J}(k). \end{equation}
The differential, restricted to a summand $B_{h_Yw_J}(k)$ in homological degree $k$, is the sum
\begin{equation} \label{def:dM} d = \bigoplus_{\substack{i \notin X\\ Xi \in M^J_{k+1}}} {}^J \startdotsign_{X}^{Xi} \oplus \bigoplus_{\substack{j \in X\\X \setminus j \in M^J_{k+1}}} {}^J \finaldotsign_{X}^{X \setminus j}. \end{equation}
Recall from Definition \ref{defn:dotsign} that we have signed dot maps $\startdotsign_X^Y$ or $\finaldotsign_X^Y$ from $B_X$ to $B_Y$ whenever $Y$ differs from $X$ by a single index. The maps ${}^J \startdotsign_X^Y$ and ${}^J \finaldotsign_X^Y$ in \eqref{defn:M} are obtained by taking the corresponding signed dot $B_X \to B_Y$, tensoring on the right with the identity map of $B_J$, precomposing with the inclusion $B_{h_X w_J} \into B_X B_J$, and postcomposing with the projection $B_Y B_J \onto B_{h_Y w_J}$.
\end{defn}

Thus, to finish the proof of Theorem \ref{thm:tensorBI}, it remains to construct an isomorphism of pseudocomplexes between $\NC_I$ and $\MC_J$. This is the content of the current chapter. The result is not terribly difficult, but we feel it is worthwhile being as explicit as possible.

\subsection{The crossover isomorphism in the finite case}

It was proven in Lemma \ref{lem:descent} that there is a bijection between $\{X \subsetneq S \mid J \subset X\}$ and $\{Y \subsetneq S \mid I \subset Y\}$ such that, whenever $X$ and $Y$
match under this bijection, we have $w_I h_X = h_Y w_J$. Thus the summands in $\NC_I$ and those in $\MC_J$ are in natural bijection. In our description of $\NC_I$ we described $B_{w_I
h_X}$ as the image of a particular idempotent endomorphism of $B_I B_X$, and this idempotent is part of the definition of ${}^I \startdotsign_X^Y$. Meanwhile, $B_{h_Y w_J}$ is defined as
the image of an idempotent endomorphism of $B_Y B_J$. Our first task is to give explicitly the isomorphism between the image of these two idempotents.

We begin in the finite setting: let $T = \{1, \ldots, N\}$ and $I = \{1, \ldots, N-1\}$ and $J = \{2, \ldots, N\}$. There are two choices of $X$: $X = T$ and $X = J$. The corresponding choices of $Y$ are $Y = T$ and $Y = I$, respectively. Consider the following maps $B_I B_X \to B_Y B_J$, which we call \emph{crossover maps}. For sake of simplicity we draw the case $N = 4$. (Colors: $\{1,2,3,4\}$ is blue, red, green, purple; $I$ is thick green and $J$ is thick purple, $T$ is thick blue, $I \cap J$ is thick red.)
\begin{equation} \label{crossoverT} X = Y = T: \quad \ig{1}{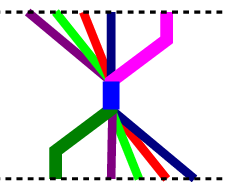} \end{equation}
\begin{equation} \label{crossoverJ} X = J, Y = I: \quad \ig{1}{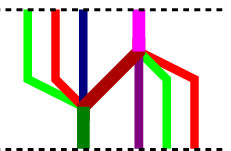} \end{equation}
Note that the colors on the bottom and top of \eqref{crossoverT} are the same, reflecting the fact that $X=Y$, while the colors on the bottom and top of \eqref{crossoverJ} are shifted.

\begin{lemma} \label{lem:crossover} The crossover maps descend to an isomorphisms $B_{w_I h_X} \to B_{h_Y w_J}$ after pre- and post-composing with the appropriate idempotents. \end{lemma} 

\begin{proof} We first show two auxiliary facts.

We have shown in Lemma \ref{lem:finitepleth} that $b_{w_I h_X}$ is smooth, and we showed in Corollary \ref{cor:ekiswhok} that the image of $e_k$ categorifies the Kazhdan-Lusztig basis element $b_{\who_k}$. When $X = T$, $w_I h_X = w_T = \who_N$, so the image of $e_N$ should be isomorphic to $B_{w_T}$. In particular, one expects the following equality.
\begin{equation} \label{XTeN} {
\labellist
\small\hair 2pt
 \pinlabel {$=$} [ ] at 73 31
 \pinlabel {$N$} [ ] at 34 12
\endlabellist
\centering
\ig{1}{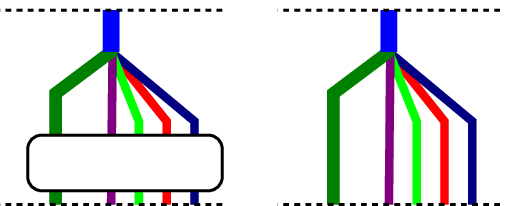}
} \end{equation}
One can prove \eqref{XTeN} using the recursive definition \eqref{eq:emrecursion}: if at each step the other idempotent $\idemtwo$ is zero when plugged in to the diagram, then plugging in $e_N$ is the same as plugging in $e_1$, which is the identity map. But plugging in $\idemtwo$ gives a subdiagram which looks like
\begin{equation} \ig{1}{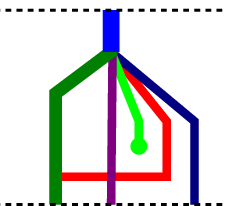}. \end{equation}
But this is zero, since the thick splitters and mergers are killed by pitchforks, c.f. \cite[Proposition 3.39]{EThick}.

When defining the idempotent $e_k$ it is only important that $I$ contains the elements $\{N+1-k, \ldots, N-1\}$. One should expect that defining $e_k$ with respect to a large enough parabolic subgroup $I' \subset I$, and then merging $I'$ into $I$, should give the same result as defining $e_k$ with respect to $I$.
\begin{equation} \label{ekslide} {
\labellist
\small\hair 2pt
 \pinlabel {$=$} [ ] at 78 38
 \pinlabel {$k$} [ ] at 42 18
 \pinlabel {$k$} [ ] at 137 53
\endlabellist
\centering
\ig{1}{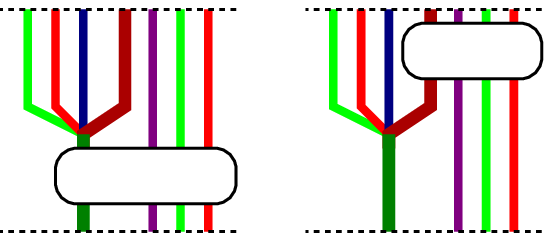}
} \end{equation}
Given the recursive definition of $e_k$, this will follow if thick trivalent vertices for colors in $I'$ can split across the splitter from $I$ to $I'$. This in turn follows directly from the definition of the thick trivalent vertex, see \eqref{thicktridefn}.

Using \eqref{XTeN} we deduce that the $X=T$ crossover from \eqref{crossoverT} is unchanged after pre- and post-composition with $e_N$. Using \eqref{XTeN} and \eqref{thicktridefn}, one
can deduce that the $X=J$ crossover from \eqref{crossoverJ} is unchanged after pre- and post-composition with $e_{N-1}$. Thus these maps between $B_I B_X$ and $B_Y B_J$ are zero on all
the direct summands of these objects except the top summands $B_{w_I h_X}$ and $B_{h_Y w_J}$. Identical statements can be made about the vertical flips of these morphsims, which give
maps in the other direction.

We claim that crossover is the inverse of the flipped crossover, giving isomorphisms between $B_{w_I h_X}$ and $B_{h_Y w_J}$. There are various ways to prove this directly and
diagrammatically. Let us resort to a cheaper trick.

In the algebraic category of Soergel bimodules, any indecomposable Soergel bimodule $B_w$ has a one-dimensional space of elements in the lowest degree $-\ell(w)$, and a particular element in this space called the \emph{1-tensor}, introduced in \cite{EKho} and discussed further in \cite{EWGr4sb}. To define this element, note that the Bott-Samelson bimodule $\BS(\un{w})$ has an element $1 \ot 1 \ot \cdots \ot 1$ called the 1-tensor, and this element lives in the top degree summand for degree reasons. The 1-tensors in Bott-Samelson bimodules are preserved by $2m$-valent vertices and enddots. Similarly, the thick splitters and mergers send the 1-tensor to the 1-tensor, effectively by definition. Thus two morphisms which are known to be colinear will be equal if they agree (and are non-zero) on the 1-tensor.

\begin{remark} This argument involving 1-tensors only works in the algebraic category of Soergel bimodules, but there is an alternate argument which works for the diagrammatics. Using
diagrammatic localization, as developed in \cite{EWGr4sb}, one can examine how a morphism between Bott-Samelsons acts on various standard submodules, and especially, on the unique
top standard submodule $R_w$. Again, $2m$-valent vertices and splitters and mergers will act as the identity on this submodule, and colinear morphisms can be measured by their action on
this summand. \end{remark}

By Soergel's Hom formula, we know that the space of degree zero endomorphisms of an indecomposable Soergel bimodule is one dimensional. Thus any endomorphism fixing the 1-tensor is the
identity map. Each of the thick splitters and mergers sends the 1-tensor to the 1-tensor, and thus so does the crossover, so the crossover and its flip must compose to the identity.
\end{proof}

\begin{lemma} The following equality holds, as morphisms from $B_{w_I h_T}$ to $B_{h_I w_J}$.
\begin{equation} \label{dotcrossover} {
\labellist
\small\hair 2pt
 \pinlabel {$=$} [ ] at 94 56
 \pinlabel {$N$} [ ] at 39 14
 \pinlabel {$N$} [ ] at 175 37
 \pinlabel {$N$} [ ] at 30 90
 \pinlabel {$N-1$} [ ] at 151 86
\endlabellist
\centering
\ig{1}{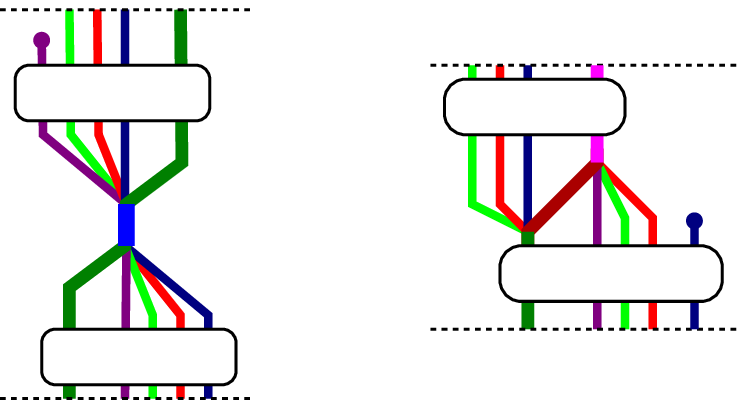}
} \end{equation}
\end{lemma}

\begin{remark} Because of \eqref{XTeN} most of these oval idempotents are redundant, except for the copy of $e_N$ in the bottom of the RHS. \end{remark}

\begin{proof} This is another result which can be proven directly with diagrammatics. Let us resort to the same cheap trick. Because both $w_I h_T$ and $h_I w_T$ are smooth, and
$\ell(w_I h_T) = \ell(h_I w_J) + 1$, the dimension of the space of degree one morphisms $B_{w_I h_T} \to B_{h_I w_T}$ is one. Thus these two maps are automatically equal up to a scalar.
Their equality follows from the fact that they both preserve the 1-tensor. \end{proof}

\begin{remark} The argument involving diagrammatic localization for this lemma is slightly more technical, as we must examine the standard module $R_{h_I w_J}$ which is not the top
standard inside $B_{w_I h_T}$. It is still not difficult, so let us omit the details. \end{remark}

\subsection{The crossover isomorphism in the affine case}

Now we bootstrap these results to the affine case. Let $J \subset X$, and let $Y$ be the corresponding set with $I \subset Y$. Let us recall from \S\ref{subsec-descentaffine} and elaborate upon the way in which one proves that $w_I h_X = h_Y w_J$. We rewrite $w_I h_X$ using \eqref{eq:rewriteaffine2} as
\begin{equation} \label{foobar1} w_I h_X = h_{I \cap A_0} (w_{I \cap A_d} h_{X \cap A_d}) \cdots (w_{I \cap A_1} h_{X \cap A_1}) (w_{\mathring{A_0}} h_{X \cap A_0}). \end{equation}
In exactly the same way we rewrite $h_Y w_J$ as
\begin{equation} \label{foobar2} h_Y w_J = (h_{Y \cap A_0} w_{\mathring{A_0}}) (h_{Y \cap A_d} w_{J \cap A_d}) \cdots (h_{Y \cap A_1} w_{J \cap A_1}) h_{J \cap A_0}. \end{equation}
By the finite case,
\begin{equation} w_{I \cap A_i} h_{X \cap A_i} = h_{Y \cap A_i} w_{J \cap A_i} \end{equation}
for each $1 \le i \le d$. Thus we need only compare what happens in the $\tau$-component $A_0$.
Since $X$ does not contain all of $A_0$ by notational assumption, $X \cap A_0 = J \cap A_0$. Similarly, $Y \cap A_0 = I \cap A_0$.  The element $w_{\mathring{A_0}}$ commutes with everything in the middle of the expression. This proves that \eqref{foobar1} and \eqref{foobar2} are equal. Finally, we should remember that
\begin{equation} \label{foobar3} h_{I \cap A_0} w_{\mathring{A_0}} = w_{I \cap A_0}, \quad w_{\mathring{A_0}} h_{J \cap A_0} = w_{J \cap A_0}. \end{equation}
This fact \eqref{foobar3} was a reminder of how the expressions in \eqref{foobar1} or \eqref{foobar2} came from expressions which started with $w_I$ or ended with $w_J$.

The equality in \eqref{foobar1}, in the categorification, are most succinctly stated in the description \eqref{style3} of the idempotent $e_{I,X,X}$. Using the corresponding description of the idempotent $e_{Y,Y,J}$ for $h_Y w_J$, we need only link these idempotents by using crossovers inside each $\tau$-component $A_i$. These crossovers render most of the ovals in \eqref{style3} redundant, and we have the following result. The proof is not worth saying aloud.

\begin{lemma} The following diagram, the \emph{mega-crossover}, is preserved by precomposition with $e_{I,X,X}$ and postcomposition with $e_{Y,Y,J}$, and gives an isomorphism between their images.
\begin{equation} \label{megacrossover} \ig{1}{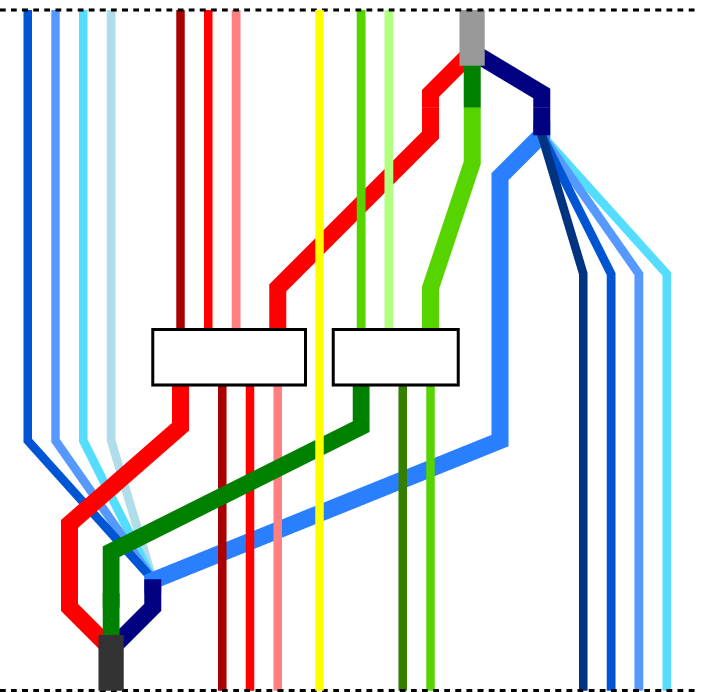} \end{equation}
Each rectangle is supposed to represent a crossover map as in \eqref{crossoverT} or \eqref{crossoverJ}. \end{lemma}

Let's tear into this picture, which is a typical example. A possible setting for this example is $n = 12$, $I = \{1,2,3,4, 6, 7, 10, 11\}$, $J = \{2, 3, 4, 5, 7, 8, 11, 0\}$, $X = \{2, 3, 4, 5, 7, 8,  9, 10, 11, 0\}$ and $Y = \{1, 2, 3, 4,6,  7, 8, 9, 10, 11, 0\}$.

Firstly, $I$ and $J$ are represented by graying thick strands, which get split into their intersections with each $\tau$-component. The reddish, yellowish, greenish, and blueish strands
represent the $\tau$-components $A_3 = \{10,11,0\}$, $A_2 = \{9\}$, $A_1 = \{6,7,8\}$, $A_0 = \{1,2,3,4,5\}$ respectively in this example. Note the yellow strand doing nothing: we have
$I \cap A_2 = \mt = J \cap A_2$, so this $\tau$-component has a single index $i$, and $i \in X \cap Y$. Of course, the idempotent $e_1$ is the identity map, so this yellow strand is
never involved in anything. In the red component we have $X = Y = A_3$ in this example, and the rectangle represents \eqref{crossoverT}. In the green component we have $X = J \cap A_1 \ne I
\cap A_1 = Y$, and the rectangle represents \eqref{crossoverJ}. We make this distinction to point out that the green strands do shift color in the crossover, while the red strands do
not. Finally, the most interesting stuff happens in the blue component $A_0$. Here $X = J \cap A_0 \ne I \cap A_0 = Y$, and the blue strands do shift color from bottom to top. The blue
splitter should be thought of as half of \eqref{crossoverT} for $T = I \cap A_0$, and the blue merger should be thought of as half of \eqref{crossoverT} for $T = J \cap A_0$, see
\eqref{foobar3}.

\begin{lemma} \label{lem:dotmegacrossover} The mega-crossover commutes with dots on strands not in $J$. That is, pick some index $x \in X \setminus J$, and put an $x$-colored startdot on the bottom of
\eqref{megacrossover}. The result is a mega-crossover for $X \setminus x$, composed with a startdot colored by $y \in Y \setminus I$. Here $y$ is the correspondent of $x$, the unique
element of $Y \setminus I$ in the same $\tau$-component as $x$. \end{lemma}

\begin{proof} Any element $x \in X \setminus J$ is in a $\tau$-component $A_i$ for $i \ne 0$. Now the result follows from \eqref{dotcrossover}. \end{proof}
	
\subsection{Comparing left and right}

Having done this work we are now capable of comparing $\NC_I$ with $\MC_J$. These complexes are not the same, but only because the signs on the differentials disagree.

Consider the differential from $B_{w_I h_X}$ to $B_{w_I h_{X'}}$, where $J \subset X,X' \subsetneq S$ differ by a single element $x \notin J$. The differential between these summands in
$\NC_I$ is a signed dot, and the sign is determined by the number of strands which come to the left of $x$. Letting $I \subset Y, Y' \subsetneq S$ be the correspondents of $X, X'$
respectively, they differ by some $y \notin I$ in the same $\tau$-component as $x$. The differential from $B_{h_Y w_J}$ to $B_{h_{Y'} w_J}$ in $\MC_J$ is again a signed dot, but the
sign is determined by the number of strands which come to the left of $y$. It was shown in Lemma \ref{lem:dotmegacrossover} that the $x$-colored dot matches the $y$-colored dot, but the
signs may not match, as evidenced in \eqref{dotcrossover}.

We can be completely precise. Suppose $x \in X$ and $X' = X \setminus x$. The element $x$ will always equal $\min(A_i)$ for some $\tau$-component $A_i$, and $y$ will equal $\max(A_i)$.
Moreover, $X \cap A_i = A_i$. Thus there are $|A_i|-1$ strands which come to the left of $x$ but not to the left of $y$, and the differentials are off by the sign $(-1)^{|A_i| - 1}$.

To remedy this, consider the map $\NC_I \to \MC_J$ which multiplies each term $B_{w_I h_X}$ by the sign $(-1)^{a(X)}$ and then applies the megacrossover, where
\begin{equation} a(X) = \sum_{A_i \subset X} (|A_i| - 1). \end{equation}
If $|A_i|$ is odd then it makes no difference where $X \cap A_i = A_i$ or $X \cap A_i = J \cap A_i$, but if $|A_i|$ is even then it does. The difference between $a(X)$ and $a(X \setminus x)$ is exactly $|A_i| - 1$, so this map will intertwine the differentials.

\begin{prop} The pseudocomplexes $\NC_I$ and $\MC_J$ are isomorphic via signed crossover maps $B_{w_I h_X} \to B_{h_Y w_J}$, where the sign is given by $(-1)^{a(X)}$. \end{prop}

\begin{remark} We believe the most satisfactory sign convention for $\NC_I$ would say that the sign on a dot is given by $(-1)$ to the number of $\tau$-components to the left of the
dot. This would produce a description of $\NC_I$ most analogous to the original description of $\FC$, where objects are determined by proper subsets of the set of $\tau$-components,
rather than by proper subsets of $S$ which contain $J$. We did not want to introduce this complication when we first defined $\NC_I$ though. \end{remark}

%%%%%%%%%%%%%%%%%%%%%%%%%
\section{The twisting isomorphism}
\label{sec-tensorBIcomps}
%%%%%%%%%%%%%%%%%%%%%%%%%

In the previous chapter we proved that $B_I \FC$ and $\FC B_{\tau(I)}$ are each homotopy equivalent to $\NC_I$, and composing these we get a homotopy equivalence
\begin{equation}\rotisom_I \co B_I \FC \to \FC B_{\tau(I)}.\end{equation} To write down this chain map explicitly can take a lot of unraveling, as each homotopy equivalence is a
complicated Gaussian elimination, and some signs are involved in the crossover. When $I = \{s\}$ it is more tractable, as the process involved a single simultaneous Gaussian elimination, but the unraveling is still somewhat annoying.

In fact, we have already proven in Proposition \ref{prop:dimhomBI} that the dimension of $\Hom(B_s \FC, \FC B_t)$ is one. If we can write down any nonzero chain map, it will agree with
the one coming from Gaussian elimination up to scalar. The minimal complex $\NC_s$ is only well-defined up to isomorphism anyway, which in this case is a scalar multiple of the identity
by Lemma \ref{lem:multfreejustid}. So before imposing additional constraints, $\rotisom_s$ is only defined up to scalar.

The additional constraint we impose is that $\square_f$ commutes for the startdot and the enddot, c.f. \eqref{startdotsquare} and \eqref{enddotsquare}. The consequence is that the
induced isomorphism $(-)\ot \FC \to \FC \ot \tau(-)$ is functorial. So we need only provide a chain map $\phi_s \co B_s \FC \to \FC B_t$ which makes these squares commute, and we are
done.

The results in this section are mostly of the form ``this is a chain map'' and ``this homotopy makes this diagram commute.'' Their proofs are completely straightforward diagrammatic
busywork. We establish some notational aids (a block matrix description of differentials and homotopies) which make the computation easier to perform, and record enough details that the reader should have an easier time reproducing the result.

\begin{remark} \label{offbyasign} We make a warning once and for all. The tensor product complex $B_s \FC$ does not satisfy the strand-counting sign rule! Instead, the differentials are
off by a sign from the strand-counting sign rule, because the first strand $B_s$ is not counted. This problem would be alleviated by using $B_s \sot \FC$, but this would create a host of
other problems instead. It is essential to use the correct monoidal structure when discussing questions of monoidal commutativity. The tensor product $\FC B_{\tau(s)}$ does satisfy the
strand-counting sign rule. \end{remark}

%===========
\subsection{Block matrices}
\label{subsec-blocks}
%===========

Let $s \in S$ and $t = \tau(s)$. We will split the proper subsets $Z \subsetneq S$ into blocks based on their intersection $X = Z \cap (S \setminus \{s,t\})$.

Let $X \subsetneq S$ with $s, t \notin X$. We write $sX$ for $X \cup \{s\}$ and $Xt$ for $X \cup \{t\}$. Whenever $X \cup \{t,s\} \ne S$, we write $X \cup \{t,s\} = X_1 t s X_2$, where $X = X_1 \sqcup X_2$ and the expression $X_1 t s X_2$ is in anticyclic order.

\begin{remark} We are careful to write all terms in anticyclic order. Since $t \notin X \cup \{s\}$, there is an anticyclic order where $s$ comes first, so $h_{sX} = s h_X$. This
justifies writing $sX$ for $X \cup \{s\}$. Similarly, $t$ can come last in $X \cup \{t\}$, so we have $h_{Xt} = h_X t$. Meanwhile, for $X \cup \{t,s\}$ it must be the case that $t$ comes
before $s$ in anticyclic order, but it may not be possible to place them at the beginning or the end of $X$. There is a decomposition $X = X_1 \sqcup X_2$ and a cyclic order on $X \cup
\{t,s\}$ such that for any $x_1 \in X_1$ and $x_2 \in X_2$ we have $x_1 > t > s > x_2$. In particular, $x_1$ commutes with $s$ and with $x_2$, and $t$ commutes with $x_2$. \end{remark}

Let us loosely write $\FC$ as a direct sum indexed by subsets $X$ for which $s, t \notin X$: \begin{equation} \label{eq:FCloose} \FC \approx \bigoplus_{s,t \notin X} B_X \oplus B_{sX}
\oplus B_{Xt} \oplus B_{X_1 ts X_2}. \end{equation} This is a very loose description: $B_X$ appears many times in $\FC$; also, for parity reasons $B_X$ and $B_{sX}$ never appear in the
same homological degree. To be more precise, in each homological degree $k$ we can write \begin{equation} \FC^k(-k) \cong \bigoplus_{s,t \notin X} B_X^{\oplus ?} \oplus
B_{sX}^{\oplus ?} \oplus B_{Xt}^{\oplus ?} \oplus B_{X_1 ts X_2}^{\oplus ?}. \end{equation} The notation $B_{sX}^{\oplus ?}$ means $B_{sX}$ if $sX \in P_k$ so that $B_{sX}$ could be a
summand of $\FC^k$, and zero if $sX \notin P_k$. But \eqref{eq:FCloose} will be more useful, as our computations will not care particularly which homological degree we are in, or which
summands are actually present.

We will write certain morphisms (e.g. differentials, chain maps, homotopies, etc) using \emph{block matrix notation}. For example, suppose we want to describe the differential in the complex $\FC$. We can give a $4 \times 4$ matrix $d_X^X$ which describes how the differential acts on the four summands
\begin{equation} \label{eq:foursummands} B_X \oplus B_{sX} \oplus B_{Xt} \oplus B_{X_1 ts X_2}. \end{equation}
We refer to these four summands as a \emph{block}, the block of $X$.

For sake of an example, here is $d_X^X$; the diagrammatic notation will be explained in the next section.
\begin{equation*} d_X^X = \left( \begin{array}{cccc}
0 & \finaldotred \rainbow & (-1)^X \rainbow \finaldotblue & 0 \\
\startdotred \rainbow & 0 & 0 & \difffinalblue \\
(-1)^X \rainbow \startdotblue & 0 & 0 & (-1)^{X+1} \difffinalred \\
0 & \diffstartblue & (-1)^{X+1} \diffstartred & 0 \end{array} \right) \end{equation*}
One should interpret this matrix as follows. Whenever $B_X(k)$ appears in degree $k$ and $B_{sX}(k+1)$ appears in degree $k+1$, the appropriate entry of the matrix (namely $\startdotred \rainbow$) gives the summand of the differential from $B_X(k) \to B_{sX}(k+1)$. If $B_X(k)$ appears in degree $k$ but $B_{sX}(k+1)$ does not appear in degree $k+1$, then there is no differential to a non-existent summand, so just don't worry about the fact that the matrix has a nonzero entry.

To fully describe the differential, we should also give a $4 \times 4$ matrix $d_X^Y$ which describes how the differential sends the block of $X$ to the block of $Y$, for $X \ne Y$. What makes this notation effective is that the matrix $d_X^X$ does not depend in any interesting way on $X$ (or, by nature of the notation, on the homological degree in which the summand
$B_X$ appears). Similarly, $d_X^Y$ will not depend on $X$ and $Y$ beyond general features, like whether $X$ and $Y$ differ by a single element of $S$.

We will also use block notation for the complexes $B_s \FC$ and $\FC B_t$. For $B_s \FC$, the \emph{block} of $X$ will refer to the four summands
\begin{equation} B_s B_X \oplus B_s B_{sX} \oplus B_s B_{Xt} \oplus B_s B_{X_1 ts X_2}. \end{equation}
To describe a chain map $g \co \FC \to B_s \FC$ one might use block notation, where $g^Y_X$ would be a $4 \times 4$ matrix indicating the chain map between summands
\begin{equation} B_X \oplus B_{sX} \oplus B_{Xt} \oplus B_{X_1 ts X_2} \to B_s B_Y \oplus B_s B_{sY} \oplus B_s B_{Yt} \oplus B_s B_{Y_1 ts Y_2}. \end{equation}

Our isomorphism $\rotisom_s$ will be block diagonal, i.e. $\rotisom_X^Y = 0$ for $X \ne Y$. This is the expected behavior, as the Gaussian elimination procedure from \S\ref{subsec-GEobjects} was also performed blockwise.

%===========
\subsection{Rewriting $\FC$}
\label{subsec-rewriting}
%===========

Let us discuss the differential on $\FC$ using block matrix notation. The differential sends $B_Z$ to $B_{Z'}$ where $Z$ and $Z'$ differ by a single element $u$.
When $Z \in \{X,sX,Xt,X_1 ts X_2\}$ is in the block of $X$, and $u \in \{s,t\}$, then $Z'$ is also in the block of $X$. If $u \notin \{s,t\}$ then $Z'$ is in the block of some $Y$ which differs from $X$ by $u$.

We draw $s$ as red and $t$ as blue. We draw $X$ (or rather, $B_X$) as a thick rainbow, and $X_1$ and $X_2$ as thinner rainbows. The remaining conventions will be discussed shortly. Here is the differential on the block diagonal.
\begin{equation} \label{eq:dXX} d_X^X = \left( \begin{array}{cccc}
0 & \finaldotred \rainbow & (-1)^X \rainbow \finaldotblue & 0 \\
\startdotred \rainbow & 0 & 0 & \difffinalblue \\
(-1)^X \rainbow \startdotblue & 0 & 0 & (-1)^{X+1} \difffinalred \\
0 & \diffstartblue & (-1)^{X+1} \diffstartred & 0 \end{array} \right) \end{equation}

Let us explain our diagrams above. We write $(-1)^X$ for $(-1)^{|X|}$. The signs which appear are governed by the strand-counting sign rule. Thus the $(-1)^X$ on the $\finaldotblue$
appears because it comes after $\rainbow$, which is really a Bott-Samelson $B_X$ with $|X|$ strands.

The most complicated differentials come to and from the term $B_{X_1 t s X_2}$. The split of an $X$ rainbow into thinner $X_1$ and $X_2$ rainbows is actually just the identity map, and represents a picture like this.
\begin{equation} {
\labellist
\small\hair 2pt
 \pinlabel {$=$} [ ] at 47 17
\endlabellist
\centering
\ig{1}{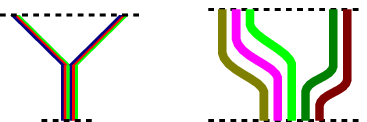}
} \end{equation}
The strands of $B_X$ that are in $X_1$ just get grouped to the left, and the strands in $X_2$ get grouped to the right. The red strand crossing over a thin $X_1$ rainbow is meant to represent a bunch of $4$-valent vertices,
and makes sense since $s$ commutes with $X_1$. Similarly, $t$ commutes with $X_2$ so the blue strand can cross over a thin $X_2$ rainbow. These crossings occur as the reduced expression
for $h_Z$ changes, and by the sign rules of \S\ref{subsec-signs}, these 4-valent vertices should all come with signs! Thus whenever a red strand crosses over $X_1$ in a differential
there should also be a sign $(-1)^{X_1}$. Thus the sign $+1$ on $\difffinalblue$ comes because there is a sign $(-1)^{X_1}$ from the red crossing, and a sign $(-1)^{X_1}$ from the
strand-counting rule for the blue dot, which cancel out. The sign $(-1)^{X+1}$ on $\difffinalred$ comes because there is a sign $(-1)^{X_2}$ from the blue crossing, and a sign $(-1)^{X_1
+ 1}$ from the strand-counting rule for the red dot, which combine to $(-1)^{X+1}$.

Note that there are two nonzero entries in each column of $d_X^X$, one corresponding to $u=s$ and one corresponding to $u=t$.

Now let us draw the off-diagonal blocks of the differential, which send terms indexed by $X$ to terms indexed by $Y$, where $Y = Xu$ or $Y = X \setminus u$ for some $u \notin \{s,t\}$.
\begin{equation} \label{eq:dXY} d_X^Y = \left( \begin{array}{cccc} 
\dotdiff & 0 & 0 & 0 \\
0 & - \linered \dotdiff & 0 & 0 \\
0 & 0 & \dotdiff \lineblue & 0 \\
0 & 0 & 0 & \thindotdiff \lineblue \linered \thindotdiff
\end{array} \right) \end{equation}
This matrix is diagonal, i.e. changing an index $u \notin \{s,t\}$ it sends $B_{sX}$ to $B_{sY}$ and $B_{tX}$ to $B_{tY}$, etcetera. We write $\dotdiff$ for either $\startdotsign_X^Y$ or $\finaldotsign_X^Y$ with its built-in sign, whenever appropriate. Thus the extra sign in the second row comes from the strand-counting sign rule, because of the extra red strand compared to the natural sign in $\dotdiff$. In the last row, the changed index $u$ is either in $X_1$ or in $X_2$. If $u \in X_1$ then this diagram is supposed to represent the identity on $B_{X_2}$ and the appropriate differential $B_{X_1} \to B_{Y_1}$, and similarly if $u \in X_2$. Note that the sign which appears on the dot is the same sign $\ep_{X,Y}$ that appears in $\dotdiff$; the addition of two extra strands $t$ and $s$ between $X_1$ and $X_2$ will not change the parity of the strand count for any index in $X_2$.

Finally, from Remark \ref{offbyasign} we remind the reader that the differential in either $B_s \FC$ (resp. $\FC B_t$) has the same sign as the differentials above, just with an extra red strand on the left (resp. an extra blue strand on the right). The extra red strand on the left is not counted for the strand-counting sign rule.

%===========
\subsection{The twisting isomorphism for a simple reflection}
\label{subsec-rotisoms}
%===========

We will now describe $\rotisom = \rotisom_s$ using block matrix notation. As mentioned previously, $\rotisom$ is block diagonal, so that $\rotisom_X^Y=0$ for $X \ne Y$.

\begin{defn} Define a block diagonal map $\rotisom$ by the block matrix
\begin{equation} \label{eq:rotisoms} \rotisom_X^X = \left( \begin{array}{cccc} 
0 & 0 & 0 & 0 \\
0 & \mergered \rainbow \startdotblue & (-1)^{X+1} \linered \rainbow	\lineblue & 0 \\
0 & (-1)^{X+1} \cupcapbow & \finaldotred \rainbow \splitblue & 0 \\
0 & 0 & 0 & \thatterm \end{array}
\right) \end{equation}
\end{defn}

\begin{remark} For any given $k$, note that a given term, such as $B_s B_{X_1 t s X_2}$, appears as a direct summand of $(B_s \FC)^k$ if and only if the corresponding term, such as
$B_{X_1 t s X_2} B_t$, appears as a direct summand of $(\FC B_t)^k$. In particular, any nonzero entry of the matrix \eqref{eq:rotisoms} always makes sense: if the source does appear in
$(B_s \FC)^k$, then the target does appear in $(\FC B_t)^k$. \end{remark}

\begin{thm} The block-diagonal morphism $\rotisom_s$, defined on blocks by $\rotisom_X^X$, is a chain map $B_s \FC \to \FC B_t$. It is not nulhomotopic. \end{thm}

\begin{proof} We just need to check that $d \rotisom = \rotisom d$. Let us check this on the summands of the form $B_s B_{sX}$, and leave the rest as an exercise to the reader.

Let us compute $d_X^X \rotisom_X^X$ applied to $B_s B_{sX}$, i.e. only look at the terms in the target with the same indexing set $X$. We obtain the column
\begin{equation} \left( \begin{array}{c} \finaldotred \rainbow \lineblue \circ \mergered \rainbow \startdotblue + (-1)^X \rainbow \finaldotblue \lineblue \circ (-1)^{X+1} \cupcapbow \\ 0 \\ 0 \\ \diffstartblue \lineblue \circ \mergered \rainbow \startdotblue + (-1)^{X+1} \diffstartred \lineblue \circ (-1)^{X+1} \cupcapbow \end{array} \right). \end{equation}
The first entry of the column adds to zero.

Meanwhile, applying $\rotisom_X^X d_X^X$ to $B_s B_{sX}$ we obtain the column
\begin{equation} \left( \begin{array}{c} 0 \\ 0 \\ 0 \\ \thatterm \circ \linered \diffstartblue \end{array} \right). \end{equation} Because $\rotisom_s$ is zero on terms of the form $B_s B_X$, only the term of the differential sending $B_s B_{sX}$ to $B_s B_{X_1 ts X_2}$ will contribute. Comparing these two columns, they are equal by the usual relation for a dot on a $6$-valent vertex, \cite[page 8]{EWGr4sb}.

There is one exceptional circumstance, which is when $B_{sX}$ appears in homological degree $k$, but $B_{X_1 ts X_2}$ does not appear in degree $k+1$. The only difference this makes is that the last row disappears from both columns, and the columns remain equal.

Now apply $d_X^Y \rotisom_X^X$ to $B_s B_{sX}$, i.e. only look at the terms in the target with the indexing set $Y$, which we assume differs from $X$ by a single element. We obtain the column
\begin{equation} \left( \begin{array}{c}
0 \\
-\linered \dotdiff \lineblue \circ \mergered \rainbow \startdotblue \\
\dotdiff \lineblue \lineblue \circ (-1)^{X+1} \cupcapbow \\
0
\end{array} \right). \end{equation}

Meanwhile, applying $\rotisom_Y^Y d_X^Y$ to $B_s B_{sX}$, we obtain the column
\begin{equation} \left( \begin{array}{c}
0 \\ \mergered \rainbow \startdotblue \circ - \linered \linered \dotdiff \\ (-1)^{Y+1} \cupcapbow \circ - \linered \linered \dotdiff \\ 0 \end{array} \right). \end{equation}
Comparing these two columns, they are equal merely by isotopy, and the observation that $(-1)^Y = - (-1)^X$.

Thus we have confirmed that $d \rotisom = \rotisom d$ when applied to any term of the form $B_s B_{sX}$. The other three terms are very similar and no more complicated, and we leave
them to the reader, who is also welcome to contact the author for details.

In fact, there are no nonzero homotopies $B_s \FC \to \FC B_t[-1]$, simply by analysis of the terms in each degree. \end{proof}

%===========
\subsection{Commutation with dots}
\label{subsec-rotisomdot}
%===========

We now examine in detail $\square_{\startdotred}$, \eqref{startdotsquare}.

\begin{thm} \label{thm:dotcommute} The chain maps $\FC \to \FC B_t(1)$ given by $\id \ot \startdotblue$ and $\rotisom_s \circ (\startdotred \ot \id)$ are homotopic. The homotopy is unique, up to homotopy of
homotopies. \end{thm}

\begin{proof} We begin by noting that there are possible homotopy of homotopies, a map $K \co \FC \to \FC B_t(1)[-2]$. By examination of the summands which appear, one can deduce that
$K$ is a block diagonal matrix, whose diagonal term is
\begin{equation} \label{eq:KXX} K_X^X = \left( \begin{array}{cccc} 0 & 0 & 0 & 0 \\ 0 & 0 & 0 & 0 \\ 0 & 0 & a_X \rainbow \splitblue & 0 \\ 0 & 0 & 0 & 0 \end{array} \right). \end{equation}
The one nonzero term represents a map $B_{Xt} \to B_{Xt} B_t(-1)$, and $a_X$ is some scalar.

Taking $dK - Kd$, we obtain the following map $\FC \to \FC B_t(1)[-1]$, which is a homotopy for the zero chain map.
\begin{equation} (dK - Kd)_X^X = \left( \begin{array}{cccc} 0 & 0 & a_X (-1)^X \rainbow \lineblue & 0 \\ 0 & 0 & 0 & 0 \\ a_X (-1)^{X+1} \rainbow \cupblue & 0 & 0 & a_X (-1)^X \foobara \\ 0 & 0 & a_X (-1)^{X+1} \foobarb & 0 \end{array} \right) \end{equation}
\begin{equation} (dK - Kd)_X^Y = \left( \begin{array}{cccc} 0 & 0 & 0 & 0 \\ 0 & 0 & 0 & 0 \\ 0 & 0 & (a_X - a_Y) \dotdiff \splitblue & 0 \\ 0 & 0 & 0 & 0  \end{array} \right) \end{equation}
	
The general form of a homotopy $h \co \FC \to \FC B_t(1)[-1]$ is the same, but with more general scalars.
\begin{equation} h_X^X = \left( \begin{array}{cccc} 0 & 0 & b_X \rainbow \lineblue & 0 \\ 0 & 0 & 0 & 0 \\ c_X \rainbow \cupblue & 0 & 0 & d_X \foobara \\ 0 & 0 & e_X \foobarb & 0 \end{array} \right) \end{equation}
\begin{equation} h_X^Y = \left( \begin{array}{cccc} 0 & 0 & 0 & 0 \\ 0 & 0 & 0 & 0 \\ 0 & 0 & f_{YX} \dotdiff \splitblue & 0 \\ 0 & 0 & 0 & 0  \end{array} \right) \end{equation}
Thus, modulo nulhomotopies of the form $dK - Kd$, for carefully chosen values $a_X$, we can assume that one of the coeficients $b_X$, $c_X$, $d_X$, or $e_X$ is zero.

It is a cumbersome and straightforward exercise to compute $dh + hd$, as well as \begin{equation}g = \id \ot \startdotblue - \rotisom_s \circ (\startdotred \ot \id).\end{equation} Let us summarize the results. In the matrices below, we have omitted the diagrams and remembered only the coefficients; the reader can verify that the diagrams are the same on both sides. (Aside: in the lower right entry of these matrices, one must use the relation for a dot on a $6$-valent vertex, \cite[page 8]{EWGr4sb}, to compare the diagrams in these two matrices.) We leave the most complicated entry for later, and write $*$ for now.
\begin{equation} (dh + hd)_X^X = \left( \begin{array}{cccc}
(-1)^X (b_X + c_X) & 0 & 0 & (-1)^X (d_X - b_X) \\ 0 & 0 & (b_X + e_X) & 0 \\ 0 & (c_X + d_X) & * & 0 \\ (-1)^X(e_X - c_X) & 0 & 0 & (-1)^{X+1}(d_X + e_X) \end{array} \right) \end{equation}
\begin{equation} g_X^X = \left( \begin{array}{cccc}
1 & 0 & 0 & 0 \\ 0 & 0 & (-1)^X & 0 \\ 0 & (-1)^X & * & 0 \\ 0 & 0 & 0 & -1 \end{array} \right) \end{equation}
 For these two matrices to be equal we must have
\begin{equation}\label{eq:bcde} b_X = d_X, \quad c_X = e_X, \quad b_X + c_X = (-1)^X. \end{equation}

Now observe that $g$ is block diagonal, so that $g_X^Y$ is the zero matrix. Meanwhile one computes that $(dh + hd)_X^Y$ has the following nonzero coefficients, when $X$ and $Y$ differ by a single element.
\begin{subequations}
\begin{equation} \label{eq:bXbYreln} b_X + b_Y + (-1)^Y \ep_{YX} f_{YX} \end{equation}
\begin{equation} c_X + c_Y + (-1)^X \ep_{YX} f_{YX} \end{equation}
\begin{equation} d_X + d_Y + (-1)^{X+1} \ep_{YX} f_{YX} \end{equation}
\begin{equation} e_X + e_Y + (-1)^{Y+1} \ep_{YX} f_{YX} \end{equation}
\end{subequations}
We should ensure that all of these are zero. It is enough to ensure that \eqref{eq:bXbYreln} to zero; the rest are zero using \eqref{eq:bcde}, and observing that $(-1)^Y = - (-1)^X$.

Finally, $h$ has a term going from the $X$ block to the $Y$ block where they differ by one element, and $d$ has a term going from the $Y$ block to the $Z$ block, where they differ by one element. Thus $hd+dh$ has a term going from the $X$ block to the $Z$ block. When $Z = X$, this contributes to the nasty $*$ term from before. When $X \ne Z$, there are exactly two possible $Y$ that this could factor through. The resulting diagram is $\dotdiff \circ \dotdiff$ next to $\splitblue$, and the coefficient is
\begin{equation} \label{eq:aeiou} \sum_Y f_{YX} \ep_{YX} + f_{ZY} \ep_{ZY}. \end{equation}

Let us now simplify by assuming $b_X = 0$ for all $X$, which we can do by adding $dK - Kd$ to $h$, as noted above. Under this assumption, \eqref{eq:bcde} implies that
\begin{equation} b_X = 0, \quad c_X = (-1)^X, \quad d_X = 0, \quad e_X = (-1)^X. \end{equation}
Then using \eqref{eq:bXbYreln} we see that
\begin{equation} f_{YX} = 0 \end{equation}
for all $X$ and $Y$. Consequently \eqref{eq:aeiou} holds, as desired.

Under these simplifying assumptions, let us describe the missing entry $*$ of $(dh+hd)_X^X$ and $g_X^X$, which is a map $B_{Xt} \to B_{Xt} B_t$. On this entry, $dh+hd$ acts as
\begin{equation} \label{eq:dhhddash} - \ig{.5}{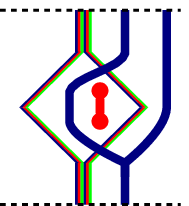} + \rainbow \finaldotblue \cupblue \end{equation}
while $g$ acts as
\begin{equation} \label{eq:gdash} \rainbow \lineblue \startdotblue - \al_s \rainbow \splitblue \end{equation}
To simplify the second term in \eqref{eq:gdash}, pull the polynomial $\al_s$ across $X_1$ and pull a blue strand across $X_2$ until the diagram looks like
\begin{equation} \ig{.5}{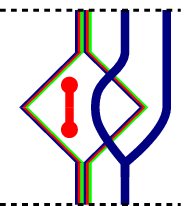}. \end{equation}
Now apply \eqref{eq:barbellforce} to force $\al_s$ across the blue line. One term will cancel the first term from \eqref{eq:dhhddash}, and the remaining terms only have blue diagrams, so they can pull back past $X_2$. In the end, one is left with diagrams which all have $\rainbow$ on the left, and equality between \eqref{eq:dhhddash} and \eqref{eq:gdash} will follow if
\begin{equation} \lineblue \startdotblue + \startdotblue \lineblue - (\lineblue \barbblue \lineblue) \circ \splitblue - \finaldotblue \cupblue = 0. \end{equation}
But this equation follows from \eqref{eq:barbellforce}, see also \cite[Remark 6.5]{EWGr4sb}.

Thus we have checked that $dh + hd = g$, and moreover, that modulo the possible homotopies of homotopies $K$, this is the unique solution for $h$. \end{proof}

We will not publically repeat this tedious process for $\square_{\finaldotred}$, \eqref{enddotsquare}. Let us make the following note for the reader. By combining the Dynkin diagram automorphisms $\si$ and $\tau$, one can reproduce the reflection on the affine Dynkin diagram that swaps $s$ and $t$. Now consider the (contravariant, monoidally contravariant) symmetry of $\Diag$ which rotates a diagram by 180 degrees and then applies this Dynkin diagram automorphism. This symmetry swaps $\square_{\finaldotred}$ and $\square_{\startdotred}$, and it fixes the matrix $\rotisom_X^X$.

For example, in the computation of $\square_{\finaldotred}$, the double homotopies $K$ are maps $B_s \FC \to \FC(1)[-2]$, and are described by
\begin{equation} \label{eq:KXXredux} K_X^X = \left( \begin{array}{cccc} 0 & 0 & 0 & 0 \\ 0 & a_X \mergered \rainbow & 0 & 0 \\ 0 & 0 & 0 & 0 \\ 0 & 0 & 0 & 0 \end{array} \right). \end{equation}
The one nonzero term represents a map $B_s B_{sX} \to B_{sX}(-1)$.

%===========
\subsection{Commutation with polynomials}
\label{subsec-rotisompoly}
%===========

We work with the standard affine realization, see Definition \ref{defn:standardaffine}. We study $\square_{f}$ for a polynomial $f$. This last computation has nothing to do with the
morphism $\rotisom_s$, since only $\rotisom_{\one} = \id_{\FC}$ is used in $\square_{f}$. Since the polynomial ring is generated over the simple roots by $x_1$, if we can show that
$\square_{x_1}$ commutes up to homotopy, then so does $\square_f$ for any polynomial $f$.

\begin{thm} The map $x_1^l - x_2^r \co \FC \to \FC(2)$, given by left multiplication by $x_1$ minus right multiplication by $x_2$, is nulhomotopic. \end{thm}

\begin{remark} We expect the homotopy to be unique up to homotopy of homotopies, just as in Theorem \ref{thm:dotcommute}. In this case there is a very high dimensional space of homotopies $\FC \to \FC(1)[-1]$ and homotopies of homotopies $\FC \to \FC(2)[-2]$, so we did not bother to prove this additional expectation. \end{remark}

\begin{proof} Let us assume that $1 \notin X \subsetneq S$, and write all summands of $\FC$ as either $B_X$ or $B_{X1}$. The crux of this proof is the realization that $x_1^l - x_2^r$
acting on $B_X$ is effectively the same as the double differential $B_X \to B_{X1} \to B_X$ (i.e. the composition of two signed dot maps), and acting on $B_{X1}$ is effectively the same
as the double differential $B_{X1} \to B_X \to B_{X1}$. We refer to this as the \emph{$1$-indexed double differential}. That is, the $1$-indexed double differential is either $\finaldotsign_{X1}^X \startdotsign_X^{X1}$, a $1$-colored barbell $\barbred$, or $\startdotsign_{X}^{X1} \finaldotsign_{X1}^X$, a $1$-colored broken line $\brokenred$. Note that any signs must cancel.

Let red denote the index $1$. We have
\begin{subequations} \label{yay}
\begin{equation} \brokenred = \poly{x_1} \linered - \linered \poly{x_2}, \end{equation}
\begin{equation} \barbred = \poly{x_1} - \poly{x_2}. \end{equation}
\end{subequations}
Thus justifies the statement that $x_1^l - x_2^r$ is the $1$-indexed double differential acting on $B_{\mt}$ or on $B_1$.

Now suppose that $X1$ is a proper subset of $S$, and choose some $\notme \notin X1$ to induce a cyclic order. Write $X1 = X_1 1 X_2$ where $s > 1 > t$ for all $s \in X_1$ and $t \in
X_2$. In particular, the polynomial $x_1$ is fixed by $s$ for all $s \in X_1$, and $x_2$ is fixed by $t$ for all $t \in X_2$. Thus $x_1^l - x_2^r$ acting on $B_{X1} \cong B_{X_1} B_1
B_{X_2}$ is the equal to $B_{X_1} \ot (x_1 B_1 - B_1 x_2) \ot B_{X_2}$. Similarly, $x_1^l - x_2^r$ acting on $B_X \cong B_{X_1} B_{X_2}$ is equal to $B_{X_1} \ot (x_1 - x_2) \ot B_{X_2}$. Thus we can use \eqref{yay} to confirm that $x_1^l - x_2^r$ is the $1$-indexed double differential on $B_X$ and $B_{X1}$ for all such $X$.

There is one case in which the $1$-indexed double differential actually does not make sense, namely when $X = S \setminus \{1\}$ so that $X1$ is not a valid summand of $\FC$. Let us instead reinterpret the $1$-indexed double differential using Lemma \ref{lem:same}. This lemma stated that, for any $Y$ with size $\le |S| - 2$, the sum of all the
$i$-indexed double differentials $B_Y \to B_{Y \pm i} \to B_Y$, over all $i \in S$, is zero modulo $\delta$. In particular, the $1$-indexed double differential is equal to minus the sum
of all the $i$-indexed double differentials for $i \in S \setminus \{1\}$. In other words, modulo $\delta$,
\begin{equation} \label{eq:yayopposite} x_1^l - x_2^r = - \sum_{i \in S \setminus Y, i \ne 1} \finaldotsign_{Yi}^{Y} \startdotsign_{Y}^{Yi} + \sum_{j \in Y,j \ne 1} \startdotsign_{Y \setminus j}^{Y} \finaldotsign_{Y}^{Y \setminus j} \end{equation}
acting on $B_Y$, whenever $Y$ has size $\le |S| - 2$. We call this sum (without the minus sign) the \emph{all-but-$1$-indexed double differential}. In fact, \eqref{eq:yayopposite} also holds (modulo $\delta$) for $Y = S \setminus \{1\}$, the only set of size $|S|-1$ for which the expression makes sense. This is an exercise in the polynomial forcing relation \cite[equation (5.2)]{EWGr4sb}.

With this observation in hand, our goal is to define a homotopy $H \co \FC \to \FC(1)[-1]$ in such a way that $dH + Hd$ is either the $1$-indexed double differential, or minus the all-but-$1$-indexed double differential. Each nonzero term of $H$ will be either $\startdotsign$ or $\finaldotsign$ (possibly with a sign), and the combinatorics will make the sum $dH+Hd$ work.

We define $H$ as follows. On summands of the form $B_{X1}$, the homotopy will be $H_{X1} = \finaldotsign_{X1}^X$; note that whenever $B_{X1}$ appears in homological degree $k$, $B_X$ does appear in degree $k-1$. On summands of the form $B_X$, the homotopy $H_X$ will typically be zero. However, when $B_X$ appears in its maximal homological degree $n-1-|X|$, then we let
\begin{equation} \label{eq:Hdefn} H_X = - \sum_Y \startdotsign_X^Y, \end{equation}
where the sum is over all ways to add an index to $X$ to obtain $Y$, except that we also assume $1 \notin Y$.

Let us restate the definition of this homotopy. The set in $P = \cup_k P_k$ of all summands of $\FC$ can be partitioned into subsets as follows. The subsets of $S \setminus 1$ form a
cube, each appearing once in their maximal homological degree, which we call the \emph{ultimate cube}. What remains can be partitioned into pairs $\{X,X1\}$ appearing in homological
degrees $k-1$ and $k$ respectively for some $k$. For every differential which stays within one of these subsets, the homotopy $H$ will do that differential backwards, and will add a sign
for everything in the ultimate cube.

We first assert that the diagonal terms of $(dH + Hd)_Y^Y$, from $B_Y$ to $B_Y$, agree with either the $1$-indexed double differential or minus the all-but-$1$-indexed double
differential. \begin{itemize} \item If $1 \notin X$ and $X$ is not in the ultimate cube, then $H$ is zero on $B_X$, and $(Hd)_X^X$ is the $1$-indexed double differential. \item
Similarly, $B_{X1}$ is not in the image of $H$, and $(dH)_{X1}^{X1}$ is the $1$-indexed double differential. \item If $1 \notin X$ is in the ultimate cube, then $(dH + Hd)_X^X$ is built from all ways to traverse the ultimate cube back and forth, and is equal to minus the all-but-$1$-indexed double differential. \end{itemize}

Now we need only check that the off-diagonal terms $(dH+Hd)^Z_Y \co B_Y \to B_Z$ are zero. For this purpose we pick two elements $Y$ and $Z$ of $P$ in the same homological degree, with $Y \ne Z$. Both
$H$ and $d$ can change the content of a set by at most one element, so $dH + Hd$ can not map $B_Y$ to $B_Z$ unless \begin{itemize}
\item $Y = Zst$, or
\item $Z = Yst$, or
\item $Ys = Zt$, \end{itemize}
for some simple reflections $s, t \in S$. We assume that one of these three conditions holds.

Suppose $1 \notin Y$ and $1 \notin Z$. If $Y$ is not in the ultimate cube, $H_Y = 0$. If $Z$ is not in the ultimate cube, then the
only summand which maps to $Z$ under $H$ is $Z1$. There is no differential $B_Y \to B_{Z1}$ since $Y$ and $Z1$ differ by more than one element. Thus if neither $Y$ nor $Z$ is in the
ultimate cube, then $dH + Hd$ can not map $B_Y$ to $B_Z$.

If $Y$ is not in the ultimate cube but $Z$ is, we must examine $Hd \co B_Y \to B_Z$. In this case we must have $|Z| = |Y|+2$, so in particular $Z = Yst$ for two indices $s, t \in S$,
neither equal to $1$. By the usual argument, $Hd$ will be a sum of two terms, one factoring through $Ys$ and one factoring through $Yt$; both terms will consist of two startdots, and the
signs will cancel, for the same reason that the signs cancelled in Lemma \ref{lem:notsame}. Thus $dH  + Hd = 0$ from $B_Y$ to $B_Z$.

If $Z$ is not in the ultimate cube but $Y$ is, we must examine $dH \co B_Y \to B_Z$. In this case we must have $|Z| = |Y|-2$, so in particular $Zst = Y$ for two indices $s, t \in S$,
neither equal to $1$. By the usual argument, $Hd$ will be a sum of two terms, one factoring through $Zs$ and one factoring through $Zt$; both terms will consist of two enddots, and the
signs will cancel.

If both $Z$ and $Y$ are in the ultimate cube, then $|Y|=|Z|$, so in particular $Ys = Zt$ for two indices $s, t \in S$, neither equal to $1$. By the usual argument, $Hd+dH$ will be a sum
of two terms, one factoring through $Ys$ and one factoring through $Y \setminus t$; both terms will consist of one startdot and one enddot, and the signs will cancel.

Suppose $Y = X1$ and $1 \notin Z$. In particular, $Z$ is obtained from $X$ by either adding or subtracting some other index $s$. Then $dH$ will send $B_Y \to B_X \to B_Z$. If $Z$ is not
in the ultimate cube, then $Hd$ will send $B_Y \to B_{Z1} \to B_Z$. Again, these two terms cancel. Meanwhile, if $Z$ is in the ultimate cube, then $dH$ will also send $B_Y \to B_X \to
B_Z$, and the signs cancel because of the sign in \eqref{eq:Hdefn}.

Suppose $1 \notin Y$ and $Z = X1$. The argument here is similar, and has two cases, based on whether $Y$ is in the ultimate cube or not.

Suppose that $Y = X1$ and $Z = X'1$. Then $Hd$ from $B_Y$ to $B_Z$ is zero because $Z$ is not in the image of $H$, while $dH$ is zero because there is no differential from $X$ to $Z$. 

Having analyzed every case, we see that $dH + Hd$ is diagonal. This concludes the proof. \end{proof}

\appendix

%%%%%%%%%%%%%%%%%%%%%%%%%
\section{Computations}
\label{sec:computations}
%%%%%%%%%%%%%%%%%%%%%%%%%

For posterity, we record the answers to a number of computations. These can all be justified by brute force, such as by confirming that the maps we have written are actually chain maps.
As such, we have not bothered to record the proofs; the reader is welcome to contact the author for more details.

We will be computing with morphism spaces, and \begin{equation} \label{eq:removeomforfun} \Hom(\Om^k \MC, \Om^k \MC') \cong \Hom(\MC,\MC')\end{equation} for complexes $\MC, \MC'$. Thus
every computation involving $\VC$ has an analog for $\FC$. We opt to work with $\FC$, and ignore the copies of $\Om$. However, this does affect one major feature, which is the bimodule
action of polynomials in $R$ on morphisms. Since $\VC$ is central, the right and left actions of $R$ on $\VC$ agree, so one may neglect to say which side a polynomial acts on. However,
the identification \eqref{eq:removeomforfun} intertwines the right action of $R$ but not the left action. Polynomials do not act the same on both sides of $\FC$. Thus we must be careful
to always perform multiplication by polynomials on the right!

Finally, the reader should either be comfortable with pseudocomplexes from \S\ref{sec-pseudo}, or work in the quotient where $\delta = 0$. For example, when $n=2$, $\al_0 \equiv -
\al_1$ modulo $\delta$, a fact which is often used when confirming the computations below. In either setting, $\delta$ acts by zero on any morphism space in the homotopy category.

%=================
\subsection{Notational shorthand} \label{subsec:shorthand}
%=================

The complex $\VC$ is color-rotation invariant, or in other words, it is invariant under $\tau$. This must be the case if it is to be in the Drinfeld center, since it must commute with
$\Om \in \DC_{\ext}$. Similarly, any morphism in $\ZC(\HC_{\ext})$ is color-rotation invariant. As our complexes get larger, we can use color-rotation invariance to simplify our
notation, keeping track of $\tau$-orbits of objects and morphisms rather than the objects themselves. The shorthand below is really useful for compactly describing complexes and chain
maps in the Drinfeld center. It is not useful for discussing Wakimoto filtrations, index-by-index computations, etcetera, which is why it has yet been used until now in this paper.

\begin{notation} (We assume $n=3$ in order to discuss examples.) The symbol $\BBB$ will indicate a direct sum over a $\tau$-orbit of objects. For example, we have
\begin{subequations}
\begin{equation} \BBB_2 := B_0 \oplus B_1 \oplus B_2, \end{equation}
\begin{equation} \BBB_{10} := B_{10} \oplus B_{21} \oplus B_{02}, \end{equation}
\begin{equation} \label{eq:blather} \BBB_{020} := B_{020} \oplus B_{101} \oplus B_{212}. \end{equation}
We do not bother with any additional notation for the orbit $R$, so that the symbol $\BBB$ is always assumed to have $n$ terms. Of course $\BBB_2 \cong \BBB_1 \cong \BBB_0$, so do not read anything into our choice of index. In \eqref{eq:blather}, $B_{020}$ represented the indecomposable Soergel bimodule $B_{s_0 s_2 s_0}$ rather than the corresponding Bott-Samelson bimodule.

Morphisms between $\BBB$-symbols will always be assumed to be $\tau$-invariant. When we discuss the morphism $R \to \BBB_2$ given by
$\startdotblue$, we actually refer to the morphism $R \to B_0 \oplus B_1 \oplus B_2$ given by the (transpose of the) matrix \begin{equation} \left( \begin{array}{ccc} \startdotgreen &
\startdotred & \startdotblue \end{array} \right). \end{equation} We may also use $\startdotred$ for this same morphism.

A morphism $\BBB_2 \to \BBB_2(2)$ should actually correspond to some $3 \times 3$ matrix. For example, $\brokenblue$ would indicate a diagonal matrix with entries $(\brokengreen,
\brokenred, \brokenblue)$. Meanwhile, the map $\bluetored$ would also have three nonzero entries, one for a map $B_2 \to B_1$, one for a map $B_1 \to B_0$, and one for a map $B_0
\to B_2$. If we wish to add these two matrices together, we will write the result as $\brokenblue \oplus \bluetored$. We prefer $\oplus$ to $+$, since one should not add diagrams for which the source and target do not match, and this helps us distinguish between matrix entries in disjoint parts of the matrix.

The tensor product $\BBB_0 \BBB_0$ is actually a direct sum of $n^2$ objects, so that
\begin{equation} \BBB_0 \BBB_0 \cong \BBB_{00} \oplus \BBB_{01} \oplus \BBB_{02}. \end{equation}
Sometimes it helps to use this simplification, but usually we will leave the tensor product $\BBB_0 \BBB_0$ as it is; in either case, morphisms from or to such an object are large matrices and will be described using the $\oplus$ notation of the previous paragraph. We just want the reader to internalize that $\BBB_0 \BBB_0$ is bigger than just $\BBB_{00}$, so the fact that $0$ appears twice is somewhat misleading.

Finally, if we describe a morphism $\BBB_0 \BBB_0 \oplus R \to \BBB_0 \oplus \BBB_{020}$, we will use a $2 \times 2$ matrix, whose entries will be $\tau$-invariant morphisms described
with the above notation. This $2 \times 2$ matrix really has entries which are themselves matrices of morphisms, and it can be expanded to a large $2n \times (n^2 + 1)$ matrix of
morphisms if desired. \end{subequations} \end{notation}

With this notation, let us rewrite the twisted standard complexes when $n=2$ and $n=3$.

When $n=2$,
\begin{equation} \label{eq:FC2redux} \FC = \left(
\begin{tikzpicture}
\node (a) at (0,0) {$R(-1)$};
\node (b) at (2,0) {$\BBB_1$};
\node (c) at (4,0) {$R(1)$};
\path
	(a) edge node[descr] {$\startdotblue$} (b)
	(b) edge node[descr] {$\finaldotblue$} (c);
\end{tikzpicture} \right)
\end{equation}

When $n=3$,
\begin{equation}\label{eq:FC3redux} \FC = \left(
\begin{tikzpicture}
\node (a) at (0,0) {$R(-2)$};
\node (b) at (2.5,0) {$\BBB_2(-1)$};
\node (c1) at (6,-1) {$R$};
\node (c2) at (6,2) {$\BBB_{21}$};
\node (d) at (9.5,0) {$\BBB_2(1)$};
\node (e) at (12,0) {$R(2)$};
\path
	(a) edge node[descr] {$\startdotblue$} (b)
	(b) edge node[descr] {$\finaldotblue$} (c1)
	(b) edge node[descr] {$\startdotblue \linered \oplus - \lineblue \startdotred$} (c2)
	(c2) edge node[descr] {$\finaldotblue \linered \oplus - \lineblue \finaldotred$} (d)
	(c1) edge node[descr] {$\startdotblue$} (d)
	(d) edge node[descr] {$\finaldotblue$} (e);
\end{tikzpicture} \right)
\end{equation}

The efficiency of this notation increases dramatically when working with tensor products. Writing $\FC$ with $6$ summands rather than $12$ is nice, but writing $\FC \ot \FC$ with $36$ instead of $144$ is almost essential to getting anything done.

\begin{notation} We can safely ignore most grading shifts, since they will always be perverse. However, in a tensor product of complexes, we will often keep track of the grading shifts to help disambiguate terms. For example, $R$ appears both in homological degree $1$ and $-1$ inside $\FC$ for $n=2$, so that in degree zero in the tensor product $\FC \ot \FC$ there is a summand $R(-1) \ot R(1)$ and a summand $R(1) \ot R(-1)$. We shorten these to $R(-1)R(1)$ and $R(1)R(-1)$. \end{notation}

%=================
\subsection{Webs: $n=2$}
%=================

It may help to review \S\ref{subsec-symmetries}.

The category $\Webs_2$ is generated by a trivalent vertex sending $2 \to 1 \ot 1$, and its images under the various symmetries of the category (these symmetries were discussed in
\S\ref{subsec-symmetries}). Note that $\VC_{\fund_2} \cong \Om^2$ is the determinant representation, so that $\Om^2$ is an invertible tensor factor of both sides. Ignoring this tensor
factor, we seek a morphism $R \to \FC \ot \FC$.

We see no reason to work with the complex $\FC \ot \FC$, when we can replace it with the isomorphic complex $\FC \sot \FC$ (see \S\ref{subsec-signconvention}) which is easier to work
with because it obeys the strand-counting sign rule. Our favorite isomorphism between $\FC \ot \FC$ and $\FC \sot \FC$ will just multiply all terms of the form $\FC^i \ot \FC^j$ by $(-1)^{i+1}$. In particular, this isomorphism is the identity map on the underlying vector space, because all copies of $R$ in $\FC$ appear in odd homological degree. Hence any statements we make about the underlying vector space still hold if one works with $\FC \ot \FC$ instead of $\FC \sot \FC$.

The complex $\FC \sot \FC$ is as follows. 
\begin{equation} \FC \sot \FC = \left(
\begin{tikzpicture}
\node (a) at (0,0) {$R(-1)R(-1)$};
\node (b1) at (3,3) {$R(-1)\BBB_1$};
\node (b2) at (3,-3) {$\BBB_1 R(-1)$};
\node (c1) at (6,4) {$R(-1)R(1)$};
\node (c2) at (6,0) {$\BBB_1 \BBB_1$};
\node (c3) at (6,-4) {$R(1)R(-1)$};
\node (d1) at (9,3) {$\BBB_1 R(1)$};
\node (d2) at (9,-3) {$R(1) \BBB_1$};
\node (e) at (12,0) {$R(1)R(1)$};
\path
	(a) edge node[descr] {$\startdotblue$} (b1)
	(a) edge node[descr] {$\startdotblue$} (b2)
	(b1) edge node[descr] {$\finaldotblue$} (c1)
	(b2) edge node[descr] {$\finaldotblue$} (c3)
	(b1) edge node[descr] {$\startdotblue \lineblue \oplus \startdotred \lineblue$} (c2)
	(b2) edge node[descr] {$- \lineblue \startdotblue \oplus -\lineblue \startdotred$} (c2)
	(c2) edge node[descr] {$\finaldotblue \lineblue \oplus \finaldotred \lineblue$} (d2)
	(c2) edge node[descr] {$- \lineblue \finaldotblue \oplus -\lineblue \finaldotred$} (d1)
	(c1) edge node[descr] {$\startdotblue$} (d1)
	(c3) edge node[descr] {$\startdotblue$} (d2)
	(d1) edge node[descr] {$\finaldotblue$} (e)
	(d2) edge node[descr] {$\finaldotblue$} (e);
\end{tikzpicture} \right)
\end{equation}

The trivalent vertex is sent by $\GC$ to the map from $R$ to $(\FC\FC)^0 = R(-1)R(1) \oplus \BBB_1\BBB_1 \oplus R(1) R(-1)$ which is given by the (transpose of the) matrix
\begin{equation}\label{eq:cupn2} \left( \begin{array}{cccccc} -1 & -\cupblue & +1 \end{array} \right). \end{equation}
To unravel the notation just once, $-\cupblue$ in \eqref{eq:cupn2} is actually a $1 \times 4$ matrix, indicating that the map from $R$ to $B_0 B_0$ is $-\cupred$, the map from $R$ to $B_1 B_1$ is $-\cupblue$, the map from $R$ to $B_0 B_1$ is zero, and the map from $R$ to $B_1 B_0$ is zero. We leave the reader to verify that this is a chain map.

Let $\{e_1, e_2\}$ correspond to a basis for the underlying vector space of $\FC$, that is, $e_2$ corresponds to $R(-1)$ and $e_1$ corresponds to $R(+1)$. Let $x$ be the copy of $R(0)$
inside $\Om^{-2} \VC_{\fund_2}$. As indicated by \eqref{eq:cupn2}, our map sends $x \mapsto e_1 e_2 - e_2 e_1$ on the underlying vector space, which is exactly the morphism of $\gl_2$
representations $\Lambda^2 V \to V \ot V$ corresponding to the trivalent vertex, as desired.

We define $\GC$ on the other generators of $\Webs_2$ by symmetry. It is easy to verify by hand that the $\gl_2$ web relations are satisfied.

There is another compatibility which should be checked: that with geometric Satake, see \S\ref{subsec:intro_satake}. This can be checked by seeing what happens after tensor product with
$B_0$. The complex $(\FC \sot \FC) B_0$ can be Gaussian eliminated until all that remains is a single copy of $B_0 B_1 B_0$ in homological degree zero. It is an exercise to keep track of
what this Gaussian elimination does to the map $B_0 \to (\FC \sot \FC) B_0$, and show that the result is a $\ig{.5}{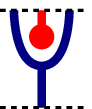}$. (Hint: the third entry $+1$ of \eqref{eq:cupn2}
tensored with $\lineblue$, followed by a signed homotopy $- \splitblue$, and then the differential $-\lineblue \startdotred$ tensored with $\lineblue$, altogether yields the ``zigzag''
which becomes $\ig{.5}{pitchforkopp}$. The remaining zigzags vanish.) On the other hand, this is also the map induced by geometric Satake (this is the non-singular version of the map
given in \cite[Definition 2.12, but upside-down]{EQuantumI}).

%=================
\subsection{The irreducible representations: $n=2$}
%=================

We let $S = \{s,t\}$ be the simple reflections of the infinite dihedral group.

Fix $\lambda \in \Z_{\ge 0}$. Let
\begin{equation} W^\lambda_k = \{ w \in W_{\aff} \mid \ell(w) = m \text{ where } 0 \le m \le \lambda -|k| \text{ and } \lambda-k-m \text{ is even} \}. \end{equation}

For example, if $\lambda$ is odd then $W^{\lambda}_0 = \{s, t, sts, tst, ststs, tstst, \ldots\}$ where the sequence stops at length $\lambda$, and $W^{\lambda}_1 = W^{\lambda}_{-1} = \{1,
st, ts, stst, tsts, \ldots\}$ where the sequence stops at length $\ell-1$. Note that $W^{\lambda}_{\pm \lambda} = \{1\}$ and $W^{\lambda}_k$ is empty for $|k| > \lambda$.

A thick calculus for the (finite and infinite) dihedral groups will appear soon, and this makes it easy to explicitly define the differentials in the complex below. Let us note for now
that there are always canonical morphisms $B_w \to B_x$ and $B_x \to B_w$ of degree $\ell(w) - \ell(x)$ whenever $x < w$ in the Bruhat order (which in the infinite dihedral group, just
means that $\ell(x) < \ell(w)$). These two maps are dual under $\DM$, and the map $B_w \to B_x$ preserves the lowest degree basis element $c_{\bot}$ (see \cite{EWHodge} for more on $c_{\bot}$).

\begin{defn} Fix $\lambda \in \Z_{\ge 0}$. Let $\VC_{\lambda}$ denote the following pseudocomplex. The chain object in homological degree $k$ is
\begin{equation} \Om^{\lambda} \bigoplus_{w \in W^{\lambda}_k} B_w(k). \end{equation}
The differential from $B_w(k)$ to $B_x(k+1)$ is a signed canonical morphism when $\ell(w) = \ell(x) \pm 1$, and is zero otherwise.
\end{defn}

\begin{ex} Here is the case $\lambda = 2$, where all the indecomposables $B_w$ which appear are actually Bott-Samelsons, so we can be explicit without needing thick calculus.
\begin{equation} \label{eq:VClambda2} \VC_\lambda = \Om^2 \left(
\begin{tikzpicture}
\node (a) at (0,0) {$R(-2)$};
\node (b) at (2.5,0) {$\BBB_1(-1)$};
\node (c1) at (6,-1) {$R$};
\node (c2) at (6,2) {$\BBB_{01}$};
\node (d) at (9.5,0) {$\BBB_1(1)$};
\node (e) at (12,0) {$R(2)$};
\path
	(a) edge node[descr] {$\startdotblue$} (b)
	(b) edge node[descr] {$\finaldotblue$} (c1)
	(b) edge node[descr] {$\startdotblue \linered \oplus - \lineblue \startdotred$} (c2)
	(c2) edge node[descr] {$\finaldotblue \linered \oplus - \lineblue \finaldotred$} (d)
	(c1) edge node[descr] {$\startdotblue$} (d)
	(d) edge node[descr] {$\finaldotblue$} (e);
\end{tikzpicture} \right)
\end{equation}
There is a striking similarity between \eqref{eq:VClambda2} and \eqref{eq:FC3redux}, but it is misleading and does not generalize. Of course, the meaning of $\BBB_1$ is different in the two contexts; here it has two summands, while in \eqref{eq:FC3redux} it has three. \end{ex}

%=================
\subsection{The generating webs: $n=3$}
%=================

The category $\Webs^+_3$ is generated by two morphisms and their images under the various symmetries: the trivalent vertex $1 \ot 1 \to 2$ and the trivalent vertex $3 \to 1 \ot 2$. We describe their images under $\GC$. After applying the appropriate factor of $\Om$, and considering more convenient replacement complexes, we are looking for morphisms $\FC \sot \FC \to \si(\FC)$ and $R \to \FC \sot \si(\FC)$ respectively.

\begin{remark} For $n=3$, $\FC \ot (-)$ and $\FC \sot (-)$ are isomorphic, since $\FC$ is an even BS complex. We write $\sot$ just to emphasize that the complexes in question satisfy the strand-counting sign rule. \end{remark}

We begin with the map $3 \to 1 \ot 2$, i.e. $R \to \FC \sot \si(\FC)$. In homological degree $0$, $\FC \sot \si(\FC)$ is a large direct sum with 36 terms. Some of these terms, like $B_{21} \ot B_{20}$, do not admit any morphisms of degree $0$ from $R$. The
relevant terms are: \begin{equation} \label{eq:V1V2oplus} R(-2) R(2) \oplus  \boxed{B_0(-1) B_0(+1)} \oplus R(0)R(0) \oplus \boxed{B_{21}(0) B_{12}(0)}  \oplus 
\boxed{B_0(1) B_0(-1)} \oplus R(2) R(-2). \end{equation} The chain map sends $R$ to $R(-2)R(2)$ and $R(2) R(-2)$ by multiplication by $+1$, to $R(0) R(0)$ by multiplication by $-1$, to $\boxed{B_0(-1) B_0(+1)}$ and $\boxed{B_0(1) B_0(-1)}$ by minus the cup, and to $\boxed{B_{21}(0) B_{12}(0)}$ by the doublecup.

\begin{remark} In general, we expect the chain map $R \to \FC \sot \si(\FC)$ to be given in homological degree zero by a signed sum of $X$-multicups. The $X$-multicup is a nested cup $R \to B_{h_X} B_{h_X\inv}$, and is the only diagram of degree $0$ in this hom space. In particular, every summand of this morphism which could be nonzero is nonzero. \end{remark}

Now we discuss the morphism $\FC \sot \FC \to \si(\FC)$ corresponding to the trivalent vertex $1 \ot 1 \to 2$. In homological degrees not between $-2$ and $+2$ it is zero, since
$\si(\FC)$ is zero in these degrees. There are some possible homotopies in this map, namely $\mergeblue$ viewed as a map \[\BBB_0(-1)\BBB_0(+1) \to \BBB_0(-1), \quad \BBB_0(+1)\BBB_0(-1) \to \BBB_0(-1), \quad \textrm{or} \;\; \BBB_0(+1)\BBB_0(+1) \to \BBB_0(+1).\] This gives some freedom in the description of this chain map; here is one representative.

In homological degree $-2$ it is
\begin{equation*} R(-2) \BBB_{10}(0) \oplus R(-2)R(0) \oplus \BBB_0(-1)\BBB_0(-1) \oplus \BBB_{10}(0)R(-2) \oplus R(0)R(-2) \to R(-2), \end{equation*} \begin{equation} \left( \begin{array}{ccccc} 0 & 1 & -\capblue & 0 & -1 \end{array} \right). \end{equation}
In homological degree $-1$ it is
\begin{equation*} R(-2)\BBB_0(-1) \oplus \BBB_0(-1)\BBB_{10}(0) \oplus \BBB_0(-1)R(0) \oplus \BBB_{10}(0) \BBB_0(-1) \end{equation*}
\begin{equation*} \quad \quad \oplus R(0) \BBB_0(-1) \oplus \BBB_0(1) R(-2) \to \BBB_0(-1), \end{equation*}
\begin{equation}
\left( \begin{array}{cccccc}
\lineblue & \mergeblue \finaldotred \oplus - \capred \lineblue & 0 & \finaldotred \mergeblue \oplus - \lineblue \capgreen & 0 & -\lineblue \end{array} \right). \end{equation}
In homological degree $0$ it is
\begin{equation*} R(-2)R(2) \oplus \BBB_0(-1)\BBB_0(+1) \oplus \BBB_{10}(0)\BBB_{10}(0) \oplus \BBB_{10}(0)R(0) \oplus R(0)\BBB_{10}(0) \oplus R(0)R(0) \end{equation*}
\begin{equation*} \quad \quad \oplus \BBB_0(1) \BBB_0(-1) \oplus R(2)R(-2) \to \BBB_{01}(0) \oplus R(0), \end{equation*}
\begin{equation} \label{eq:112n3deg0}
\left( \begin{array}{cccccccc}
0 & \linered \lineblue \oplus - \mergered \startdotblue & - \linered \capgreen \lineblue \oplus - \ig{.5}{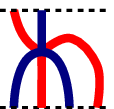} & 0 & 0 & 0 & - \linered \lineblue \oplus \startdotblue \mergegreen & 0 \\
1 & 0 & 0 & 0 & 0 & 0 & 0 & -1
\end{array} \right). \end{equation}
In homological degree $+1$ it is
\begin{equation*}\BBB_0(-1)R(2) \oplus \BBB_{10}(0)\BBB_0(1) \oplus R(0)\BBB_0(1) \oplus  \BBB_0(1)\BBB_{10}(0) \end{equation*}
\begin{equation*} \quad \quad \oplus  \BBB_0(1)R(0) \oplus  R(2)\BBB_0(-1) \to \BBB_0(+1), \end{equation*}
\begin{equation}
\left( \begin{array}{cccccc}
\lineblue & \pitchfork \oplus - \lineblue \capgreen & 0 & \pitchfork \oplus - \capred \lineblue & 0 & -\lineblue \end{array} \right). \end{equation}
In homological degree $+2$ it is
\begin{equation*} \BBB_{10}(0)R(2)  \oplus R(0)R(2) \oplus \BBB_0(1)\BBB_0(-1) \oplus R(2)\BBB_{10}(0) \oplus R(2)R(0) \to R(2), \end{equation*} \begin{equation} \left( \begin{array}{ccccc} 0 & 1 & -\capblue & 0 & -1 \end{array} \right). \end{equation}

% funkyhuh is a 2mvalent vertex with bottom boundary 1010 and top 01.

We leave the reader to confirm that these are chain maps, and that they match the webs on the underlying vector space. For example, the sixth column of \eqref{eq:112n3deg0} is zero,
which corresponds to the fact that $e_2 e_2 \mapsto 0$ under the map $V \ot V \to \Lambda^2 V$. Confirming the compatibility with geometric Satake can be done by Gaussian elimination,
but it is a very hard exercise.

%=================
\subsection{The complex for $S^2 V$: $n=3$}
%=================

Here is the complex without the differentials.

\begin{equation} \label{eq:S2Vn3} S^2 \VC = \Om^2 \left(
\begin{tikzpicture}
\node (a) at (0,0) {$R(-4)$};
\node (b) at (1.5,0) {$\BBB_1(-3)$};
\node (c1) at (3,-1) {$R(-2)$};
\node (c2) at (3,0) {$\BBB_{12}(-2)$};
\node (c3) at (3,1) {$\BBB_{21}(-2)$};
\node (d1) at (4.5,-1.5) {$\BBB_1(-1)$};
\node (d2) at (4.5,-.5) {$\BBB_1(-1)$};
\node (d3) at (4.5,.5) {$\BBB_{121}(-1)$};
\node (d4) at (4.5,1.5) {$\BBB_{210}(-1)$};
\node (e1) at (6,-2) {$R(0)$};
\node (e2) at (6,-1) {$R(0)$};
\node (e3) at (6,0) {$\BBB_{12}(0)$};
\node (e4) at (6,1) {$\BBB_{21}(0)$};
\node (e5) at (6,2) {$\BBB_{2102}(0)$};
\node (f1) at (7.5,-1.5) {$\BBB_1(1)$};
\node (f2) at (7.5,-.5) {$\BBB_1(1)$};
\node (f3) at (7.5,.5) {$\BBB_{121}(1)$};
\node (f4) at (7.5,1.5) {$\BBB_{210}(1)$};
\node (g1) at (9,-1) {$R(2)$};
\node (g2) at (9,0) {$\BBB_{12}(2)$};
\node (g3) at (9,1) {$\BBB_{21}(2)$};
\node (h) at (10.5,0) {$\BBB_1(3)$};
\node (i) at (12,0) {$R(4)$};
\end{tikzpicture} \right) \end{equation}

Let us use even more shorthand to write the differential. All the differentials will be one of the following ($\tau$-orbits of) morphisms or their vertical flips.
\begin{equation} \startdotblue \quad -\lineblue \startdotred \quad -\lineblue \startdotgreen \quad \startdotred \lineblue \quad \startdotgreen \lineblue \quad \dotfirsttoteal \quad \ig{.5}{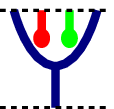} \end{equation}
%thelongone is the degree 1 map t to tsut.
Instead of drawing these morphisms, we just keep track of their coefficients. Recall that teal represented the indecomposable Soergel bimodule $B_{121}$, see \S\ref{subsec-stdrep3} for more details.

The differential from degree $-4$ to $-3$ is
\begin{equation} \left( \begin{array}{c} 1 \end{array} \right). \end{equation}
The differential from degree $-3$ to $-2$ is
\begin{equation} \left( \begin{array}{c} 1 \\ 1 \\ 1 \end{array} \right). \end{equation}
The differential from degree $-2$ to $-1$ is
\begin{equation} \left( \begin{array}{ccc} 1 & 1 & 0 \\ 1 & 1 & 0 \\ 1 & 1 & 2 \\ 0 & 1 & 1 \end{array} \right). \end{equation}
The differential from degree $-1$ to $0$ is
\begin{equation} \left( \begin{array}{cccc} 1 & 1 & 1 & -2 \\ 1 & 1 & 1 & 0 \\ 1 & 1 & 1 & -1 \\ 0 & 0 & 1 & -1 \\ 0 & 0 & 1 & 0 \end{array} \right). \end{equation}
The differential from degree $0$ to $+1$ is
\begin{equation} \left( \begin{array}{ccccc} 1 & 1 & 2 & 0 & 0 \\ 1 & 1 & 2 & 0 & 0 \\ 1 & 1 & 2 & 2 & 1 \\ -1 & 0 & -1 & -1 & 0 \end{array} \right). \end{equation}
The differential from degree $+1$ to $+2$ is
\begin{equation} \left( \begin{array}{cccc} 1 & 1 & 1 & 0 \\ 1 & 1 & 1 & 2 \\ 0 & 0 & 1 & 1 \end{array} \right). \end{equation}
The differential from degree $+2$ to $+3$ is
\begin{equation} \left( \begin{array}{ccc} 1 & 1 & 2 \end{array} \right). \end{equation}
The differential from degree $+3$ to $+4$ is
\begin{equation} \left( \begin{array}{c} 1 \end{array} \right). \end{equation}

On the underlying vector space, this description of the complex corresponds to the basis of $S^2 V \sumset V \ot V$ given by
\begin{equation} \{e_1 \ot e_1, e_1 \ot e_2 + e_2 \ot e_1 , e_1 \ot e_3 + e_3 \ot e_1, e_2 \ot e_2, e_2 \ot e_3 + e_3 \ot e_2, e_3 \ot e_3\}. \end{equation}

%=================
\subsection{Braiding: $n=2$} \label{subsec:braidingn2}
%=================

We abandon the notation $\BBB_w$ henceforth, because we will be describing morphisms (some of them $\tau$-invariant) using methods which are not obviously $\tau$-invariant. (The notation continues to work in this particular section, but we still abandon it for reasons of sanity.)

Recall that 
\begin{equation*} \FC = \left(
\begin{tikzpicture}
\node (a) at (0,0) {$R(-1)$};
\node (a0) at (2,1) {$B_1$};
\node (a1) at (2,-1) {$B_0$};
\node (b) at (4,0) {$R(1)$};
\path
	(a) edge node[descr] {$\startdotblue$} (a0)
	(a) edge node[descr] {$\startdotred$} (a1)
	(a0) edge node[descr] {$\finaldotblue$} (b)
	(a1) edge node[descr] {$\finaldotred$} (b);
%	\draw [brown] (current bounding box.south west) rectangle (current bounding box.north east);
\end{tikzpicture} \right)
\end{equation*}

The reader should reread \S\ref{ssec:intro_centrality}, as we are about to follow the process discussed therein for computing the commutation map $\phi_{\MC} \co \MC \ot \VC \to \VC \ot
\MC$ in the case $\MC = \VC$. Analogously, we can compute with the natural transformation $\phi \co (-) \ot \FC \to \FC \ot \tau(-)$.

We have already defined a homotopy equivalence $\phi_B \co B \ot \FC \to \FC \ot \tau(B)$ for all objects $B$ of $\Diag$. If there exists a unique chain map $\phi_{\MC}$ which preserves
the homological filtration on $\MC$ and acts on the associated graded as $\phi_{\MC^k[-k]}$ for each degree $k$, then $\phi_{\MC}$ is the commutation map.

The chain map $\phi_{B_1} \co B_1 \FC \to \FC B_0$ was computed in \S\ref{subsec-rotisoms}, which we now recall for the special case of $n=2$. In homological degrees $-1$ and $+1$ it is the zero map. In homological degree zero it is the map $B_1 B_1 \oplus B_1 B_0 \to B_1 B_0 \oplus B_0 B_0$ given by the matrix
\begin{equation} \label{eq:phiB1redux} \phi_{B_1}^0 = \left( \begin{array}{cc} \mergered \startdotblue & - \linered \lineblue \\ - \capredcupblue & \finaldotred \splitblue \end{array} \right). \end{equation}
The chain map $B_0 \FC \to \FC B_1$ is the color-rotation of this. The chain map $R \FC \to \FC R$ is the identity map.

\begin{remark} When comparing this with \S\ref{subsec-rotisoms}, the reader should note that $X = \mt$ is the only valid choice of block, and that $X \cup \{0,1\} = S$ so this term never appears. \end{remark}

Now we consider the ``braiding map'' $\phi_{\FC} \co \FC \FC \to \FC \FC$, and its homological filtration. It is fairly important that we work with $\FC \ot \FC$ and not $\FC \sot \FC$, so one must be careful when computing the signs on the differential. The source is viewed as
\begin{equation}
\left(
\begin{tikzpicture}
\node (a) at (0,0) {${\color{purple} R(-1)}$};
\node (a0) at (1,.5) {${\color{orange} B_1}$};
\node (a1) at (1,-.5) {${\color{pink} B_0}$};
\node (b) at (2,0) {${\color{teal} R(1)}$};
%	\draw [brown] (current bounding box.south west) rectangle (current bounding box.north east);
\end{tikzpicture} \right) \ot \FC, \end{equation}
so that, for instance, the source in homological degree zero is
\begin{equation} {\color{teal} R} \oplus {\color{orange} B_1 B_1 \oplus B_1 B_0} \oplus {\color{pink} B_0 B_1 \oplus B_0 B_0} \oplus {\color{purple} R}. \end{equation}
Here, we have colored the summands based on which part of the filtration they appear in; the teal summands form a subcomplex, and the purple summands a quotient complex. 
The target is viewed as
\begin{equation}
\FC \ot \left(
\begin{tikzpicture}
\node (a) at (0,0) {${\color{purple} R(-1)}$};
\node (a0) at (1,.5) {${\color{orange} B_0}$};
\node (a1) at (1,-.5) {${\color{pink} B_1}$};
\node (b) at (2,0) {${\color{teal} R(1)}$};
%	\draw [brown] (current bounding box.south west) rectangle (current bounding box.north east);
\end{tikzpicture} \right), \end{equation}
so that, for instance, the target in homological degree zero is
\begin{equation} {\color{teal} R} \oplus {\color{orange} B_1 B_0 \oplus B_0 B_0} \oplus {\color{pink} B_1 B_1 \oplus B_0 B_1} \oplus {\color{purple} R}. \end{equation}

The braiding map must preserve this filtration, sending teal to teal, pink to teal or pink, orange to teal or orange, and purple to anything. Thus it is a block upper-triangular matrix,
which acts by $\phi_{R(\mp 1)[\pm 1]}$ or $\phi_{B_i}$ on the block diagonal. In \S\ref{ssec:intro_braiding} we discussed that, while the map $\phi_R$ must be the identity map,
$\phi_{R[1]}$ is not the identity map. Instead, it is the chain map which multiplies the $k$-th homological degree in $\FC$ by $(-1)^k$.

Here is the unique way to fill in the upper-triangular pieces to obtain a chain map. We ignore all the grading shifts. In homological degree $-2$,
\begin{equation} \phi_{\FC}^{-2} \co {\color{purple} R} \to {\color{purple} R}, \quad \phi_{\FC}^{-2} = \left( \begin{array}{c} -1 \end{array} \right). \end{equation}
In homological degree $-1$,
\begin{equation} \phi_{\FC}^{-1} \co {\color{orange} B_1} \oplus {\color{pink} B_0} \oplus {\color{purple} B_1} \oplus {\color{purple} B_0} \to {\color{orange} B_0} \oplus {\color{pink} B_1} \oplus {\color{purple} B_1} \oplus {\color{purple} B_0} , \quad \phi_{\FC}^{-1} = \left( \begin{array}{cccc} 0 & 0 & -\linered & 0 \\ 0 & 0 & 0 & -\lineblue \\ 0 & 0 & \linered & 0 \\ 0 & 0 & 0 & \lineblue \end{array} \right). \end{equation}
In homological degree $0$,
\begin{equation*} \phi_{\FC}^{0} \co {\color{teal} R} \oplus {\color{orange} B_1 B_1 \oplus B_1 B_0} \oplus {\color{pink} B_0 B_1 \oplus B_0 B_0} \oplus {\color{purple} R} \to {\color{teal} R} \oplus {\color{orange} B_1 B_0 \oplus B_0 B_0} \oplus {\color{pink} B_1 B_1 \oplus B_0 B_1} \oplus {\color{purple} R}, \end{equation*}
\begin{equation} \label{braiding0} \phi_{\FC}^{0} = \left( \begin{array}{cccccc}
-1 & \capred & 0 & 0 & \capblue & 0 \\
0 & \mergered \startdotblue & - \linered \lineblue & 0 & 0 & 0 \\
0 & - \capredcupblue & \finaldotred \splitblue & 0 & 0 & - \cupblue \\
0 & 0 & 0 & \finaldotblue \splitred & - \capbluecupred & - \cupred \\
0 & 0 & 0 & - \lineblue \linered & \mergeblue \startdotred & 0 \\
0 & 0 & 0 & 0 & 0 & -1
 \end{array} \right). \end{equation}
In homological degree $+1$,
\begin{equation} \phi_{\FC}^{+1} \co {\color{teal} B_1} \oplus {\color{teal} B_0} \oplus {\color{orange} B_1} \oplus {\color{pink} B_0} \to {\color{teal} B_1} \oplus {\color{teal} B_0} \oplus {\color{orange} B_0} \oplus {\color{pink} B_1}, \quad \phi_{\FC}^{+1} = \left( \begin{array}{cccc} \linered & 0 & -\linered & 0 \\ 0 & \lineblue & 0 & -\lineblue \\ 0 & 0 & 0 & 0 \\ 0 & 0 & 0 & 0 \end{array} \right). \end{equation}
In homological degree $+2$,
\begin{equation} \phi_{\FC}^{+2} \co {\color{teal} R} \to {\color{teal} R}, \quad \phi_{\FC}^{+2} = \left( \begin{array}{c} -1 \end{array} \right). \end{equation}

Note that $\phi_{\FC}$ is concentrated on the diagonal after passage to the underlying vector space. That is, the upper right corner of the matrix \eqref{braiding0} is zero. Consequently, the braiding matches the Satake braiding $e_i e_j \mapsto -e_j e_i$ on the underlying vector space, as discussed in \S\ref{ssec:intro_braiding}.

%=================
\subsection{The endomorphism $\chi$: $n=2,3$} \label{subsec:chi}
%=================

Let us give some quick motivation from the conjectures of Gorsky-Negut-Rasmussen \cite{GNR}. Recall that $\Hilb_n$ is the (isospectral) Hilbert scheme of $n$ (ordered) points in $\C^2$
which are all concentrated on the locus $Y=0$. The tautological sheaf $\TC$ is filtered by line bundles $\LC_i$, which essentially describe the location of the $i$-th point. Now $\TC$
comes naturally equipped with two commuting operators $X$ and $Y$ which preserve the filtration, where $Y$ is nilpotent.

The corresponding complex $\EC \in \ZC(\HC_{\fin})$ has a filtration with subquotients Jucys-Murphy complexes $J_i$, each appearing once. According to their conjecture, $\EC$ should also
have two commuting operators $X$ and $Y$, where $X$ has degree $[0](2)$ and $Y$ has degree $[2](-2)$. Now $\EC$ also has an action of the polynomial ring $R = \C[x_1, \ldots, x_n]$, and
the action of $X$ on $J_i$ in the associated graded should agree with $x_i$.

We expect that these operators $X$ and $Y$ lift to operators $\chi$ and $\mu$ on $\VC$ before flattening, where $\mu$ is the monodromy map. In order to have an action of $x_i$ in the
affine setting, we work with the standard affine realization, see Definition \ref{defn:standardaffine}. Our goal is to define an endomorphism $\chi \co \VC \to \VC(2)$ which acts (in
the associated graded) by multiplication by $x_i$ on the Wakimoto complex $Y_i$. In fact, $\chi$ will be the unique filtered chain map with the correct associated graded. We will
describe the corresponding endomorphisms of $\FC$ rather than $\VC$, being careful to right multiply by $x_i$.

Again, recall that 
\begin{equation*} \FC = \left(
\begin{tikzpicture}
\node (a) at (0,0) {$R(-1)$};
\node (a0) at (2,1) {$B_1$};
\node (a1) at (2,-1) {$B_0$};
\node (b) at (4,0) {$R(1)$};
\path
	(a) edge node[descr] {$\startdotblue$} (a0)
	(a) edge node[descr] {$\startdotred$} (a1)
	(a0) edge node[descr] {$\finaldotblue$} (b)
	(a1) edge node[descr] {$\finaldotred$} (b);
%	\draw [brown] (current bounding box.south west) rectangle (current bounding box.north east);
\end{tikzpicture} \right)
\end{equation*}
The endomorphism $\chi$ is multiplication by $x_2$ in homological degree $-1$, multiplication by $x_1$ in homological degree $+1$, and is given by the following matrix in homological degree $0$.
\begin{equation} \chi^0 = \left( \begin{array}{cc} \lineblue \poly{x_1} & \redtoblue \\ 0 & \linered \poly{x_2} \end{array} \right). \end{equation}
Every matrix coefficient in $\chi$ is determined by the fact that it acts by $x_i$ on the associated graded piece $Y_i$, except for the upper right corner $\redtoblue$ of the matrix $\chi^0$. This matrix coefficient is uniquely determined by the requirement that $\chi$ is a chain map.

Though it is not immediately transparent, it is a good exercise to see that $\chi$ and $\tau(\chi)$ are homotopic. By adding $\chi$ and $\tau(\chi)$ we can give the following
description of a morphism $\chi'$ homotopic to $2 \chi$ which is more obviously color-rotation invariant (modulo $\delta$): in degrees $-1$ and $+1$ it is
multiplication by $x_1 + x_2$, and in degree $0$ it is given by the matrix \begin{equation} (\chi')^0 = \left( \begin{array}{cc} 2 \lineblue \poly{x_1}  & \redtoblue \\ \bluetored & 
2 \linered \poly{x_2} \end{array} \right). \end{equation} Note however that $\chi'$ does not preserve the Wakimoto filtration.

\begin{remark} In the main body of the paper we were considering endomorphisms of perverse complexes of degree $0$, which have no homotopies. Now we examine endomorphisms of degree $2$,
which may have homotopies! With the addition of grading shifts, all hell breaks loose. \end{remark}

\begin{remark} To confirm that $\chi$ is a morphism in the Drinfeld center, one must check that $\chi$ commutes with the natural transformation $\phi$. This can be checked on the generating objects $\Om$ and $B_i$ (and one need only check $B_i$ for one simple reflection). That $\chi$ commutes with $\phi_{\Om}$ is equivalent to it being fixed by $\tau$. That $\chi$ commutes with $\phi_{B_i}$ is an exercise. \end{remark}

For $n=3$, recall that
\begin{equation*}
\FC = \left(
\begin{tikzpicture}
\node (a) at (0,0) {$R(-2)$};
\node at (2,1) {$B_2(-1)$};
\node (b) at (2,0) {$B_1(-1)$}; 
\node at (2,-1) {$B_0(-1)$};

\node at (4,.5) {$R$};
\node at (4,1.5) {$B_{21}$};
\node at (4,-1.5) {$B_{10}$};
\node at (4,-.5) {$B_{02}$};
\node (c) at (4,0) {$\;\;\;$};

\node  at (6,1) {$B_2(1)$};
\node (d) at (6,0) {$B_1(1)$};
\node at (6,-1) {$B_0(1)$};
\node (e) at (8,0) {$R(2)$};
\path
	(a) edge (b)
	(b) edge (c)
	(c) edge (d)
	(d) edge (e);	
\end{tikzpicture} \right)
\end{equation*}
The endomorphism $\chi$ is multiplication by $x_3$ in homological degree $-2$. In homological degree $-1$ it is given by
\begin{equation} \chi^{-1} = \left( \begin{array}{ccc} \lineblue \poly{x_2} & \redtoblue & \greentoblue \\ 0 & \linered \poly{x_3} & 0 \\ 0 & 0 & \linegreen \poly{x_3} \end{array} \right). \end{equation}
In homological degree $0$ it is given by
\begin{equation} \chi^0 = \left( \begin{array}{cccc} \lineblue \linered \poly{x_1} & - \startdotblue \startdotred & - \finaldotgreen \lineblue \startdotred & - \startdotblue \linered \finaldotgreen \\ 0 & \poly{x_2} & 0 & 0 \\ 0 & 0 & \linegreen \lineblue \poly{x_2} & - \finaldotred \linegreen \startdotblue \\ 0 & 0 & 0 & \linered \linegreen  \poly{x_3} \end{array} \right). \end{equation}
In homological degree $+1$ it is given by
\begin{equation} \chi^{+1} = \left( \begin{array}{ccc} \lineblue \poly{x_1} & 0 & 0 \\ 0 & \linered \poly{x_1} & \greentored \\ 0 & 0 & \linegreen \poly{x_2} \end{array} \right). \end{equation}
In homological degree $+2$ it is multiplication by $x_1$. Again, all but the (block-)upper triangular part is determined on the associated graded, and the upper triangular part is unique under the requirement that $\chi$ is a chain map.

%=================
\subsection{The complete endomorphism ring: $n=2$} \label{subsec:endosn2}
%=================

First we compute the bigraded endomorphism ring $\END_{\HC_{\ext}}(\FC)$. It is an easy but nontrivial computation that, modulo homotopy and multiplication by $\delta$, there are no
morphisms $\FC \to \FC(k)[c]$ for any $k \in \Z$ when $c \ne 0, 2$. Of course we already know about the identity map $\id$ and the monodromy map $\mu$. It is also not hard to check that
$\mu \cdot f$ is nulhomotopic, for $f \in R$, if and only if $f$ is in the ideal $(\al_0,\al_1)$.

Concerning endomorphisms $\FC \to \FC(2)[0]$, we have already computed a morphism $\chi$ in the previous section. This morphism is not defined in the realization spanned by the simple
roots, so let us describe another endomorphism $\theta \co \FC \to \FC(2)$ which is. In homological degree $-1$, $\theta$ is multiplication by $\al_0$. In degree $+1$, $\theta$ is multiplication by zero. In homological degree $0$, $\theta$ is given by
\begin{equation} \theta^0 = \left( \begin{array}{cc} 0 & \redtoblue \\ 0 & \linered \poly{\al_0} \end{array} \right). \end{equation}
One can also verify directly that an $R$-linear combination $\id \cdot f_1 + \theta \cdot f_2$ is nulhomotopic if and only if $f_1$ and $f_2$ are multiples of $\delta$.
In the standard affine realization one can compute that $\chi = \id \cdot x_1 + \theta$ modulo homotopy and $\delta$, so either $\chi$ or $\theta$ could be used interchangeably as a generator.

It is easy to compute that
\begin{equation} \theta^2 = \theta \cdot \al_0, \end{equation} from which we deduce that
\begin{equation} \chi^2 = \chi \cdot (x_1 + x_2) - \id \cdot x_1 x_2 \end{equation}
modulo $\delta$. In addition,
\begin{equation} \theta \mu = \mu \theta = 0, \end{equation}
\begin{equation} \chi \mu = \mu \chi = \mu x_1 = \mu x_2. \end{equation}
In particular, $\mu$ and $\theta$ and $\chi$ commute.

Thus the endomorphism ring of $\FC$ in $\HC_{\ext}$ in any realization is the commutative ring
\begin{equation} \END(\FC) \cong R[\theta,\mu]/(\delta, \mu \al_1, \theta(\theta - \al_0), \mu \theta, \mu^2). \end{equation}
Note that an element $f \in R$ here represents the endomorphism of right multiplication by $f$, so that $f \mu = \mu f$ is the composition of $\mu$ and right multiplication by $f$.

Over the standard affine realization, using the generator $\chi$ instead, we have the presentation
\begin{equation} \label{eq:endFC} \END(\FC) \cong R[\chi,\mu]/(\delta, \mu(x_1 - x_2), (\chi - x_1)(\chi - x_2), \mu (\chi-x_1),\mu^2). \end{equation}

Note that $\tau$ fixes $\chi$ but does not fix $\theta$. Similarly, $\tau$ fixes $(x_1 + x_2)$ and $x_1 x_2$ modulo $\delta$, but it does not fix $x_1$ or $x_2$ individually. Hence the endomorphism ring of $\VC$ inside $\ZC(\HC_{\ext})$ is the subring of \eqref{eq:endFC} given by the invariants under swapping $x_1$ and $x_2$. Letting $R' = \C[x_1 + x_2, x_1 x_2]$ (the finite invariant ring) we have
\begin{equation} \END_{\ZC}(\VC) \cong R'[\chi,\mu]/( (x_1+x_2)\mu - 2 \chi \mu, x_1 x_2 \mu - \chi^2 \mu, (\chi - x_1)(\chi - x_2), \mu^2). \end{equation}
The first two relations are simplified by the ``non-relation'' $\chi \mu = x_1 \mu = x_2 \mu$.

%=================
\subsection{The flattening: $n=2$}
%=================

The flattening functor sends $\Om$ to the Rouquier complex $F_1$. It sends $B_0$ to $F_1 B_1 F_1\inv$, which is isomorphic to $B_1$. However, the isomorphism between $F_1 B_1 F_1\inv$ and $B_1$ is not functorial, rather it adds a sign to certain morphisms, see \cite{EHBraidToolkit}. Up to a global choice of sign, we can assume that this isomorphism sends the startdot to the startdot, and the enddot to minus the enddot.

Thus the flattening of $\FC$ is the complex
\begin{equation} \flat(\FC) \cong \left(
\begin{tikzpicture}
\node (a) at (0,0) {$R(-1)$};
\node (a0) at (2,1) {$B_1$};
\node (a1) at (2,-1) {$B_1$};
\node (b) at (4,0) {$R(1)$};
\path
	(a) edge node[descr] {$\startdotblue$} (a0)
	(a) edge node[descr] {$\startdotblue$} (a1)
	(a0) edge node[descr] {$\finaldotblue$} (b)
	(a1) edge node[descr] {$-\finaldotblue$} (b);
%	\draw [brown] (current bounding box.south west) rectangle (current bounding box.north east);
\end{tikzpicture} \right).
\end{equation}
By changing basis in homological degree zero, this complex is isomorphic to
\begin{equation} \flat(\FC) \cong \left(
\begin{tikzpicture}
\node (a) at (0,0) {$R(-1)$};
\node (a0) at (2,1) {$B_1$};
\node (a1) at (2,-1) {$B_1$};
\node (b) at (4,0) {$R(1)$};
\path
	(a) edge node[descr] {$\startdotblue$} (a0)
	(a) edge node[descr] {$0$} (a1)
	(a0) edge node[descr] {$0$} (b)
	(a1) edge node[descr] {$\finaldotblue$} (b);
%	\draw [brown] (current bounding box.south west) rectangle (current bounding box.north east);
\end{tikzpicture} \right),
\end{equation}
or in other words $\flat(\FC) \cong F_s \oplus F_s\inv$.

Thus $\flat(\VC) \cong F_s F_s \oplus \one$.

%%%%%%%%%%%%%%%%%%%%%%%%%
\section{Proofs related to extended diagrammatics}
\label{sec:extendeddiagproofs}
%%%%%%%%%%%%%%%%%%%%%%%%%

%=================
\subsection{A proof from topology}
%=================

Let $\BC \subset \EDiag$ denote the following (non-full, non-graded, non-additive) monoidal subcategory. The objects are generated by $\Om^+$ and $\Om^-$. The morphisms are generated by
the oriented cups and caps of \eqref{eq:blackcupscaps}, modulo the relations \eqref{eq:cupcaprelations}. Unlike the category $\EDiag$, we do not allow linear combinations of diagrams or
polynomials, but only allow genuine diagrams as morphisms. (Note that the relations we apply do not involve any linear combinations.) Similarly, let $\BC_n \subset \EDiag_n$ denote the
non-full monoidal subcategory, extending $\BC$, which also allows the black $n$-valent vertex of \eqref{eq:blacknv}, modulo the relations \eqref{eq:nvalentrelns}. We refer to either $\BC$
or $\BC_n$ (depending on context) as the \emph{black subcategory}.

To prove that the diagrammatic categories $\EDiag$ and $\EDiag_n$ have the expected size, we will first show that the subcategories $\BC$ and $\BC_n$ have the expected size, which is the
most interesting part of the argument. To accomplish this we will use the same topological arguments as were found in \cite{EWFenn}. Let us briefly recall the main story, and refer the
reader to \cite{EWFenn} for more details and additional references.

Given a group $G$, its $2$-groupoid is the monoidal category $\KM_G$ whose objects are elements of $G$, where monoidal composition is multiplication in $G$. The only morphisms are
identity maps, and thus it is clear how the monoidal structure is defined on morphisms. We will prove the following theorem.

\begin{thm} \label{thm:2groupoid} The category $\BC$ is equivalent to $\KM_G$ for $G = \Z$. The category $\BC_n$ is equivalent to $\KM_G$ for $G = \Z/ n\Z$. \end{thm}

This is expected from Lemmas \ref{lem:noblackendo} and \ref{lem:noblackendo2}, since any two morphisms with the same source and target are equal, and these morphisms exist if and only if
the source and target represent equal elements of the corresponding group $G$.

Let $X$ be a topological space. One can view $\pi_2(X)$ as a 2-category: the objects are points of $X$, the 1-morphisms are paths from source to target, and the 2-morphisms are disks
glued on to the corresponding boundary, modulo homotopy. Suppose that $X$ is the classifying space $BG = K(G,1)$. Fixing a basepoint, the corresponding monoidal subcategory of $\pi_2(X)$
is equivalent to $\KM_G$.

There is a lovely diagrammatic calculus, consisting of \emph{Igusa diagrams}\footnote{It seems these diagrammatics were due independently to Igusa and Fenn.}, that takes a sufficiently
combinatorial CW complex $X$ and encodes elements of $\pi_2(X)$ (and of relative $\pi_2$ with respect to an element of $\pi_1(X)$) as planar diagrams. The diagrammatic calculus only
depends on the 3-skeleton of $X$, as one might expect since $\pi_2(X)$ only depends on the 3-skeleton. This can be packaged as a 2-category given by generators and relations, which is a
combinatorial rigidification of the category $\pi_2(X)$ above. It has one object for each 0-cell of $X$, a generating invertible 1-morphism for each 1-cell of $X$, and a generating
invertible 2-morphism for each 2-cell of $X$. In addition to a number of ``generic'' relations, each 3-cell of $X$ provides an additional relation.

In \cite[Chapter 2]{EWFenn} we recall this diagrammatic calculus in the style presented by Fenn's book \cite{Fenn}. We then use this to address the question of how one can extend a
presentation of $G$ by generators and relations to a diagrammatic presentation of $\KM_G$. A presentation of a group is the same data as a 2-dimensional CW-complex $X_2$ with one 0-cell,
where the generators give 1-cells and the relations give 2-cells. Note that $\pi_1(X_2) = G$. If we can determine which 3-cells to glue onto $X_2$ to produce a CW-complex $X$ for which
$\pi_2(X)=0$, then applying the Igusa diagram construction will produce a presentation for $\KM_G$. In \cite{EWFenn} we accomplish this task for the usual Coxeter presentation of a Coxeter
group, and conjecturally for braid groups as well. In this appendix we do the much easier groups $G = \Z$ and $G = \Z/n\Z$. The $n=2$ case is the same as the Coxeter group of type $A_1$,
which is done in great detail in \cite[Chapter 3]{EWFenn}. The general case is largely the same.

To be more precise, Igusa diagrams are associated to marked CW-complexes. Each $2$-cell has a marked point on its boundary, which is glued to some point in the $0$-skeleton.

\begin{proof}[Proof of Theorem \ref{thm:2groupoid}] We leave the reader to verify that the category $\BC$ is exactly the diagrammatic category for the circle $S^1$, coming from the presentation of $\Z$ with one generator and
no relations. This particular category is discussed in Fenn's book \cite{Fenn}; the relations \eqref{eq:cupcaprelations} appear there verbatim. Since $\pi_2(S^1)$ is trivial, $\BC \cong
\KM_G$ for $G = \Z$.

Now fix $n \ge 2$, and consider the following CW-complex $X_2$. It has \begin{itemize} \item one 0-cell, \item one 1-cell, which we call $\Om$, \item one 2-cell, glued in along $\Om^n$.
\end{itemize} Clearly $\pi_1(X_2) = \Z/n\Z$. The corresponding diagrammatic category is not quite $\BC_n$. It does contain $\BC$: these diagrams correspond to $S^2 \to X_2$ and
nulhomotopies which take place entirely within the 1-skeleton (which is again $S^1$). It also contains an additional pair of generators corresponding to the $2$-cell: these are drawn as $n$-valent vertices (one with orientations pointing in, the other with orientations pointing out). However, unlike the $n$-valent vertex of $\BC_n$, this $n$-valent vertex also has a \emph{tag} indicating where the marked point on the boundary of the $2$-cell went. Consequently, this $n$-valent vertex is not rotationally invariant: rotating by $360/n$ degrees will alter the placement of the tag. The morphism is still cyclic (rotation by $360$ degrees returns the tag to its original
location) so planar diagrams are still the appropriate tool.

Keeping track of the tag, one has the following relations for the Igusa diagram version of $\pi_2(X_2)$, instead of \eqref{eq:nvalentrelns}.
\begin{subequations} \label{eq:nvalentrelnswtag}
\begin{equation} \label{eq:blackbreaknwtag} {
\labellist
\small\hair 2pt
 \pinlabel {$=$} [ ] at 40 18
\endlabellist
\centering
\ig{1.5}{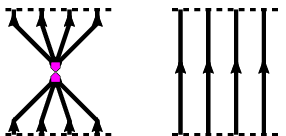}
} \end{equation}
\begin{equation} \label{eq:blackncirclewtag} {
\labellist
\small\hair 2pt
 \pinlabel {$=$} [ ] at 40 18
\endlabellist
\centering
\ig{1.5}{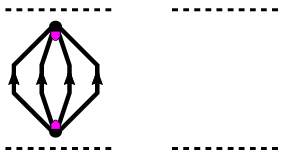}
} \end{equation}
\begin{equation} \label{eq:blackncyclicwtag} {
\labellist
\small\hair 2pt
 \pinlabel {$=$} [ ] at 56 21
\endlabellist
\centering
\ig{1.5}{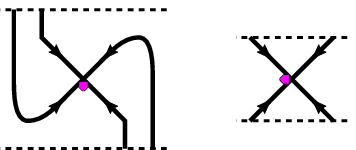}
} \end{equation}
\end{subequations}
Note that \eqref{eq:blackbreaknwtag} and \eqref{eq:blackncirclewtag} only apply when the tags are in matching regions.

Now $\pi_2(X_2)$ is nontrivial. A non-nulhomotopic function $S^2 \to X^2$ maps the northern hemisphere along the 2-cell, and the southern hemisphere along the orientation-reversed 2-cell, but twisted once by $\Om$ (via the conjugation action of $\pi_1$ on relative $\pi_2$). The corresponding Igusa diagram of this map $S^2 \to X^2$ is as follows.
\begin{equation} \label{eq:blackncirclewtagnonzero} \ig{1.5}{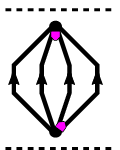} \end{equation}

So let us attach a 3-cell along the attaching map \eqref{eq:blackncirclewtagnonzero}, to produce $X = X_3$. This makes the map of \eqref{eq:blackncirclewtagnonzero} nulhomotopic in $X$,
so it becomes equal to the empty diagram. Any relation with empty boundary gives rise to a relation with boundary obtained by slicing the sphere into northern and southern hemispheres, as
explained in \cite{Fenn}. By slicing \eqref{eq:blackncirclewtagnonzero}, we deduce that the $n$-valent vertex with tag is now equal to its rotation, and the location of the tag is no
longer relevant. After forgetting the positioning of the tag, and the relations \eqref{eq:nvalentrelnswtag} become \eqref{eq:nvalentrelns}. Thus the Igusa category for $X$ is equivalent
(in fact, isomorphic) to $\BC_n$.

Thus we need only prove that $\pi_2(X)$ is trivial. Instead of performing this healthy topological exercise, let us make the argument in a slightly different language: where one unravels
the presentation to a description of $EG$ instead of $BG$. This train of thought is discussed at length in \cite[Chapter 2.4]{EWFenn}. Consider instead of the following CW-complex
$\tilde{X}$ which admits a free action of $G = \Z/ n\Z$: it has \begin{itemize} \item $n$ 0-cells, labeled by the elements of $G$, \item $n$ 1-cells, connecting the 0-cell $i$ with $i+1$
for each $i \in G$, \item $n$ 2-cells, glued in to the big loop $i \to i+1 \to \ldots \to n-1 \to 0 \to 1 \to \ldots \to i$ in the 1-skeleton, for each $i \in G$, \item $n$ 3-cells, which
are attached along \eqref{eq:blackncirclewtagnonzero} but with the tag based at $i$ for $i \in G$. \end{itemize} It is clear that the quotient of $\tilde{X}$ by the (free) action of (the
discrete group) $G$ is $X$. We need only prove that $\pi_2(\tilde{X})$ is also trivial.

Your picture when $n=4$ should go as follows, and the case of general $n$ is analogous. The 1-skeleton of $\tilde{X}$ is a square. Then one attaches four disks to this square (via the
same attaching map, if one forgets the base point), producing a stack of four pancakes stuck along their rims, or some bubbles attached to a square frame. Then one attaches four 3-cells.
Three of these fill in the gaps between the pancakes/bubbles, making one solid mass of pancake flesh homeomorphic to $D^3$. The fourth 3-cell attaches around the outside, along a map
which is now nulhomotopic (after the first three 3-cells are attached), and hence producing a space homeomorphic to $S^3$. Regardless, $\pi_2$ is trivial. \end{proof}

Again, to see a very similar argument done in more detail, see \cite[Chapter 3]{EWFenn}.

%=================
\subsection{Straightening}
%=================

In the rest of this chapter we prove several results about $\EDiag$. The corresponding results for $\EDiag_n$ follow by entirely analogous arguments, and there is no additional difficulty dealing with the black $n$-valent vertex.

Fix $B, C \in \Diag$ and $k, \ell \in \Z$, and consider the objects $X
= \Om^k B$ and $Y = \Om^\ell C$ in $\EDiag$. There is certainly a linear map
\[ \rho \co \Hom_{\BC}(\Om^k,\Om^\ell) \sqot \Hom_{\Diag}(B,C) \to \Hom_{\EDiag}(X,Y) \]
given by concatenating diagrams. This map is well-defined because the relations of $\BC$ and of $\Diag$ are also present in $\EDiag$. We call any morphism in the image of $\rho$ \emph{straight}.  Conversely, let $\phi$ be a diagram representing a morphism in $\Hom_{\EDiag}(X,Y)$. Using \eqref{eq:pullapart} and
\eqref{eq:mixedpoly} and \eqref{eq:pullnvalent}, we can pull apart the colored portion of the diagram from the black portion (including all polynomials in the colored portion). We call
the resulting diagram the \emph{straightening} of $\phi$. The next definition demonstrates that the straightening process is unambiguous, giving an inverse to $\rho$.

\begin{defn} Given a diagram $\phi$ representing a morphism in $\EDiag$, let $b(\phi)$ denote the \emph{black portion} of $\phi$, obtained by ignoring all colored strands and polynomials.
It is a diagram representing a morphism in $\BC$. Similarly, let $c(\phi)$ denote the \emph{colored portion} of $\phi$, obtained as follows: \begin{itemize} \item Take each colored edge
in the graph, and alter the color by applying $\tau\inv$ a number of times equal to the number of black strands to its right. \item Take each polynomial and apply $\tau\inv$ a number of
times equal to the number of black strands to its right. \item Now remove all black strands. \end{itemize} The result is a diagram representing a morphism in $\Diag$. \end{defn}

\begin{lemma} The above definition makes sense. That is, $b(\phi)$ and $c(\phi)$ are diagrams in $\BC$ and $\Diag$ respectively. \end{lemma}
	
\begin{proof} This can be checked easily for each of the defining vertices of $\EDiag$. Showing that this process glues monoidally is also straightforward. \end{proof}

\begin{lemma} \label{lem:phiisom} If $\phi$ is a morphism in $\Hom_{\EDiag}(\Om^k B, \Om^\ell C)$ with notation as above, then 
\begin{equation} \label{eq:phiandrho} \phi = \rho(b(\phi) \ot c(\phi)).\end{equation} Consequently, $\rho$ is an isomorphism.
\end{lemma}

\begin{proof} The equality \eqref{eq:phiandrho} is a straightforward application of \eqref{eq:pullapart} and \eqref{eq:mixedpoly} and \eqref{eq:pullnvalent}. Each of these relations does
not change the colored portion or the black portion of a diagram. Thus $\rho$ is surjective. It is also injective, since $b(\rho(\xi \ot \psi)) = \xi$ and $c(\rho(\xi \ot \psi)) = \psi$.
\end{proof}

\begin{proof}[Remainder of proof of Lemmas \ref{lem:pullapart} and \ref{lem:noblackendo}] That $\Hom(\Om^k B, \Om^\ell C) = 0$ unless $k = \ell$ is clear, because if $k \ne \ell$ then $\Hom_{\BC}(\Om^k,\Om^\ell)=0$.

Let us deduce that $\Hom_{\EDiag}(\Om^k,\Om^k) \cong R_{\tau^k}$ as $R$-bimodules. We know that $\Hom_{\Diag}(\one,\one) \cong R$. Also, $\Hom_{\BC}(\Om^k,\Om^k) = \{\id\}$ by Theorem
\ref{thm:2groupoid}. Thus by Lemma \ref{lem:phiisom}, we see that $\Hom_{\EDiag}(\Om^k,\Om^k) \cong R$ as vector spaces. Now using \eqref{eq:mixedpoly}, its $R$-bimodule structure is
clear. From this, Lemma \ref{lem:phiisom} implies \eqref{eq:homOmB}. We also see that any morphism in $\Hom_{\EDiag}(\Om^k,\Om^k)$ can be uncolored.

To argue that any morphism in $\Hom_{\EDiag}(B,C)$ can be unblackened, observe again that $\Hom_{\BC}(\one,\one) = \{\id\}$ by Theorem \ref{thm:2groupoid}. Thus the isomorphism of
\eqref{eq:homOmB} implies that any morphism in $\Hom_{\EDiag}(B,C)$ is in the image of $\Hom_{\Diag}(B,C)$, as desired. \end{proof}

%=================
\subsection{Remaining loose ends}
%=================

\begin{proof}[Proof of Lemma \ref{lem:onlyoneiso}] We know that any diagrams built from black cups and caps and from mixed crossings is an isomorphism, because these generators are isomorphisms. Let $\phi$ and $\psi$ be two such diagrams, from $X$ to $\Om^k B$. Recall the duality functor $\DM$ from \S\ref{subsec-duality}. Then the composition $\phi \circ \DM(\psi)$ is an automorphism of $\Om^k B$, so it is in the image of $\rho$. Moreover, applying $b$ or $c$ to a mixed crossing yields an identity map, while applying $b$ to a black cup/cap gives the black cup/cap, and $c$ gives the empty diagram. Thus $\phi \circ \DM(\psi)$ straightens to the identity map, so it is the identity map. Thus $\DM(\psi)$ is the inverse of $\phi$. But this was true for arbitrary such diagrams; in particular, $\DM(\phi)$ is also the inverse of $\phi$, and thus $\psi = \phi$. \end{proof}

\begin{proof}[Proof of Theorem \ref{thm:EDCQsize}] At this point, the only thing which remains unproven is the calculation of the Grothendieck group. Let $k, \ell \in \Z$ and $B, C \in
\EDiag$. We have already seen that any object in $\EDiag$ is isomorphic to one of the form $\Om^k \ot B$. Moreover, if $B$ and $C$ are non-isomorphic indecomposable objects in
$\Kar(\Diag)$, or if $k \ne \ell$, then $\Om^k \ot B$ and $\Om^\ell \ot C$ can not be isomorphic, by \eqref{eq:homOmB}. Thus the indecomposable objects of the Karoubi envelope, up to
grading shift and isomorphism, are enumerated by $\{\Om^k \ot B_w \}_{k \in \Z, w \in W}$. This gives one a basis for the Grothendieck group, after which point it is easy to verify the
statements of Theorem \ref{thm:EDCQsize}. \end{proof}

\begin{remark} We have omitted the argument that $\Kar(\EDiag)$ is Krull-Schmidt. That $\Kar(\Diag)$ is Krull-Schmidt was proven in \cite{EWGr4sb} using the fact that it is an
object-adapted cellular category. In fact, $\EDiag$ is also an object-adapted cellular category, having a double leaves basis. The proof of this is trivial given the corresponding facts
for $\Diag$, but recalling all this technology seems like overkill. \end{remark}

\bibliographystyle{plain}
\bibliography{mastercopy}

\begin{thebibliography}{10}

\bibitem{BFO}
Roman Bezrukavnikov, Michael Finkelberg, and Victor Ostrik.
\newblock Character {$D$}-modules via {D}rinfeld center of {H}arish-{C}handra
  bimodules.
\newblock {\em Invent. Math.}, 188(3):589--620, 2012.

\bibitem{BezYunMonodromy}
Roman Bezrukavnikov and Zhiwei Yun.
\newblock On {K}oszul duality for {K}ac-{M}oody groups.
\newblock {\em Represent. Theory}, 17:1--98, 2013.

\bibitem{CKM}
Sabin Cautis, Joel Kamnitzer, and Scott Morrison.
\newblock Webs and quantum skew {H}owe duality.
\newblock {\em Math. Ann.}, 360(1-2):351--390, 2014.

\bibitem{GNR}
A.~Negut E.~Gorsky and J.~Rasmussen.
\newblock Flag {H}ilbert schemes, colored projectors and {K}hovanov-{R}ozansky
  homology.
\newblock Preprint, 2016.
\newblock arXiv:1608.07308.

\bibitem{EFlat}
Ben Elias.
\newblock Flattening in type {A}.
\newblock In preparation.

\bibitem{ELLCC}
Ben Elias.
\newblock Light ladders and clasp conjectures.
\newblock Preprint, 2015.
\newblock arXiv 1510.06840.

\bibitem{EThick}
Ben Elias.
\newblock Thicker {S}oergel calculus in type {$A$}.
\newblock {\em Proc. Lond. Math. Soc. (3)}, 112(5):924--978, 2016.

\bibitem{ECathedral}
Ben Elias.
\newblock The two-color {S}oergel calculus.
\newblock {\em Compos. Math.}, 152(2):327--398, 2016.

\bibitem{EQuantumI}
Ben Elias.
\newblock Quantum {S}atake in type {A}: part {I}.
\newblock {\em J. Comb. Algebra}, 1(1):63--125, 2017.
\newblock arXiv:1403.5570.

\bibitem{EHDrinfeld}
Ben Elias and Matthew Hogancamp.
\newblock Extending isomorphisms of functors to the homotopy category.
\newblock In preparation.

\bibitem{EHBraidToolkit}
Ben Elias and Matthew Hogancamp.
\newblock A toolkit for conjugation by {R}ouquier complexes.
\newblock In preparation.

\bibitem{EHDiag}
Ben Elias and Matthew Hogancamp.
\newblock Categorical diagonalization.
\newblock Preprint, 2017.
\newblock arXiv:1707.04349.

\bibitem{EHDiag2}
Ben Elias and Matthew Hogancamp.
\newblock Categorical diagonalization of full twists.
\newblock Preprint, 2017.
\newblock arXiv:1801.00191.

\bibitem{ElHog16a}
Ben Elias and Matthew Hogancamp.
\newblock On the computation of torus link homology.
\newblock {\em Compositio Mathematica}, To appear.
\newblock arXiv:1603.00407.

\bibitem{EKho}
Ben Elias and Mikhail Khovanov.
\newblock Diagrammatics for {S}oergel categories.
\newblock {\em Int. J. Math. Math. Sci.}, pages Art. ID 978635, 58, 2010.

\bibitem{EKra}
Ben Elias and Dan Krasner.
\newblock Rouquier complexes are functorial over braid cobordisms.
\newblock {\em Homology, Homotopy Appl.}, 12(2):109--146, 2010.

\bibitem{EWHodge}
Ben Elias and Geordie Williamson.
\newblock The {H}odge theory of {S}oergel bimodules.
\newblock {\em Ann. of Math. (2)}, 180(3):1089--1136, 2014.

\bibitem{EWGr4sb}
Ben Elias and Geordie Williamson.
\newblock Soergel calculus.
\newblock {\em Represent. Theory}, 20:295--374, 2016.
\newblock arXiv:1309.0865.

\bibitem{EWFenn}
Ben Elias and Geordie Williamson.
\newblock Diagrammatics for {C}oxeter groups and their braid groups.
\newblock {\em Quantum Topol.}, to appear.
\newblock arXiv:0902.4700.

\bibitem{Fenn}
Roger~A. Fenn.
\newblock {\em Techniques of geometric topology}, volume~57 of {\em London
  Mathematical Society Lecture Note Series}.
\newblock Cambridge University Press, Cambridge, 1983.

\bibitem{Gait01}
D.~Gaitsgory.
\newblock Construction of central elements in the affine {H}ecke algebra via
  nearby cycles.
\newblock {\em Invent. Math.}, 144(2):253--280, 2001.

\bibitem{GORS}
Eugene Gorsky, Alexei Oblomkov, Jacob Rasmussen, and Vivek Shende.
\newblock Torus knots and the rational {DAHA}.
\newblock {\em Duke Math. J.}, 163(14):2709--2794, 2014.

\bibitem{Harterich}
M.~H\"arterich.
\newblock {\em Kazhdan-{L}usztig-{B}asen, unzerlegbare {B}imoduln und die
  {T}opologie der {F}ahnenmannigfaltigkeit einer {K}ac-{M}oody-{G}ruppe}.
\newblock PhD thesis, Albert-Ludwigs-Universit\"at Freiburg, 1999.

\bibitem{HogGor}
Matt Hogancamp and Eugene Gorsky.
\newblock Hilbert schemes and $y$-ification of {K}hovanov-{R}ozansky homology.
\newblock Preprint, 2017.
\newblock arXiv:1712.03938.

\bibitem{HumpCox}
James~E. Humphreys.
\newblock {\em Reflection groups and {C}oxeter groups}, volume~29 of {\em
  Cambridge Studies in Advanced Mathematics}.
\newblock Cambridge University Press, Cambridge, 1990.

\bibitem{JMW14}
Daniel Juteau, Carl Mautner, and Geordie Williamson.
\newblock Parity sheaves.
\newblock {\em J. Amer. Math. Soc.}, 27(4):1169--1212, 2014.

\bibitem{Khov07}
Mikhail Khovanov.
\newblock Triply-graded link homology and {H}ochschild homology of {S}oergel
  bimodules.
\newblock {\em Internat. J. Math.}, 18(8):869--885, 2007.

\bibitem{Kupe}
Greg Kuperberg.
\newblock Spiders for rank {$2$} {L}ie algebras.
\newblock {\em Comm. Math. Phys.}, 180(1):109--151, 1996.

\bibitem{LuszCharactersheaves85}
George Lusztig.
\newblock Character sheaves. {I}.
\newblock {\em Adv. in Math.}, 56(3):193--237, 1985.

\bibitem{MacThi}
Marco Mackaay and Anne-Laure Thiel.
\newblock Categorifications of the extended affine {H}ecke algebra and the
  affine {$q$}-{S}chur algebra {$\widehat{\bold S}(n,r)$} for {$3\leq r<n$}.
\newblock {\em Quantum Topol.}, 8(1):113--203, 2017.

\bibitem{MackyMonodromy}
S.~Makisumi.
\newblock Modular {K}oszul duality for {S}oergel bimodules.
\newblock Preprint, 2017.
\newblock arXiv:1703.01576.

\bibitem{Mellit}
Anton Mellit.
\newblock Homology of torus knots.
\newblock Preprint, 2017.
\newblock arXiv:1704:07630.

\bibitem{ORS}
Alexei Oblomkov, Jacob Rasmussen, and Vivek Shende.
\newblock The {H}ilbert scheme of a plane curve singularity and the {HOMFLY}
  homology of its link.
\newblock {\em Geom. Topol.}, 22(2):645--691, 2018.
\newblock With an appendix by Eugene Gorsky.

\bibitem{OblRoz}
Alexei Oblomkov and Lev Rozansky.
\newblock Knot homology and sheaves on the {H}ilbert scheme of points on the
  plane.
\newblock {\em Selecta Math. (N.S.)}, 24(3):2351--2454, 2018.

\bibitem{OkoVer}
Andrei Okounkov and Anatoly Vershik.
\newblock A new approach to representation theory of symmetric groups.
\newblock {\em Selecta Math. (N.S.)}, 2(4):581--605, 1996.

\bibitem{AMRW}
S.~Riche P.~Achar, S.~Makisumi and G.~Williamson.
\newblock Free-monodromic mixed tilting sheaves on flag varieties.
\newblock Preprint, 2017.
\newblock arXiv:1703.05843.

\bibitem{Rasmussen}
Jacob Rasmussen.
\newblock Some differentials on {K}hovanov-{R}ozansky homology.
\newblock {\em Geom. Topol.}, 19(6):3031--3104, 2015.

\bibitem{RouqBraid-pp}
Raphael Rouquier.
\newblock Categorification of the braid groups.
\newblock Preprint, 2004.
\newblock arXiv:math/0409593.

\bibitem{Soer90}
Wolfgang Soergel.
\newblock Kategorie {$\mathcal{O}$}, perverse {G}arben und {M}oduln \"uber den
  {K}oinvarianten zur {W}eylgruppe.
\newblock {\em J. Amer. Math. Soc.}, 3(2):421--445, 1990.

\bibitem{Soer07}
Wolfgang Soergel.
\newblock Kazhdan-{L}usztig-{P}olynome und unzerlegbare {B}imoduln \"uber
  {P}olynomringen.
\newblock {\em J. Inst. Math. Jussieu}, 6(3):501--525, 2007.

\bibitem{WillSingular}
Geordie Williamson.
\newblock Singular {S}oergel bimodules.
\newblock {\em Int. Math. Res. Not. IMRN}, (20):4555--4632, 2011.

\end{thebibliography}

\end{document}